\documentclass{article}
\usepackage[a4paper, margin=3cm, bottom=4cm]{geometry}
\usepackage[utf8]{inputenc}
\usepackage[T1]{fontenc}
\usepackage{lmodern}

\usepackage[
    backend=biber,
    style=alphabetic, maxalphanames=5, % write all initials, e.g "FBHW21" instead of "Foo+21"
    maxbibnames=100, maxsortnames=100, % disable "et al." in \printbibliography
    sorting=nyt,
    natbib=true,
    url=false,doi=false,isbn=false,eprint=false
]{biblatex}
\usepackage[pdftex]{graphicx}
\usepackage{subcaption}
\graphicspath{ {images/} }
\usepackage[ruled,vlined]{algorithm2e} 
\DontPrintSemicolon 
\usepackage{enumitem}
\usepackage[normalem]{ulem}
\usepackage[toc,page]{appendix}

\usepackage[bookmarks, bookmarksnumbered,
			colorlinks=true,
            linkcolor = blue,
            urlcolor  = blue,
            citecolor = teal]{hyperref}

\newcommand\fnsurl[1]{{\footnotesize\url{#1}}}

\usepackage{myquicksetup}
\numberwithin{definition}{section}
\numberwithin{theorem}{section}
\numberwithin{corollary}{section}
\numberwithin{proposition}{section}
\numberwithin{lemma}{section}
\numberwithin{claim}{section}
\numberwithin{fact}{section}
\numberwithin{remark}{section}
\numberwithin{example}{section}
\numberwithin{equation}{section}
\mathtoolsset{showonlyrefs} % only display equation numbers for equations that are referenced to. https://tex.stackexchange.com/a/4729

\makeatletter
\newcommand\Autoref[1]{\@first@ref#1,@}
\def\@throw@dot#1.#2@{#1}% discard everything after the dot
\def\@set@refname#1{%    % set \@refname to autoefname+s using \getrefbykeydefault
    \edef\@tmp{\getrefbykeydefault{#1}{anchor}{}}%
    \xdef\@tmp{\expandafter\@throw@dot\@tmp.@}%
    \ltx@IfUndefined{\@tmp autorefnameplural}%
         {\def\@refname{\@nameuse{\@tmp autorefname}s}}%
         {\def\@refname{\@nameuse{\@tmp autorefnameplural}}}%
}
\def\@first@ref#1,#2{%
  \ifx#2@\autoref{#1}\let\@nextref\@gobble% only one ref, revert to normal \autoref
  \else%
    \@set@refname{#1}%  set \@refname to autoref name
    \@refname~\ref{#1}% add autoefname and first reference
    \let\@nextref\@next@ref% push processing to \@next@ref
  \fi%
  \@nextref#2%
}
\def\@next@ref#1,#2{%
   \ifx#2@ and~\ref{#1}\let\@nextref\@gobble% at end: print and+\ref and stop
   \else, \ref{#1}% print  ,+\ref and continue
   \fi%
   \@nextref#2%
}
\makeatother

\usepackage{mylivemacros}
\usepackage{todonotes}

\newcommand{\rqedhere}{\tag*{\hspace{-1em}\qedhere}}  % https://groups.google.com/g/comp.text.tex/c/hNqGXm3krSY

\usepackage{bilinearminmax} % project-specific macros

\usepackage{subfiles} % Best loaded last in the preamble

\newcommand{\AdditionalMaterial}{Additional Material}
\newif\ifextended  % extended version: include all appendices
\extendedtrue
\title{An Exponentially Converging Particle Method for the Mixed Nash Equilibrium of Continuous Games}
\date{\today}
\author{Guillaume Wang\thanks{Institute of Mathematics, École polytechnique fédérale de Lausanne (EPFL), Station Z, CH-1015 Lausanne} ,\qquad Lénaïc Chizat\footnotemark[1]}

\hypersetup{
    pdftitle={An Exponentially Converging Particle Method for the Mixed Nash Equilibrium of Continuous Games},
    pdfauthor={Guillaume Wang, Lénaïc Chizat}
}

\begin{document}
\maketitle

% !TEX root = ../main.tex
% \documentclass[../main]{subfiles}
%\begin{document}

\abstract{
	We consider the problem of computing mixed Nash equilibria of two-player zero-sum games with continuous sets of pure strategies and with first-order access to the payoff function.
	This problem arises for example in game-theory-inspired machine learning applications, such as distributionally-robust learning.
	In those applications, the strategy sets are high-dimensional and thus methods based on discretisation cannot tractably return high-accuracy solutions.
	
	In this paper, we introduce and analyze a particle-based method that enjoys guaranteed local convergence for this problem.
    This method consists in parametrizing the mixed strategies as atomic measures and applying proximal point updates to both the atoms' weights and positions. It can be interpreted as an implicit time discretization of the ``interacting'' Wasserstein-Fisher-Rao gradient flow. 

    We prove that, under non-degeneracy assumptions, this method converges at an exponential rate to the exact mixed Nash equilibrium 
    from any initialization
    satisfying a natural notion of closeness to optimality.
    We illustrate our results with numerical experiments and discuss applications to max-margin and distributionally-robust classification using two-layer neural networks, where our method has a natural interpretation as a simultaneous training of the network's weights and of 
    the adversarial distribution.
}

% \end{document}

% !TEX root = ../main.tex
% \documentclass[../main]{subfiles}
%\begin{document}

\section{Introduction} \label{sec:intro}

Consider the min-max, or saddle-point, optimization problem
\begin{equation} \label{eq:intro:SP}
	\min_{\mu \in \PPP(\XXX)} \max_{\nu \in \PPP(\YYY)} 
%	\EE_{x \sim \mu, y \sim \nu}[f(x, y)]
	\int_\XXX \int_\YYY f(x, y) d\mu(x) d\nu(y)
\end{equation}
where $\PPP(\XXX)$ and $\PPP(\YYY)$ are the sets of probability distributions over the sets of \emph{pure strategies} $\XXX$ and $\YYY$, and $f: \XXX \times \YYY \to \RR$ is the \emph{payoff function}.
In the language of game theory, $\PPP(\XXX)$ and $\PPP(\YYY)$ are the sets of \emph{mixed strategies} and
solutions $(\mu^*, \nu^*)$ of \eqref{eq:intro:SP} are the \emph{mixed Nash equilibria} (MNEs) of the two-player zero-sum game $(f, \XXX, \YYY)$.
The conditions for the existence of a MNE are well-known since the 1950s.
In particular, by Glicksberg's theorem~\cite{glicksberg_further_1952} --- later generalized as Sion's minimax theorem \cite{sion_general_1958} --- a MNE always exists if $\XXX$ and $\YYY$ are finite, or if $\XXX$ and $\YYY$ are compact and $f$ is continuous.

Many methods have been proposed to compute MNEs given zeroth-order access to $f$, including in noisy, online or decentralized settings \cite{mertikopoulos_learning_2022, cen_fast_2021}.
Those methods are typically derived and studied for finite games (i.e., with finite strategy sets). 
When $\XXX$ and $\YYY$ are continuous (say, differentiable manifolds) and we additionally have access to the gradients of $f$, it is still possible to reduce the game to a finite one by discretization, but this not only wastes the gradient information, it may also incur a prohibitively high discretization cost when $\XXX$ and $\YYY$ have high dimension.
In this paper, we thus study how to efficiently compute a MNE of $(f, \XXX, \YYY)$ to a high accuracy, when the strategy sets $\XXX, \YYY$ are continuous and given first-order access to $f$.

\paragraph{Particle methods.}
In this setting, a possible strategy is to parametrize the unknowns via
\begin{equation*} 
	\mu = \sum_{i=1}^n a_i \delta_{x_i},
    \qquad
	\nu = \sum_{j=1}^m b_j \delta_{y_j},
\end{equation*}
and to use gradient methods to solve the reparametrized min-max problem
\begin{equation} \label{eq:intro:SP_Fnm}
	\min_{\substack{a \in \Delta_n\\ x \in \XXX^n}} \max_{\substack{b \in \Delta_m\\ y \in \YYY^m}}~ 
    \left\lbrace 
        F_{n,m} \left( (a, x), (b, y) \right)
        ~\coloneqq~ \sum_{i=1}^n \sum_{j=1}^m a_i b_j f(x_i, y_j)
    \right\rbrace,
\end{equation}
where $\Delta_n \coloneqq \left\lbrace a \in \RR^n_+ ; \sum_{i=1}^n a_i=1 \right\rbrace$ is the $n$-simplex.
This approach takes its inspiration from the recent guarantees obtained for ``weighted particle methods'' for convex minimization on the space of measures. In particular, adapting the Conic%
\footnote{The term ``conic'' refers to the particular geometry on the space of couples $(a_i,x_i)$ that leads to multiplicative updates on $a_i$ and additive updates on $x_i$ in \eqref{eq:intro:CP-GDA}.}
Particle Gradient Descent (CP-GD) method~\cite{chizat_sparse_2021} to the constrained min-max context, leads to the \emph{Conic Particle Mirror Descent-Ascent} (CP-MDA) method
which defines a sequence of iterates $(a^k,x^k,b^k,y^k)_{k\geq 0}$ by the update rule
\begin{align}\label{eq:intro:CP-GDA}
	&\begin{cases}
		a^{k+1}_i \propto a^k_i e^{-\eta \frac{\partial}{\partial a_i} F_{n,m}(a^k, x^k, b^k, y^k)} \\
		x^{k+1}_i = x^k_i - \sigma \frac{1}{a^k_i} \frac{\partial}{\partial x_i} F_{n,m}(a^k, x^k, b^k, y^k)
	\end{cases}
	&
	&\begin{cases}
		b^{k+1}_j \propto b^k_j e^{\eta \frac{\partial}{\partial b_j} F_{n,m}(a^k, x^k, b^k, y^k)} \\
		y^{k+1}_j = y^k_j + \sigma \frac{1}{b^k_j} \frac{\partial}{\partial y_j} F_{n,m}(a^k, x^k, b^k, y^k).
	\end{cases}
\end{align}
Here $\eta, \sigma>0$ are step-sizes to be chosen and $a^k$ and $b^k$ are normalized to sum to $1$ at each step.
(The equations above are for the case where $\XXX$, $\YYY$ are Euclidean and without boundaries; in general a retraction step is needed for the update of $x^k$, $y^k$.)
% In the limit where $\eta, \sigma \to 0$, we obtain a continuous-time dynamics which was introduced in~\cite{domingo-enrich_mean-field_2021} under the name ``interacting Wasserstein-Fisher-Rao gradient flow''. There, it is shown in particular that this dynamics admits a mean-field limit, i.e., if the initial parameters are randomly and independently sampled, then there is a well-posed limiting dynamics as $n,m\to \infty$.
In the limit where $\eta, \sigma \to 0$, we obtain a continuous-time dynamics studied by~\cite{domingo-enrich_mean-field_2021} under the name ``interacting Wasserstein-Fisher-Rao gradient flow'', which they show admits a mean-field limiting dynamics when $n, m \to \infty$ and the initial iterates are randomly independently sampled.

\paragraph{Convergence to mixed Nash equilibria.}
The reparametrized saddle-point objective $F_{n,m}$ is finite-dimensional, but is unfortunately not convex-concave in general, and there is no known convergence guarantee for CP-MDA. 
In fact, taking $\sigma=0$, we recover the Mirror Descent-Ascent algorithm on the finite game $(f, \{x^0_i\}_i, \{y^0_j\}_j)$ (a.k.a.\ multiplicative weight updates),
% and it is known that the Bregman divergence of the iterates to saddle points is non-decreasing for this algorithm~\cite{bailey_multiplicative_2018}! 
and it is known that the Bregman divergence of the iterates to the MNE is then non-decreasing~\cite{bailey_multiplicative_2018}! 
For finite games, this non-convergence issue can be resolved by considering instead the implicit version of the same algorithm, or other methods which can be interpreted as tractable approximations
of it \cite{daskalakis_last-iterate_2020, mokhtari_unified_2019}.

In this paper, we propose the implicit version of CP-MDA, which we call the \emph{Conic Particle Proximal Point} algorithm (CP-PP).
\begin{itemize}
    \item We show that, if \eqref{eq:intro:SP} admits a unique and non-degenerate sparse saddle point $(\mu^*,\nu^*)$,
    % admits a unique, sparse and non-degenerate saddle point -> better?
    % (see \autoref{subsec:main_res:pb_stmt}), 
    and if CP-PP is initialized close enough to optimality, then it converges to $(\mu^*,\nu^*)$ at an exponential rate. 
    Note that one can always find such an initialization by sampling sufficiently many particles,
    setting $\sigma=0$ and taking the averaged iterate
    in an initial warm-up phase, see \autoref{subsec:main_res:main_res}.
    The convergence is established both for the Nikaido-Isoda error (the natural~measure of optimality for min-max problems) and for the Wasserstein-Fisher-Rao distance to~$(\mu^*, \nu^*)$.
    \item While CP-PP itself is not directly implementable, we also prove in a simplified setting that a computationally efficient approximation of CP-PP, the \emph{Conic Particle Mirror Prox} algorithm (CP-MP), also converges to $(\mu^*, \nu^*)$ under the same conditions and with the same rate. We observe experimentally that its convergence behavior is the same in the general setting.
    \item We illustrate our work with numerical experiments, including examples of applications to max-$\FFF_1$-margin and distributionally-robust classification with two-layer neural networks.
    
    We observe experimentally that the explicit method CP-MDA does not always converge (although it does in some cases), so that using an implicit time discretization of the dynamics, like CP-PP, is necessary for convergence in general.
\end{itemize}

\subsection{Related work} \label{subsec:intro:related_work}

\paragraph{Infinite-dimensional Mirror Descent-Ascent.}
For finite strategy sets $\XXX$ and $\YYY$, the min-max problem \eqref{eq:intro:SP} is finite-dimensional, and the classical Mirror Descent-Ascent algorithm can be applied to obtain convergence of the averaged iterate to a MNE \cite{bubeck_orf523_2013-1}.
In \cite{hsieh_finding_2018}, the authors use sampling by Langevin dynamics to formulate an implementable version of Mirror Descent-Ascent for continuous strategy sets, which coincides with the bona fide infinite-dimensional Mirror Descent-Ascent algorithm in expectation.
% Accordingly, they give an ergodic convergence rate which holds in expectation, but with no high-probability guarantee.

\paragraph{Computing approximate MNEs via regularization.}
% \cite{domingo-enrich_simultaneous_2022} --> they withdrew the paper
The authors in~\cite{domingo-enrich_mean-field_2021, ma_provably_2022, lu2022two} propose and analyze methods that solve an entropy-regularized variant of \eqref{eq:intro:SP}:
\begin{equation*}
	\min_{\mu \in \PPP(\XXX)} \max_{\nu \in \PPP(\YYY)} 
%	\EE_{x \sim \mu, y \sim \nu}[f(x, y)] 
	\int_\XXX \int_\YYY f(x, y) d\mu(x) d\nu(y)
	+ \beta^{-1} H(\mu) - \beta^{-1} H(\nu)
\end{equation*}
where $H(\mu) = \int_\XXX \log \frac{d\mu}{dx}(x) d\mu(x)$ denotes negative differential entropy
% (and $H(\mu) = +\infty$ if $\mu$ is not absolutely continuous with respect to the reference measure $dx$),
(and $H(\mu) = +\infty$ if $\mu$ is not absolutely continuous with respect to Lebesgue measure),
and $\beta$ is a fixed regularization parameter.
The methods analyzed in these papers correspond to the continuous-time dynamics called ``entropy-regularized interacting Wasserstein gradient flow'' in \cite{domingo-enrich_mean-field_2021}.
Qualitative convergence properties are shown in \cite{domingo-enrich_mean-field_2021}, \cite{ma_provably_2022} proves convergence to an approximate MNE in the quasi-static regime --- i.e., 
when the step-size used to update $\mu$ is infinitely smaller than the one for $\nu$ ---,
and \cite{lu2022two} proves convergence in a regime with finite timescale separation.
The continuous-time guarantees of the aforementioned works can be translated to discrete-time algorithms, thanks to the general framework developed in \cite{karimi2023dynamical} and \cite{karimi2023stochastic}.

% \paragraph{Computing MNEs of some games without discretization and without gradients.}
% %\TODO{I merged two paragraphs into this one.}
% The methods mentioned above, as well as the one analyzed in this work, can be applied to arbitrary differentiable payoff functions. When more structure is known, other methods can be derived. 
% 
% For polynomial games or more generally for separable games, i.e., when
% the payoff function is a finite sum $f(x, y) = \sum_{kl} c_{kl} g_k(x) h_l(y)$ for some known continuous $g_k: \XXX \to \RR, h_l: \YYY \to \RR$ and $c_{kl} \in \RR$,
% \cite{parrilo_polynomial_2006} and \cite{stein_separable_2008} 
% show how to reduce the min-max problem \eqref{eq:intro:SP} to a finite-dimensional maximization problem which is tractable in some cases.
% 
% In \cite{adam_double_2020} and the follow-up work \cite{kroupa_multiple_2021}, the authors propose generalizing the double oracle algorithm --- first introduced by \cite{mcmahan_planning_2003} for finite games --- to compute the MNEs of continuous games. 
% Their work assumes that we have access to an oracle which, given a value of $\nu$, returns a global minimizer of $F\nu: x \mapsto \int_\YYY f(x, y) d\nu(y)$, as well as to a similar oracle for when the roles of $\mu$ and $\nu$ are inverted.

\paragraph{Last-iterate convergence of proximal point methods for min-max optimization.}
In the optimization and learning community, there has been much interest in general 
convex-concave min-max problems $\min_x \max_y G(x, y)$.
Some works focus on ergodic convergence, that is, convergence of the averaged iterate: 
$\left( \frac{1}{T} \sum_{k=1}^T x^k, \frac{1}{T} \sum_{k=1}^T y^k \right)$.
For instance, the Mirror Prox and Proximal Point methods were introduced in \cite{nemirovski_prox-method_2004} to attain $O(1/T)$ ergodic convergence for convex-concave $C^{1,1}$ functions, instead of $O(1/\sqrt{T})$ using Mirror Descent-Ascent \cite{bubeck_orf523_2013-1, bubeck_orf523_2013-2}.

Recent works showed that, while Mirror Descent-Ascent may not in fact converge in the last-iterate sense \cite{mertikopoulos_cycles_2017}, Proximal Point and related methods (e.g.\ Mirror Prox, Optimistic Mirror Descent-Ascent) do exhibit last-iterate convergence \cite{liang_interaction_2019}.
In the special case of finding MNEs of finite two-player zero-sum games, convergence rates for the last iterate have been derived by \cite{daskalakis_last-iterate_2020} and by \cite{wei_linear_2021} for Optimistic Mirror Descent-Ascent, under the assumption that the MNE is unique.
%This last work is closest to ours in spirit.

It should be noted that, when using particle methods for the problem \eqref{eq:intro:SP}, the averaged iterate $\left( \frac{1}{T} \sum_{k=1}^T \mu^k, \frac{1}{T} \sum_{k=1}^T \nu^k \right)$ consists of $(n+m) T$ atoms in general. 
This means that averaging in measure space would result in unacceptably large memory requirements in cases where the domain $\XXX$ or $\YYY$ is large, such as for mixtures of GANs \cite{domingo-enrich_mean-field_2021}.
Another option is to take the average of the $(a^k_i,x^k_i)$ directly, but it would not a priori improve upon the last iterate because $F_{n,m}$ is not convex-concave.

There also exists a growing literature on nonconvex-nonconcave min-max optimization, which focuses on the problem of finding local saddle points or even just stationary points of the gradient descent-ascent flow \cite{gidel2018variational, diakonikolas_efficient_2021}. Our result are stronger than what one could expect to achieve with techniques from that literature: We are able to find the solution of \eqref{eq:intro:SP}, which gives a global Nash equilibrium of $F_{n,m}$, instead of simply stationary points.
% \TODO{it would be interesting to determine whether there exist spurious stationary points of $F_{n,m}$...}
% -> possible argument: take n<n*, m<m*, then the dynamics converges to some non-NE stationary point (*). Then use n>n*, m>m* but initialized degenerately: dynamics converges to that same point.
% -> (*) however it's unclear that for n<n*, m<m*, the dynamics even converges. While unclear how one could show this in general, experimentally yes that seems to happen sometimes, which is sufficient to make our point

\vspace{1em}
The remainder of this paper is structured as follows.
In \autoref{sec:main_res}, we state the problem, describe the algorithm, and present the main result.
In \autoref{sec:cv_proof}, we prove the main result.
In \autoref{sec:experiments}, we discuss examples of applications and present numerical experiments.
In \autoref{sec:ccl} we conclude and proof details are deferred to the appendix.

% \end{document}

% !TEX root = ../main.tex
%\documentclass[../main]{subfiles}
%\begin{document}

\section{Main result} \label{sec:main_res}

\subsection{Problem setting: Computing MNEs of continuous games} \label{subsec:main_res:pb_stmt}

\paragraph{Preliminaries.}

Let us first recall the general convex-concave min-max optimization framework.
For convex sets $M, N$ and a convex-concave function $F: M \times N \to \RR$, a \emph{saddle point} or \emph{solution of the min-max problem}
\begin{equation*}
	\min_{\mu \in M} \max_{\nu \in N} F(\mu, \nu)
\end{equation*}
is any pair $(\mu^*, \nu^*)$ such that%
\footnote{When $F$ is not convex-concave, several notions of min-max solutions can be considered \cite[Sec.~2.1]{diakonikolas_efficient_2021}, but we will never need such considerations in this article.
To fix ideas, we will use the strongest such notion --- that of global saddle point --- and call \emph{solution of a nonconvex-nonconcave problem $\min_x \max_y G(x,y)$} any $(x^*, y^*)$ satisfying $G(x^*, y) \leq G(x, y^*)$ for all $x, y$;
but we emphasize that this choice has no impact on any of our discussion.}
\begin{equation*}
	\forall \mu \in M, \forall \nu \in N,~
	F(\mu^*, \nu) \leq F(\mu^*, \nu^*) \leq F(\mu, \nu^*).
\end{equation*}
The existence of saddle points is guaranteed for example by Sion's minimax theorem when $M$ and $N$ are compact and $F$ is continuous.
The goodness of a pair $(\hmu, \hnu)$ can be quantified by its \emph{Nikaido-Isoda error} (NI error), a.k.a.~duality gap, 
% \TODO{LC: I prefer the term duality gap (which is more informative). } -> GW: Nikaido-Isoda error is better-established in convex min-max literature. but I'm open to further discussion
defined as
\begin{equation*}
	\NI(\hmu, \hnu) = \max_{\mu, \nu} F(\hmu, \nu) - F(\mu, \hnu).
\end{equation*}
Indeed, it is easily seen that $\NI(\hmu, \hnu) \geq 0$ with equality if and only if $(\hmu, \hnu)$ is a saddle point.

As an example, the problem of finding the MNE of a two-player zero-sum game with finite strategy sets $\XXX = [n] \coloneqq \{1,...,n\}$ and $\YYY = [m]$ and payoff function $f(i, j)$ can be written as
\begin{equation*}
	\min_{a \in \Delta_n} \max_{b \in \Delta_m} ~
    \left\lbrace
    	F(a, b) = \sum_{i=1}^n \sum_{j=1}^m a_i b_j f(i, j) = a^\top \gmat b
    \right\rbrace
	% F(a, b) 
	% ~~~~ \text{where} ~~~~
	% F(a, b) = \sum_{i=1}^n \sum_{j=1}^m a_i b_j f(i, j) = a^\top \gmat b
\end{equation*}
where $\gmat_{ij} = f(i, j)$ and $\Delta_n = \left\{ a \in \RR_+^n; \sum_i a_i = 1 \right\} \simeq \PPP([n])$.
Since $\Delta_n$ and $\Delta_m$ are convex compact and $F(a, b) = a^\top \gmat b$ is convex-concave and continuous, Sion's minimax theorem applies so saddle points exist. 

\paragraph{Problem setting and assumptions.}

The min-max problem we are concerned with in this paper is that of finding a MNE of the continuous game $(f, \XXX, \YYY)$, as defined in \eqref{eq:intro:SP}:
%$\min_{\mu \in \PPP(\XXX)} \max_{\nu \in \PPP(\YYY)} 
%%\EE_{x \sim \mu, y \sim \nu}[f(x, y)]
%\int_\XXX \int_\YYY f(x, y) d\mu(x) d\nu(y)$,
\begin{equation*}
	\min_{\mu \in \PPP(\XXX)} \max_{\nu \in \PPP(\YYY)} ~
    \left\lbrace
        F(\mu, \nu) \coloneqq \int_\XXX \int_\YYY f(x, y) d\mu(x) d\nu(y)
    \right\rbrace,
	% F(\mu, \nu)
	% %= \EE_{x \sim \mu, y \sim \nu}[f(x, y)]
	% ~~~~ \text{where} ~~~~
	% F(\mu, \nu) = \int_\XXX \int_\YYY f(x, y) d\mu(x) d\nu(y),
\end{equation*}
with the following assumptions. 
% make references to these labels as: \hyperref[assum:1]{Assumption~\ref*{assum:1}}
\begin{assumptions*}
	~
	\begin{enumerate}
        \item \label{assum:1}
        The strategy sets are the $d_x$- resp.\ $d_y$-dimensional tori $\XXX = \TT^{d_x}$, $\YYY = \TT^{d_y}$. 
		\item \label{assum:2}
        The payoff function $f: \XXX \times \YYY \to \RR$ is $C^{2,1}$, i.e., it has Lipschitz-continuous second-order differentials.
		\item \label{assum:3}
        The MNE $(\mu^*, \nu^*)$ of \eqref{eq:intro:SP} is unique.
		\item \label{assum:4}
        The MNE $(\mu^*, \nu^*)$ of \eqref{eq:intro:SP} is sparse, that is,
% 		$\support(\mu^*) = \{ x^*_I, I \in [n^*] \}$ and $\support(\nu^*) = \{ y^*_J, J \in [m^*] \}$ 
		\begin{equation*}
			\support(\mu^*) = \left\{ x^*_I, I \in [n^*] \right\}
			~~~\text{and}~~~
			\support(\nu^*) = \left\{ y^*_J, J \in [m^*] \right\}
		\end{equation*}
		for some $n^*, m^* < \infty$.
	\end{enumerate}
\end{assumptions*}
Note that points~1 and~2 imply the existence of a MNE by Glicksberg's theorem.
Also note that point~4 is guaranteed to hold if the game is \emph{separable}, that is, if $f$ can be written as a finite sum of the form
$
	f(x, y) = \sum_{kl} c_{kl} g_k(x) h_l(y)
$
for some $c_{kl} \in \RR$ and $g_k: \XXX \to \RR, h_l: \YYY \to \RR$ continuous \cite[Corollary~2.10]{stein_separable_2008}.
% (see \autoref{subsec:intro:related_work}).
% \TODO{determine whether my analyticity argument is true and/or interesting, and if so add it in an appendix} --> I kind of gave up on proving it; I'll still write it up in a commented out appendix for personal reference
Moreover, point~1 could be replaced by assuming only that $\XXX$ and $\YYY$ are compact Riemannian manifolds without boundaries; our analysis could be generalized to this setting (following \cite{chizat_sparse_2021}) at the expense of more technical notation.
% \TODO{think about whether compactness is needed for our proof.}
% LC: I guess it's not, provided we assume a suitable uniform lower/upper bound on Fmu or nuF outside of a large set containing mu^*/nu^*
Point 3 is crucial to our
analysis, but also to all known last-iterate convergence analyses for MNEs in finite dimension \cite{wei_linear_2021, daskalakis_last-iterate_2020}.

Before stating the rest of our assumptions, let us remark a useful fact about the structure of the problem.
By definition the MNE $(\mu^*, \nu^*)$ is characterized by
\begin{equation*}
	\forall \mu \in \PPP(\XXX), \forall \nu \in \PPP(\YYY),~
	F(\mu^*, \nu) \leq \gval \leq F(\mu, \nu^*)
	~~~ \text{where}~ \gval \coloneqq F(\mu^*, \nu^*),
\end{equation*}
i.e.,
\begin{equation*}
	\forall \mu \in \PPP(\XXX), \int_\XXX \underbrace{ \left( \int_\YYY f(x, y) d\nu^*(y) \right) }_{\eqqcolon~ (F \nu^*)(x)} d\mu(x)
	\geq \gval
	~~~\text{and}~~~
	\forall \nu \in \PPP(\YYY), \int_\YYY \underbrace{ \left( \int_\XXX f(x, y) d\mu^*(x) \right) }_{\eqqcolon~ ((\mu^*)^\top F)(y)} d\nu(y) 
	\leq \gval.
\end{equation*}
The function $F\nu \in \CCC(\XXX)$ defined by this equation is 
% the Fréchet differential
the first variation 
of $F$ with respect to $\mu$ at any $(\mu, \nu)$ \cite{santambrogio__2016}; note that it is independent of $\mu$ thanks to bilinearity of $F$.
Since $\min_{\mu \in \PPP(\XXX)} \int g~ d\mu = \min_\XXX g$ for any $g \in \CCC(\XXX)$, 
the above inequalities are
equivalent to
% \footnote{
% 	Another way of stating the same fact is to observe that for all $(\mu', \nu')$,
% 	\begin{align*}
% 		\NI(\mu', \nu') = \max_{\mu \in \PPP(\XXX), \nu \in \PPP(\YYY)} F(\mu', \nu) - F(\mu, \nu')
% 		= \left[
% 		\max_{y \in \YYY} ((\mu')^\top F)(y) - \gval 
% 		\right] - \left[ 
% 		\min_{x \in \XXX} (F \nu')(x) - \gval
% 		\right]. 
% 	\end{align*}
% }
\begin{equation} \label{eq:main_res:charact_MNE_firstvars_ineqs}
	\forall x \in \XXX,~
	(F \nu^*)(x) \geq \gval
	~~~\text{and}~~~
	\forall y \in \YYY,~
	((\mu^*)^\top F)(y) \leq \gval.
\end{equation}
As a partial converse, we also have
\begin{equation} \label{eq:main_res:charact_MNE_firstvars_ineqs_partialconverse}
	\forall I \in [n^*],~
	(F \nu^*)(x^*_I) = \gval
	~~~\text{and}~~~
	\forall J \in [m^*],~
	((\mu^*)^\top F)(y^*_J) = \gval
\end{equation}
since if $(F \nu^*)(x^*_I) > \gval$ for some $I$ then we would have 
$F(\mu^*, \nu^*) = \int_\XXX (F \nu^*)(x) d\mu^*(x) > \gval$.

Our second set of assumptions requires the inequalities \eqref{eq:main_res:charact_MNE_firstvars_ineqs} to be strict wherever possible, and even ``strong'' locally.
\begin{samepage}
\begin{assumptions*}[Non-degeneracy]
	~
	\begin{enumerate}[start=5]
		\item \label{assum:5}
        The first variations at optimum, $F \nu^* \in \CCC(\XXX)$ resp.\ $(\mu^*)^\top F \in \CCC(\YYY)$, are equal to $\gval \coloneqq F(\mu^*, \nu^*)$ only at the $\{ x^*_I, I \in [n^*] \}$ resp.\ $\{ y^*_J, J \in [m^*] \}$.
		\item \label{assum:6} 
        The local kernels are non-degenerate, that is,
		\begin{equation*}
			\forall I \in [n^*],~ 
			% \partial_{xx} (F\nu^*)(x^*_I) \succ 0
            \nabla^2 (F\nu^*)(x^*_I) \succ 0
			~~~\text{and}~~~
			\forall J \in [m^*],~ 
			% \partial_{yy} ((\mu^*)^\top F)(y^*_J) \prec 0.
			\nabla^2 ((\mu^*)^\top F)(y^*_J) \prec 0.
		\end{equation*}
	\end{enumerate}
\end{assumptions*}
\end{samepage}

These non-degeneracy assumptions are analogous to the ones made in \cite{chizat_sparse_2021} for minimization in the space of measures. As for minimization, they generally cannot be checked a priori; but they can be checked a posteriori after computing 
% (an approximation of) 
$(\mu^*, \nu^*)$ by computing the Hessians of $F \nu^*$ on $\support(\mu^*)$ and of $(\mu^*)^\top F$ on $\support(\nu^*)$, or simply visually in the one-dimensional case as in \autoref{subsec:experiments:RFF} (\autoref{fig:random_fourier_1D:firstvars_smoothedmeas_lastiter}).
Even though we do not have any rigorous result in this direction, informally we expect the non-degeneracy assumptions to be generic in some sense. For example they turned out to be satisfied in all of our experiments with random payoff functions of the form of \autoref{subsec:experiments:RFF}, with any dimension $d_x, d_y$.
\autoref{example:experiments:synthetic} provides a case where all the assumptions can be checked analytically, including the uniqueness of the MNE.

Our work provides the first convergence guarantee for the sparse continuous-game MNE setting (Assumptions~1,~2,~4). On the other hand, Assumptions~3,~5,~6 are admittedly difficult to check a priori for any given $(f, \XXX, \YYY)$. But they are the minimal ones for our convergence analysis to go through, and we expect that they are very difficult to relax.

% \subsection{The Conic Particle Proximal Point and Conic Particle Mirror Prox algorithms} \label{subsec:main_res:algo}
\subsection{The Conic Particle Proximal Point algorithm} \label{subsec:main_res:algo}

In order to solve the saddle-point problem \eqref{eq:intro:SP}, we reparametrize the problem via $\mu = \sum_{i=1}^n a_i \delta_{x_i}$ and $\nu = \sum_{j=1}^m b_j \delta_{y_j}$
and we use a particle gradient algorithm to tackle the reparametrized problem~\eqref{eq:intro:SP_Fnm}:
$\min_{a,x} \max_{b,y} F_{n,m}((a,x),(b,y))$.
Specifically, we analyze the Conic Particle Proximal Point (CP-PP) algorithm given by the update rule
\begin{multline}
	\!\!\!\!\!
    ((a^{k+1}\!, x^{k+1}), (b^{k+1}\!, y^{k+1})) = 
	\argmin_{\substack{a \in \Delta_n\\ x \in \XXX^n}} \argmax_{\substack{b \in \Delta_m\\ y \in \YYY^m}} \,
	F_{n,m}((a, x), (b, y)) 
	+ \frac{1}{\eta} \KLdiv(a, a^k) + \frac{1}{2 \sigma} \sum_{i=1}^n a^k_i \norm{x_i-x^k_i}^2 \\
	- \frac{1}{\eta} \KLdiv(b, b^k) - \frac{1}{2 \sigma} \sum_{j=1}^m b^k_j \norm{y_j-y^k_j}^2
\label{eq:main_res:ppa_upd}
\end{multline}
where $\eta, \sigma>0$ are constant step-sizes to be chosen and $\KLdiv(w, \hw) = \sum_i w_i \log \frac{w_i}{\hw_i}$ denotes Kullback-Leibler divergence a.k.a.\ relative entropy.
% The number of particles used, $n$ and $m$, are also parameters of the algorithm, and should be chosen at least equal to $n^*$ resp.\ $m^*$ so that convergence to the solution $(\mu^*, \nu^*) = \left( \sum_{I \in [n^*]} a^*_I \delta_{x^*_I}, \sum_{J \in [m^*]} b^*_J \delta_{x^*_J} \right)$ is possible.
% -> captured by the NI << eps initial condt
Interestingly, the function $D\left((a, x), (\ha, \hx) \right) = \KLdiv(a, \ha) + \frac{\eta}{2 \sigma} \sum_i \ha_i \norm{x_i-\hx_i}^2$ is technically not a Bregman divergence (it does not satisfy the last point of \autoref{lm:aux_lemmas:bregman_divgce}), due to the ``cross-terms'' in the second term.

The following lemma, proved in \autoref{sec:proof_welldef}, justifies that the CP-PP update is well-defined.
\begin{lemma} \label{lm:main_res:ppm_upd_well_defined}
    Under Assumptions~\ref{assum:1}-\ref{assum:2},
	there exist $\eta_0, \sigma_0>0$
	(dependent only on $(f, \XXX, \YYY)$)
	such that if $\eta \leq \eta_0$ and $\sigma \leq \sigma_0$, then
	the objective function in \eqref{eq:main_res:ppa_upd} is convex-concave over a ball 
%	in $(\Delta_n \times \XXX^n) \times (\Delta_m \times \YYY^m)$ 
	centered at $((a^k, x^k), (b^k, y^k))$,
	and it has a saddle point in the interior of that ball.
\end{lemma}

\paragraph{Implementable variant: Conic Particle Mirror-Prox (CP-MP).}
% Note that the above update rule is not directly implementable, as it requires solving a min-max optimization problem at each iteration.
Note that every CP-PP update requires solving a min-max optimization problem \eqref{eq:main_res:ppa_upd} exactly.
In practice, in the spirit of \cite{nemirovski_prox-method_2004}, one may obtain an approximate solution to \eqref{eq:main_res:ppa_upd} by running an inner loop where $F_{n,m}((a, x), (b, y))$ is replaced by its first-order approximation at the current point, starting from $(\ta^0, \tx^0, \tb^0, \ty^0) = (a^k, x^k, b^k, y^k)$:
\begin{multline*}
	(\ta^{l+1}, \tx^{l+1}, \tb^{l+1}, \ty^{l+1})
	= \\
	\argmin_{\substack{a \in \Delta_n\\ x \in \XXX^n}} \argmax_{\substack{b \in \Delta_m\\ y \in \YYY^m}}
	\innerprod{\nabla F_{n,m} \left( \ta^l, \tx^l, \tb^l, \ty^l \right) }{(a, x, b, y)}
	+ \frac{1}{\eta} \KLdiv(a, a^k) + \frac{1}{2 \sigma} \sum_{i=1}^n a^k_i \norm{x_i-x^k_i}^2 \\
	- \frac{1}{\eta} \KLdiv(b, b^k) - \frac{1}{2 \sigma} \sum_{j=1}^m b^k_j \norm{y_j-y^k_j}^2 \qquad
\end{multline*}
and letting $(a^{k+1}, x^{k+1}, b^{k+1}, y^{k+1}) = (\ta^L, \tx^L, \tb^L, \ty^L)$, where $L \geq 1$ is the number of times we run the inner loop at each $k$.
Each iteration of the inner loop decomposes into four independent mirror descent updates.
Pseudocode for this implementable variant of CP-PP is given in \autoref{alg:main_res:ppm},
where to lighten the notation we use the shorthand $z^k = (a^k, x^k, b^k, y^k)$.

One can check that $L=1$ corresponds to the CP-MDA algorithm described in the introduction, and that for $L=\infty$ we recover CP-PP.
When $L=2$, we refer to this method as Conic Particle Mirror Prox (CP-MP).
Similarly as in \cite{nemirovski_prox-method_2004}, one can expect that $L=2$ actually suffices to obtain the same convergence behavior as $L=\infty$; this is confirmed in numerical experiments, and proved in a simplified setting (\autoref{prop:cv_proof:exactparam:mp_cv}).
This behavior can be explained by the fact that Proximal Point and Mirror Prox updates coincide up to order-3 terms in the step-size
--- see \Autoref{lm:mp_exactparam:mp_expan,lm:mp_exactparam:pp_expan}.

% \vspace{-1em}
\begin{algorithm}[t]
\caption{Conic Particle Proximal Point, implementable variant} \label{alg:main_res:ppm}
	\KwIn{$f: \XXX \times \YYY \to \RR$, ~
	$n, m \in \NN^*$, ~
	$\eta, \sigma>0$, ~
	$T, L \in \NN^*$}
	Initialize $z^0 = (a^0, x^0, b^0, y^0) \in \Delta_n \times \XXX^n \times \Delta_m \times \YYY^m$ \;
	\For{$k=0, ..., T-1$}{
		$\tz^0 \gets z^k$ \;
		\For{$l=0, ..., L-1$}{
			$\forall i,~ \ta^{l+1}_i \gets \ta^l_i e^{-\eta \frac{\partial}{\partial a_i} F_{n,m}(\tz^l)} / Z$ ~~~ with $Z$ such that $\ta^{l+1} \in \Delta_n$ \;
			$\forall i,~ \tx^{l+1}_i \gets \tx^l_i -  \sigma \frac{1}{a^k_i} \frac{\partial}{\partial x_i} F_{n,m}(\tz^l)$ \;
			similarly for $\tb^{l+1}$ and $\ty^{l+1}$
		}
		$z^{k+1} \gets \tz^L$
	}
	\Return{$\mu^T = \sum_{i=1}^n a^T_i \delta_{x^T_i}$, $\nu^T = \sum_{j=1}^m b^T_i \delta_{y^T_j}$}
\end{algorithm}

\paragraph{Relation to Wasserstein-Fisher-Rao (WFR) geometry.}
% \TODO{make this into a Remark? it's not really crucial (except that it justifies the choice of Lyapunov function, but that connection is already mentioned directly in the "Proof" section)}
% --> actually i think it's ok
The Wasserstein-Fisher-Rao distance, a.k.a.~Hellinger-Kantorovich distance \cite{liero_optimal_2018,kondratyev2016new, chizat_unbalanced_2019} is a distance on the set of non-negative measures which
metrizes narrow convergence and
combines features of the Fisher-Rao and of the Wasserstein distances.
This last fact is perhaps easiest to see in its dynamical formulation 
%\cite[Def.~2.3]{gallouet_jko_2018}
\cite[Thm.~8.18]{liero_optimal_2018}: 
\begin{equation} \label{eq:main_res:def_WFR22}
    \WFR_2^2(\mu_0, \mu_1) 
    = \inf_{(\mu_t, v_t, r_t) \in \AAA(\mu_0, \mu_1)} 
    \int_0^1 \int_\XXX \left( \frac{\eta}{2 \sigma} \norm{v_t(x)}^2 + \frac{1}{2} r_t(x)^2 \right) d\mu_t(x) dt
\end{equation}
where $\AAA(\mu_0, \mu_1)$ is the set of triples $(\mu_t, v_t, r_t)_{0 \leq t \leq 1}$ such that $(\mu_t)_{t\in [0,1]}$ is a weakly continuous curve in $\MMM_+(\XXX)$ the set of non-negative measures with endpoints $\mu_0$ and $\mu_1$, $v_t \in L^2_{\mu_t}(\XXX)^{d_x}$, $r_t \in L^2_{\mu_t}(\XXX)$, and satisfying the continuity equation with source
% \begin{equation}
$
    \partial_t \mu_t + \nabla \cdot (\mu_t v_t) = \mu_t r_t
$
% \end{equation}
in the sense of distributions.
(See \cite{liero_optimal_2018, chizat_unbalanced_2019} for equivalent static formulations).
Here the scalars $\eta$ and $\sigma$ trade off the Fisher-Rao (``local growth and destruction of mass'') and Wasserstein (``movement of mass'') components of the distance.

The CP-MDA, CP-PP and CP-MP algorithms are time-discretizations of the interacting WFR gradient flow \cite{domingo-enrich_mean-field_2021}, which is the finite-particle version of the WFR gradient flow in measure space \cite[Prop.~2.1]{chizat_sparse_2021}.
Because of this connection, the CP-PP iterates' WFR distance to the MNE (or rather a simpler proxy of it that is sufficient for our purpose) will play a central role in our convergence analysis.

\subsection{Main convergence result} \label{subsec:main_res:main_res}

Our main result is that the CP-PP algorithm \eqref{eq:main_res:ppa_upd} converges locally at an exponential rate.
Its proof follows from \autoref{prop:cv_proof:gencase:NI_equiv_V}, \autoref{prop:cv_proof:gencase:WFR_equiv_V} and \autoref{thm:cv_proof:gencase:loc_exp_cv}, as shown in \autoref{sec:proof_mainres}.
\begin{theorem} \label{thm:main_res:loc_exp_cv_NI}
    Fix any $\Gamma_0 \geq 1$. 
    % Under the assumptions of \autoref{subsec:main_res:pb_stmt},
    Under Assumptions~\ref{assum:1}-\ref{assum:6},
    there exist $\eta_0, \sigma_0>0$ such that for all $\eta \leq \eta_0, \sigma \leq \sigma_0$ with $\Gamma_0^{-1} \leq \frac{\sigma}{\eta} \leq \Gamma_0$,
	there exists $C, C', r_0, \kappa>0$ such that if $\NI(\mu^0, \nu^0) \leq r_0$, then the CP-PP iterates $(\mu^k, \nu^k) = \left( \sum_{i=1}^n a^k_i \delta_{x^k_i}, \sum_{j=1}^m b^k_j \delta_{y^k_j} \right)$ satisfy for all $k\in \mathbb{N}$
	\begin{align*}
	    \NI(\mu^k, \nu^k) &\leq C (1-\kappa)^k \\
	    % \intertext{and}
        \text{and}~~~~
		\WFR_2^2(\mu^k, \mu^*) + \WFR_2^2(\nu^k, \nu^*) &\leq C' (1-\kappa)^k.
	\end{align*}
	In particular, since $\WFR_2$ metrizes narrow convergence, $\mu^k, \nu^k$ converge narrowly to $\mu^*, \nu^*$, that is 
    $\forall \phi \in \mathcal{C}(\mathcal{X})$, $\int \phi d\mu^k\to \int \phi d\mu^*$ and $\forall \psi \in \mathcal{C}(\mathcal{Y})$, $\int \psi d\nu^k\to \int \psi d\nu^*$.
\end{theorem}

% \paragraph{Dependency on the problem parameters.}
The dependency of $C, C', r_0, \kappa$ on the problem data $(f, \XXX, \YYY)$ appears quite subtle, unfortunately. It can be traced back to \autoref{lm:growth_conds:error_bound} establishing an ``error bound'' type inequality which relies on uniqueness of the MNE, making it difficult to quantify. The known analyses for finite-game MNEs \cite{wei_linear_2021, daskalakis_last-iterate_2020} face the same drawback.

\paragraph{Dependency of the constants on the step-sizes $\eta, \sigma$.}
The quantities $C, C', r_0, \kappa$ appearing in the theorem depend on the step-sizes.
Our proof technique requires to take them at most of order
$r_0 \lesssim \eta^{17/4}$,
$\kappa \lesssim \eta^2$,
and
$C \gtrsim \eta^{-1/5} r_0^{2/5}$,
$C' \gtrsim \eta^{-2/5} r_0^{4/5}$
    % For the exact-param case, or more generally if we didn't have the additional factor of \sqrt{sigma} in \autoref{prop:cv_proof:gencase:NI_equiv_V}, it would be: 
    % $r_0 \lesssim \eta^{15/4}$,
    % $\kappa \lesssim \eta^2$,
    % and
    % $C \gtrsim r_0^{2/5}$,
    % $C' \gtrsim r_0^{4/5}$
--- supposing $\eta \asymp \sigma$ ---,
as one can check from 
the proof in \autoref{sec:proof_mainres} and
the statements of 
% Props.~\ref{prop:cv_proof:gencase:NI_equiv_V} and \ref{prop:cv_proof:gencase:WFR_equiv_V}
\autoref{prop:cv_proof:gencase:NI_equiv_V}, \autoref{prop:cv_proof:gencase:WFR_equiv_V} 
and \autoref{thm:cv_proof:gencase:loc_exp_cv}.
So our approach only applies for discrete-time algorithms, and would not allow to show convergence of the continuous-time flow associated to CP-PP.

Indeed if we formally let $\eta, \sigma \to 0$, the localness requirement on the initial iterate $\NI(\mu^0, \nu^0) \leq r_0 \xrightarrow{\eta, \sigma \to 0} 0$ reduces to requiring $(\mu^0, \nu^0)$ to already be the MNE.
Even ignoring the localness requirement, the bound $\left( 1- \Theta(\eta^2) \right)^k$ becomes constant when $k = \floor{\frac{t}{\eta}}$ for a fixed $t$ and $\eta \asymp \sigma \to 0$; so the bound does not vanish as $t$ increases.

Interestingly, experimentally we do observe convergence of the continuous-time flow in some (but not all) cases; see \autoref{subsec:experiments:RFF}. 
It is worth mentioning that for $\XXX$ and $\YYY$ finite, the Fisher-Rao gradient flow does not converge \cite{mertikopoulos_cycles_2017}, so the fact that we sometimes observe convergence of the flow may be specific to conic particle gradient methods.

\paragraph{Necessity of taking $\eta, \sigma>0$ (comparison with pure Fisher-Rao and pure Wasserstein particle methods).}~
\begin{itemize}
    \item If $\sigma=0$,
    the position variables ($x^k, y^k$) of CP-PP stay constant throughout the algorithm, and only the weights ($a^k, b^k$) vary. This corresponds to the Fisher-Rao gradient dynamics.
    \begin{itemize}
        \item If the initial measure variables $(\mu^0, \nu^0)$ were supported on a large number of points $\left\lbrace \hx_1,...,\hx_n \right\rbrace$ resp.\ $\left\lbrace \hy_1,...,\hy_m \right\rbrace$ covering $\XXX$ resp.\ $\YYY$ uniformly,
        one may expect CP-PP to converge locally exponentially to $(\mu^\infty, \nu^\infty)$ a MNE of the finite game $\left( f, \{\hx_i\}_i, \{\hy_j\}_j \right)$.%
        \footnote
            {In fact there is no guarantee so far that the last iterate of CP-PP with $\sigma=0$ will converge locally to a MNE of the finite game $\left( f, \{\hx_i\}_i, \{\hy_j\}_j \right)$ (although the averaged iterate is guaranteed to converge globally to one by \cite{nemirovski_prox-method_2004}).
            Indeed, 
            while \cite[App.~D]{mertikopoulos_optimistic_2019} shows qualitative convergence without a rate,
            known quantitative last-iterate convergence guarantees for finite games \cite{wei_linear_2021, daskalakis_last-iterate_2020} require the MNE to be unique, which may not be the case a priori for $\left( f, \{\hx_i\}_i, \{\hy_j\}_j \right)$.
        }
        By continuity of $f$, the ``price'' of the discretization with respect to the original game can then be bounded by $\NI(\mu^\infty, \nu^\infty) = O\left( n^{-1/d_x} + m^{-1/d_y} \right)$.
        That is, in order to achieve a NI error of $O(\eps)$, it is sufficient to let $n = (1/\eps)^{d_x}, m = (1/\eps)^{d_y}$, and to run CP-PP to convergence with $\sigma=0$;
        however the computational complexity of each iteration would then be $\Theta(nm)$, which can be prohibitively costly if $d_x, d_y$ are large.
        \item Also note that if $(\mu^0, \nu^0)$ are supported on the entire space $\XXX$ resp.\ $\YYY$, we can only expect convergence of CP-PP to the exact MNE at a rate $\asymp 1/k$ in the worst case, and not at an exponential rate \cite{chizat_convergence_2021}
        (a worst-case example can be constructed
        % similarly as in
        by using a similar idea as in Prop.~5.5, Setting~II of that paper).
        % This phenomenon is discussed in \cite{chizat_convergence_2021} in the minimization case, and a counter-example is given by $f(x, y) = \Phi(x) - \Phi(y)$ where $\Phi$ is any smooth function such that $\Phi(x) = \norm{x}^2$ when $\norm{x} \leq \frac{1}{2}$ and $\Phi(x) \geq \frac{1}{4}$ otherwise.
        % Indeed, $\mu^k \propto \mu^0 \odot \exp(-k \eta \frac{1}{2} \norm{x}_2^2)$
        % and $\nu^k \propto \nu^0 \odot \exp(-k \eta \frac{1}{2} \norm{y}_2^2)$,
        % which are essentially Gaussian distributions of variance $\frac{1}{\eta k}$,
        % and by straightforward computations $\NI(\mu^k, \nu^k) \approx \frac{2}{\eta k}$,
        % if $\mu^0, \nu^0$ are ``flat'' enough.
    \end{itemize}
    \item If $\eta=0$,
    the weight variables of CP-PP stay fixed, and only the positions vary. This corresponds to the interacting Wasserstein gradient dynamics \cite{domingo-enrich_mean-field_2021}. 
    \begin{itemize}
        \item This dynamics has a degraded behavior because it has fewer degrees of freedom. For instance, for $(\mu^0, \nu^0)$ supported on finitely many point and $a^0_i = \frac{1}{n}, b^0_j = \frac{1}{m}$, the CP-PP iterates cannot converge to $(\mu^*, \nu^*)$ exactly --- unless the solution weights $a^*_I, b^*_J$ happen to all lie in $\frac{1}{n} \ZZ$ resp.\ $\frac{1}{m} \ZZ$. 
        \item Even allowing for $(\mu^0, \nu^0)$ supported on the entire space, we are not aware of any convergence guarantee to $(\mu^*, \nu^*)$ for this dynamics, in continuous time or otherwise.
        % \item In \autoref{prop:cv_proof:gencase:NI_equiv_V} we show a two-sided bound relating NI error with (a proxy of) the WFR distance to the MNE, indicating that weighted particle methods --- and in particular CP-PP --- may have a definite edge over pure Wasserstein particle methods.    \TODO{is this nonsense?} -> maybe
    \end{itemize}
\end{itemize}

\paragraph{Agnosticity to the numbers of particles $n, m$.}
The numbers of particles $n, m$ used in the CP-PP algorithm do not appear in the theorem, nor are they hidden in the constants.
In particular the convergence rate does not deteriorate with large $n,m$ while the per-iteration cost is linear in $n+m$.
Even the condition that $n \geq n^*$ and $m \geq m^*$ (the sparsities of $(\mu^*, \nu^*)$) does not appear explicitly, but it is implied by the localness condition $\NI(\mu^0, \nu^0) \leq r_0$.
That is, $r_0$ is defined such that, if $n < n^*$ or $m < m^*$, then there simply do not exist atomic measures $\mu^0, \nu^0$ with $n$ resp.\ $m$ atoms that achieve NI error less than $r_0$.

The fact that our result is agnostic to such overparametrization should be viewed as a strength.
Indeed the convergence guarantee does not deteriorate with large $n,m$, and on the other hand allowing ourselves to use arbitrary $(n,m) \neq (n^*, m^*)$ enables simpler warm-up procedures.
In terms of the application to classification with two-layer neural networks presented in \autoref{sec:experiments}, this agnosticity means that our results apply regardless of the number of hidden neurons, as long as it is not too small.

% \paragraph{Towards a globally convergent algorithm.}
\paragraph{A possible two-phase procedure.}
While our proposed algorithm is only shown to be locally convergent,
we would like to stress that it is the only known one that can provably converge to the actual solution $(\mu^*, \nu^*)$.
This is in contrast to any algorithm relying on discretization of the strategy spaces $\XXX, \YYY$, since these algorithms can only ever output measures $(\mu^k, \nu^k)$ whose support does not even match the optimal one (unless one is extremely lucky when choosing the discretization).
So one way to take advantage of our proposed algorithm and performance guarantee, is to use it as a second ``high-accuracy'' phase, preceded by a warm-up phase with global convergence guarantees.

A simple such warm-up procedure is 
to fix $\eps_{\text{warm-up}} \leq r_0$,
to discretize $\XXX$ and $\YYY$ by 
% $(x_i)_{i \leq n}$, $(y_j)_{j \leq m}$ with 
$n = (1/\eps_{\text{warm-up}})^{d_x}$, $m = (1/\eps_{\text{warm-up}})^{d_y}$
grid-points $\{\hx_i\}_i, \{\hy_j\}_j$,
to run Fisher-Rao Proximal Point (i.e, CP-PP with $\sigma=0$) for $T_{\text{warm-up}} = 1/\eps_{\text{warm-up}}$ iterations and take the average of the iterates.
Indeed this ensures a NI error of $O(\eps_{\text{warm-up}})$ for the discretized game \cite{nemirovski_prox-method_2004}, and the NI error for the discretized game is $O(\eps_{\text{warm-up}})$-close to the NI error for the original game by the choice of $n, m$.
The per-iteration complexity of the first phase is $\Theta(nm)$ and that of the second phase is $\Theta(n+m) = \Theta((1/\eps_{\text{warm-up}})^{d_x \vee d_y})$, which could still be large.
But note that \emph{with the same per-iteration complexity}, one can then exponentially fast achieve NI error less than $\eps$, for any $\eps < \eps_{\text{warm-up}}$. This ``high-accuracy'' regime is where our method improves upon classical discretization-based algorithms.

Note however that, unfortunately, we cannot control the size $r_0$ of the neighborhood where our local result applies. Even worse: even if the quantity $r_0$ was known, this would still not suffice to certify efficiently that that neighborhood is reached, as the NI error is difficult to compute and even to upper-bound.
% Using the Lyapunov function we defined in Sec.~3 is not an option either, as it depends on the solution $(\mu^*, \nu^*)$ itself.
% An alternative attempt could be, instead of using a warm-up procedure followed by our algorithm CP-PP, to simply use CP-PP with a small step-size for the position variables: $\sigma \ll \eta$.
% One could indeed hope that the global convergence guarantees of Fisher-Rao Proximal Point for finite games would transfer to global guarantees for CP-PP with small $\sigma$.
% However, the known analyses for last-iterate convergence of Fisher-Rao Proximal Point (and related methods: Mirror Prox, Optimistic Mirror Descent-Ascent) are not adapted for this purpose, as they require the finite game to have a unique MNE, and the convergence rate itself depends on the finite game in a complicated way.
Thus we are at present unable to deduce a provably globally convergent algorithm from our work.
(If convergence rates are not desired, then it suffices to use the two-phase procedure proposed above along with some form of the doubling trick to choose $\eps_{\text{warm-up}}$.)
% but repeatedly running Fisher-Rao Proximal Point with the doubling trick would also suffice, albeit much slower.)
% In practice, 

%\end{document}

% !TEX root = ../main.tex
% \documentclass[../main]{subfiles}
%\begin{document}

\section{Convergence proof} \label{sec:cv_proof}

Throughout this section, we make 
the Assumptions~\ref{assum:1}-\ref{assum:6} described in \autoref{subsec:main_res:pb_stmt}.
Furthermore, to lighten notation, we denote by $z^k = (a^k, x^k, b^k, y^k)$ the iterates of the CP-PP algorithm.
We denote the sparse unique MNE as
\begin{equation*}
    \mu^* = \sum_{I \in [n^*]} a^*_I \delta_{x^*_I}
    ~~~~\text{and}~~~~
    \nu^* = \sum_{J \in [m^*]} b^*_J \delta_{y^*_J}
    ~~~~\text{with}~ a^*_I, b^*_J > 0.
\end{equation*}

\paragraph{The variational inequality characterizing CP-PP.}

From \autoref{lm:main_res:ppm_upd_well_defined}, we know that $z^{k+1}$ is well-defined as the saddle point of a convex-concave function in the interior of its domain.
Just by writing out the first-order optimality conditions in the argmin/argmax, we see that the CP-PP update \eqref{eq:main_res:ppa_upd} is characterized by the variational formula
\begin{align} \label{eq:cv_proof:ppa_ineq_foc}
	\forall z = (a, x, b, y),~
	\eta \tgap(z ; z^{k+1}) 
	&\leq 
	\sum_i (a_i-a^{k+1}_i) \log \frac{a^{k+1}_i}{a^k_i}
	+ \sum_j (b_j-b^{k+1}_j) \log \frac{b^{k+1}_j}{b^k_j} \\
	&\!\!\! + \frac{\eta}{\sigma} \sum_i a^k_i \innerprod{x^{k+1}_i-x^k_i}{x_i-x^{k+1}_i}
	+ \frac{\eta}{\sigma} \sum_j b^k_j \innerprod{y^{k+1}_j-y^k_j}{y_j-y^{k+1}_j}
\end{align}
(both sides are linear in $z-z^{k+1}$),
where we introduce, for all $z = (a, x, b, y)$ and $\hz = (\ha, \hx, \hb, \hy)$,
\begin{equation*}
	\tgap(z ; \hz) 
	\coloneqq \innerprod{\begin{pmatrix}
		\nabla_a \\
		\nabla_x \\
		-\nabla_b \\
		-\nabla_y
	\end{pmatrix}
	F_{n,m}(\hz)}{\begin{pmatrix}
		\ha - a \\
		\hx - x \\
		\hb - b \\
		\hy - y
	\end{pmatrix}},
    \quad
    \gap(z; \hz) \coloneqq 
    F_{n,m}((\ha, \hx), (b, y)) - F_{n,m}((a, x), (\hb, \hy)).
    % F_{n,m}\begin{pmatrix} \ha, \\ \hx, \\ b, \\ y \end{pmatrix} - F_{n,m}\begin{pmatrix} a, \\ x, \\ \hb, \\ \hy \end{pmatrix}
\end{equation*}
% Since it simply expresses the first-order optimality condition of $z^{k+1}$,  \eqref{eq:cv_proof:ppa_ineq_foc} is linear in $z-z^{k+1}$.
% This is why we can assert that 
% %there exists $z^{k+1}$ such that 
% the inequality holds for all $z \in (\Delta_n \times \XXX^n) \times (\Delta_m \times \YYY^m)$ and not just for $z$ in a neighborhood of $z^k$, as it would first seem from \autoref{lm:main_res:ppm_upd_well_defined}.
% Furthermore, 
One can check that $z^k$ is in the (relative) interior of $(\Delta_n \times \XXX^n) \times (\Delta_m \times \YYY^m)$ for all $k$ --- provided $z^0$ is ---, so \eqref{eq:cv_proof:ppa_ineq_foc} holds in fact with an equality
(but we continue to write inequalities to show the generality of our arguments).

The significance of the quantity $\tgap(z; \hz)$ comes from the fact that, if $F_{n,m}$ was convex-concave, $\hz$ would be a saddle point if and only if 
$\forall z, \tgap(z; \hz) \leq 0$,
% $\max_z \tgap(z; \hz) \leq 0$,
as a solution of the Stampacchia variational inequality \cite[Sec.~2.1]{diakonikolas_efficient_2021}.
% Technically IIUC (local) saddle points can be defined exactly defined as solutions of the Stampacchia variational inequality, but in this project we only need to define saddle points for convex-concave functions
Furthermore, $\tgap(z; \hz)$ can also be interpreted as a first-order approximation of $\gap(z; \hz)$, whose significance is that $\NI \left(\sum_i \ha_i \delta_{\hx_i}, \sum_j \hb_j \delta_{\hy_j} \right) = \max_z \gap(z; \hz)$.

\subsection{Exact-parametrization case} \label{subsec:cv_proof:exactparam}

In this subsection, we present a short proof of our result in the simpler case where we additionally assume that the number of particles $(n, m)$ is exactly equal to the sparsity of the solution $(n^*, m^*)$.
Relabel the solution particles $(a^*_I, x^*_I)$ resp.\ $(b^*_J, y^*_J)$ arbitrarily so that they are indexed by $i \in [n] = [n^*]$ resp.\ $j \in [m] = [m^*]$.

The convergence analysis relies on the Lyapunov function
$V(a, x, b, y) = V(a, x) + V(b, y)$ where
\begin{equation} \label{eq:cv_proof:exactparam:def_V}
	V(a, x) = 
    \underbrace{
        \KLdiv(a^*, a) 
    }_{\eqqcolon V_\wei(a, x)}
    ~+~ \frac{\eta}{\sigma}~
    \underbrace{
        \frac{1}{2}
        \sum_{i=1}^n a_i \norm{x^*_i-x_i}^2
    }_{\eqqcolon V_\pos(a, x)}
\end{equation}
and similarly for $V(b, y)$.
For ease of presentation, also let
$V_1(a, x) = V_\wei(a, x) + V_\pos(a, x)$,
and similarly for $V_1(b, y)$ and $V_1(a, x, b, y)$.
Note that we always have 
$(1 \wedge \sigma/\eta) V \leq V_1 \leq (1 \vee \sigma/\eta) V$.

Note that $V(a, x, b, y) \geq 0$ and that equality holds if and only if $(a, x, b, y) = (a^*, x^*, b^*, y^*)$.
We can also relate this quantity to the NI error
% , which is after all the natural measure of optimality for min-max problems, 
as follows;
in particular $V$ is arbitrarily small for $\NI$ small, and vice-versa.
The proof can be found in \autoref{subsec:rel_Lya_NI:exactparam}.
\begin{proposition} \label{prop:cv_proof:exactparam:NI_equiv_V}
	Assume that $n = n^*, m = m^*$ and define $V_1$ as in \eqref{eq:cv_proof:exactparam:def_V}.
	
	There exists a constant $C>0$ dependent only on $(f, \XXX, \YYY)$ such that,
	for any $z = (a, x, b, y)$, denoting $\mu = \sum_i a_i \delta_{x_i}$ and $\nu = \sum_j b_j \delta_{y_j}$,
	\begin{equation*}
		\NI(\mu, \nu) \leq C \sqrt{V_1(z)}.
	\end{equation*}
	
	Moreover, there exist $C', r>0$ dependent only on $(f, \XXX, \YYY)$ such that if $\NI(\mu, \nu) \leq r$, then up to permuting the labels of the solution particles (so that, for each $i \in [n]$, $x_i$ is in a neighborhood of $x^*_i$, and not necessarily of $x^*_{i'}$ for $i' \neq i$),
	\begin{equation*}
	    C' V_1(z)^{5/4} \leq \NI(\mu, \nu).
	\end{equation*}
\end{proposition}

Our choice of Lyapunov function is essentially a proxy for the squared WFR distance \eqref{eq:main_res:def_WFR22} of $\mu$ to $\mu^*$ and of $\nu$ to $\nu^*$, as shown in \cite{chizat_sparse_2021}. In our notation:
\begin{proposition}[{\cite[Lemma~D.1]{chizat_sparse_2021}}] \label{prop:cv_proof:exactparam:WFR_equiv_V}
    Assume that $n = n^*$ and define $V$ as in \eqref{eq:cv_proof:exactparam:def_V}.
    There exist constants $C, r>0$ (dependent only on $\mu^*$) such that for any $(a, x)$ with $V(a, x) \leq r$,
    denoting $\mu = \sum_i a_i \delta_{x_i}$,
    \begin{equation*}
        \WFR_2^2(\mu, \mu^*) \leq 2 V(a, x) \left( 1 + C \frac{\eta}{\sigma} \right).
    \end{equation*}
\end{proposition}

The main result of this subsection is that the CP-PP algorithm converges locally at an exponential rate, as measured by the Lyapunov function.
Convergence measured by NI error and by WFR distance (\autoref{thm:main_res:loc_exp_cv_NI} with the additional exact-parametrization assumption) could be shown by combining \autoref{thm:cv_proof:exactparam:loc_exp_cv} with \autoref{prop:cv_proof:exactparam:NI_equiv_V} and \autoref{prop:cv_proof:exactparam:WFR_equiv_V}.

\begin{theorem} \label{thm:cv_proof:exactparam:loc_exp_cv}
	Assume that $n = n^*, m = m^*$ and define $V$ as in \eqref{eq:cv_proof:exactparam:def_V}.
	
    Fix any $\Gamma_0 \geq 1$. There exist $\eta_0, \sigma_0$ such that for all $\eta \leq \eta_0, \sigma \leq \sigma_0$ with $\Gamma_0^{-1} \leq \frac{\sigma}{\eta} \leq \Gamma_0$,
	there exists $r_0>0$ such that if $V(z^0) \leq r_0$, then
    the CP-PP iterates $z^k$ satisfy
	\begin{equation*}
		\forall k,~ V(z^k) \leq V(z^0) \left( 1-\kappa \right)^k
	\end{equation*}
	for some constant $\kappa>0$.
\end{theorem}

More precisely, one can check from the last step of the proof that the rate $\kappa$ and the localness level $r_0$ can at most be chosen equal to $\eta^2$ times a constant (dependent on $(f, \XXX, \YYY)$ and $\Gamma_0$).

\begin{proof}
	Evaluate \eqref{eq:cv_proof:ppa_ineq_foc} at $z = (a^*, x^*, b^*, y^*)$.
	By the Bregman three-point identity on $h: a \mapsto a \log a - a + 1$ --- so that the Kullback-Leibler divergence $\KLdiv(\cdot, \cdot)$ is equal to the Bregman divergence of $h$ summed component-wise --- and 
% 	Pythagorean identity
	on $x \mapsto \frac{1}{2} \norm{x}^2$, we can rewrite the obtained inequality as
	\begin{align} \label{eq:cv_proof:exactparam:ppa_ineq_star}
		~ \eta \tgap(z^*; z^{k+1}) 
		&\leq V(z^k) - V(z^{k+1}) \\
		&\!\!\!\! -\, \underbrace{\! \left(
			\KLdiv(a^{k+1}\!, a^k) + \KLdiv(b^{k+1}\!, b^k)
			+ \frac{\eta}{2\sigma} \sum_i a^k_i \norm{x^{k+1}_i-x^k_i}^2 + \frac{\eta}{2\sigma} \sum_j b^k_j \norm{y^{k+1}_j-y^k_j}^2
		\right) \!}_{\eqqcolon~ D(k+1, k)} \\
		&\!\!\!\! +\, \underbrace{
            \frac{\eta}{2 \sigma} \sum_i (a^{k+1}_i-a^k_i) \norm{x^*_i-x^{k+1}_i}^2 
		      + \frac{\eta}{2 \sigma} \sum_j (b^{k+1}_j-b^k_j) \norm{y^*_j-y^{k+1}_j}^2
        }_{\eqqcolon~ \mathrm{[err]}}.
	\end{align}

	Now one can show that, 
    if $\eta \leq \eta_0, \sigma \leq \sigma_0$ and $V(z^k) \leq r_0$
    for some $\eta_0, \sigma_0, r_0$ dependent only on $(f, \XXX, \YYY$) and $\Gamma_0$,
    then both $\left[ \min_i a^k_i \wedge \min_j b^k_j \right]$ and $\left[ \min_i a^{k+1}_i \wedge \min_j b^{k+1}_j \right]$ are lower-bounded by a fixed positive constant (\autoref{lm:proof_gencase:ulb_on_weights_locally_k+1}),
    and so
	\begin{itemize}
		\item (\autoref{lm:proof_gencase:lb_tgap_olw0})
		The left-hand side is lower-bounded as
		$
			\eta \tgap(z^*; z^{k+1}) \geq 
			\eta \frac{\sigma_{\min}}{2} 
            % \left( 
			% 	\sum_i a^{k+1}_i \norm{x^*_i-x^{k+1}_i}^2 
			% 	+ \sum_j b^{k+1}_j \norm{y^*_j-y^{k+1}_j}^2 
			% \right)
            V_\pos(z)
			+ O ( \eta V(z^{k+1})^{3/2} )
		$
        for some $\sigma_{\min}>0$ dependent only on $(f, \XXX, \YYY)$.
        
		This inequality is a consequence of the ``quadratic growth'' and ``star-convexity'' properties discussed in 
        Sec.~\ref{subsubsec:cv_proof:proofingr_growthconds:quadr}, resp.\ \ref{subsubsec:cv_proof:proofingr_growthconds:starconvex}.
		\item (\autoref{lm:proof_gencase:lb_Dk+1_k})
		There exists a constant $C>0$ dependent only on $(f, \XXX, \YYY)$ and $\Gamma_0$ such that
		% the term appearing with a negative sign on the second line (the ``divergence from $(k+1)$ to $k$'') is lower-bounded by
		$
            % \MoveEqLeft
			% \KLdiv(a^{k+1}, a^k) + \KLdiv(b^{k+1}, b^k)
			% + \frac{\eta}{2\sigma} \sum_i a^k_i \norm{x^{k+1}_i-x^k_i}^2 + \frac{\eta}{2\sigma} \sum_j b^k_j \norm{y^{k+1}_j-y^k_j}^2 \\
            D(k+1, k)
            \geq C \eta^2 V(z^{k+1})
			+ O \left( \eta V(z^{k+1})^2 \right)
		$.
  
		This inequality follows from an ``error bound''-type result discussed in \autoref{subsubsec:cv_proof:proofingr_growthconds:error_bound}.
		\item (\autoref{lm:proof_gencase:errterm3__control_nonBregmanness_errterm})
		The terms on the third line, that arise due to the fact that the divergence used in the update \eqref{eq:main_res:ppa_upd} is not a Bregman divergence, are bounded as
		$
			% \frac{\eta}{2 \sigma} \sum_i (a^{k+1}_i-a^k_i) \norm{x^*_i-x^{k+1}_i}^2
            % + \frac{\eta}{2 \sigma} \sum_j (b^{k+1}_j-b^k_j) \norm{y^*_j-y^{k+1}_j}^2
            \mathrm{[err]}
			= O \left( \eta V(z^{k+1})^{3/2} \right)
		$.
	\end{itemize}
    In each of the bounds above, the $O(\cdot)$ hides a constant dependent only on $(f, \XXX, \YYY)$ and $\Gamma_0$.
    
	By plugging these bounds back into \eqref{eq:cv_proof:exactparam:ppa_ineq_star}, we obtain
	\begin{equation*}
		V(z^{k+1}) 
		\leq V(z^k) - C \eta^2 V(z^{k+1})
		+ O \left( \eta V(z^{k+1})^{3/2} \right).
	\end{equation*}
	In particular, since $V(z^{k+1})$ is bounded by a constant (due to \autoref{lm:proof_gencase:ulb_on_weights_locally_k+1}), for small enough $\eta_0$ and $\sigma_0$ we have that $V(z^{k+1}) \leq 2 V(z^k)$.
	By rearranging we get
	\begin{equation*}
		V(z^{k+1}) 
		\leq \frac{V(z^k)}{
		    1 + \eta^2 \left[ C - O \left( \frac{\sqrt{V(z^{k+1})}}{\eta} \right) \right]
	   }.
	\end{equation*}
	Hence, for appropriately small choices of $r_0$,
    we have
	$V(z^k) \leq r_0 \implies V(z^{k+1}) \leq 2r_0 \implies V(z^{k+1}) \leq V(z^k) (1-\kappa)$ 
	for some $\kappa>0$.
	The final result follows by induction.
\end{proof}

% \begin{remark}
% 	The Lyapunov function $V$ used in this section depends on the ordering of the ``true'' particles $(a^*_I, x^*_I)$ resp.\ $(b^*_J, y^*_J)$.
% 	It may seem strange that our result depends on this ordering, but our convergence analysis is only local, so it makes sense.
% \end{remark}

\subsubsection{Convergence of CP-MP}

In the exact-parametrization case it is relatively easy to extend our convergence result for CP-PP to CP-MP (the implementable variant of CP-PP, \autoref{alg:main_res:ppm}, with $L=2$).
Namely CP-MP converges under the same conditions and with the same rate as CP-PP. 
\begin{proposition} \label{prop:cv_proof:exactparam:mp_cv}
     The statement of \autoref{thm:cv_proof:exactparam:loc_exp_cv} also holds for $z^k$ being the iterates of CP-MP.
\end{proposition}
The proof of the proposition, in \autoref{sec:mp_exactparam}, essentially relies on the convergence result for CP-PP and on the fact that the CP-MP and CP-PP updates coincide up to
% terms of order 3 in the step-size times the gradient norm.
order-3 terms in the step-size.
In particular we derive order-2 expressions for the Mirror Prox
% (\autoref{lm:mp_exactparam:mp_expan})
and Proximal Point updates 
% (\autoref{lm:mp_exactparam:pp_expan}) 
under quite general assumptions
(\Autoref{lm:mp_exactparam:mp_expan,lm:mp_exactparam:pp_expan}),
which may be of independent interest.

% Moreover we show the error terms are also proportional to the projected gradient norm,
% which makes the conversion of convergence results for Proximal Point into ones for Mirror Prox very easy indeed.

\subsection{General case} \label{subsec:cv_proof:gencase}

In general, the sparsity of the solution $(n^*, m^*)$ is not known in advance, and $n \neq n^*, m \neq m^*$.
Contrary to the exact-parametrization case where the choice of Lyapunov function was relatively straightforward,
% (it was just the divergence used in the algorithm \eqref{eq:main_res:ppa_upd} evaluated at the solution on the left), 
here it must be carefully designed, due to overparametrization.
Indeed, the variables $(a, x, b, y)$ and the solution $(a^*, x^*, b^*, y^*)$ live in different spaces: $a \in \Delta_n \neq \Delta_{n^*} \ni a^*$,
so we cannot just evaluate the algorithm's characterizing inequality \eqref{eq:cv_proof:ppa_ineq_foc} at the solution.

\begin{figure}
    \begin{subfigure}[t]{0.48\textwidth}
    % \begin{minipage}{.48\textwidth}
        \centering
        \includegraphics[width=\textwidth]{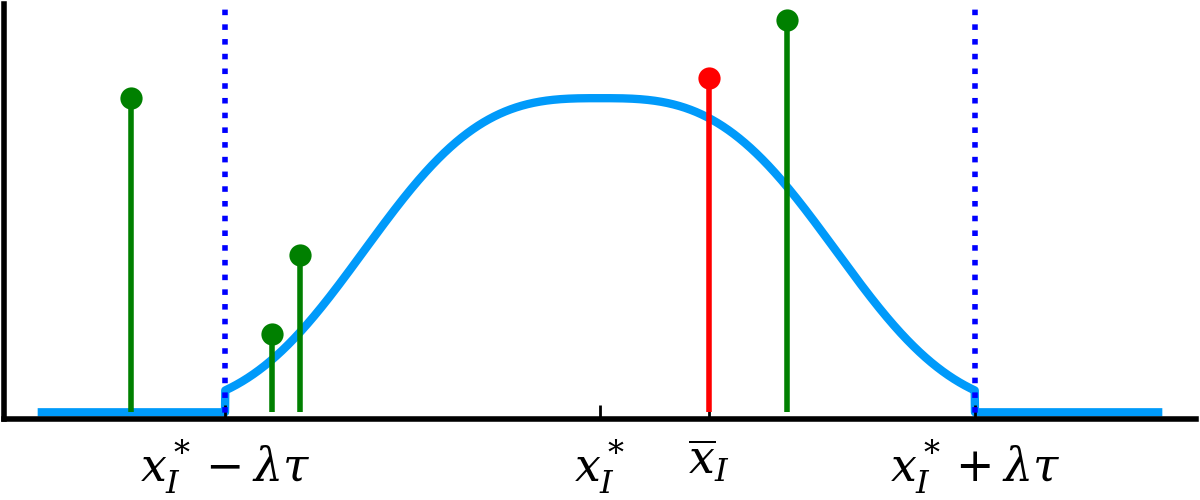}
        \caption{
            Green stems represent the particles: abscissa indicates $x_i$ and stem height indicates $a_i$
        }
        \label{fig:misc:constru_lya}
    % \end{minipage}
    \end{subfigure}
    \hfill
    \begin{subfigure}[t]{0.47\textwidth}
    % https://tex.stackexchange.com/questions/37581/latex-figures-side-by-side
    % \begin{minipage}{.47\textwidth} 
        \centering
        \includegraphics[width=\textwidth]{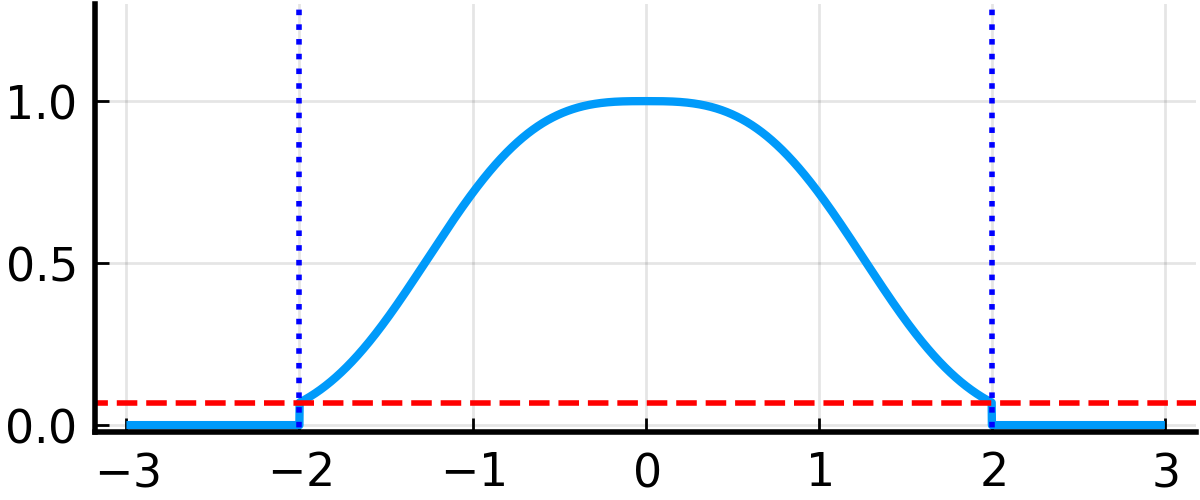}
        % \caption{$y = \ind_{\abs{x} \leq 2} \cdot \exp\left( -\frac{\abs{x}^3}{3} \right)$}
        \caption{$y = \exp\left( -\frac{\abs{x}^3}{3} \right)$ if $\abs{x} \leq 2$, $y=0$ otherwise}
        \label{fig:misc:psi}
    % \end{minipage}
    \end{subfigure}
    \caption{Illustration of the construction defining $V(a, x)$}
    \label{fig:misc}
\end{figure}

We define a Lyapunov function $V(a, x, b, y) = V(a, x) + V(b, y)$ 
%(for $a \in \Delta_n, x \in \XXX^n, b \in \Delta_m, y \in \YYY^m$) 
by the following construction, generalizing \cite[Eq.~(20)]{chizat_sparse_2021}.
See \autoref{fig:misc:constru_lya} for an illustration.
\begin{itemize}
	\item Fix $(\varphi_I)_{I \in [0, n^*]}$ a partition of unity of $\XXX$ centered at the $(x^*_I)_I$, i.e.,
	\begin{itemize}
		\item Each $\varphi_I$ is a measurable function $\XXX \to \RR$;
		\item $\forall I \in [n^*],~ \varphi_I \geq 0$ 
		and $\varphi_0 = 1 - \sum_{I=1}^{n^*} \varphi_I \geq 0$
		over $\XXX$;
		\item $\forall I \in [n^*],~ \varphi_I(x^*_I) = 1$.
	\end{itemize}
	\item For any $a \in \Delta_n, x \in \XXX^n$, define the \emph{aggregated weights}, the \emph{aggregated positions} and the \emph{local covariance matrices} of $\mu = \sum_i a_i \delta_{x_i}$ as
	\begin{align*}
		\forall I \in [0,n^*],~
		\ola_I &= \int_\XXX \varphi_I d\mu
		& &~~\text{and}~~ &
		\forall I \in [n^*],~
		\olx_I &= 
% 		\int_\XXX x~ d\left( \frac{\varphi_I \mu}{\ola_I} \right)(x) \\
		\int_\XXX x~ \frac{\varphi_I(x)}{\ola_I} d\mu(x) \\
		& & & & 
		\Sigma_I &= 
% 		\int_\XXX (x-\olx_I) (x-\olx_I)^\top d\left( \frac{\varphi_I \mu}{\ola_I} \right)(x).
		\int_\XXX (x-\olx_I) (x-\olx_I)^\top~ \frac{\varphi_I(x)}{\ola_I} d\mu(x).
	\end{align*}
	I.e., in discrete notation,
	\begin{align*}
		% \label{eq:cv_proof:gencase:localmoments}
		\forall I \in [n^*],~
		\ola_I &= \sum_i \varphi_{Ii} a_i
		&
		\olx_I &= \sum_i \frac{\varphi_{Ii} a_i}{\ola_I} x_i
		&
		\Sigma_I &= \sum_i \frac{\varphi_{Ii} a_i}{\ola_I} (x_i-\olx_I) (x_i-\olx_I)^\top
	\end{align*}
	where $\varphi_{Ii} = \varphi_I(x_i)$, 
	and $\ola_0 = 1 - \sum_I \ola_I$ is the \emph{stray weight}.
% 	(not counted in any of the aggregates).
	\item Let, for any $a \in \Delta_n, x \in \XXX^n$,
	\begin{equation} \label{eq:cv_proof:gencase:def_V}
		V(a, x) = \underbrace{
			\KLdiv(a^*, \ola)
		}_{\eqqcolon V_\wei(a, x)}
		~+~ \frac{\eta}{\sigma}~
        \underbrace{ 
            \frac{1}{2}
			\sum_I \ola_I \left( 
			\norm{x^*_I-\olx_I}^2 
			+ \trace(\Sigma_I)
			\right)
		}_{\eqqcolon V_\pos(a, x)}.
	\end{equation}
\end{itemize}
Similarly, fix $(\psi_J)_{J \in [0,m^*]}$ a partition of unity of $\YYY$ centered at the $y^*_J$,
similarly define $\olb \in \Delta_{m^*}$ and $\oly \in \YYY^{m^*}$ for any $b \in \Delta_m, y \in \YYY^m$,
and similarly define $V(b, y)$.
For ease of presentation, also let
$V_1(a, x) = V_\wei(a, x) + V_\pos(a, x)$,
and similarly for $V_1(b, y)$ and $V_1(a, x, b, y)$.
Note that we always have 
$(1 \wedge \sigma/\eta) V \leq V_1 \leq (1 \vee \sigma/\eta) V$.

The Lyapunov function $V$ depends on the choice of partitions of unity $(\varphi_I)_I$ and $(\psi_J)_J$.
They can be freely designed so as to make the proof go through, as long as they satisfy the conditions announced above (non-negative, sum to $1$, $\varphi_I(x^*_I)=1$).
For example, our analysis for the exact-parametrization case was equivalent to choosing as $\varphi_I$ the indicator function of a small ball centered at $x^*_I$
(see~\autoref{claim:growth_conds:gencase_imply_exactparam}).
Proving convergence in the general case requires a much subtler choice;
specifically, the partitions of unity we use for the proof of our main result are defined as
\begin{equation} \label{eq:cv_proof:gencase:def_varphi}
	\varphi_I(x) = \begin{cases}
		\exp\left( -\frac{\norm{x-x^*_I}^3}{3 \tau^3} \right) & ~\text{if}~ \norm{x-x^*_I} \leq \lambda\tau \\
		0 & ~\text{otherwise}
	\end{cases}
\end{equation}
for some bandwidth resp.\ cut-off parameters $\tau, \lambda > 0$ chosen as functions of $\eta$ and $\sigma$.
The cut-off parameter $\lambda$ is used to ensure that $\sum_{I' \neq I} \varphi_{I'}(x^*_I) = 0$ and so $\varphi_0(x^*_I) \geq 0$.
See \autoref{fig:misc:psi} for an illustration in one dimension with $\tau=1$, $\lambda=2$.

Note that $V(a, x, b, y) \geq 0$ with equality if and only if $(\ola, \olx, \olb, \oly) = (a^*, x^*, b^*, y^*)$ and $\Sigma_I, \Sigma_J = 0$ for all $I, J$, i.e., if and only if $(\mu, \nu) = (\mu^*, \nu^*)$.
Beyond this equivalence, similarly as for the exact-parametrization case, we can show the following relation between $V$ and NI error.
The proof of the following proposition, as well as a more quantitative version of it, can be found in \autoref{sec:rel_Lya_NI}.%
\footnote{
	The reason why we need to split \autoref{prop:cv_proof:gencase:NI_equiv_V} into two parts is that the inequality
	$V^\alpha \lesssim \NI$ (for any exponent $\alpha$)
	cannot be true for all $(a, x, b, y)$. Indeed, $\NI$ is bounded, but $V$ may be infinite due to the terms $\KLdiv(a^*, \ola)$ if $\mu$ has no mass near one of the $x^*_I$'s.
% 	(Incidentally, in this case the Lyapunov function becomes a poor approximation of Wasserstein-Fisher-Rao distance...)
} 
\begin{proposition} \label{prop:cv_proof:gencase:NI_equiv_V}
    Define $V_1$ as in \eqref{eq:cv_proof:gencase:def_V} with the partitions of unity $(\varphi_I)_I$ and $(\psi_J)_J$ as in \eqref{eq:cv_proof:gencase:def_varphi}. 
    Suppose that 
    % $\lambda, \tau$ are such that the $\support(\varphi_I)$ resp.\ $\support(\psi_J)$ do not overlap, and that   --> this is included in the following condition
    $\lambda\tau$ is less than some constant dependent on $(f, \XXX, \YYY)$.
	
	There exists a constant $C>0$ dependent only on $(f, \XXX, \YYY)$ such that, for any $z = (a, x, b, y)$,
    denoting $\mu = \sum_i a_i \delta_{x_i}$ and $\nu = \sum_j b_j \delta_{y_j}$,
	\begin{equation*}
	    \NI(\mu, \nu) \leq C \sqrt{V_1(z)}.
	\end{equation*}
	
	Moreover, there exist $C', r>0$ dependent on $(f, \XXX, \YYY)$, $\lambda$ and $\tau$ such that, if $\NI(\mu, \nu) \leq r$, then
	\begin{equation*}
		C' V_1(z)^{5/4} \leq \NI(\mu, \nu).
	\end{equation*}
    More precisely
    if $\lambda, \tau$ are chosen as functions of $\eta, \sigma$ as in \eqref{eq:proof_gencase:choice_lambda_tau}
    and $\Gamma_0^{-1} \leq \frac{\sigma}{\eta} \leq \Gamma_0$ for some $\Gamma_0\geq 1$,
	then $r$ and $C'$ can be chosen as $\sqrt{\sigma}$ times constants dependent only on $(f, \XXX, \YYY)$ and $\Gamma_0$.
\end{proposition}

The Lyapunov function $V$ is by design essentially a proxy for squared WFR distance \eqref{eq:main_res:def_WFR22} of $\mu$ to $\mu^*$ and of $\nu$ to $\nu^*$.
Indeed we followed the same construction as in \cite[Lemma~D.1]{chizat_sparse_2021}, with the nuance that $\varphi_I$ and $\psi_J$ are not necessarily indicator functions.
A simple modification of their proof shows that, in our notation:% 
\footnote{
    Namely the modification to bring to the proof of \cite[Lemma~D.1]{chizat_sparse_2021} is (in the notations of that paper)
    to use the transport plan that sends
    $(r, \theta)$ to $\left( \frac{r}{\olr_I} r^*_I, \theta^*_I \right)$ with probability $\varphi_I(\theta)$ and to $(0, \theta)$ with probability $\varphi_0(\theta)$.
    The present work uses Kullback-Leibler divergence while \cite{chizat_sparse_2021} uses squared Hellinger distance, but the difference can be controlled similarly as in \autoref{lm:aux_lemmas:KL_chi2_comparison}, thanks to \autoref{lm:growth_conds:ulb_on_weights_locally}.
    The factor $\frac{\eta}{\sigma}$ can be thought of simply as a linear rescaling of $\norm{\cdot}_\XXX^2$.
}
\begin{proposition}[{Modification of \cite[Lemma~D.1]{chizat_sparse_2021}}]  \label{prop:cv_proof:gencase:WFR_equiv_V}
	Define $V$ as in \eqref{eq:cv_proof:gencase:def_V} with the partitions of unity $(\varphi_I)_I$ and $(\psi_J)_J$ as in \eqref{eq:cv_proof:gencase:def_varphi}.
    There exist constants $C, r>0$ (dependent only on $\mu^*$) such that for any $(a, x)$ with $V(a, x) \leq r$,
    denoting $\mu = \sum_i a_i \delta_{x_i}$,
    \begin{equation*}
        \WFR_2^2(\mu, \mu^*) 
        \leq 2 V(a, x) \left( 1+ C \frac{\eta}{\sigma} (\lambda\tau)^2 \right).
    \end{equation*}
\end{proposition}

The main result of this subsection, proved in \autoref{sec:proof_gencase}, is that the CP-PP algorithm converges locally at an exponential rate, as measured by the Lyapunov function.
Convergence measured by NI error and by WFR distance (\autoref{thm:main_res:loc_exp_cv_NI}) follows by combining \autoref{thm:cv_proof:gencase:loc_exp_cv} with \autoref{prop:cv_proof:gencase:NI_equiv_V} and \autoref{prop:cv_proof:gencase:WFR_equiv_V} as shown in \autoref{sec:proof_mainres}.

\begin{theorem} \label{thm:cv_proof:gencase:loc_exp_cv}
	Define $(\varphi_I)_I, (\psi_J)_J$ as in \eqref{eq:cv_proof:gencase:def_varphi}
	and define $V$ as in \eqref{eq:cv_proof:gencase:def_V}.
	Choose $\lambda, \tau$ as functions of $\eta, \sigma$ as in \eqref{eq:proof_gencase:choice_lambda_tau}.
	
    Fix any $\Gamma_0 \geq 1$.
    There exist $\eta_0, \sigma_0$ such that for all $\eta \leq \eta_0, \sigma \leq \sigma_0$ with $\Gamma_0^{-1} \leq \frac{\sigma}{\eta} \leq \Gamma_0$,
	there exists $r_0>0$ such that if $V(z^0) \leq r_0$, then the CP-PP iterates satisfy
	\begin{equation*}
		\forall k,~ V(z^k) 
        \leq V(z^0) (1-\kappa)^k
	\end{equation*}
	for some constant $\kappa>0$.
\end{theorem}

More precisely, one can check from the last step of the proof (\autoref{subsec:proof_gencase:proof_conclusion}) that 
% $\eta_0$ and $\sigma_0$ only need to depend on $(f, \XXX, \YYY)$ and on $\Gamma_0$, and that      --> this is stupid, what else could they depend on
the rate $\kappa$ can at most be chosen equal to $\eta^2$ times a constant,
and that the localness level $r_0$ can at most be chosen equal to $\eta^3$ times a constant
(dependent on $(f, \XXX, \YYY)$ and $\Gamma_0$).

The full proof of the theorem can be found in \autoref{sec:proof_gencase}.
It has the same general structure as for the exact-parametrization case, but needs to deal with the following difficulties:
\begin{itemize}
    \item The variables $(a, x, b, y)$ and the solution $(a^*, x^*, b^*, y^*)$ live in different spaces so we cannot just evaluate the characterizing inequality \eqref{eq:cv_proof:ppa_ineq_foc} at the solution particles. 
    Instead we identify a notion of ``proxy solution particles'' $(\aps, \xps, \bps, \yps) \in \Delta_n \times \XXX^n \times \Delta_m \times \YYY^m$, namely
    \begin{equation} \label{eq:cv_proof:gencase:choice_of_proxy*}
    	\xps_i \coloneqq x^{k+1}_i + \sum_{I \in [n^*]} \varphi^{k+1}_{Ii} (x^*_I - x^{k+1}_i)
    	~~~\text{and}~~~
    	\aps_i \coloneqq \sum_{I \in [n^*]} a^*_I \frac{\varphi^{k+1}_{Ii} a^{k+1}_i}{\ola^{k+1}_I}
    \end{equation}
    where $\varphi^{k+1}_{Ii} = \varphi_I(x^{k+1}_i)$,
    and similarly for $\bps, \yps$.
    We show that evaluating \eqref{eq:cv_proof:ppa_ineq_foc} at $(\aps, \xps, \bps, \yps)$ makes the Lyapunov function emerge naturally, yielding a general-case equivalent of \eqref{eq:cv_proof:exactparam:ppa_ineq_star}.
    \item Compared to the exact-parametrization case, several additional error terms appear, which are much more delicate to control.
    For example, we need to bound
	% $\frac{1}{2} \sum_I \sum_i a^k_i \left( \varphi^{k+1}_{Ii}-\varphi^k_{Ii} \right) \norm{x^*_I-x^k_i}^2$,
	$\frac{1}{2} \sum_I \sum_i a^k_i ( \varphi^{k+1}_{Ii}-\varphi^k_{Ii} ) \norm{x^*_I-x^k_i}^2$,
	a term which did not appear in \eqref{eq:cv_proof:exactparam:ppa_ineq_star}.
	This requires some technical work, and it is here that we benefit from choosing the partitions of unity as \eqref{eq:cv_proof:gencase:def_varphi}; the choice of the parameters $\lambda$ and $\tau$ also requires care.
	\item The stray-weight and variance terms of the Lyapunov function do not appear in the error-bound-type inequality of \autoref{subsubsec:cv_proof:proofingr_growthconds:error_bound}, so we do need to combine that result with the quadratic growth and star-convexity properties, contrary to the exact-parametrization case.
\end{itemize}

\subsection{A crucial proof ingredient: lower growth properties} \label{subsec:cv_proof:proofingr_growthconds}

For unconstrained min-max optimization of a smooth convex-concave objective $G(x, y)$, 
the proof of convergence of the (Euclidean) Proximal Point method
essentially reduces to three steps:
\begin{enumerate}
    \item Letting $\gap(z; \hz) = G(\hx, y) - G(x, \hy)$ and 
    $\tgap(z; \hz) = \innerprod{
        \begin{pmatrix}
            \nabla_x \\
            -\nabla_y
        \end{pmatrix} G(\hz)
    }{
        \hz-z
    }$,
    notice that 
    \begin{equation*}
        \tgap(z; \hz) = \gap(z; \hz) + D_{G(\cdot, \hy)}(x, \hx) - D_{G(\hx, \cdot)}(y, \hy)
        \geq \gap(z; \hz)
    \end{equation*}
    where $D$ denotes a Bregman divergence, by convexity-concavity.
    \item The Proximal Point update is characterized by the variational inequality, analogous to \eqref{eq:cv_proof:ppa_ineq_foc},
    % (here we have in fact equality due to the absence of constraints)
    \begin{equation*}
        % \forall z \in \RR^{d_x} \times \RR^{d_y},~
        \forall z,~
        \eta \tgap(z; z^{k+1}) 
        \leq 
        \innerprod{z-z^{k+1}}{z^{k+1}-z^k}
        = \frac{1}{2} \norm{z-z^k}^2 
        - \frac{1}{2} \norm{z-z^{k+1}}^2
        - \frac{1}{2} \norm{z^{k+1}-z^k}^2.
    \end{equation*}
    \item In particular evaluating at the saddle point $z^*$ assumed unique for simplicity,
    \begin{align*}
        % \gap(z^*; z^{k+1}) 
        % \leq \tgap(z^*; z^{k+1}) 
        % \leq \frac{1}{2\eta} \norm{z^*-z^k}^2 
        % - \frac{1}{2\eta} \norm{z^*-z^{k+1}}^2
        % - \frac{1}{2\eta} \norm{z^{k+1}-z^k}^2 \\
        \frac{1}{2} \norm{z^*-z^{k+1}}^2 
        \leq \frac{1}{2} \norm{z^*-z^k}^2 
        &- \eta \gap(z^*; z^{k+1}) \\
        &- \eta \left[ D_{G(\cdot, y^{k+1})}(x^*, x^{k+1}) - D_{G(x^{k+1}, \cdot)}(y^*, y^{k+1}) \right] \\
        &- \frac{1}{2} \norm{z^{k+1}-z^k}^2.
    \end{align*}

    Depending on the properties of $G$, we lower-bound one of the three terms appearing with a
    negative sign on the right-hand side. For example,
    \begin{itemize}
        \item If $G$ satisfies (a min-max analog of) the quadratic growth property \cite[Def.~5.1]{guille-escuret_study_2020}, i.e., if there exists $C>0$ such that
        \begin{equation*}
           \forall z,~ \gap(z^*; z) \geq  \frac{C}{2} \norm{z^*-z}^2,
        \end{equation*}
        then we can directly conclude
        % to linear convergence with a rate $\asymp \frac{C}{2} \eta$.
        to exponential decrease of the Lyapunov function $V(z) = \frac{1}{2} \norm{z^*-z}^2$
        with a rate at least $\frac{C}{2} \eta$.
        % since
        % \begin{align*}
        %     \frac{1}{2} \norm{z^*-z^{k+1}}^2 
        %     &\leq \frac{1}{2} \norm{z^*-z^{k+1}}^2 
        %     - \frac{\eta C}{2} \norm{z^*-z^{k+1}}^2 \\
        %     \text{so}~~~
        %     \frac{1}{2} \norm{z^*-z^{k+1}}^2
        %     &\leq \frac{1}{1 + \eta C/2} \cdot \frac{1}{2} \norm{z^*-z^k}^2
        %     \leq \left( \frac{1}{1 + \eta C/2} \right)^{k+1} \cdot \frac{1}{2} \norm{z^*-z^0}^2.
        % \end{align*}
        % \item If $G$ is $\mu$-strongly-convex-strongly-concave, in particular it satisfies quadratic growth with $C=\mu$, and moreover
        % % $
        % %     % \forall z, \hz,~ 
        % %     D_{G(\cdot, \hy)}(x, \hx) - D_{G(\hx, \cdot)}(y, \hy)
        % %     \geq 
        % %     \frac{\mu}{2} \norm{z-\hz}^2
        % % $
        % % for all $z, \hz$.
        % \begin{equation*}
        %     \forall z,
        %     D_{G(\cdot, y)}(x^*, x) - D_{G(x, \cdot)}(y^*, y)
        %     \geq 
        %     \frac{\mu}{2} \norm{z^*-z}^2.
        % \end{equation*}
        % So we directly get exponential convergence with a rate at least $\mu \eta$.
        \item If $G$ satisfies the error bound property with constant $C'>0$ \cite{tseng_linear_1995, hsieh_explore_2020},
        % \cite{wei_linear_2021}, 
        i.e., if
        \begin{equation*}
            \forall z,~ \norm{\begin{pmatrix}
                \nabla_x G(z) \\
                -\nabla_y G(z)
            \end{pmatrix}} \geq C' \norm{z-z^*},
        \end{equation*}
        then since
        % one can check that
        $z^{k+1}-z^k = \eta \begin{pmatrix}
            -\nabla_x G(z^{k+1}) \\
            \nabla_y G(z^{k+1})
        \end{pmatrix}$, 
        % then $\frac{1}{2} \norm{z^{k+1}-z^k}^2 \geq \frac{(C' \eta)^2}{2} \norm{z^*-z^{k+1}}^2$, and 
        we can conclude to exponential convergence with a rate at least $(C' \eta)^2$.
        \item We note that convexity-concavity of $G$ is not essential. (Suppose existence and uniqueness of the saddle point $z^*$ is guaranteed by some other property that convexity-concavity.)
        Indeed, to lower-bound the second term on the right-hand side, it suffices to have $G$ ($\mu'$-strongly) star-convex-concave \cite[Def.~5.1]{guille-escuret_study_2020}, that is
        \begin{equation*}
            \forall z,~
            D_{G(\cdot, y)}(x^*, x) - D_{G(x, \cdot)}(y^*, y)
            \geq 
            \frac{\mu'}{2} \norm{z^*-z}^2.
        \end{equation*}
    \end{itemize}
\end{enumerate}
In total, if $G$ is ($\mu'$-strongly) star-convex-concave, satisfies quadratic growth with constant $C$, and error bound with constant $C'$, then 
we can conclude to exponential convergence with a rate at least $\frac{C + \mu'}{2} \eta + (C' \eta)^2$.
% If $C$ or $C' > 0$, $\mu'$ need not be positive.
% (Smoothness is required to prove convergence for Gradient Descent-Ascent, but not for Proximal Point.)

In our case --- constrained min-max optimization of 
% the nonconvex-nonconcave objective 
the overparametrized objective
$F_{n,m}$ using the divergence 
$D\left((a, x), (\ha, \hx) \right) = \KLdiv(a, \ha) + \frac{\eta}{2\sigma} \sum_i \ha_i \norm{x_i-\hx_i}^2$
which is non-Euclidean and not even a Bregman divergence --- the analysis is significantly more technical, but it involves the same basic ingredients.

\subsubsection{``Quadratic growth'' with respect to the position and stray weight variables} \label{subsubsec:cv_proof:proofingr_growthconds:quadr}

We establish a quadratic growth property for $F_{n,m}$ involving only some of the desired terms in the lower bound.
Analogously to the analysis of \cite{chizat_sparse_2021} for minimization, the proof relies on the non-degeneracy Assumptions~\ref{assum:5}-\ref{assum:6};
a crucial difference however, is that we do not have quadratic growth in the weight variables.
To be precise, compared to the assumption~(A5) of \cite{chizat_sparse_2021}, our non-degeneracy assumption concerns only the so-called local kernels $H_I, H_J$, and the min-max analog of the global kernel ($K$ in that paper's notations) is necessarily zero due to the bilinearity of $F(\mu, \nu)$.

The precise statement of our result is given in \autoref{subsec:growth_conds:quadr}; here we state a simplified version to give the intuition.

\begin{lemma}[``Quadratic growth'', simplified]
    There exist constants $r, C>0$ only dependent on $(f, \XXX, \YYY)$ (and on $\lambda, \tau$ in the general case) such that, for any $z = (a, x, b, y)$ with $V_1(z) \leq r$, then
    \begin{equation*}
        F(\mu, \nu^*) - F(\mu^*, \nu) 
        \geq
		C \left( V_\pos(z) + \ola_0 + \olb_0 \right)
    \end{equation*}
    where $\mu = \sum_{i=1}^n a_i \delta_{x_i}$ and $\nu = \sum_{j=1}^m b_j \delta_{y_j}$
    (and $\ola_0=\olb_0=0$ by definition in the exact-\hspace{0pt}parametrization case).
\end{lemma}

Note that $\norm{\ola - a^*}_1$ does not appear in this inequality, but it will appear in the error-bound-type inequality discussed in the next paragraph. 
% Interestingly, 
Conversely $\ola_0$ and $\sum_I \ola_I \trace(\Sigma_I)$ 
% --- the variance part of $V_\wei$ --- 
appear in the inequality of this paragraph but not of the next.

\subsubsection{``Error bound'' with respect to the weight and aggregated position variables} \label{subsubsec:cv_proof:proofingr_growthconds:error_bound}

It is well-known that the error bound property holds 
for strongly-convex-strongly-concave and smooth min-max objectives, or 
for bilinear objectives constrained to a product of polytopes \cite{tseng_linear_1995}.
In our case, the reparametrized objective $F_{n,m}(a, x, b, y)$ 
% --- whose gradient descent-ascent flow we morally follow --- 
is bilinear in the weight components $(a,b)$,
and intuitively it possesses some local strong convexity-concavity in the position components $(x,y)$ thanks to Assumption~\ref{assum:6}.
But these two facts are not enough to directly show an error bound inequality, 
because the constant $C$ for the components $(a, b)$ may depend arbitrarily badly on $(x, y)$ a priori. % not super rigorous... but should be fine
Instead we use an argument, inspired by \cite[Lemma~14]{wei_linear_2021}, that also exploits the Assumption~\ref{assum:3} of uniqueness of the MNE.

Again, the precise statement of our result is deferred to the appendix \autoref{subsec:growth_conds:error_bound};
here we state an informal version to give the intuition.
We also refer to the second paragraph of that appendix for an interpretation of the quantity appearing on the left-hand side.

\begin{lemma}[``Error bound'', informal]
    Consider any $\hz = (\ha, \hx, \hb, \hy) \in \Delta_n \times \XXX^n \times \Delta_m \times \YYY^m$.
    For any $(A_I)_{I \in [n^*]}, (B_J)_{J \in [m^*]}$ (and $A_0=B_0=0$) and $(X_I)_{I \in [n^*]}, (Y_J)_{J \in [m^*]}$,
    define ``proxy particles'' from $(A, X, B, Y)$ and $(\ha, \hx, \hb, \hy)$ analogously to \eqref{eq:cv_proof:gencase:choice_of_proxy*}, and denote them by $z$.
    
	There exist $r, C > 0$ only dependent on $(f, \XXX, \YYY)$ such that if $V_1(\hz) \leq r$, then
	\begin{align*}
		\max_{A,X,B,Y} \tgap(z; \hz) 
		&\geq 
		C \sqrt{
			\sum_I d_h(a^*_I, \hola_I) + \sum_J d_h(b^*_J, \holb_J)
			+ \sum_I \hola_I \norm{\holx_I-x^*_I}^2 + \sum_J \holb_J \norm{\holy_J-y^*_J}^2
		} \\
		&~~ + O \left( V_1(\hz) \right).
	\end{align*}
\end{lemma}

\subsubsection{Local ``star-convexity-concavity'' (strong with respect to the position variables)}
\label{subsubsec:cv_proof:proofingr_growthconds:starconvex} 

Note that $\tgap(z; \hz)- \gap(z; \hz) = D_{F_{n,m}(\cdot, \hy)}(x, \hx) - D_{F_{n,m}(\hx, \cdot)}(y, \hy)$
where $D$ denotes a Bregman divergence.
Intuitively, $F_{n,m}$ is bilinear in the weight variables and, in a neighborhood of the MNE, it possesses some local strong convexity-concavity with respect to the position variables thanks to Assumption~\ref{assum:6}.
And indeed, by Taylor expansions, one can obtain lower-bounds on $\tgap(z; \hz)- \gap(z; \hz)$ consisting of positive terms and of error terms in $V(z)$, $V(\hz)$.
Note however that (for the general case) due to the overparametrization, there are many ways to write Taylor expansions.

Specifically, we will use the following bound.
Again the precise statement of the result is deferred to the appendix \autoref{subsec:growth_conds:starconvex}; here we state an informal version to give the intuition.

\begin{lemma}[``Local star-convexity-concavity'', informal]
    Consider any $\hz = (\ha, \hx, \hb, \hy) \in \Delta_n \times \XXX^n \times \Delta_m \times \YYY^m$,
    and let $\zps = (\aps, \xps, \bps, \yps)$ the ``proxy solution particles'' defined as in \eqref{eq:cv_proof:gencase:choice_of_proxy*}. 
    Denote $\hmu = \sum_{i=1}^n \ha_i \delta_{\hx_i}$ and $\hnu = \sum_{j=1}^m \hb_j \delta_{\hy_j}$.
    
	There exist $r, C > 0$ only dependent on $(f, \XXX, \YYY)$ such that if $V_1(\hz) \leq r$, 
    then for an appropriate choice of the partitions of unity $\varphi_I, \psi_J$,
    % (only dependent on $(f, \XXX, \YYY)$),
    \begin{equation*}
    	\tgap(\zps; \hz)
    	\geq F(\hmu, \nu^*) - F(\mu^*, \hnu) 
        % \ha^\top \gmaths b^* - (a^*)^\top \gmatsh \hb
		+ C V_\pos(\hz)
    	+ O \left( V_1(\hz)^{3/2} \right).
    \end{equation*}
\end{lemma}

% \end{document}

% !TEX root = ../main.tex
% \documentclass[../main]{subfiles}
%\begin{document}

\section{Numerical experiments} \label{sec:experiments}

In this section, we illustrate the CP-PP algorithm and its convergence properties on simple examples of applications.
As discussed in \autoref{subsec:main_res:algo}, in experiments we actually run the CP-MP algorithm since the CP-PP update cannot be computed exactly;
but based on \Autoref{lm:mp_exactparam:mp_expan,lm:mp_exactparam:pp_expan} we strongly expect the same convergence behavior for these two algorithms, as proved in the exact-parametrization case in \autoref{prop:cv_proof:exactparam:mp_cv}.
% behavior up to higher order terms.

Julia code to reproduce the experiments is publicly available online at
\url{https://github.com/guillaumew16/particle-MNE}.

\subsection{Payoff drawn from a Gaussian process} 
\label{subsec:experiments:RFF}

We start by an application of our method on a toy example where the payoff function is drawn from a Gaussian process.
More precisely we apply CP-MP on the function $f: \TT^{d_x} \times \TT^{d_y} \to \RR$ defined by
\begin{equation*}
	f(x, y) = \Re \sum_{\abs{k} \leq K} \sum_{\abs{l} \leq L} c_{k,l} 
% 	e^{2\pi i \left( \bk \cdot \bx + \bl \cdot \by \right)}
    e^{2\pi i \left( \innerprod{k}{x} + \innerprod{l}{y} \right)}
\end{equation*}
where the $c_{k,l} \in \CC$ are drawn randomly, namely $\Re[c_{k,l}], \Im[c_{k,l}]$ are drawn i.i.d.\ from the standard normal distribution.
The orders $K$ and $L$ control the smoothness of the function.
% since we do not make the coefficients decay.
Remark that the game is separable, 
i.e., $f$ can be written as a finite sum of the form $f(x, y) = \sum_{k,l} c'_{k,l} g_k(x) h_l(y)$
% for some continuous $g_k: \TT^{d_x} \to \RR, h_l: \TT^{d_y} \to \RR$ and $c'_{k,l} \in \RR$,
for some $c'_{k,l} \in \RR$ and continuous $g_k, h_l$,
since
${
	f(x, y) = \sum_{k,l} \abs{c_{k,l}} \cos \left( 
        2\pi \innerprod{k}{x} + 2\pi \innerprod{l}{y} 
        + \mathrm{arg}(c_{k,l})
    \right)
}$
and $\cos(a+b) = \cos a \cos b - \sin a \sin b$;
so we are guaranteed that a sparse MNE exists \cite[Corollary~2.10]{stein_separable_2008}.

\begin{figure}[t]
    \begin{subfigure}[t]{0.5\textwidth}
        \centering
        \includegraphics[width=0.8\textwidth]{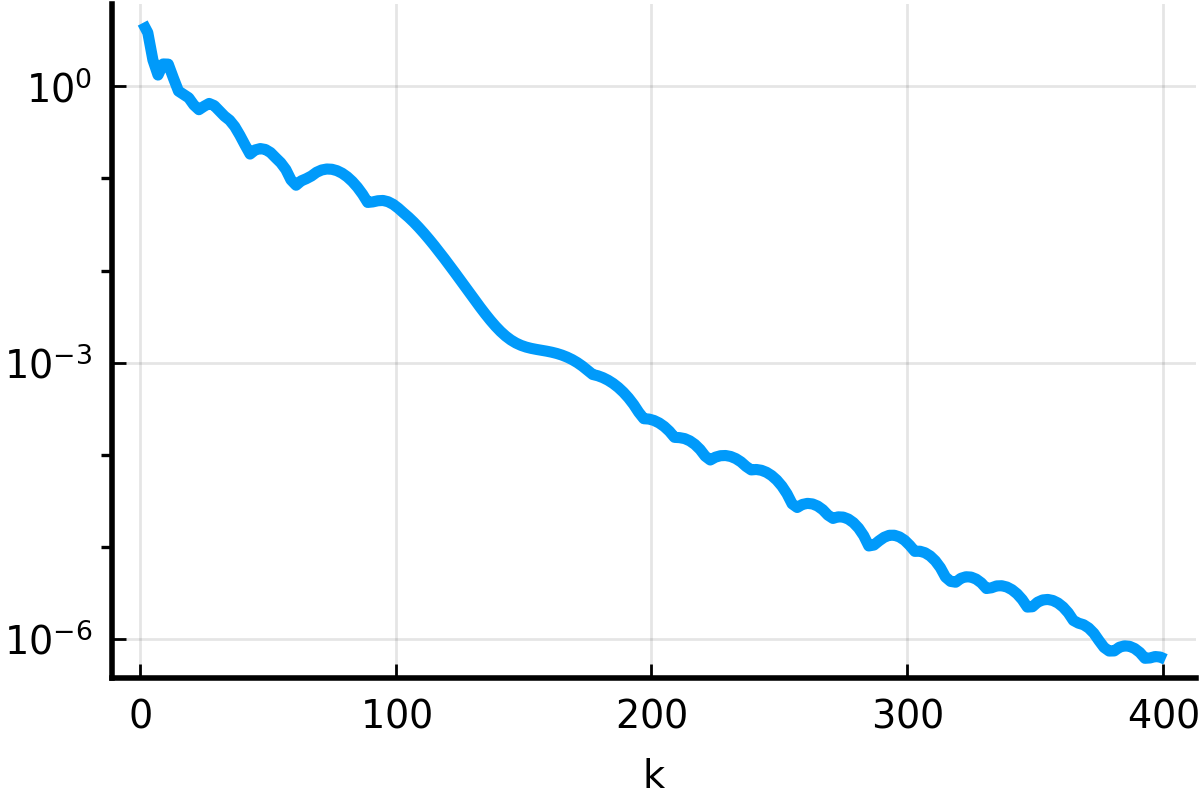}
        \caption{NI errors $\NI(\mu^k, \nu^k)$ (log-linear scale)}
        \label{fig:random_fourier_1D:NI_error}
    \end{subfigure}
    \hfill
    \begin{subfigure}[t]{0.5\textwidth}
        \centering
        \includegraphics[width=0.8\textwidth]{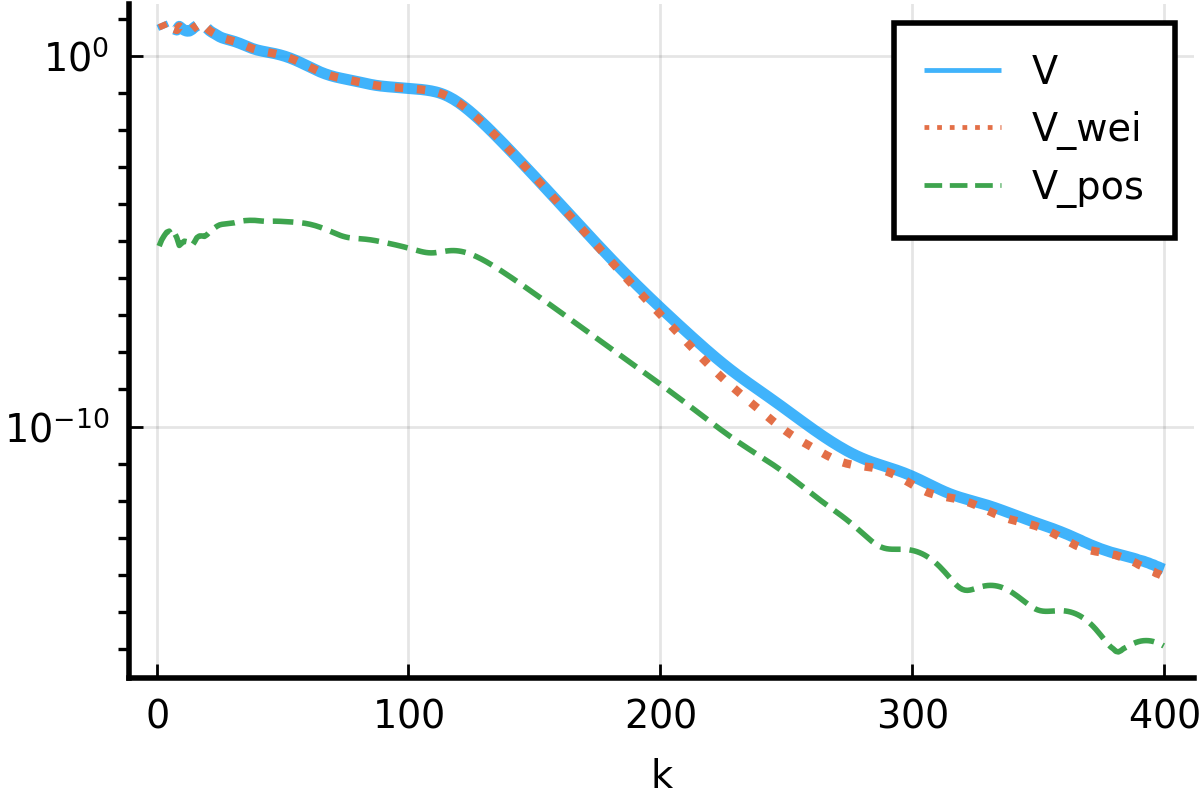}
        \caption{Lyapunov potentials $V(z^k)$ (log-linear scale)}
        \label{fig:random_fourier_1D:lyapunov_xy}
    \end{subfigure}
	\caption{Optimality metrics of CP-MP iterates for Gaussian process payoff}
    \label{fig:random_fourier_1D:optimality_metrics}
% \end{figure}
    \vspace{1em}
% \begin{figure}
    \begin{subfigure}[t]{0.37\textwidth}
        \centering
        \includegraphics[width=\textwidth]{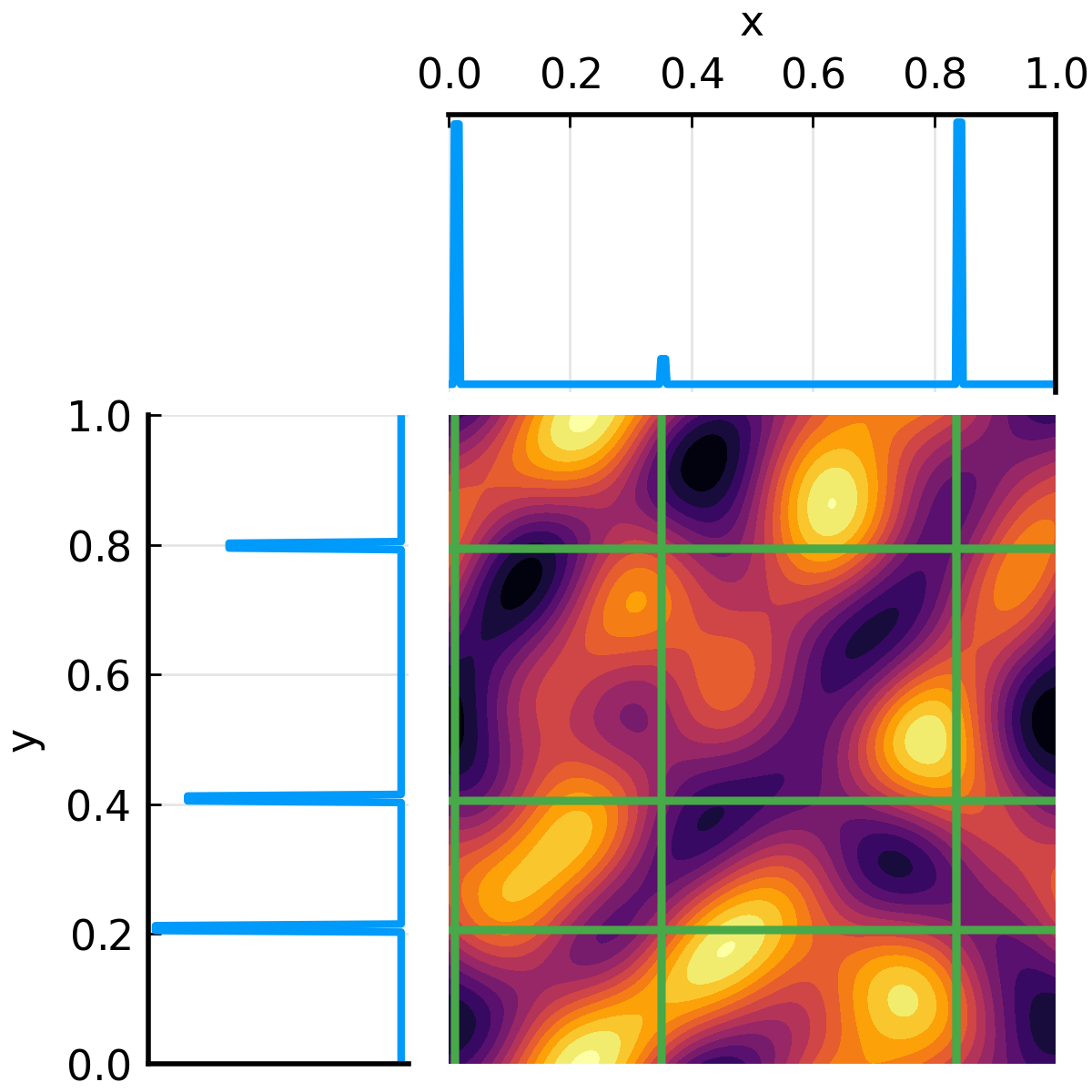}
        \caption{The measures smoothed by convolution with $h = \ind_{[0, 0.01]}$. \\
        Middle: contour plot of $f(x, y)$. \\
        Top: $(\mu^T * h)(x)$.
        Left: $(\nu^T * h)(y)$}
        \label{fig:random_fourier_1D:last_iter}
    \end{subfigure}
    \hfill
    \begin{subfigure}[t]{0.3\textwidth}
        \centering
        \includegraphics[width=\textwidth]{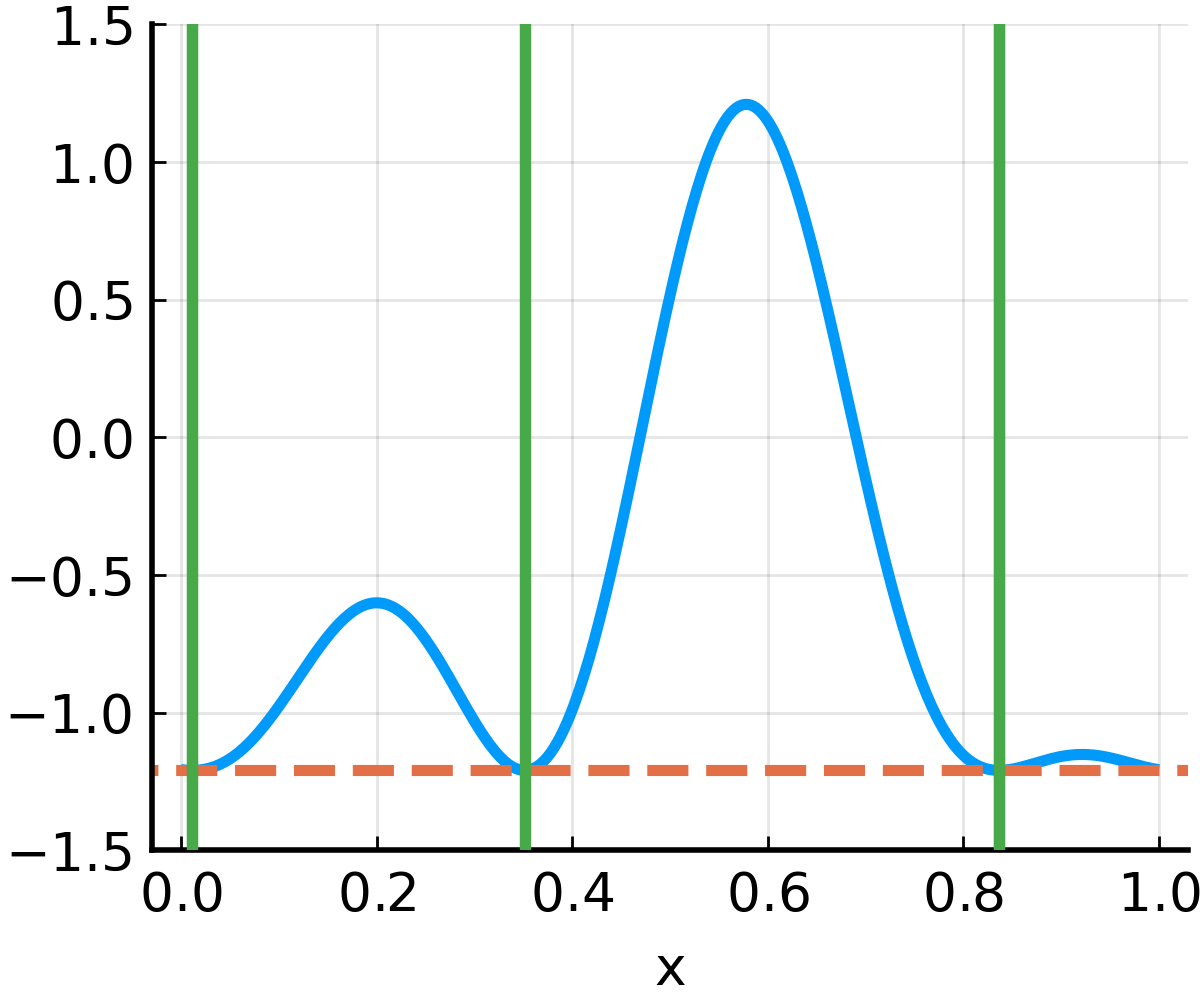}
        \caption{First variation $(F \nu^T)(x)$ at the last iterate. \\
        Orange dashed line: payoff at equilibrium $F(\mu^*, \nu^*)$}
        \label{fig:random_fourier_1D:firstvar_x}
    \end{subfigure}
    \hfill
    \begin{subfigure}[t]{0.3\textwidth}
        \centering
        \includegraphics[width=\textwidth]{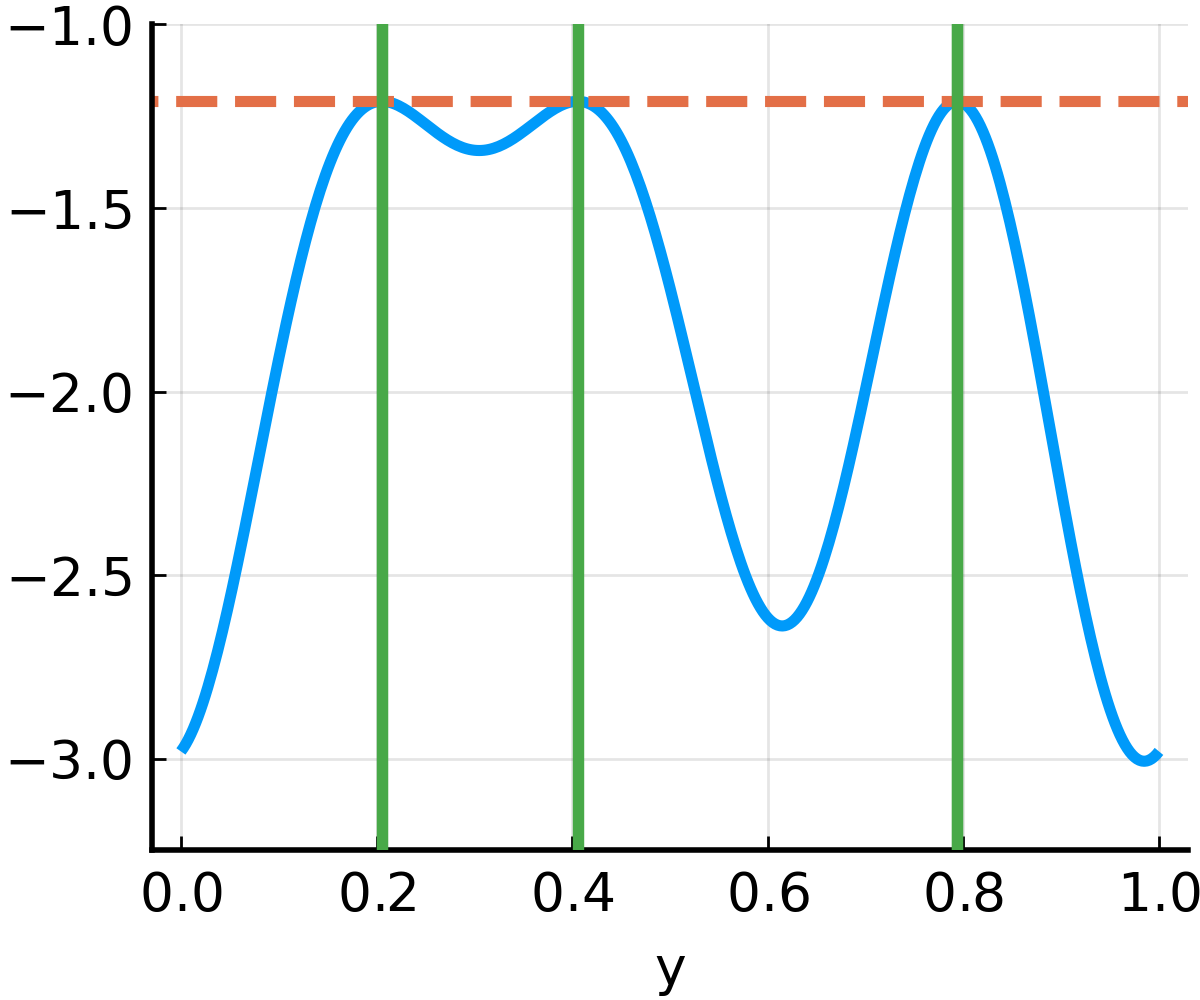}
        \caption{First variation $((\mu^T)^\top F)(y)$ at the last iterate. \\
        Orange dashed line: payoff at equilibrium $F(\mu^*, \nu^*)$}
        \label{fig:random_fourier_1D:firstvar_y}
    \end{subfigure}
    \caption{Smoothed measures and first variations at the last iterate}
    \label{fig:random_fourier_1D:firstvars_smoothedmeas_lastiter}
\end{figure}

We illustrate the behavior of the CP-MP algorithm on such a payoff function with $d_x=d_y=1$, $K=L=3$ and $n=m=15$.
A contour plot of $f$ is contained in \autoref{fig:random_fourier_1D:last_iter}.

In \autoref{fig:random_fourier_1D:optimality_metrics}, we plot the NI errors and Lyapunov potentials of the iterates for $\eta=0.04$ and $\sigma=0.001$, up to $T=400$.
Those values decrease exponentially as expected from our upper bounds.
In order to compute the NI errors, we computed $\max_\nu F(\mu^k, \nu) = \max_\nu \int_\YYY ((\mu^k)^\top F) d\nu = \max_{y \in \YYY} ((\mu^k)^\top F)(y)$ simply by discretization of $\YYY = \TT^1$.
In order to compute the Lyapunov potentials $V(z^k)$ defined in \eqref{eq:cv_proof:gencase:def_V}, we used an estimation of the $(a^*_I, x^*_I)_I, (b^*_J, y^*_J)_J$ obtained by clustering the particles of $\mu^{2T}, \nu^{2T}$.
% Interestingly the computed Lyapunov potential is not monotonous; this may be due to the fact that CP-MP was used instead of CP-PP, but we believe it is more likely due to estimation errors and inappropriate choice of $\lambda, \tau$. --> fixed, I had taken a bad scaling for choice of lambda

In \autoref{fig:random_fourier_1D:firstvars_smoothedmeas_lastiter}, we plot (a smoothed version of) the measures $(\mu^T, \nu^T)$ as well as the first variations $(F \nu^T)(x)$ and $((\mu^T)^\top F)(y)$ at the last iterate.
On all three subfigures, green lines indicate the support of $\mu^T$ and $\nu^T$.
The iterates converge to a sparse measure (here $n^* = m^* = 3$), as expected.
The first variations visibly satisfy the inequalities \eqref{eq:main_res:charact_MNE_firstvars_ineqs}, which characterize the MNE, as well as the non-degeneracy Assumptions~\ref{assum:5}-\ref{assum:6}.

\paragraph{Convergence of the continuous-time flow.}
Interestingly, for payoff functions of this form, we observe experimentally that the continuous-time flow corresponding to CP-MP typically also converges to the MNE. This behavior is not captured by our upper bound.
Indeed, we observe that the slopes of the lines in the log-linear plots of \autoref{fig:random_fourier_1D:optimality_metrics} scale as $\eta$ instead of $\eta^2$.
Moreover, experimentally, the explicit time-discretization CP-MDA (i.e.\ \autoref{alg:main_res:ppm} with $L=1$) also converges exponentially to the solution in this setting.

This phenomenon is specific to CP-MP, as it does not arise for Mirror Prox in finite games.%
\footnote{For finite games with a unique MNE, Mirror Descent-Ascent diverges, unless the MNE consists of two Dirac deltas (i.e.\ there exists a pure-strategy Nash equilibrium) \cite{bailey_multiplicative_2018}.
Moreover, in all our experiments with random payoff matrices, we observed that the convergence rate of Mirror Prox scaled as $\eta^2$.}
% In order to understand its origin, it may be helpful to consider the exact-parametrization case, as the aforementioned experimental observations also hold there. 
% 
% We emphasize that it is not always the case that the continuous-time flow of CP-MP converges. Experimentally, it does not converge for the examples described in the next two subsections, nor for the synthetic example described below.
We emphasize that it is not always the case that the continuous-time flow of CP-MP converges, as shown experimentally for the synthetic example below.

The mechanism behind this phenomenon is explained in depth in the follow-up work \cite{wang2023local} (subsequent to the completion of this work), where in particular the conditions for convergence of the continuous-time flow are described precisely, in the exact-parametrization case.

\begin{example}[The continuous-time flow may not converge] \label{example:experiments:synthetic}
    Take $d_x=d_y=1$ and
    % \begin{equation*}
    %     c_{2,0} = -i,~~~~ 
    %     c_{0,2} = -i,~~~~ 
    %     c_{1,1} = 2,~~~~
    %     c_{k,l} = 0 ~~\text{otherwise},
    % \end{equation*}
    $c_{2,0} = -i$, $c_{0,2} = -i$, $c_{1,1} = 2$, $c_{k,l} = 0$ otherwise, i.e.,
    \begin{equation*}
        % f(\frac{x}{2\pi}, \frac{y}{2\pi}) = \sin(2x) + \sin(2y) + 2 \cos(x+y),
        % f(x, y) = \sin(2\cdot 2\pi x) ~+~ \sin(2\cdot 2\pi y) ~+~ 2 \cos(2\pi x + 2\pi y).
        f(x, y) = \sin(4\pi x) ~+~ \sin(4\pi y) ~+~ 2 \cos(2\pi x + 2\pi y).
    \end{equation*}
    For this payoff function, the MNE is unique and equal to 
    $(\mu^*, \nu^*) = \left( 
        \frac{1}{2} \delta_{\frac{3}{8}} + \frac{1}{2} \delta_{\frac{7}{8}},\frac{1}{2} \delta_{\frac{1}{8}} + \frac{1}{2} \delta_{\frac{5}{8}}
    \right)$.
    Indeed,
    \begin{enumerate}
        \item One can check that
        $\int_{[0,1]} f(x,y) d\mu^*(x) = \sin(4\pi y)-1$ and $\int_{[0,1]} f(x,y) d\nu^*(y) = \sin(4\pi x)+1$.
        So since $-1 \leq \sin \leq 1$, this $(\mu^*, \nu^*)$ is a MNE, and the value at optimum is $\rho=0$.
        % Using that $-1 \leq \sin \leq 1$ one can check that this $(\mu^*, \nu^*)$ is a MNE, and so the value at optimum is $\rho=0$.
        \item Suppose by contradiction that there exists a MNE $(\mu', \nu')$ such that $\EE_{\nu'}[\sin (4\pi y) -1] <0$, and pose $x_+ = \frac{3}{8}, x_- = \frac{7}{8}$.
        Using that 
        $\cos(2\pi x_+ + 2\pi y) + \cos(2\pi x_- + 2\pi y) = 0$ for all $y$,
        we have either 
        % $\EE_{\delta_{x_+}, \nu'}[2 \cos(2\pi x + 2\pi y)] \leq 0$ or
        % $\EE_{\delta_{x_-}, \nu'}[2 \cos(2\pi x + 2\pi y)] \leq 0$.
        $\EE_{\nu'}[2 \cos(2\pi x^+ + 2\pi y)] \leq 0$ or
        $\EE_{\nu'}[2 \cos(2\pi x^- + 2\pi y)] \leq 0$.
        So $F(\delta_{x_+}, \nu') < 0=\rho$ or $F(\delta_{x_-}, \nu') < \rho$, contradicting optimality of $\nu'$.
        \item By the previous point and the symmetric argument for $\mu'$, any MNE $(\mu', \nu')$ must satisfy 
        $\EE_{\mu'}[\sin (4\pi x) +1] = 0$ and
        $\EE_{\nu'}[\sin (4\pi y) -1] = 0$, i.e., must be of the form
        $\mu' = a \delta_{\frac{3}{8}} + (1-a) \delta_{\frac{7}{8}}$,
        $\nu' = b \delta_{\frac{1}{8}} + (1-b) \delta_{\frac{5}{8}}$.
        By explicit calculations, one can show that necessarily $a=b=\frac{1}{2}$.
    \end{enumerate}
    On the other hand, 
    experimentally we observe that CP-MDA does not converge, while CP-MP converges with an exponential rate
    that scales as $\eta^2$.
\end{example}

\subsection{Max-\texorpdfstring{$\FFF_1$}{F1}-margin classification with two-layer neural networks}

A well-known machine learning task which uses the min-max framework is max-margin classification. In particular when using a two-layer neural network as the classifier, the training task is exactly of the form \eqref{eq:intro:SP}.
Indeed, a two-layer network with non-decreasing positive-homogeneous activation $\sigma$ (without bias terms) can be represented as a signed measure $\nu_\pm$ on the space of normalized hidden neurons 
$\Theta = \bbS^{d-1} = \left\lbrace \theta \in \RR^d; \norm{\theta}_2 = 1 \right\rbrace$ 
via
\begin{equation*}
    \mathrm{NN}(x; \nu_\pm) = \int_\Theta \sigma(\theta^\top x) d\nu_\pm(\theta),
\end{equation*}
or equivalently as a non-negative measure $\nu$ on the space $\Theta_+ \sqcup \Theta_-$, the disjoint union of two copies of $\Theta$,
% \cite[Appendix~A]{chizat_sparse_2021},
via
\begin{equation*}
    \mathrm{NN}(x; \nu) = \int_{\Theta_+} \sigma(\theta_+^\top x) d\nu(\theta_+) - \int_{\Theta_-} \sigma(\theta_-^\top x) d\nu(\theta_-).
\end{equation*}
Two-layer networks with bias terms can be represented in the same way, by appending a constant component $1$ to the input vector $x$ and taking $\Theta = \bbS^d$ instead of $\bbS^{d-1}$.
One can define the $\FFF_1$ norm of a function $f:\mathbb{R}^d\to \mathbb{R}$ as the infimum of $\nu(\Theta_+ \sqcup \Theta_-)$ over all $\nu$ such that $f=\mathrm{NN}(\cdot; \nu)$. Balls for this norm admit advantageous estimation/approximation trade-offs in a supervised learning task 
\cite{bach_breaking_2017}.

Consider a supervised classification task with covariates $x \in \RR^d$ and labels $y \in \{-1,1\}$.
Given $N$ observations $(x_i, y_i)_{1 \leq i \leq N}$, the max-$\FFF_1$-margin classification task consists in finding $\nu$ that maximizes the following problem
\begin{align*}
    & \max_{\substack{
        \nu \in \MMM_+(\Theta_+ \sqcup \Theta_-) \\
        \nu(\Theta_+ \sqcup \Theta_-) = 1
    }}~
    \min_{1 \leq i \leq N}~
    y_i ~ \mathrm{NN}(x_i; \nu) \\
    \equiv~ & \max_{\nu \in \PPP(\Theta_+ \sqcup \Theta_-)} ~ \min_{a \in \PPP([N])}~
    \sum_{i=1}^N a_i y_i \left( \int_{\Theta_+} \sigma(\theta_+^\top x_i) d\nu(\theta_+) - \int_{\Theta_-} \sigma(\theta_-^\top x_i) d\nu(\theta_-) \right).
\end{align*}
This problem is indeed
an instance of \eqref{eq:intro:SP} when we set
$\XXX = [N]$,
$\YYY = \Theta_+ \sqcup \Theta_-$,
and
$f(i, \theta) = \begin{cases}
    y_i \sigma(\theta^\top x_i) &\text{if}~ \theta \in \Theta_+ \\
    -y_i \sigma(\theta^\top x_i) &\text{if}~ \theta \in \Theta_-.
\end{cases}$
% $f(i, \theta) = y_i \sigma(\theta^\top x_i)$ if $\theta \in \Theta_+$ and $-y_i \sigma(\theta^\top x_i)$ if $\theta \in \Theta_-$.
% on which we may thus apply 
% a semi-discrete version of 
% the CP-PP algorithm.

\paragraph{Numerical results.}
As detailed above, the max-$\FFF_1$-margin classification problem can be written as
\begin{equation*}
    \max_{\nu \in \PPP(\Theta_+ \sqcup \Theta_-)}~
    \min_{a \in \PPP([N])}~
    \sum_{i=1}^N \int_{\Theta_+ \sqcup \Theta_-} a_i f(i,\theta) d\nu(\theta).
\end{equation*}
It is straightforward to adapt the CP-MP algorithm to this setting.
Namely, choose $m' = 2m$ a number of neurons, reparametrize by 
% $\nu = \sum_{j=1}^m b_j \delta_{\theta_j} + \sum_{j=m+1}^{2m} b_j \delta_{\theta_j}$ 
$\nu = \sum_{j=1}^{2m} b_j \delta_{\theta_j}$ for $\theta_1, ..., \theta_m \in \Theta_+$ and $\theta_{m+1}, ..., \theta_{2m} \in \Theta_-$, 
and consider the reparametrized problem
\begin{equation*}
    \min_{a \in \Delta_N}~ \max_{\substack{b \in \Delta_{2m} \\ \theta \in (\bbS^{d-1})^{2m}}}~
    \left\lbrace
        \sum_{i=1}^N \sum_{j=1}^{2m} a_i b_j \cdot 
        \left( \bmone_{[j \leq m]} - \bmone_{[j > m]} \right) \cdot
        y_i \sigma(\theta_j^\top x_i)
        ~\eqqcolon~ F_{n,m}(a, (b, \theta))
    \right\rbrace.
\end{equation*}
We can then apply \autoref{alg:main_res:ppm} with $y=\theta$ and with $x_i$ kept constant equal to $i$ for all $i$.

\begin{figure}[t!]
    \begin{subfigure}[t]{0.3\textwidth}
        \centering
        \includegraphics[width=\textwidth]{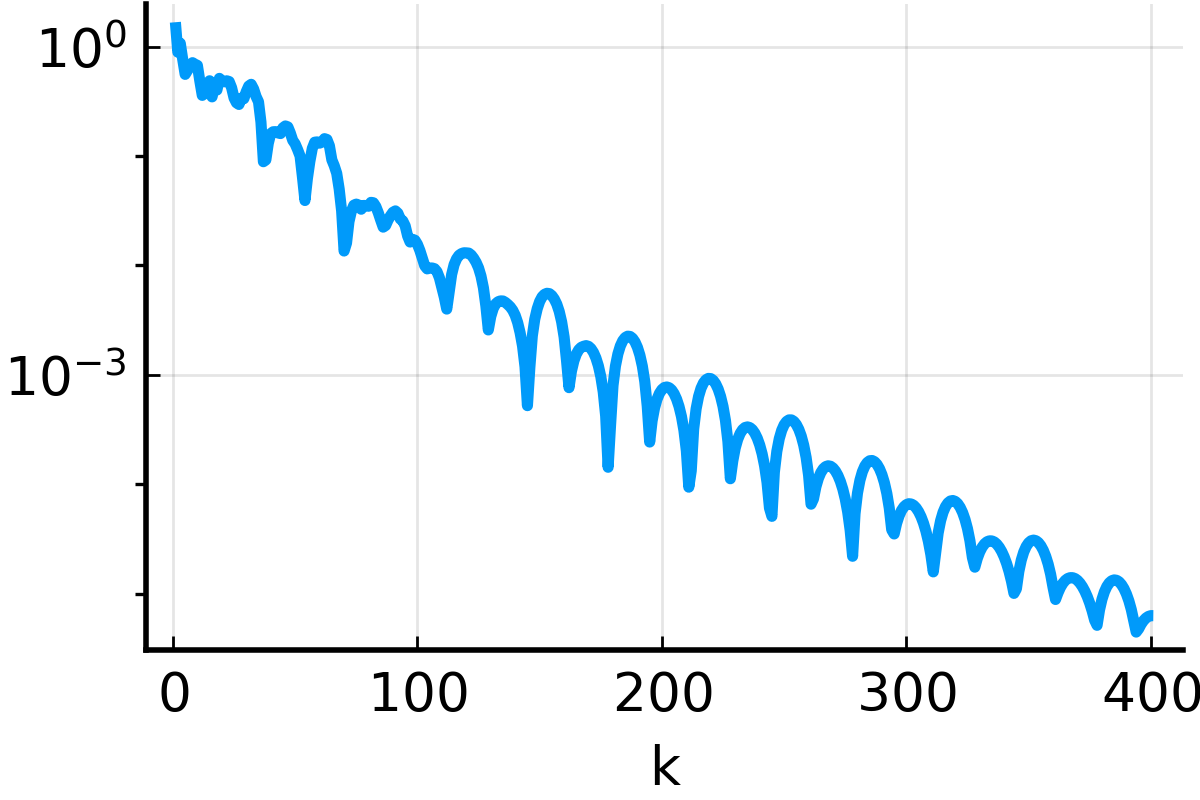}
        \caption{NI errors $\NI(a^k, \nu^k)$ (log-linear scale)}
        \label{fig:maxF1margin:NI_error}
    \end{subfigure}
    \hfill
    \begin{subfigure}[t]{0.35\textwidth}
        \centering
        \includegraphics[width=\textwidth]{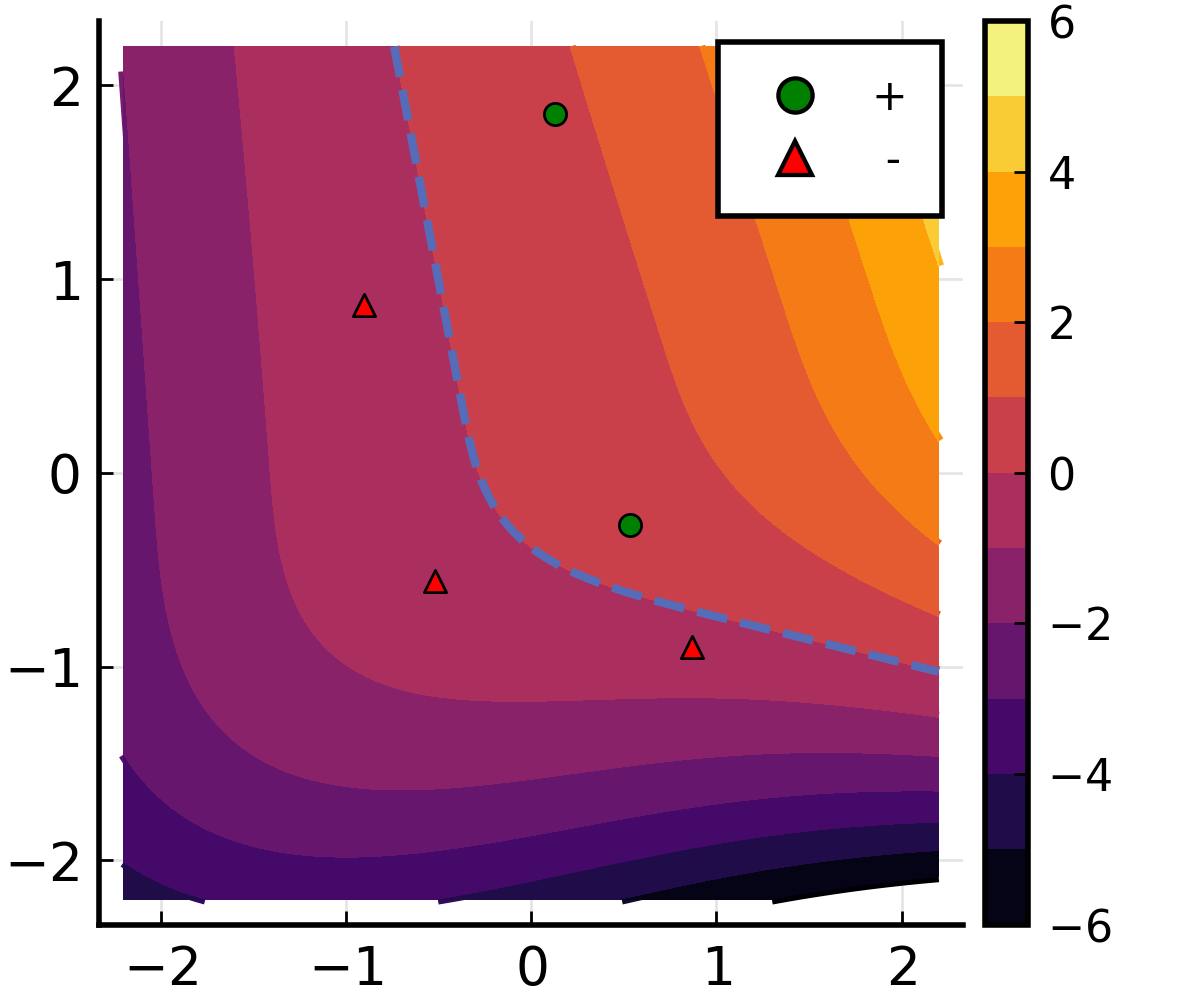}
        \caption{Soft decision regions (logits) at the last iteration.
        The dashed blue line indicates the decision boundary.}
        \label{fig:maxF1margin:contour_soft_lastiter}
    \end{subfigure}
    \hfill
    \begin{subfigure}[t]{0.3\textwidth}
        \centering
        \includegraphics[width=\textwidth]{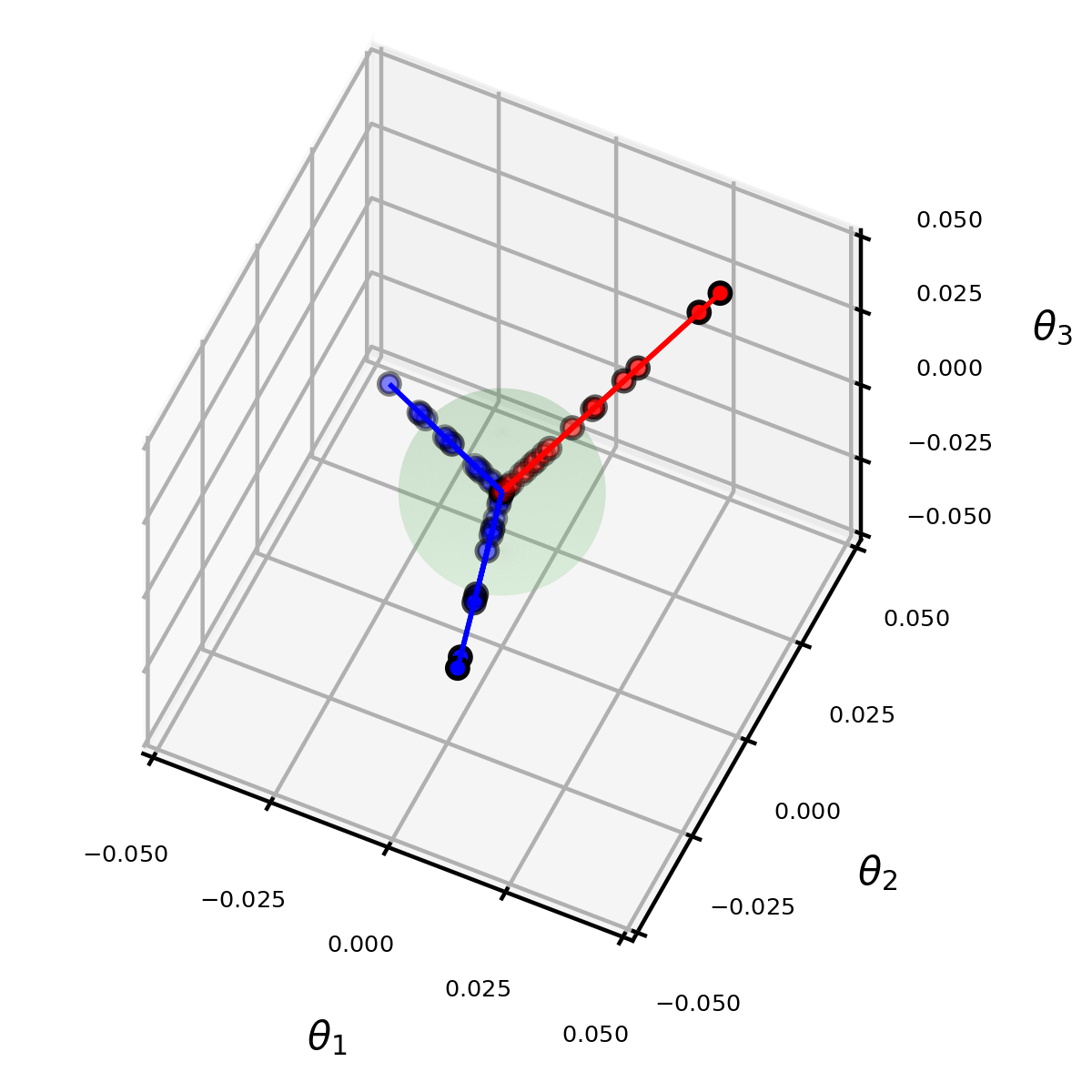}
        \caption{Neurons at the last iteration}
        \label{fig:maxF1margin:neurons_lastiter}
    \end{subfigure}
	\caption{Results for the max-$\FFF_1$-margin classification experiment}
    \label{fig:maxF1margin}
% \end{figure}
% 
% \begin{figure}[t]
    \begin{subfigure}[t]{0.3\textwidth}
        \centering
        \includegraphics[width=\textwidth]{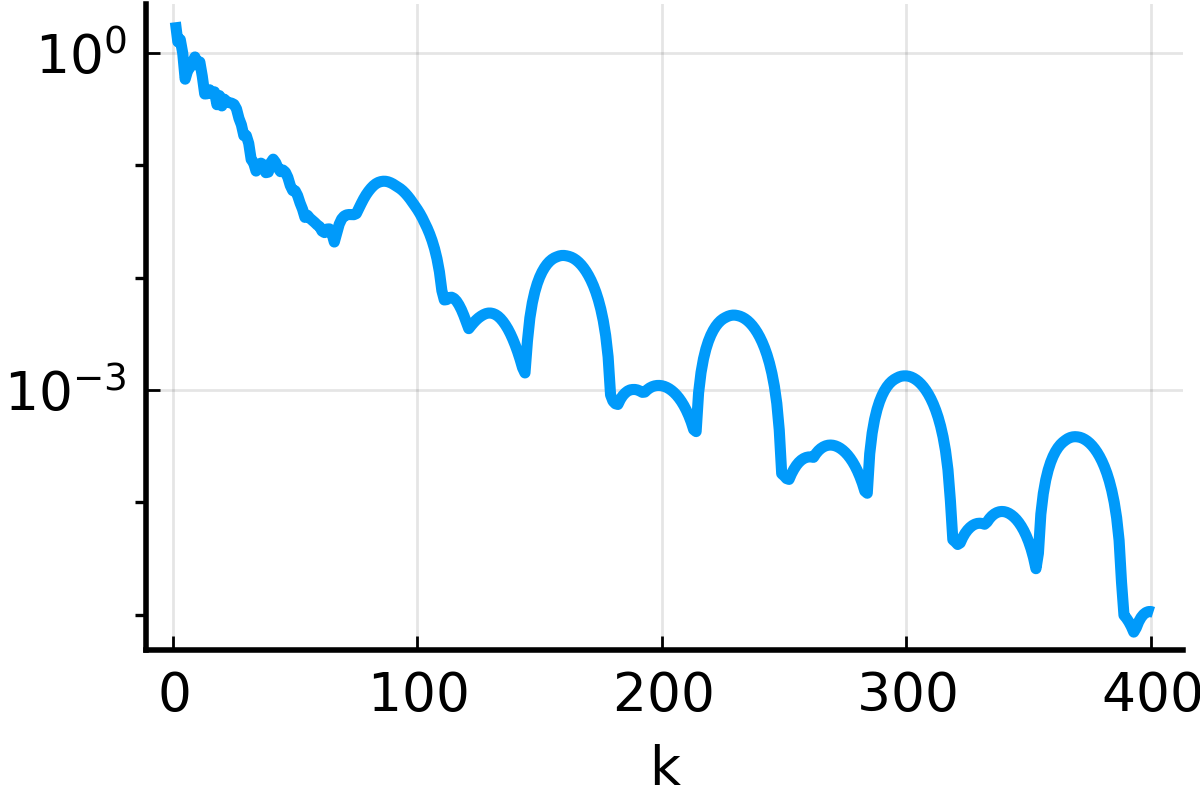}
        \caption{NI errors $\NI(\mu^k, \nu^k)$ (log-linear scale)}
        \label{fig:distrib_rob_0.2:NI_error}
    \end{subfigure}
    \hfill
    \begin{subfigure}[t]{0.35\textwidth}
        \centering
        \includegraphics[width=\textwidth]{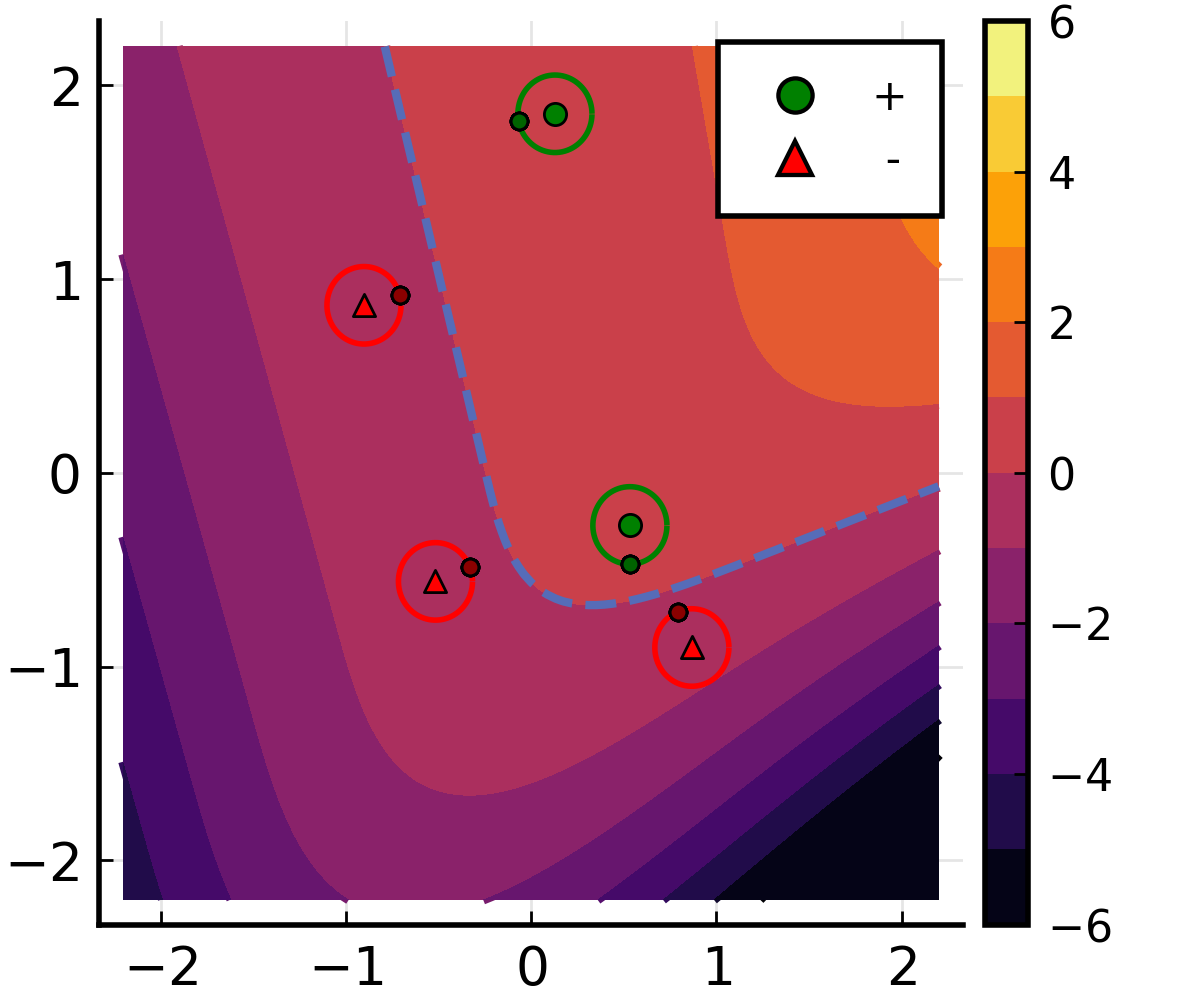}
        \caption{Decision regions at the last iteration. See text for further description.}
        % (used max.(xxx, -6) to artificially fill in the blank at the bottom right)
        \label{fig:distrib_rob_0.2:contour_soft_lastiter}
    \end{subfigure}
    \hfill
    \begin{subfigure}[t]{0.3\textwidth}
        \centering
        \includegraphics[width=\textwidth]{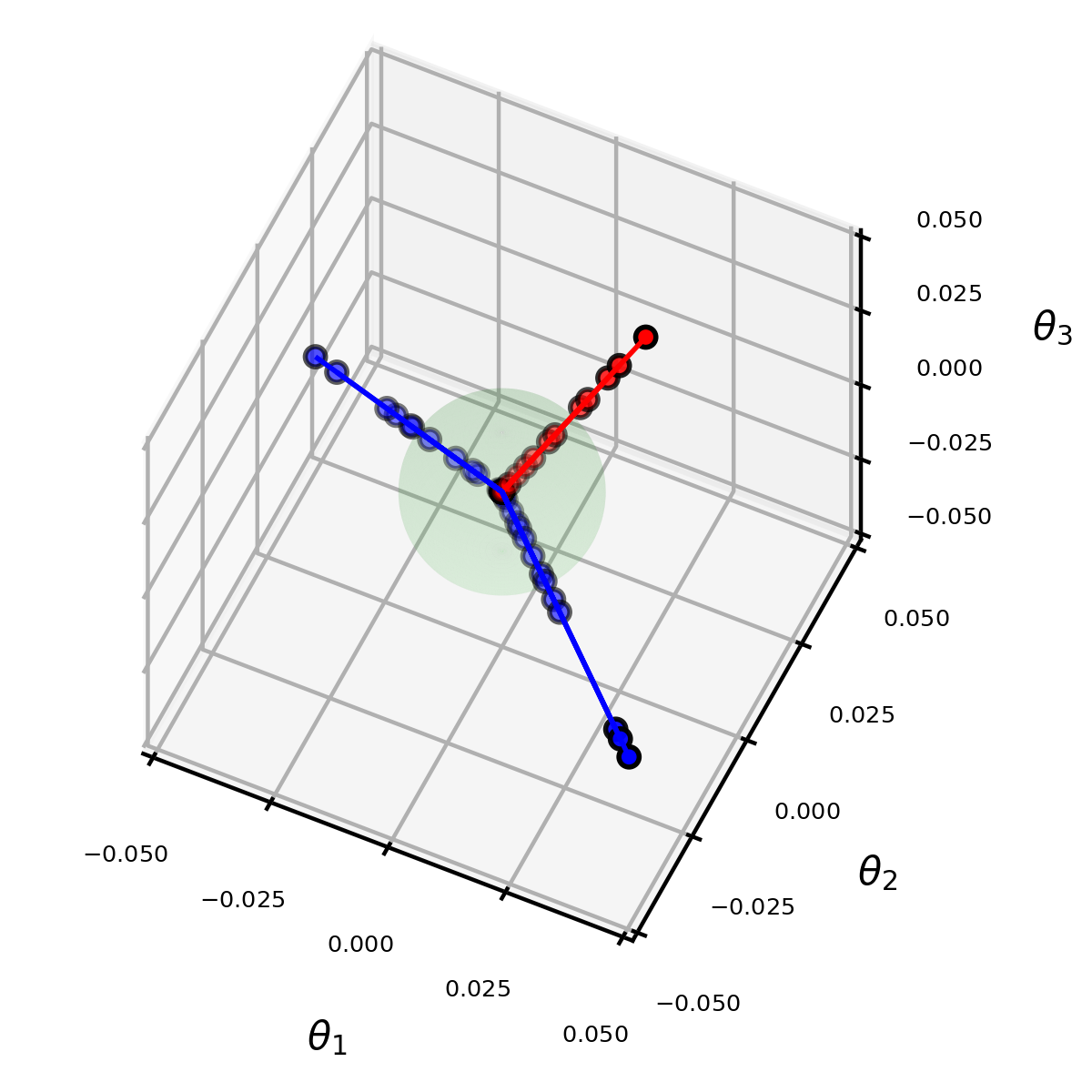}
        \caption{Neurons at the last iteration
        %\TODO{The ticks take way too much space but I didn't find a way to fix it. I don't see another way than to adjust it a posteriori by photoshopping it. Is it ok to leave like this for an arxiv version?}
        }
        \label{fig:distrib_rob_0.2:neurons_lastiter}
    \end{subfigure}
	\caption{Results for the distributionally-robust classification experiment with $r=0.2$}
    \label{fig:distrib_rob_0.2} 
% \end{figure}
% https://tex.stackexchange.com/a/290249
% \begin{figure}[t]
    \begin{subfigure}[t]{0.3\textwidth}
        \centering
        \includegraphics[width=\textwidth]{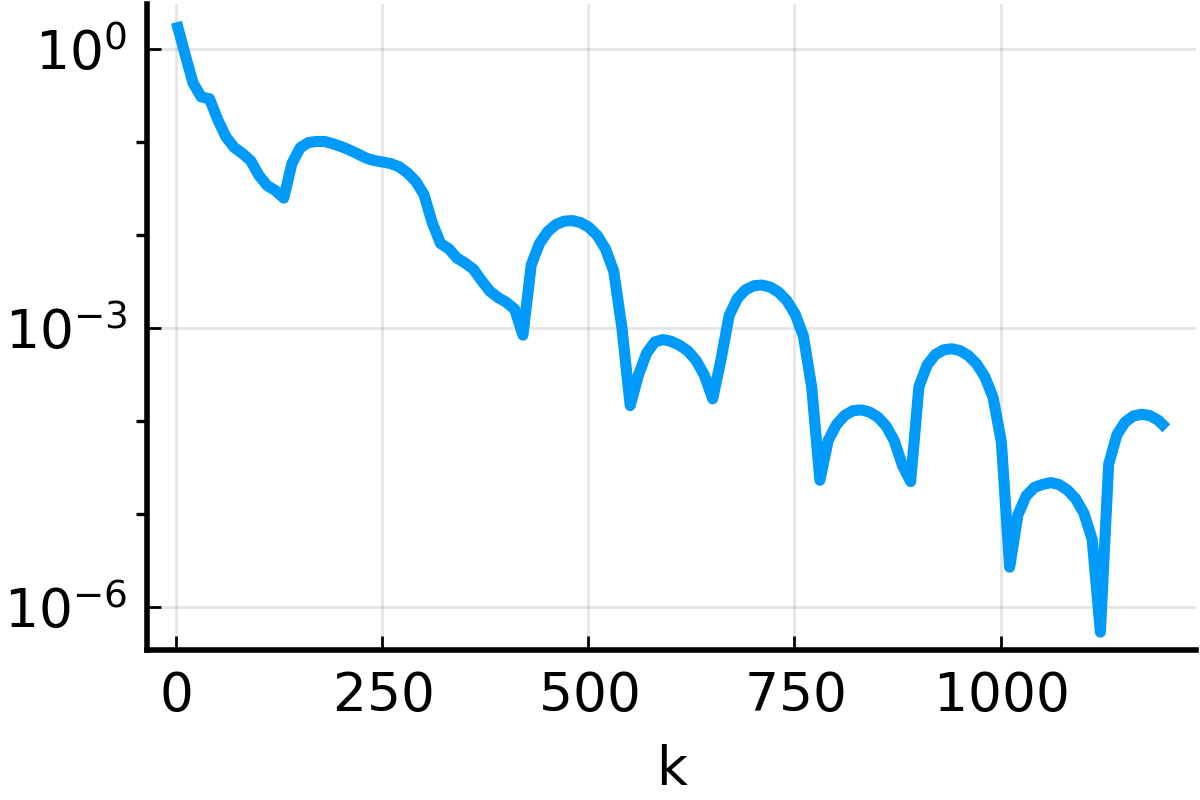}
        \caption{NI errors $\NI(\mu^k, \nu^k)$ (log-linear scale)}
        \label{fig:distrib_rob_0.3:NI_error}
    \end{subfigure}
    \hfill
    \begin{subfigure}[t]{0.35\textwidth}
        \centering
        \includegraphics[width=\textwidth]{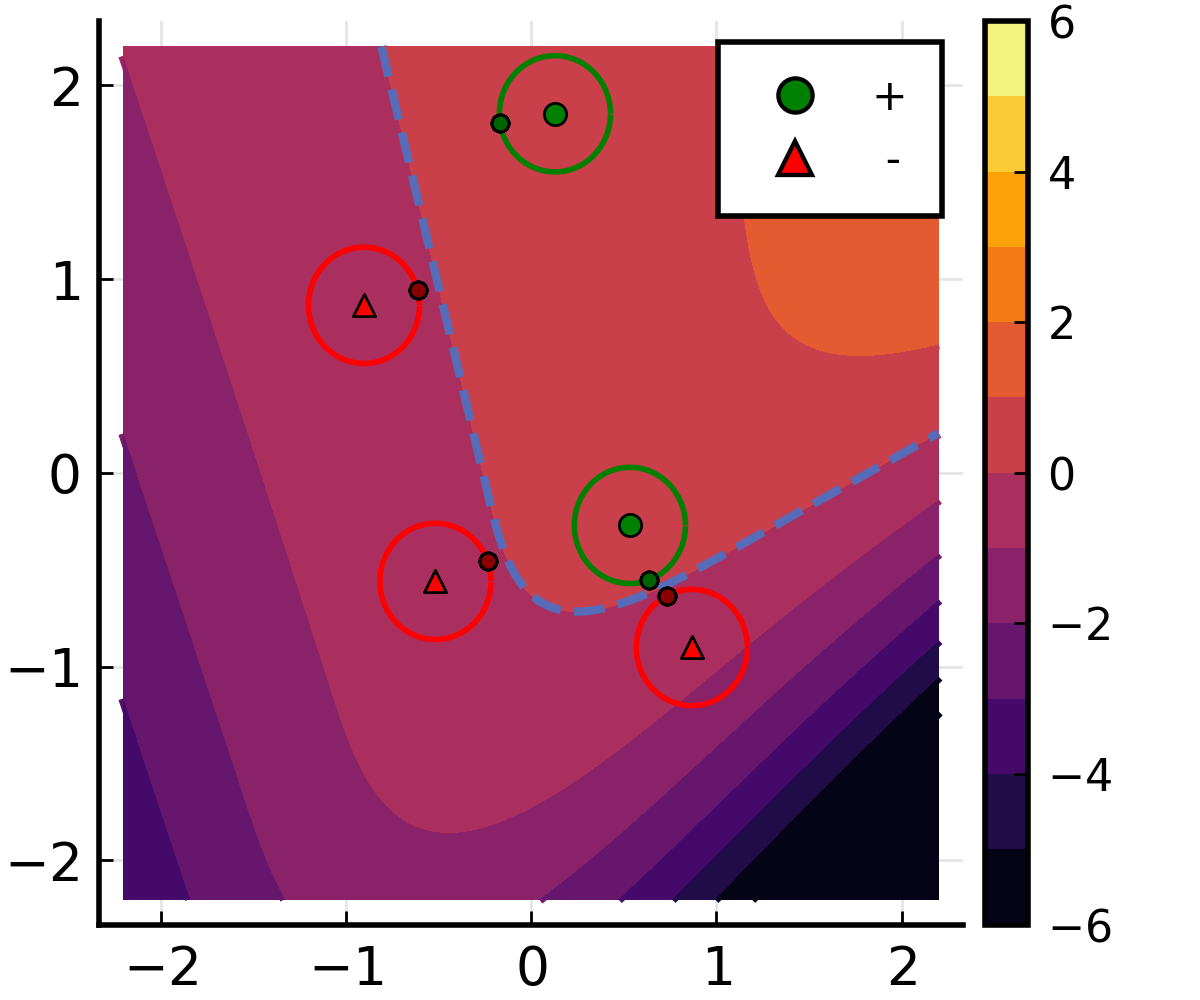}
        \caption{Decision regions at the last iteration. See text for further description.}
        \label{fig:distrib_rob_0.3:contour_soft_lastiter}
    \end{subfigure}
    \hfill
    \begin{subfigure}[t]{0.3\textwidth}
        \centering
        \includegraphics[width=\textwidth]{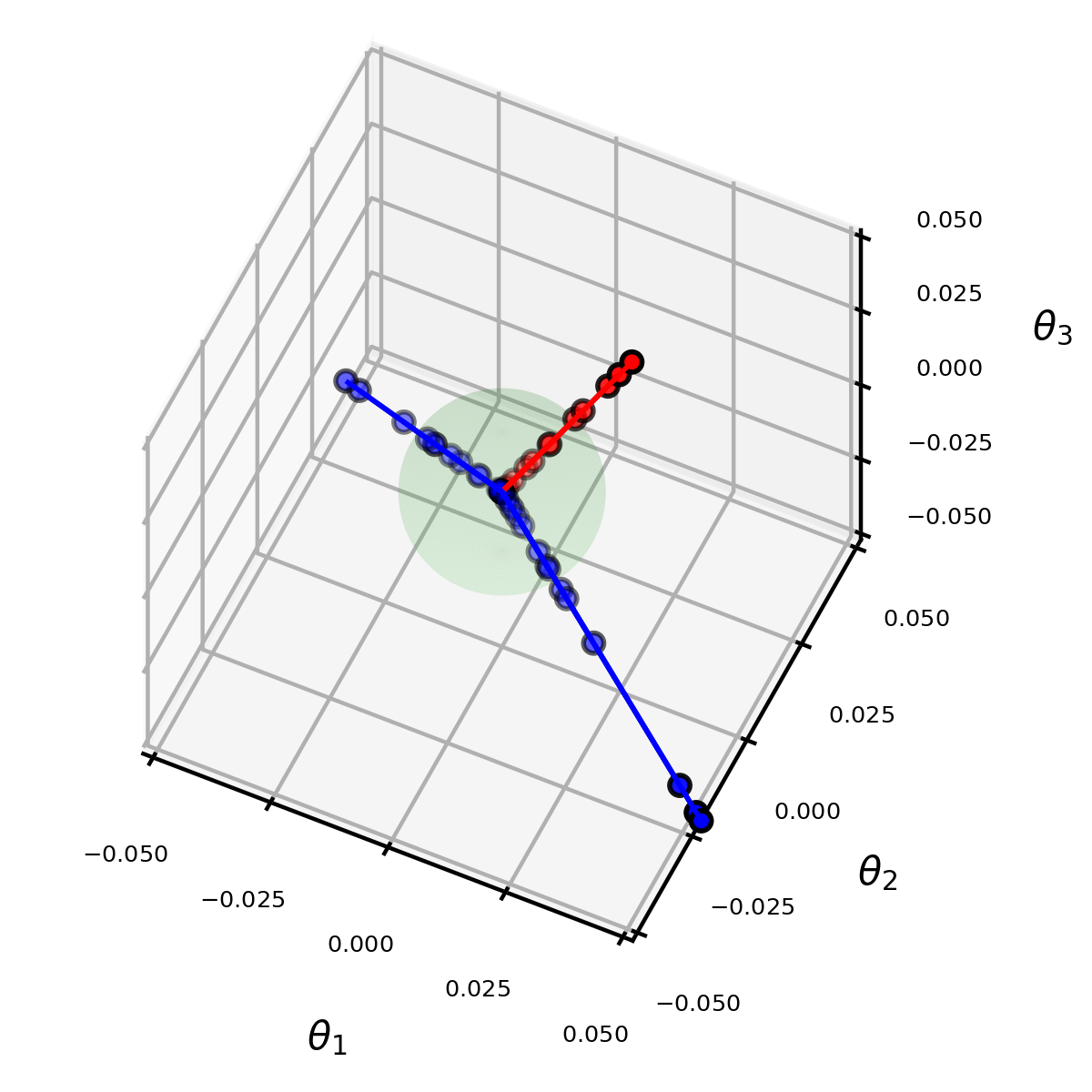}
        \caption{Neurons at the last iteration}
        \label{fig:distrib_rob_0.3:neurons_lastiter}
    \end{subfigure}
	\caption{Results for the distributionally-robust classification experiment with $r=0.3$}
    \label{fig:distrib_rob_0.3}
\end{figure}

In \autoref{fig:maxF1margin} we present the results of an experiment with $N=5$ samples, two positively labeled and three negatively labeled, and $2m=2*50$ neurons and activation $\sigma(s) = \max(0, s)^3$.
The dimensionality of the problem is $d=3$ with each sample having $1$ as the last coordinate, meaning that the last component of $\theta$ acts as a bias term.
Our analysis does not cover this case, strictly speaking, since one strategy space is discrete, and there is no guarantee that the MNE is unique;
yet the experimental results indicate a similar behavior as for 
% smooth
continuous games.
\begin{itemize}
    \item \autoref{fig:maxF1margin:NI_error} shows that the NI error, here $\NI(a^k\!, \nu^k) = \max_\theta \!\abs{\sum_{i=1}^N a_i y_i \sigma(\theta^\top \!x_i)} - \min_i y_i \mathrm{NN}(x_i; \nu^k)$, decreases exponentially to $0$.
    \item In particular the margin is non-negative at optimum so all points are classified correctly, as expected from the universality of two-layer neural networks \cite{pinkus_approximation_1999}.
    This can also be seen from the decision regions shown in \autoref{fig:maxF1margin:contour_soft_lastiter}.
    \item The solution found ($\nu^T$) turns out to be sparse, as shown by the plots in \autoref{fig:maxF1margin:neurons_lastiter}, where blue dots correspond to positively weighted neurons and red dots to negatively weighted neurons, and the distance from the origin represents the associated magnitude $b^T_j$. A green sphere of radius $\frac{1}{m}$ was added for scale.
    \item While they are not represented in the figure so as not to overload it, the variables $a$ are also of interest as a measure of each sample's importance. For example in this experiment, $a^T_i$ is close to zero for the topmost sample ($x_i \approx (0, 2)$ and $y_i=+1$), and non-zero for all other samples.
    In particular, removing the topmost sample from the dataset does not modify the learned network.
\end{itemize}

The activation function $\sigma(s) = \max(0, s)^3$ chosen for this experiment has locally Lipschitz-continuous second derivative, so our results' smoothness assumption on the payoff is verified. 
Interestingly, when using the ReLU activation $\sigma(s) = \max(0, s)$ for the same toy dataset, we observe that the NI error first decreases at an exponential rate and then oscillates around a value of about $10^{-3}$, even for large $m$. % The max I tried is m=200 
For $\sigma(s) = \max(0, s)^2$, in all our experiments we observed that the NI error vanishes exponentially.

\subsection{Distributionally-robust classification with two-layer neural networks} \label{subsec:experiments:distrib_robust}

Consider again a supervised classification task with covariates $x \in \RR^d$ and labels $y \in \{-1,1\}$.
Consider a dataset of $N$ observations $(\hx_k, \hy_k)_{1 \leq k \leq N}$ and let $\hmu = \frac{1}{N} \sum_{k=1}^N \delta_{(\hx_k, \hy_k)}$ the empirical distribution.
Let $W_\infty$ denote the $L^\infty$-Wasserstein distance on $\PPP(\RR^d \times \{-1,1\})$, defined by
\begin{align*}
    % W_\infty(\mu, \mu') 
    % = \inf_{\gamma \in \Pi(\mu, \mu')} \max_{((x,y), (x',y')) \in \support(\gamma)}
    % = \begin{cases}
    %     \norm{x-x'} & ~\text{if}~ y=y' \\
    %     +\infty & ~\text{otherwise}
    % \end{cases}
    % 
    \,W_\infty(\mu, \mu') 
    &=\!\! \inf_{\gamma \in \Pi(\mu, \mu')} \max_{\substack{((x,y), (x',y')) \\ \in \support(\gamma)}}
    d((x, y), (x', y')) 
	& &\text{\small where} &
    d((x, y), (x', y')) &=\! \begin{cases}
        \norm{x-x'}_2  ~\text{\small if $y=y'$} \\
        +\infty \qquad \text{\small otherwise}
    \end{cases}
\end{align*}
where $\Pi(\mu, \mu')$ is the set of couplings of $\mu$ and $\mu'$.
Fix a ``robustness level'' $r>0$.
In the spirit of \cite{esfahani_data-driven_2017},
the distributionally-robust classification task with respect to $W_\infty$, using two-layer neural networks $\mathrm{NN}(\cdot; \nu)$, is to find $\nu$ that maximizes
\begin{align*}
    & \max_{\substack{
        \nu \in \MMM_+(\Theta_+ \sqcup \Theta_-) \\
        \nu(\Theta_+ \sqcup \Theta_-) = 1
    }}~
    \min_{\substack{
        \mu \in \PPP(\RR^d \times \{-1,1\}) \\
        W_\infty(\mu, \hmu) \leq r
    }}~
    \int_{\RR^d \times \{-1,1\}} y~ \mathrm{NN}(x; \nu) d\mu(x, y) \\
    \equiv ~
    & \max_{\nu \in \PPP(\Theta_+ \sqcup \Theta_-)}~
    % \min_{\mu \in \PPP\left( \support(\hmu) + B_r \right)}~
    \min_{\substack{
        \mu \in \PPP(\RR^d \times \{-1,1\}) \\
        W_\infty(\mu, \hmu) \leq r
    }}~
    \int_{\RR^d \times \{-1,1\}} \int_{\Theta_+ \sqcup \Theta_-} f\left( (x,y), \theta \right) d\nu(\theta) d\mu(x, y)
\end{align*}
with
$f\left( (x,y), \theta \right) = \begin{cases}
    y \sigma(\theta^\top x) &~\text{if}~ \theta \in \Theta_+ \\
    -y \sigma(\theta^\top x) &~\text{if}~ \theta \in \Theta_-
\end{cases}$~
a ``payoff'' function over $\left( \RR^d \times \{-1,1\} \right) \times (\Theta_+ \sqcup \Theta_-)$.
More concretely, since $\hmu = \frac{1}{N} \sum_{k=1}^N \delta_{\hx_k, \hy_k}$, then $W_\infty(\mu, \hmu) \leq r$ means that 
\begin{equation*}
    \support(\mu) 
    % \subset \support(\hmu) + r B_d
    \subset \left\lbrace
        (x, y) ;~
        \exists k,
        d\left( (x, y), (\hx_k, \hy_k) \right) \leq r
    \right\rbrace
    =
    \bigcup_{k \in [N]} \left( \hx_k + r \BB \right) \times \{ \hy_k \}
\end{equation*}
where $\BB$ denotes the unit Euclidean ball.
In the language of adversarial robustness, the inner minimization means that we train the model $\mathrm{NN}(\cdot; \nu)$ to correctly classify potential adversarial examples that are within a distance of $r$ from an instance present in the dataset.

\paragraph{Numerical results.}

We showed how the task of distributionally-robust classification can be rewritten as
\begin{equation*}
    \max_{\nu \in \PPP(\Theta_+ \sqcup \Theta_-)}~
    \min_{\mu \in \PPP\left( 
        \bigcup_{k \in [N]} \left( \hx_k + r \BB \right) \times \{ \hy_k \}
    \right)}~
    \int_{\RR^d \times \{-1,1\}} \int_{\Theta_+ \sqcup \Theta_-} f\left( (x,y), \theta \right) d\nu(\theta) d\mu(x, y).
\end{equation*}
Let us adapt the CP-MP algorithm to this setting.
\begin{itemize}
    \item For the classifier ($\nu$), similarly to the previous example,
    % of max-$\FFF_1$-margin classification, 
    choose $m'=2m$ a number of neurons and 
    let $\nu = \sum_{j=1}^{2m} b_j \delta_{\theta_j}$ with 
    $\theta_1,...,\theta_m \in \Theta_+$ and 
    $\theta_1,...,\theta_m \in \Theta_-$.
    \item For the adversary ($\mu$), choose $n' = Nn$ a number of particles --- $n$ per sample ---, and 
    let $\mu = \sum_{k=1}^N \sum_{i=1}^n a_{ki} \delta_{\left( \hx_k+u_{ki}, \hy_k \right)}$
    with $\norm{u_{ki}}_2 \leq r$.
    To deal with the constraint on the $u_{ki}$'s, we project those variables back to $r \BB$ after each gradient step.
\end{itemize}
We obtain the reparametrized problem
\begin{equation*}
    \min_{\substack{a \in \Delta_{Nn} \\ u \in (r \BB)^{Nn}}} \,
    \max_{\substack{b \in \Delta_{2m} \\ \theta \in (\bbS^{d-1})^{2m} }} \!
    \left\lbrace
        \sum_{k=1}^N \sum_{i=1}^n \sum_{j=1}^{2m} a_{ki} b_j \cdot
        \left( \bmone_{[j \leq m]} \!-\! \bmone_{[j > m]} \right) \cdot
        \hy_k \sigma\left( \theta_j^\top (\hx_k + u_{ki}) \right)
        \eqqcolon F_{n,m}((a, u), (b, \theta))
    \! \right\rbrace
\end{equation*}
on which we can apply \autoref{alg:main_res:ppm} with $x=u$, $y=\theta$, modified with a projection step to ensure $u \in r \BB$.

In \autoref{fig:distrib_rob_0.2} (resp.\ \autoref{fig:distrib_rob_0.3}), 
we show the results of experiments with the same dataset and the same network architecture as in the previous subsection, with robustness level $r=0.2$ (resp.\ $r=0.3$), and using $n=10$ particles per datapoint. 
The bias terms are taken into account, i.e., each $u_{ki}$ has $0$ as the last coordinate.

\begin{itemize}
    \item In both experiments, similar to the previous subsection,
    the NI error decreases exponentially to $0$ (\autoref{fig:distrib_rob_0.2:NI_error}, \autoref{fig:distrib_rob_0.3:NI_error}).
    \item In particular the robust margin
    $\min_k \min_{x \in \hx_k+r\BB} \hy_k \mathrm{NN}(x; \nu^T)$
    is non-negative at optimum. In other words, the disks of radius $r$ around the sample covariates
    are classified correctly, 
    as can be seen in \autoref{fig:distrib_rob_0.2:contour_soft_lastiter}, \autoref{fig:distrib_rob_0.3:contour_soft_lastiter}, where the disks' boundaries are shown by green and red circles.
    \item In \autoref{fig:distrib_rob_0.2:contour_soft_lastiter}, \autoref{fig:distrib_rob_0.3:contour_soft_lastiter}
    we also represented the adversary's support points ($\hx_k+u_{ki}^T$) by slightly darker marks. We see that they are concentrated on the points of the disks that are closest to the decision boundary (dashed blue line).
    \item Just like in the max-$\FFF_1$-margin experiment of the previous subsection, the learned network ($\nu^T$) is sparse, as shown in \autoref{fig:distrib_rob_0.2:neurons_lastiter}, \autoref{fig:distrib_rob_0.3:neurons_lastiter}.
    In fact, max-$\FFF_1$-margin can be seen as an instance of distributionally-robust classification with level $r=0$, and increasing $r$ seems to perturb the learned neurons in a continuous way.
    \item Again, the variables $a$ are not represented in the figures to avoid overloading them. In both experiments, $\sum_{i=1}^n a^T_{ki}$ is close to zero for the topmost sample ($\hx_k \approx (0, 2)$ and $\hy_k=+1$) and non-zero for all other samples, just like in the max-$\FFF_1$-margin experiment.
\end{itemize}

%\end{document}

% !TEX root = ../main.tex
% \documentclass[../main]{subfiles}
%\begin{document}

\section{Conclusion} \label{sec:ccl}

In this paper, we showed that weighted particle methods can be successfully used to compute the MNE of continuous games. Specifically, we prove local exponential convergence of Conic Particle Proximal Point (CP-PP) under non-degeneracy assumptions. This algorithm is easily implementable as a descent-ascent method on a reparametrized finite-dimensional (but nonconvex-nonconcave) objective, and corresponds to an implicit time-discretization of the Wasserstein-Fisher-Rao gradient flow.
Applied to max-margin and distributionally-robust classification, our result indicates --- and our numerical experiments confirm --- that training the classifier and the adversary simultaneously is sufficient for convergence, with no need for timescale separation nor for any reformulation as in \cite{esfahani_data-driven_2017}.

An interesting question for further research would be to relax the assumption that the step-sizes for the weight ($\eta$) and position variables ($\sigma$) are of the same order, as this would allow a direct comparison with the convergence behavior of pure Fisher-Rao or pure Wasserstein gradient methods.
Another open direction is to adapt our algorithm and analysis to the case where only noisy access to the payoff function or its derivatives is available.
Finally, it could be interesting to extend our study of distributionally-robust classification (\autoref{subsec:experiments:distrib_robust}) to regression tasks, or to classification using the logistic loss.

\section*{Acknowledgments}
We would like to thank Praneeth Netrapalli for insightful discussions.

% \end{document}

% \newpage
% \nocite{*}
\printbibliography
\addcontentsline{toc}{section}{\refname} % Add References to (ToC and) bookmarks

%%%%%%%%%%%%%%%%%%%%%%%%%%%%%%%%%%%%%%%%%%%%%%%%%%%%%%
\newpage
\appendix
%\phantomsection
%\addcontentsline{toc}{section}{APPENDIX}
% \addcontentsline{toc}{chapter}{APPENDIX}

% \newpage

\section*{Appendix}
The appendix is structured as follows.
\begin{itemize}
    \item In \autoref{sec:notations_for_proofs} we introduce notations used throughout the proofs in the appendix.
    \item In \autoref{sec:proof_welldef} we prove \autoref{lm:main_res:ppm_upd_well_defined} stating that the CP-PP update is well-defined.
    \item In \autoref{sec:growth_conds} we formalize and prove the lower growth properties discussed in \autoref{subsec:cv_proof:proofingr_growthconds}.
    
    This section, and in particular the ``steepness'' result \autoref{claim:growth_conds:steepness_of_reduced_game}, represents the crux of our analysis.
    In particular much of \autoref{sec:rel_Lya_NI} will rely on the same proof ideas, and \autoref{sec:proof_gencase} will make crucial use of the results from this section.
    \item In \autoref{sec:rel_Lya_NI} we prove 
    the bounds of 
    \autoref{prop:cv_proof:exactparam:NI_equiv_V} and \autoref{prop:cv_proof:gencase:NI_equiv_V} relating NI error and our Lyapunov potential.
    % This proof relies on the same ideas as for the lower growth properties.
    \item In \autoref{sec:proof_gencase} we present the complete convergence analysis of CP-PP for the general case, proving \autoref{thm:cv_proof:gencase:loc_exp_cv}.
    % This proof makes use of the lower growth properties.
    
    % After evaluating the variational inequality characterizing CP-PP at ``proxy solution particles'', the proof essentially consists of the same steps as in the exact-parametrization case (proof of \autoref{thm:cv_proof:exactparam:loc_exp_cv}), with several additional error terms appearing.
    It is instructive to see how the aggregated weights and positions naturally appear in the derivations, so that the steps of the proof almost perfectly match the ones for the exact-parametrization case (proof of \autoref{thm:cv_proof:exactparam:loc_exp_cv}).
    The manipulations required to deal with the additional error terms are purely technical however,
    % , which can easily be identified by using the informal rule of thumb that they are $o\left( (\eta \vee \sigma)^2 V(z^k) \right)$; 
    and they are deferred to 
\ifextended%
    the last subsection.
    % its own subsection, \autoref{subsec:proof_gencase:delayed_proofs}.
\else%
    the arXiv version of this paper as \AdditionalMaterial.
\fi
    \item \autoref{sec:aux_lemmas} collects elementary auxiliary facts used in some of the above sections.
\ifextended%
    \item \autoref{sec:tedious_calcs} contains some delayed calculatory proofs for the above sections.
\fi
    \item In \autoref{sec:mp_exactparam} we prove \autoref{prop:cv_proof:exactparam:mp_cv} stating that (in the exact-parametrization case) CP-MP has the same convergence behavior as CP-PP.
    
    The proof consists in deriving generically applicable approximate expressions for the Mirror Prox and Proximal Point updates, up to order-3 terms in the step-size.
    % , that may be of independent interest. 
    In particular we show and exploit the fact that the error terms are also proportional to the projected gradient norm.
    \item In \autoref{sec:proof_mainres} we show in detail how our main result \autoref{thm:main_res:loc_exp_cv_NI} follows from combining \autoref{prop:cv_proof:gencase:NI_equiv_V}, \autoref{prop:cv_proof:gencase:WFR_equiv_V} and \autoref{thm:cv_proof:gencase:loc_exp_cv}.
\end{itemize}

% !TEX root = ../main.tex
% \documentclass[../main]{subfiles}
% \begin{document}

\section{Notations used in the proofs} \label{sec:notations_for_proofs}
% {\color{red} LC: since the appendix is quite long, can you include un sommaire?} -> done (cf main.tex)
In this section we collect notations used throughout the proof. 
Most of them are natural, except perhaps our use of the $O(\cdot)$ notation (last paragraph).

\paragraph{Relative entropy.}

% Let $h: [0,1] \to \RR$ defined by % we will use strong convexity of $h$ restricted to $[0,1]$, but writing this is just weird
Let $h: \RR_+ \to \RR$ defined by 
$h(s) = s \log s - s + 1$. 
$h$ is convex and its Bregman divergence is
\begin{equation*}
	d_h(s, s') = \begin{cases}
	    +\infty &~\text{if}~ s'=0, s>0 \\
	    s \log \frac{s}{s'} - s + s' &~\text{otherwise}.
	\end{cases}
\end{equation*}
The Kullback Leibler (KL) divergence between $w$ and $\hw \in [0,1]^n$ is given by
$\KLdiv(w, \hw) = \sum_i d_h(w_i, w'_i)$.

%Define the divergence
%\begin{equation*}
%	D(wp, \hwp) = \KLdiv(w, \hw) + \frac{\eta}{\sigma} \sum_i w_i \norm{p_i-\hp_i}^2.
%\end{equation*}

\paragraph{Indexing.}

\begin{itemize}
	\item We use $I \in [n^*]$ resp.\ $J \in [m^*]$ to index the ``true'' particles, i.e., 
	the unique MNE $(\mu^*, \nu^*)$ is
	\begin{align*}
		\mu^* &= \sum_{I \in [n^*]} a^*_I \delta_{x^*_I}
		~~~~ (a^*_I>0)
		&
		\nu^* &= \sum_{J \in [m^*]} b^*_J \delta_{y^*_J}
		~~~~ (b^*_J>0).
	\end{align*}
	
	We use $i \in [n]$ resp.\ $j \in [m]$ to index the particles used by the algorithm.
	\item Let by convention $a^*_0 = b^*_0 = 0$. 
	In particular we can write that
	\begin{equation*}
		\forall \ola \in \Delta_{[0, n^*]},~
		\KLdiv(a^*, \ola) = \sum_{I \in [0, n^*]} a^*_I \log \frac{a^*_I}{\ola_I}.
	\end{equation*}
	Unless specified, summations over $I$ refer to $I \in [n^*]$ (excluding index $0$).
%	: $\sum_I \equiv \sum_{I \in [n^*]}$.
	
%	To be clear (although it won't really matter), 
% 	To lift any ambiguity, $\norm{\ola-a^*}_1$ refers to $\ell_1$-norm for $\Delta_{[0, n^*]}$:
% 	${
% 	    \norm{\ola-a^*}_1 = \sum_{I \in [n^*]} \abs{\ola-a^*} + \ola_0
% 	}$.
% 	    --> I changed my mind.
	To lift any ambiguity, $\norm{\ola-a^*}_1$ refers to $\ell_1$-norm for $\Delta_{[n^*]}$:
	${
	    \norm{\ola-a^*}_1 = \sum_{I \in [n^*]} \abs{\ola-a^*}
	}$,
	even when $\ola \in \Delta_{[0, n^*]}$.
% 	(Actually 
%     $0 \leq \ola_0 = \sum_{I \in [n^*]} a^*_I - \ola_I \leq \norm{\ola-a^*}_1$, so
%     $\norm{\ola-a^*}_1 \leq \sum_{I \in [0, n^*]} \abs{\ola-a^*} \leq 2\norm{\ola-a^*}_1$,
%     so it doesn't really matter.)
	\item For $\ola \in \Delta_{[0, n^*]}$, $\Delta \ola_I = \ola_I - a^*_I$,
	and for $\olx \in \XXX^{n^*}$, $\Delta \olx_I = \olx_I - x^*_I$.
%	
%	In \autoref{subsubsec:lb_max_WP} the convention will be slightly different: We will write $\Delta \ola$ for the vector in $\RR^{[n^*]}$.
%	--> nah should be clear from context
% 
% 	\item With slight abuse of notation,
% 	for $\olw = (\ola, \olb) \in \Delta_{[0, n^*]} \times \Delta_{[0, m^*]}$, 
% 	we will sometimes write
% 	\begin{align*}
% 		&\sum_{I \in [0, n^*]} \text{expr}(\olw_I)
% 		& &\text{to mean}&
% 		&\sum_{I \in [0, n^*]} \text{expr}(\ola_I)
% 		+ \sum_{J \in [0, m^*]} \text{expr}(\olb_J).
% 	\end{align*}
    \item Generically denote the joint weight resp.\ position variables by
    $w = \begin{pmatrix}
        a \\
        b
    \end{pmatrix}
    \in \Delta_n \times \Delta_m$,
    resp.\  
    $p = \begin{pmatrix}
        x \\
        y
    \end{pmatrix}
    \in \XXX^n \times \YYY^m$.
    Summations over $w_i$ will implicitly be over $[n] \sqcup [m]$, that is,
    $\sum_i f(w_i) = \sum_{i=1}^n f(a_i) + \sum_{j=1}^m f(b_j)$.
    Likewise for 
    $\olw = \begin{pmatrix}
        \ola \\
        \olb
    \end{pmatrix}
    \in \Delta_{[n^*]} \times \Delta_{[m^*]}$
    and
    $\olp = \begin{pmatrix}
        \olx \\
        \oly
    \end{pmatrix}
    \in \XXX^{n^*} \times \YYY^{m^*}$,
    for which summations will implicitly be over $[n^*] \sqcup [m^*]$ (excluding the two indices $0$).
    Finally, $\olw_0 = \ola_0 + \olb_0$.
    % 
    % This shorthand is sometimes extremely handy, but we will use it sparingly to minimize the possibility of confusion.
\end{itemize}

\paragraph{Local payoff matrices.}

\begin{itemize}
	\item Generically let, for $(\ha, \hx, \hb, \hy) \in (\Delta_n \times \XXX^n) \times (\Delta_m \times \YYY^m)$, 
	\begin{align*}
		\gmathh_{ij} &= f(\hx_i,\hy_j)
		&
		\gmaths_{iJ} &= f(\hx_i, y^*_J) \\
		\gmatsh_{Ij} &= f(x^*_I, \hy_j)
		&
		\gmatss_{IJ} &= f(x^*_I, y^*_J)
	\end{align*}
	as well as 
	$\partial_x \gmathh_{ij} = \partial_x f(\hx_i,\hy_j)$,
	$\partial_x \gmatss_{IJ} = \partial_x f(x^*_I, y^*_J)$,
	etc.,
	and similarly for $\partial_y$, $\partial^2_{xx}$, $\partial^2_{xy}$, and $\partial^2_{yy}$.
	For example, we have the Taylor expansion 
	for all $i, j, I$
	\begin{equation*}
	    \gmathh_{ij} = \gmatsh_{Ij} + (\hx_i-x^*_I)^\top \partial_x \gmatsh_{Ij} + O(\norm{\hx_i-x^*_I}^2).
	\end{equation*}
	
	\item Let
	\begin{equation*}
		\forall I \in [n^*],~ H_I = \sum_J \partial_{xx}^2 \gmatss_{IJ} b^*_J
		~~~~~~ \text{resp.} ~~~~~~
		\forall J\in [m^*],~ H_J = - \sum_I a^*_I \partial_{yy}^2 \gmatss_{IJ}
	\end{equation*}
	the local kernels; that is,
	$H_I = \partial^2_{xx} (F \nu^*)(x^*_I)$
	and $H_J = -\partial^2_{yy} ((\mu^*)^\top F)(y^*_J)$.
	
	Let $H_x$ the $n^* d_x$ by $n^* d_x$ block-diagonal matrix with blocks $(H_I)_{I \in [n^*]}$,
	and likewise let $H_y$ the block-diagonal matrix with blocks $(H_J)_{J \in [m^*]}$.
% 	\item We will typically drop transpose signs $\bullet^\top$ to lighten notation, for example
% 	\begin{align*}
% 		a M b
% 		& \quad\quad \text{should be read as} \quad\quad
% 		a^\top M b \\
% 		\text{and} \quad\quad
% 		a ~ \delta x ~ \partial_x M b
% 		& \quad\quad \text{should be read as} \quad\quad
% 		\sum_{ij} a_i (\delta x_i)^\top \partial_x M_{ij} b_j.
% 	\end{align*}
    \item We will use the usual matrix-vector product notation, so that for example $H_I = \partial_{xx}^2 \gmatss_{I\bullet} b^*$.
    In addition we introduce the following shorthands:
    For $\ola, \olx, \olb, \oly = (\Delta_{n^*} \times \XXX^{n^*}) \times (\Delta_{m^*} \times \YYY^{m^*})$, 
    \begin{itemize}
        \item Even though $\olx$ is a vector in $(\RR^{d_x})^{n^*}$ and $\ola$ is a vector in $\RR^{n^*}$, we denote by $\ola \odot \olx$ the vector in $(\RR^{d_x})^{n^*}$ such that $(\ola \odot \olx)_I = \ola_I \olx_I$.
        \item We denote by $\ola H_x$ the block-diagonal matrix such that $\left[ \ola H_x \right]_{I I'} = \bmone_{I=I'} \ola_I H_I$. Likewise,
        \item $\ola \partial_x \gmatss$ is the matrix such that $\left[ \ola \partial_x \gmatss \right]_{IJ} = \ola_I \partial_x \gmatss_{IJ}$, and
        \item $\partial_y \gmatss \olb$ is the matrix such that $\left[ \partial_y \gmatss \olb \right]_{IJ} = \olb_J \partial_y \gmatss_{IJ}$, and
        \item $\ola \partial_{xy}^2 \gmatss \olb$ is the matrix such that $\left[ \ola \partial_{xy}^2 \gmatss \olb \right]_{IJ} = \ola_I \olb_J \partial_{xy}^2 \gmatss_{IJ}$.
    \end{itemize}
%     \item \TODO{(I dropped transpose signs $\bullet^\top$ in a lot of places; put them back)}
% 	\item We will write sums of the following form symbolically as vector-matrix-vector products: \TODO{is it ok to drop the $\bullet^\top$ in this already very symbolic notation?}
% 	\begin{align*}
% 		\sum_{ij} a_i (\delta x_i)^\top \partial_x \gmat_{ij} b_j
% 		& ~~~\text{will be written as}~~~
% 		a^\top (\delta x)^\top \partial_x \gmat b \\
% 		\sum_{ij} a_i \cdot (\delta x_i)^\top \partial^2_{xy} \gmat_{ij} \delta y_j \cdot b_j
% 		& ~~~\text{will be written as}~~~
% 		a^\top (\delta x^2)^\top \partial^2_{xx} \gmat b \\
% 		\sum_{ij} a_i \cdot (\delta x_i, \delta x_i)^\top \partial^2_{xx} \gmat_{ij} \cdot b_j
% 		& ~~~\text{will be written as}~~~
% 		a^\top (\delta x)^\top \partial^2_{xy} \gmat \delta y b.
% 	\end{align*}
% 	This will considerably lighten notation and enhance clarity when we write Taylor expansions containing many such sums.
    Finally, we use $\id$ to denote the identity matrix, and its size will be clear from context.
\end{itemize}

\paragraph{Norms and dual norms on $\XXX$ and $\YYY$.}

We assume that $\XXX$ and $\YYY$ are the $d_x$- resp.\ $d_y$-dimensional tori, that is, $\XXX = \TT^{d_x} = (\RR / \ZZ)^{d_x}$ and
\begin{equation*}
    % \forall x, x' \in \XXX,~ \norm{x-x'}_\XXX = \inf_{k \in \ZZ^{d_x}} \norm{x-x'+k}_{\ell_2^{d_x}}.
    \forall x, x' \in \XXX,~ \norm{x-x'}_\XXX = \inf_{k \in \ZZ^{d_x}} \norm{x-x'+k}_2.
\end{equation*}
In particular $\XXX$ is a compact Riemannian manifold,
the tangent space at any point is isometric to $\RR^{d_x}$,
and the norm of a tangent vector is
\begin{equation*}
    % \forall v \in T_x \XXX,~ \norm{v}_x = \norm{v}_{\ell_2^{d_x}}.
    \forall v \in T_x \XXX,~ \norm{v}_x = \norm{v}_2.
\end{equation*}
The same considerations apply for $\YYY = \TT^{d_y}$.

To lighten notation, we use $\norm{\cdot}$ to denote the norm over $\XXX$ or $\YYY$ or $T_x \XXX$ or $T_y \YYY$; which one is meant in each situation will be clear from context.

\paragraph{The quantities arising from the Assumptions~\ref{assum:1}-\ref{assum:6}.}

\begin{itemize}
	\item Since $\XXX$ and $\YYY$ are compact Riemannian manifolds, let $\diamXY = \mathrm{diameter}(\XXX) \vee \mathrm{diameter}(\YYY) < \infty$.
	\item 
% 	Since $f \in \CCC^{2,1}(\XXX \times \YYY)$, % $\CCC^{2,1}$ means C2 + 1-Hoelder second-order derivatives
    Since $f$ has bounded differentials of order up to $3$,
	let $\partial_x, \partial_y$ the partial derivative operators and 
	denote the smoothness constants of $f$ as
	\begin{align*}
		\smoothnessf_0 &= \sup_{\XXX \times \YYY} f - \inf_{\XXX \times \YYY} f, &
		\smoothnessf_1 &= \sup_{\XXX \times \YYY} \norm{\partial_x f} \vee \norm{\partial_y f},  &
		\smoothnessf_2 &= \sup_{\XXX \times \YYY} \norm{\partial^2_{xx} f} \vee \norm{\partial^2_{xy} f} \vee \norm{\partial^2_{yy} f}
	\end{align*}
% 	and $\smoothnessf_3 = \sup_{\XXX \times \YYY} \norm{\partial^3_{xxx} f} \vee \norm{\partial^3_{xxy} f} \vee \norm{\partial^3_{xyy} f} \vee \norm{\partial^3_{yyy} f}$,
    and $\smoothnessf_3$ such that $\partial^2_{xx} f$, $\partial^2_{xy} f$, $\partial^2_{xx} f$ are $\smoothnessf_3$-Lipschitz-continuous,
	and $\smoothnessf = \smoothnessf_0 \vee \smoothnessf_1 \vee \smoothnessf_2 \vee \smoothnessf_3$.
% 	Standard notation is: $\smoothnessf = \norm{f}_{\CCC^{2,1}}$, but I don't think we need it
	\item By definition of MNE, the local kernels $H_I, H_J \succeq 0$, 
	and by non-degeneracy assumption, $H_I, H_J \succ 0$.
	Denote $\sigma_{\min} = (\min_I \sigma_{\min}(H_I)) \wedge (\min_J \sigma_{\min}(H_J)) > 0$ the least eigenvalue.
\end{itemize}

\paragraph{Shorthands for partitions of unity.}

We recall the following notations, already introduced in our construction of the Lyapunov function in \autoref{subsec:cv_proof:gencase}.
\begin{itemize}
    \item For each $I \in [n^*]$, $\varphi_I: \XXX \to \RR$ is the function defined in \eqref{eq:cv_proof:gencase:def_varphi}, and $\varphi_0 = 1 - \sum_I \varphi_I$.
    \item Generically denote, for any $a \in \Delta_n, x \in \XXX^n$: ~
    ${
        \forall I \in [0, n^*], \forall i \in [n],~
        \varphi_{Ii} = \varphi_I(x_i)
    }$,
    and
    \begin{align*}
    % 	\label{eq:notations_for_proofs:localmoments} --> don't want to risk confusing with \label{eq:cv_proof:gencase:def_varphi}
    	\forall I \in [n^*],~
    	\ola_I &= \sum_i \varphi_{Ii} a_i
    	&
    	\olx_I &= \sum_i \frac{\varphi_{Ii} a_i}{\ola_I} x_i
    	&
    	\Sigma_I &= \sum_i \frac{\varphi_{Ii} a_i}{\ola_I} (x_i-\olx_I) (x_i-\olx_I)^\top
    \end{align*}
    as well as $\ola_0 = 1 - \sum_I \ola_I$.
    We refer to $\ola_I$ as the aggregated weights, to $\olx_I$ as the aggregated positions, to $\Sigma_I$ as the local covariance matrices, and to $\ola_0$ as the stray weight.
    
    For example, the iterates at $k$ have aggregated weights 
    $\ola^k_I = \sum_i \varphi^k_{Ii} a^k_i$.
    % $\ola^k_I = \sum_i \varphi_I(x^k_i) a^k_i$.
    % and $\olx^k_I = \sum_i \frac{\varphi_I(x^k_i) a^k_i}{\ola^k_I} x^k_i$.
    \item For any $J \in [0, m^*]$, then $\psi_J: \YYY \to \RR$ is defined similarly.
    
    For any $b \in \Delta_m$ and $y \in \YYY^m$, then
    $\olb \in \Delta_{m^*}$ and $\oly \in \YYY^{m^*}$ are defined similarly.
\end{itemize}

In addition, we let $\eps = e^{-\lambda^3/3}$ be the value of $\varphi_I$ and $\psi_J$ at the cut-off.

\paragraph{What we hide in the $O(\cdot)$'s.}
% \paragraph{Big-O shorthand notation.}

\begin{itemize}
% 	\item Fix large constants $H_0, \Sigma_0, \Gamma_0$ and restrict attention to choices of step-sizes such that 
% 	$\eta \leq H_0$, $\sigma \leq \Sigma_0$ and $\Gamma_0^{-1} \leq \frac{\eta}{\sigma} \leq \Gamma_0$.%
% 	\footnote{$H_0$ and $\Sigma_0$ are introduced only to make it clear from the start that $\eta$ and $\sigma$ are bounded; in the course of the proof we eventually need to assume that they are sufficiently small anyway. The restriction that $\Gamma_0^{-1} \leq \frac{\eta}{\sigma} \leq \Gamma_0$ does remain though, and our convergence rate depends on that constant a priori.}
% % 	\TODO{actually introducing $H_0, \Sigma_0$ seems a bit unnecessary, they are almost never used... maybe even strictly never used (to check).} -> bah it's fine. We do use it only once: in the proof of \autoref{lm:proof_gencase:errterm3__control_nonBregmanness_errterm}. The alternative would be to write $\eta, \sigma \leq 1000$ but that would probably be even more confusing
% % 	which is fixed a priori; in fact the $\gamma_0$ in the theorem statements is simply $\Gamma_0^{-1}$.
    \item Fix an arbitrary constant $\Gamma_0 \geq 1$ and restrict attention to choices of step-sizes such that 
    $\eta, \sigma \leq 1000$ and
    $\Gamma_0^{-1} \leq \frac{\sigma}{\eta} \leq \Gamma_0$.
% 
%	Likewise fix large constants $\alpha_0, \lambda_0, \tau_0$ and restrict attention to choices of degree $\alpha \leq \alpha_0$, $\lambda \leq \lambda_0$, $\tau \leq \tau_0$. \TODO{check that adding this $\alpha_0$ is WLOG (the rest is clearly ok).}
%	Likewise fix a large constant $\alpha_0$ and restrict attention to choices of degree $3 \leq \alpha \leq \alpha_0$. 
%	\red{In fact, I think $\alpha=3$ works out.}
    \item Let
    \begin{equation*}
        \wulb = \frac{a^*_{\min} \wedge b^*_{\min}}{4} = \frac{\left( \min_I a^*_I \right) \wedge \left( \min_J b^*_J \right)}{4}.
    \end{equation*}
    We will justify in
    \autoref{lm:proof_gencase:ulb_on_weights_locally_k+1}
    that, locally, this quantity is a uniform lower bound on the iterates' aggregated weights:
    $\min_I \ola^k_I, \min_J \olb^k_J, \min_I \ola^{k+1}_I, \min_J \olb^{k+1}_J \geq \wulb$.
	\item We will use $O(\cdot)$ and $\lesssim$ and $\asymp$ to hide only 
	constants dependent on $(f, \XXX, \YYY)$ 
	(such as $\wulb, \diamXY, \smoothnessf, n^*, m^*$...)
%	(note that these are constants only dependent on $(f, \XXX, \YYY)$),
    and on $\Gamma_0$; in particular these constants are independent of $n$ and $m$.
	That is,
	\begin{itemize}
		\item $a = O(b)$ means that there exists a constant $C$ only dependent on those quantities, such that $\abs{a} \leq C \abs{b}$.
		\item $a \lesssim b$ means that there exists a constant $C$ only dependent on those quantities, such that $a \leq C b$.
		\item $a \asymp b$ means that $a \lesssim b$ and $a \gtrsim b$.
	\end{itemize} 
	For example, we have $\eta, \sigma = O(1)$ and $\eta \asymp \sigma$.
% 	\red{I really think this is an OK notation: There is no risk of confusion since the analysis contains no asymptotics whatsoever. An alternative would be to write "there exists a constants $C$ dependent only etc.\ etc."\ every time, which would be exactly the same i.t.o meaning and much more verbose}
% %	The dependency will be specified more precisely every time in the lemma statements.
% %	--> nah too much work
% 
    \item Likewise, by ``for $\eta$ sufficiently small'' we mean that a property holds for all $\eta \leq \eta_0$ for some $\eta_0$ dependent only on those quantities.
\end{itemize}

% \end{document}

% !TEX root = ../main.tex
% \documentclass[../main]{subfiles}
% \begin{document}

\section{Proof of \autoref{lm:main_res:ppm_upd_well_defined}} \label{sec:proof_welldef}

Fix $k$ and denote the objective function in \eqref{eq:main_res:ppa_upd} as
\begin{align*}
    G((a, x), (b, y)) \coloneqq
    F_{n,m}((a, x), (b, y)) 
	&+ \frac{1}{\eta} \KLdiv(a, a^k) + \frac{1}{2 \sigma} \sum_{i=1}^n a^k_i \norm{x_i-x^k_i}^2 \\
	&- \frac{1}{\eta} \KLdiv(b, b^k) - \frac{1}{2 \sigma} \sum_{j=1}^m b^k_j \norm{y_j-y^k_j}^2.
\end{align*}
Recall that $\smoothnessf_1$ denotes the Lipschitz constant of $f$ and $\smoothnessf_2$ its smoothness constant.
We prove a quantitative version of \autoref{lm:main_res:ppm_upd_well_defined}.
\begin{lemma}
    $G$ is convex-concave over $(A^k \times \XXX^n) \times (B^k \times \YYY^n)$ where
    \begin{align*}
		A^k &= \left\lbrace
    		a \in \Delta_n ; \forall i,~ 
    		c_1 a^k_i \leq a_i \leq c_2 a^k_i
		\right\rbrace \\
		\text{and}~~~
		B^k &= \left\lbrace
    		b \in \Delta_m ; \forall j,~ 
    		c_1 b^k_j \leq b_j \leq c_2 b^k_j
		\right\rbrace
    \end{align*}
	for any $c_1, c_2$ such that
	\begin{align*}
		c_1 &\leq 1 \leq c_2
		& &\text{and} &
% 		\frac{c_1}{c_2} &\geq \frac{\eta \smoothnessf_1}{2}
		\frac{c_2}{c_1} &\leq \frac{2}{\eta \smoothnessf_1}
		& &\text{and} &
		c_2 &\leq \frac{1}{( \smoothnessf_1 + \smoothnessf_2 ) \sigma}.
	\end{align*}
	
	Such $c_1, c_2$ exist
% 	and so $A^k$ and $B^k$ are non-empty
	if and only if
	$\eta \leq \frac{2}{\smoothnessf_1}$ and $\sigma \leq \frac{1}{\smoothnessf_1 + \smoothnessf_2}$.
	In particular, if $\eta \leq \frac{1}{\smoothnessf_1}$ and $\sigma \leq \frac{1}{2 (\smoothnessf_1 + \smoothnessf_2)}$, then we can take $c_1 = 0.75$ and $c_2 = 1.5$.
	
	Furthermore, let $((a^*, x^*), (b^*, y^*))$ denote a saddle point of $G$ over $(A^k \times \XXX^n) \times (B^k \times \YYY^n)$.
	If $\eta \leq \frac{1}{\smoothnessf_1} \frac{c_1}{c_2}$,
	then we have $\KLdiv(a^*, a^k) + \KLdiv(b^*, b^k) \leq O(\eta)$.
	In particular, for $\eta$ small enough, $a^*$ resp.\ $b^*$ belong to the interior of $A^k$ resp.\ $B^k$.
\end{lemma}

\begin{proof}
    Fix $(\hb, \hy) \in B^k \times \YYY^m$ and let us show that $G(\cdot, (\hb, \hy))$ is convex over $A^k \times \XXX^n$. For this, it suffices to show that
    its Bregman divergence is non-negative, i.e., that
    \begin{align*}
        &\forall (a, x), (\ha, \hx) \in A^k \times \XXX^n,~ \\
        &D_{G(\cdot, (\hb, \hy))}((a, x), (\ha, \hx)) \coloneqq
        G((a, x), (\hb, \hy)) - G((\ha, \hx), (\hb, \hy)) 
        - \innerprod{\!
            \begin{pmatrix}
                \nabla_a \\
                \nabla_x
            \end{pmatrix}\!
            G((\ha, \hx), (\hb, \hy))
        }{
            \begin{pmatrix}
                a-\ha \\
                x-\hx
            \end{pmatrix}
        \!} 
        \geq 0.
    \end{align*}
    By straightforward calculations summarized in \autoref{lm:aux_lemmas:bregman_divgce}, this quantity 
    % (which is just the Bregman divergence of $G(\cdot, (\hb, \hy))$)
    is equal to
    \begin{equation*}
        D_{G(\cdot, (\hb, \hy))}((a, x), (\ha, \hx)) 
        = D_{F_{n,m}(\cdot, (\hb, \hy))}((a, x), (\ha, \hx))
        + \frac{1}{\eta} \KLdiv(a, \ha) + \frac{1}{2\sigma} \sum_i a^k_i \norm{x_i-\hx_i}^2.
    \end{equation*}
    % where 
    % $
    %     D_{F_{n,m}(\cdot, (\hb, \hy))}((a, x), (\ha, \hx)) =
    %     F_{n,m}((a, x), (\hb, \hy)) - F_{n,m}((\ha, \hx), (\hb, \hy)) 
    %     - \innerprod{
    %         \begin{pmatrix}
    %             \nabla_a \\
    %             \nabla_x
    %         \end{pmatrix}
    %         F_{n,m}((\ha, \hx), (\hb, \hy))
    %     }{
    %         \begin{pmatrix}
    %             a-\ha \\
    %             x-\hx
    %         \end{pmatrix}
    %     }
    % $.
    
    Let us now estimate the term in $D_{F_{n,m}(\cdot, (\hb, \hy))}$.
    Using the shorthands for the local payoff matrices,
	% \begin{align*}
	% 	\gmatnh_{i,j} &= f(x_i,\hy_j)
	% 	&
	% 	\gmathh_{i,j} &= f(\hx_i,\hy_j)
	% 	&
	% 	[\partial_x \gmathh]_{i,j} &= \partial_x f(\hx_i,\hy_j)
	% \end{align*}
	% so that
	% \begin{align*}
	% 	F_{n,m}((\ha, \hx), (\hb, \hy)) &= \ha^\top \gmathh \hb,
	% 	&
	% 	\nabla_a F_{n,m}((\ha, \hx), (\hb, \hy)) &= \gmathh \hb,
	% 	&
	% 	\nabla_{x_i} F_{n,m}((\ha, \hx), (\hb, \hy)) &= \ha_i \sum_j [\partial_x \gmathh]_{i,j} \hb_j.
	% \end{align*}
	% With these notations,
	\begin{align*}
		D_{F_{n,m}(\cdot, (\hb, \hy))}((a, x), (\ha, \hx))
		&= F_{n,m}((a, x), (\hb, \hy)) - F_{n,m}((\ha, \hx), (\hb, \hy)) \\
		&~~~~ - \innerprod{\nabla_{(a, x)} F_{n,m}((\ha, \hx), (\hb, \hy))}{(a, x) - (\ha, \hx)} \\
		&= a^\top \gmatnh \hb - \ha^\top \gmathh \hb
		- (a-\ha)^\top \gmathh \hb
		- \sum_{i,j} \ha_i (x_i-\hx_i) \cdot [\partial_x \gmathh]_{i,j} \hb_j \\
		&= a^\top \left( \gmatnh - \gmathh \right) \hb
		- \ha^\top \left[ \Diag(x-\hx) \partial_x \gmathh \right] \hb \\
		&= \underbrace{
			\ha^\top \left( \gmatnh - \gmathh - \left[ \Diag(x-\hx) \partial_x \gmathh \right] \right) \hb
		}_{\eqqcolon S_1}
		+ \underbrace{
			(a-\ha)^\top \left( \gmatnh - \gmathh \right) \hb
		}_{\eqqcolon S_2}.
	\end{align*}
	For the first term: For all $i, j$, by $\smoothnessf_2$-smoothness of $f(\cdot, \hy_j)$,
	\begin{align*}
		\abs{
			\left( \gmatnh - \gmathh - \left[ \Diag(x-\hx) \partial_x \gmathh \right] \right)_{i,j}
		}
		&= \abs{
			f(x_i, \hy_j) - f(\hx_i, \hy_j) - (x_i-\hx_i) \partial_x f(\hx_i, \hy_j)
		} \\
		&\leq \frac{\smoothnessf_2}{2} \norm{x_i-\hx_i}^2 \\
		\text{so}~~~
%		\abs{
%			\ha^\top \left( \gmatnh - \gmathh - \left[ \Diag(x-\hx) \partial_x \gmathh \right] \right) \hb
%		}
		\abs{S_1}
		&\leq \frac{\smoothnessf_2}{2} \sum_i \ha_i \norm{x_i-\hx_i}^2.
	\end{align*}
	For the second term: For all $i, j$, by $\smoothnessf_1$-Lipschitz-continuity of $f(\cdot, \hy_j)$,
	\begin{align*}
		\abs{\left( \gmatnh - \gmathh \right)_{i,j}}
		= \abs{f(x_i, \hy_j) - f(\hx_i, \hy_j)}
		&\leq \smoothnessf_1 \norm{x_i-\hx_i} \\
		\text{so}~~~
%		\abs{
%			(a-\ha)^\top \left[ \Diag(x-\hx) \partial_x \gmathh \right] \hb
%		}
		\abs{S_2}
		&\leq \smoothnessf_1 \sum_i \abs{a_i-\ha_i} \norm{x_i-\hx_i} \\
		&= \smoothnessf_1 \sum_i \frac{\abs{a_i-\ha_i}}{\sqrt{\ha_i}} \cdot \sqrt{\ha_i} \norm{x_i-\hx_i} \\
		&\leq \frac{\smoothnessf_1}{2} \left( 
		    \sum_i \frac{(a_i-\ha_i)^2}{\ha_i}
		    + \sum_i \ha_i \norm{x_i-\hx_i}^2
		\right).
	\end{align*}
	
	Thus we have that for all $(a, x), (\ha, \hx)$,
	\begin{align*}
        D_{G(\cdot, (\hb, \hy))}((a, x), (\ha, \hx)) 
%         &\geq 
% 	    -\frac{\smoothnessf_1}{2} \left( 
% 	        \sum_i \frac{(a_i-\ha_i)^2}{\ha_i}
% 		    + \sum_i \ha_i \norm{x_i-\hx_i}^2 
% 		\right)
%         + \frac{1}{\eta} \KLdiv(a, \ha) + \frac{1}{2\sigma} \sum_i a^k_i \norm{x_i-\hx_i}^2 \\
        &\geq 
	    -\frac{\smoothnessf_1}{2} \sum_i \frac{(a_i-\ha_i)^2}{\ha_i}
	    -\frac{\smoothnessf_1+\smoothnessf_2}{2} \sum_i \ha_i \norm{x_i-\hx_i}^2 \\
        &~~~~ + \frac{1}{\eta} \KLdiv(a, \ha) + \frac{1}{2\sigma} \sum_i a^k_i \norm{x_i-\hx_i}^2 \\
        &= \underbrace{
            \frac{1}{\eta} \KLdiv(a, \ha) - \frac{\smoothnessf_1}{2} \chi^2(a, \ha)
        }
        ~+~
        \underbrace{
            \sum_i \left( \frac{1}{2\sigma} a^k_i - \frac{\smoothnessf_1+\smoothnessf_2}{2} \ha_i \right) \norm{x_i-\hx_i}^2.
        }
	\end{align*}
	\begin{itemize}
        \item By \autoref{lm:aux_lemmas:KL_chi2_comparison}, if $\max_i \frac{a_i}{\ha_i} \leq \frac{2}{\eta \smoothnessf_1}$, then $\chi^2(a, \ha) \leq \frac{2}{\eta \smoothnessf_1} \KLdiv(a, \ha)$, and so the first underbrace is non-negative.
        \item If furthermore $\frac{1}{2\sigma} a^k_i - \frac{\smoothnessf_1+\smoothnessf_2}{2} \ha_i \geq 0$ for all $i$, then the second underbrace is non-negative.
    \end{itemize}
    Both of these conditions can be ensured by imposing $a, \ha \in A^k$ with $c_1, c_2$ as defined in the lemma, since $c_1 a^k_i \leq a_i, \ha_i \leq c_2 a^k_i \implies \frac{a_i}{\ha_i} \leq \frac{c_2}{c_1}$.
    
    The conditions on $\eta, \sigma$ and the possible choices of $c_1, c_2$ are straightforward to check. 
    
    Finally, fix some admissible $c_1, c_2$ and let us prove the last part of the lemma. 
    Let $((a^*, x^*), (b^*, y^*))$ denote a saddle point of $G$ over $(A^k \times \XXX^n) \times (B^k \times \YYY^n)$, i.e., such that
    \begin{equation*}
        \forall ((a, x), (b, y)) \in (A^k \times \XXX^n) \times (B^k \times \YYY^n),~ 
        G((a^*, x^*), (b, y)) \leq G((a^*, x^*), (b^*, y^*)) \leq G((a, x), (b^*, y^*)).
    \end{equation*}
    We have shown above that for all $(a, x), (\ha, \hx) \in A^k \times \XXX^n$ and $(b, y), (\hb, \hy) \in B^k \times \YYY^n$,
    \begin{align*}
        D_{G(\cdot, (\hb, \hy))}((a, x), (\ha, \hx)) 
        &\geq \frac{1}{\eta} \KLdiv(a, \ha) - \frac{\smoothnessf_1}{2} \chi^2(a, \ha)
        + \sum_i \left( \frac{1}{2\sigma} a^k_i - \frac{\smoothnessf_1+\smoothnessf_2}{2} \ha_i \right) \norm{x_i-\hx_i}^2 \\
        % &\geq \frac{1}{\eta} \KLdiv(a, \ha) - \frac{\smoothnessf_1}{2} \left( \max_i \frac{a_i}{\ha_i} \right) \KLdiv(a, \ha) \\
        &\geq \left( \frac{1}{\eta} - \frac{\smoothnessf_1}{2} \frac{c_2}{c_1} \right) \KLdiv(a, \ha)
    \end{align*}
    and symmetrically,
    \begin{equation*}
        D_{G((\ha, \hx), \cdot)}((b, y), (\hb, \hy)) 
        \leq - \left( \frac{1}{\eta} - \frac{\smoothnessf_1}{2} \frac{c_2}{c_1} \right) \KLdiv(b, \hb).
    \end{equation*}
    In particular, for $(a, x) = (a^*, x^*)$, $(\ha, \hx) = (a^k, x^k)$ and $(b, y) = (b^*, y^*)$, $(\hb, \hy) = (b^k, y^k)$, the difference of the left-hand sides reads
    \begin{align*}
        & D_{G(\cdot, (b^k, y^k))}((a^*, x^*), (a^k, x^k)) 
        - D_{G((a^k, x^k), \cdot)}((b^*, y^*), (b^k, y^k)) \\
        &= G((a^*, x^*), (b^k, y^k)) - G((a^k, x^k), (b^k, y^k))
        - \innerprod{
            \begin{pmatrix}
                \nabla_a \\
                \nabla_x
            \end{pmatrix} G(a^k, x^k, b^k, y^k)
        }{
            \begin{pmatrix}
                a^* - a^k \\
                x^* - x^k
            \end{pmatrix}
        } \\
        &~~ - \left(
            G((a^k, x^k), (b^*, y^*)) - G((a^k, x^k), (b^k, y^k))
            - \innerprod{
                \begin{pmatrix}
                    \nabla_b \\
                    \nabla_y
                \end{pmatrix} G(a^k, x^k, b^k, y^k)
            }{
                \begin{pmatrix}
                    b^* - b^k \\
                    y^* - y^k
                \end{pmatrix}
            }
        \right) \\
        &= G((a^*, x^*), (b^k, y^k)) - G((a^k, x^k), (b^*, y^*))
        - \innerprod{
            \begin{pmatrix}
                \nabla_a \\
                \nabla_x \\
                -\nabla_b \\
                -\nabla_y
            \end{pmatrix} G(a^k, x^k, b^k, y^k)
        }{
            \begin{pmatrix}
                a^* - a^k \\
                x^* - x^k \\
                b^* - b^k \\
                y^* - y^k
            \end{pmatrix}
        }.
    \end{align*}
    Now $G((a^*, x^*), (b^k, y^k)) - G((a^k, x^k), (b^*, y^*)) \leq 0$ by definition of the saddle point, so we have
    \begin{equation*}
        \left( \frac{1}{\eta} - \frac{\smoothnessf_1}{2} \frac{c_2}{c_1} \right) \left( \KLdiv(a^*, a^k) + \KLdiv(b^*, b^k) \right)
        \leq 
        - \innerprod{
            \begin{pmatrix}
                \nabla_a \\
                \nabla_x \\
                -\nabla_b \\
                -\nabla_y
            \end{pmatrix} G(a^k, x^k, b^k, y^k)
        }{
            \begin{pmatrix}
                a^* - a^k \\
                x^* - x^k \\
                b^* - b^k \\
                y^* - y^k
            \end{pmatrix}
        }.
    \end{equation*}
    Since $\restr{ \nabla_a \KLdiv(a, a^k) }{a=a^k} = 0$ and $\restr{ \nabla_{x_i} \frac{1}{2} \norm{x_i-x^k_i}^2 }{x_i=x^k_i} = 0$, the right-hand side is equal to
    \begin{equation*}
        -\innerprod{
            \begin{pmatrix}
                \nabla_a \\
                \nabla_x \\
                -\nabla_b \\
                -\nabla_y
            \end{pmatrix} F_{n,m}(a^k, x^k, b^k, y^k)
        }{
            \begin{pmatrix}
                a^* - a^k \\
                x^* - x^k \\
                b^* - b^k \\
                y^* - y^k
            \end{pmatrix}
        }.
    \end{equation*}
    By proceeding similarly as for our bound of 
    % $D_{F_{n,m}(\cdot, (\hb, \hy))}((a, x), (\ha, \hx))$,
    $D_{F_{n,m}(\cdot, (\hb, \hy))}$,
    one can show that it is bounded by 
    $C \coloneqq 2 (\smoothnessf_0 \vee \smoothnessf_1) \left( 1 + \diamXY \right)$
    % where $\diamXY$ is an upper bound on the diameter of $\XXX$ and $\YYY$ 
    (in particular the bound does not depend on $n$, $m$).
    Thus we have as announced
    \begin{equation*}
        \KLdiv(a^*, a^k) + \KLdiv(b^*, b^k)
        \leq \frac{C}{\left( \frac{1}{\eta} - \frac{\smoothnessf_1}{2} \frac{c_2}{c_1} \right) }
        = O(\eta).
    \qedhere
    \end{equation*}
\end{proof}

\section{Proof of the lower growth properties} \label{sec:growth_conds}

Here we give the proof of the lower growth properties, which are crucial ingredients in our convergence analysis of CP-PP, and which were discussed in \autoref{subsec:cv_proof:proofingr_growthconds}.

% \subsection{Proof of ``quadratic growth'' with respect to the position and stray weight variables} \label{subsec:growth_conds:quadr}
\subsection{Proof of ``quadratic growth''} \label{subsec:growth_conds:quadr}

\begin{lemma}[``Quadratic growth'', general case] \label{lm:growth_conds:quadr_gencase}
    Consider the Lyapunov function $V$ as in \eqref{eq:cv_proof:gencase:def_V} with the partitions of unity $(\varphi_I)_I$ and $(\psi_J)_J$ as in \eqref{eq:cv_proof:gencase:def_varphi}.
    % We make no assumption on $\lambda, \tau$ other than that the $\support(\varphi_I)$, resp.\ $\support(\psi_J)$, do not overlap.
    
    Then for any $(\ha, \hx, \hb, \hy) \in \Delta_n \times \XXX^n \times \Delta_m \times \YYY^m$, 
	\begin{equation*}
	    F(\hmu, \nu^*) - F(\mu^*, \hnu) 
	   % = \ha \gmaths b^* - a^* \gmatsh \hb
		\geq 
		\max \left\{
    		\left[
    		    \frac{\sigma_{\min}}{2} \wedge \frac{2 \gvsep}{(\lambda\tau)^2}
    		\right]
    		V_\pos(\ha, \hx, \hb, \hy), ~~
    		\left[ 
    			\frac{\sigma_{\min}}{4} \frac{3 (\lambda\tau)^2}{\lambda^3} \wedge \gvsep
    		\right]
    		(\hola_0 + \holb_0)
    	\right\}
	\end{equation*}
	for some constant $\gvsep>0$ only dependent on $(f, \XXX, \YYY)$.
\end{lemma}

\begin{proof}
    We have
	\begin{equation*}
	    F(\hmu, \nu^*) - F(\mu^*, \hnu) 
	    = \innerprod{\hmu}{F\nu^* - \gval}_{\CCC(\XXX)} + \innerprod{\gval - (\mu^*)^\top F}{\hnu}_{\CCC(\YYY)}.
	\end{equation*}
	Focus on the first term. 
	By the non-degeneracy Assumptions~\ref{assum:5} and \ref{assum:6}, $F \nu^*$ grows quadratically as $\frac{1}{2} \sigma_{\min}(H_I)$ on the neighborhood of $x^*_I$ for each $I$, and is lower-bounded by a constant everywhere else. 
	In symbols, there exists a constant $\gvsep>0$ such that
	\begin{equation} \label{eq:growth_conds:nondegengrowth_Fnustar}
		\forall x \in \XXX,~
		(F \nu^*)(x) - \gval
		\geq \left( \frac{\sigma_{\min}}{4} \min_I \norm{x-x^*_I}^2 \right) \wedge \gvsep.
	\end{equation}
	
	This directly implies a lower bound in terms of the position variables. Indeed,
    \begin{equation*}
    	\forall x \in \support(\varphi_I) = B_{x^*_I, \lambda\tau},~ (F \nu^*)(x)- \gval \geq \left[ \frac{\sigma_{\min}}{4} \wedge \frac{\gvsep}{(\lambda\tau)^2} \right] \norm{x-x^*_I}^2,
    \end{equation*}
    and $(F \nu^*)(x) - \gval \geq 0$ on all of $\XXX$, so by decomposing $\hmu = \sum_I \varphi_I \hmu + \varphi_0 \hmu$,
    \begin{equation*}
    	\innerprod{\hmu}{F\nu^* - \gval}_{\CCC(\XXX)}
    	\geq \left[ \frac{\sigma_{\min}}{4} \wedge \frac{\gvsep}{(\lambda\tau)^2} \right] 
    	\cdot
    	\underbrace{
    	    \sum_{I \in [n^*]} \sum_{i=1}^n \hvarphi_{Ii} \ha_i \norm{\hx_i-x^*_I}^2
    	}_{= 2 V_\pos(\ha, \hx)}.
    \end{equation*}
    
	We can also get a lower bound in terms of the ``stray weights'' $\hola_0$ and $\holb_0$. Indeed,
	\begin{equation*}
		\forall r \leq \lambda,~
		1 - e^{-\frac{r^3}{3}}
		\leq \frac{r^3}{3}
		\leq \frac{\lambda}{3} \cdot r^2
	\end{equation*}
	so that for each $I$,
	\begin{align*}
		\forall x \in B_{x^*_I, \lambda\tau},~
		\varphi_0(x) &= 1 - \exp\left( -\frac{\norm{x-x^*_I}^3}{3 \tau^3} \right)
		\leq \frac{\lambda}{3} \cdot \frac{\norm{x-x^*_I}^2}{\tau^2} \\
		\text{and so}~~~~
		\forall x \in \XXX,~
		\varphi_0(x) &\leq \left( \frac{\lambda}{3 \tau^2} \min_I \norm{x-x^*_I}^2 \right) \wedge 1.
	\end{align*}
	Hence, \eqref{eq:growth_conds:nondegengrowth_Fnustar} implies
	\begin{align*}
		\forall x \in \XXX, \,
		(F \nu^* - \gval)(x)
		\geq \left( \frac{\sigma_{\min}}{4} \min_I \norm{x-x^*_I}^2 \right) \wedge \gvsep
		&\geq 
		\left[ 
			\frac{\sigma_{\min}}{4} \frac{3 \tau^2}{\lambda} \wedge \gvsep
		\right] \!\cdot\!
		\left[ 
			\left( \frac{\lambda}{3 \tau^2} \min_I \norm{x-x^*_I}^2 \right) \wedge 1
		\right] \\
		&\geq 
		\left[ 
			\frac{\sigma_{\min}}{4} \frac{3 \tau^2}{\lambda} \wedge \gvsep
		\right] \cdot
		\varphi_0(x) \\
        \text{and so finally} \qquad\qquad
		\innerprod{\hmu}{F\nu^* - \gval}_{\CCC(\XXX)}
		&\geq 
		\left[ 
			\frac{\sigma_{\min}}{4} \frac{3 \tau^2}{\lambda} \wedge \gvsep
		\right]
		\cdot \underbrace{ \int_\XXX \varphi_0 d\hmu }_{=\hola_0}.
    \rqedhere
	\end{align*}
\end{proof}

For the exact-parametrization case, we can actually reuse the result for the general case, using that our two Lyapunov functions ``coincide locally''. We summarize that fact in the following easily checked claim, which will also be useful 
% in \autoref{subsec:growth_conds:error_bound} and \autoref{subsec:rel_Lya_NI:exactparam}.
in Sections~\ref{subsec:growth_conds:error_bound} and \ref{subsec:rel_Lya_NI:exactparam}.

\begin{claim}
\label{claim:growth_conds:gencase_imply_exactparam}
    Consider the exact-parametrization case.
    Denote
    \begin{itemize}
        \item $d_* = \min_{i' \neq i''} \norm{x^*_{i'} - x^*_{i''}} \wedge \min_{j' \neq j''} \norm{y^*_{j'} - y^*_{j''}}$.
        Fix any $\uld_* \leq d_*$.
        \item $V_e$ the Lyapunov function designed for the exact-parametrization case \eqref{eq:cv_proof:exactparam:def_V}.
        \item $V_g$ the Lyapunov function designed for the general case  \eqref{eq:cv_proof:gencase:def_V} with $\varphi_I$ being the indicator function of $B_{x^*_I, \uld_*/4}$ (well-defined since those balls do not intersect).
        \item $V_{\lambda, \tau}$ the Lyapunov function designed for the general case  \eqref{eq:cv_proof:gencase:def_V} with the parametrized partitions of unity $(\varphi_I)_I$ and $(\psi_J)_J$ from \eqref{eq:cv_proof:gencase:def_varphi}.
    \end{itemize}
    
    There exists $r>0$ such that if $V_e(\hz) \leq r$, then $\ha_i \geq \frac{a^*_i}{2}$ and $\norm{\hx_i-x^*_i} \leq \uld_*/4$ for all $i \in [n]$.
    That is, for all $i$, $\hx_i \in B_{x^*_i, \uld_*/4}$.
    For such $\hz$, we have $V_e(\hz) = V_g(\hz)$.
    
    Moreover, $V_g$ is the point-wise limit of
    $V_{\lambda, \tau}$
    when $\lambda\tau$ is held constant equal to $\uld_*/4$ and $\tau \to \infty$.
\end{claim}

% \begin{proof}
%     Let $d_* = \min_{i' \neq i''} \norm{x^*_{i'} - x^*_{i''}} \wedge \min_{j' \neq j''} \norm{y^*_{j'} - y^*_{j''}}$.
%     % It is not hard to show that
%     There exists $r>0$ such that if $(\eta \vee \sigma) V(\hz) \leq r$, then $\ha_i \geq \frac{a^*_i}{2}$ and $\norm{\hx_i-x^*_i} \leq d_*/4$ for all $i \in [n]$.
%     That is, for all $i$, $\hx_i \in B_{x^*_i, d_*/4}$ and those balls do not intersect.
%     For such $\hz$, the Lyapunov function $V$ from \eqref{eq:cv_proof:exactparam:def_V} (the one designed for the exact-parametrization case) is equivalent to the one from \eqref{eq:cv_proof:gencase:def_V} (the one designed for the general case), with $\varphi_I$ being the indicator function of $B_{x^*_I, d_*/4}$.
%     The lemma then follows as a limit case of the next lemma for $\lambda\tau$ held constant equal to $d_*/4$, and $\tau \to \infty$.
% \end{proof}

The result for the exact-parametrization case now follows immediately from the above lemma and claim:

\begin{lemma}[``Quadratic growth'', exact-parametrization case] \label{lm:growth_conds:quadr_exactparam}
    Assume $n=n^*, m=m^*$, and consider the Lyapunov function $V$ defined in \eqref{eq:cv_proof:exactparam:def_V}.
    There exist constants $r, C>0$ only dependent on $(f, \XXX, \YYY)$ such that, for any $\hz = (\ha, \hx, \hb, \hy)$ with $V(\hz) \leq r$, then
    \begin{equation*}
        F(\hmu, \nu^*) - F(\mu^*, \hnu) 
        \geq
		C V_\pos(\ha, \hx, \hb, \hy).
    \end{equation*}
\end{lemma}

% \subsection{Proof of ``error bound'' with respect to the weight and aggregated position variables} \label{subsec:growth_conds:error_bound}
\subsection{Proof of ``error bound''} \label{subsec:growth_conds:error_bound}

Our error-bound-type result is contained in the following lemma.
It is stated for the general case, i.e., $V_1$ refers to the Lyapunov function defined in \eqref{eq:cv_proof:gencase:def_V}.
But since we make no assumption on the partitions of unity $(\varphi_I)_I, (\psi_J)_J$, the conclusion of the lemma (``if $V_1(\hz) \leq r$ then we have this error bound inequality'') is also true for the exact-parametrization case with $V_1$ referring to the Lyapunov function from \eqref{eq:cv_proof:exactparam:def_V}, as one can deduce a posteriori thanks to \autoref{claim:growth_conds:gencase_imply_exactparam}.

\begin{lemma} \label{lm:growth_conds:error_bound}
    Consider any $\hz = (\ha, \hx, \hb, \hy) \in \Delta_n \times \XXX^n \times \Delta_m \times \YYY^m$.
    For any $(A_I)_{I \in [n^*]}, (B_J)_{J \in [m^*]}$ (and $A_0=B_0=0$) and $(X_I)_{I \in [n^*]}, (Y_J)_{J \in [m^*]}$,
    let analogously to \eqref{eq:cv_proof:gencase:choice_of_proxy*} the ``proxy particles''
    \begin{equation} \label{eq:growth_conds:encode_WP_hat}
    	x_i = \hx_i + \sum_I \hvarphi_{Ii} (X_I - \hx_i)
    	~~~\text{and}~~~
    	a_i = \sum_I A_I \frac{\hvarphi_{Ii} \ha_i}{\hola_I}
    \end{equation}
    and similarly for $b, y$, and $z = (a, x, b, y)$.
    
	There exist $r, C > 0$ only dependent on $(f, \XXX, \YYY)$ such that if $V_1(\hz) \leq r$, then
	\begin{equation*}
		\max_{A,X,B,Y} \tgap(z; \hz) 
		\geq 
		C \sqrt{
			\sum_I d_h(w^*_I, \holw_I)
			+ \sum_I \holw_I \norm{\Delta \holp_I}^2
		}
		+ O \left( V_1(\hz) \right).
	\end{equation*}
\end{lemma}

Informally, $\max_{A,X,B,Y} \tgap(z; \hz)$ can be interpreted as a lower bound (up to a constant) on
$\max_{\norm{\delta z}_* \leq 1} \innerprod{
    \begin{pmatrix}
        \nabla_a \\ \nabla_x \\ -\nabla_b \\ -\nabla_y
    \end{pmatrix}
    F_{n,m}(\hz)
}{
    \begin{pmatrix}
        \delta a \\ \delta x \\ \delta b \\ \delta y
    \end{pmatrix}
}
= \norm{\nabla F_{n,m}(\hz)}$.
Note that in the latter expression, $\delta x$ has $n$ degrees of freedom, whereas in $\max_{A,X,B,Y} \tgap(z; \hz)$, $X$ has only $n^*$ degrees of freedom (and similarly for $A, B$ and $Y$).
Thus, $\max_{A,X,B,Y} \tgap(z; \hz)$ represents a ``norm'' of $\nabla F_{n,m}(\hz)$ --- which justifies why we refer to \autoref{lm:growth_conds:error_bound} as an error bound ---, for a notion of norm that is adapted to the geometry of the algorithm and of the problem at hand.

The remainder of this subsection is dedicated to proving the above lemma.
To lighten notation, we continue to leave the dependence of $z = (a, x, b, y)$ on $A, X, B$ and $Y$ implicit.

% Furthermore, for concision we denote as usual
% $w = \begin{pmatrix}
%     a \\
%     b
% \end{pmatrix}$,
% $p = \begin{pmatrix}
%     x \\
%     y
% \end{pmatrix}$,
% $W = \begin{pmatrix}
%     A \\
%     B
% \end{pmatrix}$, and
% $P = \begin{pmatrix}
%     X \\
%     Y
% \end{pmatrix}$.

Let us start by showing that locally (i.e.\ for $z$'s with small enough Lyapunov potential), we have a constant lower bound on the aggregated weights $\ola_I, \olb_J$.
This fact will be used repeatedly throughout this appendix and the next ones.

\begin{lemma} \label{lm:growth_conds:ulb_on_weights_locally}
% 	There exists $r>0$ (only dependent on $a^*$) such that if $V_\wei(a, x) \leq r$, then
% 	\begin{equation*}
% 		\forall I \in [n^*], \ola_I \geq \frac{a^*_{\min}}{2}.
% 	\end{equation*}
% 
	There exists $r>0$ (only dependent on $a^*, b^*$) such that if $V_\wei(a, x, b, y) \leq r$, then
	\begin{equation*}
		\left( \min_{I \neq 0} \ola_I \right) \wedge \left( \min_{J \neq 0} \olb_J \right)
        \geq \frac{a^*_{\min} \wedge b^*_{\min}}{2}.
	\end{equation*}
\end{lemma}

\begin{proof}
	$h: [0, 1] \to \RR, s \mapsto s \log s - s + 1$ is $1$-strongly convex
    (as $h'' \geq h''(1)=1$),
    % (just bound $h''$), 
    so
	for any $I$, 
	\begin{equation*}
		(a^*_I-\ola_I)^2 \leq 2 d_h(a^*_I, \ola_I)
		\leq 2 V_\wei(a, x)
		\leq 2 r.
	\end{equation*}
    Similarly, for any $J$, $(b^*_J - \olb_J)^2 \leq 2r$.
    So by choosing $r$ small enough, we can ensure that
    % \begin{equation*}
    %     \min_I \ola_I \geq \frac{\min_I a^*_I}{2}
    %     \qquad\text{and}\qquad
    %     \min_J \olb_J \geq \frac{\min_J b^*_J}{2},
    % \end{equation*}
    $\min_I \ola_I \geq \frac{\min_I a^*_I}{2}$ and 
    $\min_J \olb_J \geq \frac{\min_J b^*_J}{2}$,
    and the lemma follows by combining these two inequalities.
\end{proof}

\subsubsection{Step 1: Reduce to a bilinear functional in the aggregated weights and positions}

\begin{claim} \label{claim:growth_conds:tgaplocs_proxyz}
    For any $\hz$ and $A, X, B, Y$,
	$\tgap(z; \hz)$ is approximately given by
	\begin{align*} 
% \label{eq:gencase:tgaplocs}
		\tgaplocs(z; \hz)
		&= - \Delta A^\top \gmatss \tDelta \holb ~+~ \tDelta \hola^\top \gmatss \Delta B \\
		&~~~~ + \tDelta \hola^\top \partial_y \gmatss (\Delta Y \odot b^*)
		~-~ (a^* \odot \Delta X)^\top \partial_x \gmatss \tDelta \holb \\
		&~~~~ - \Delta A^\top \partial_y \gmatss (\Delta \holy \odot b^*)
		~+~ (a^* \odot \Delta \holx)^\top \partial_x \gmatss \Delta B \\
		&~~~~
% 		- \sum_I w^*_I \innerprod{\Delta P_I}{\Delta \holp_I}_{H_I} \\
		- \sum_I a^*_I \innerprod{\Delta X_I}{\Delta \holx_I}_{H_I}
		- \sum_J b^*_J \innerprod{\Delta Y_I}{\Delta \holy_J}_{H_J} \\
		&~~~~ - (a^* \odot \Delta X)^\top \partial_{xy}^2 \gmatss (\Delta \holy \odot b^*)
		~+~ (a^* \odot \Delta \holx)^\top \partial_{xy}^2 \gmatss (\Delta Y \odot b^*)
	\end{align*}
    where $\tDelta \hola = \Delta \hola + \hola_0 a^*$ and $\tDelta \holb = \Delta \holb + \holb_0 b^*$,
	and more precisely
	\begin{equation*}
		\tgap(z; \hz)
		= \tgaplocs(z; \hz)
		+ O \left( (\min_I \holw_I)^{-1} V_1(\hz) \right).
	\end{equation*}
\end{claim}

The proof consists of simple but tedious calculations,
\ifextended% 
    which we defer to 
    \autoref{subsec:tedious_calcs:proof__claim--growth_conds--tgaplocs_proxyz}.
\else%
    which can be found as \AdditionalMaterial\ in Section~G.1 of the arXiv version of this paper.
\fi
Essentially we just write out the expression of $\tgap(z; \hz)$ and do Taylor expansions of $f(\hx_i, \hy_j)$ around $(x^*_I, y^*_J)$.

We see that $\tgaplocs(z; \hz)$ has a bilinear structure; in matrix form, denoting $[a^* H_x]_{I I'} = \ind_{I=I'} a^*_I H_I$
and $[a^* \partial_x \gmatss]_{IJ} = a^*_I \partial_x \gmatss_{IJ}$, etc., 
\begin{equation} \label{eq:growth_conds:tgaplocs_symbolic}
	\tgaplocs(z; \hz) 
	= -\begin{pmatrix}
		\Delta A \\
		\Delta X \\
		\Delta B \\
		\Delta Y
	\end{pmatrix}^\top
	\begin{bmatrix}
		0 & 0 & \gmatss & \partial_y \gmatss b^* \\
		0 & a^* H_x & a^* \partial_x \gmatss & a^* \partial_{xy}^2 \gmatss b^* \\
		-(\gmatss)^\top & -(a^* \partial_x \gmatss)^\top & 0 & 0 \\
		-(\partial_y \gmatss b^*)^\top & -(a^* \partial_{xy}^2 \gmatss b^*)^\top & 0 & b^* H_y
	\end{bmatrix}
	\begin{pmatrix}
		\tDelta \hola \\
		\Delta \holx \\
		\tDelta \holb \\
		\Delta \holy
	\end{pmatrix}.
\end{equation}
Here the first component $\Delta A$/$\tDelta \hola$ is a vector in $\RR^{n^*}$, and the second component $\Delta X$/$\Delta \holx$ is in $\XXX^{n^*}$.
Note that on the right, the first component $\tDelta \hola = \Delta \hola + \hola_0 a^*$ sums to $0$ by definition.
For concision and for clarity of the argument to follow, introduce some notation for the remainder of this section:
\begin{itemize}
	\item Let $\HHH$ denote the above block matrix.
	\item Let 
	\begin{equation*} %\label{eq:growth_conds:def_ZZZ}
        \ZZZ = \Delta_{n^*} \times \XXX^{n^*} \times \Delta_{m^*} \times \YYY^{m^*}
        ~~~~\text{and}~~~~
		Z = \begin{pmatrix}
			A \\
			X \\
			B \\
			Y
		\end{pmatrix},~
		Z^* = \begin{pmatrix}
			a^* \\
			x^* \\
			b^* \\
			y^*
		\end{pmatrix}
        \in \ZZZ.
	\end{equation*}
    Also let 
    $\hZ = \begin{pmatrix}
        \hola \\
        \holx \\
        \holb \\
        \holy
    \end{pmatrix}$, which does not belong to $\ZZZ$ (unless $\hola_0=\holb_0=0$),
    and 
    $\tDelta \hZ = \begin{pmatrix}
        \tDelta \hola \\
        \Delta \holx \\
        \tDelta \holb \\
        \Delta \holy
    \end{pmatrix}$
    so that $\tDelta \hZ + Z^* \in \ZZZ$.
	\item For any $\tZ \in \ZZZ$, denote $\Delta \tZ = \tZ - Z^*$ and let
	\begin{equation} \label{eq:growth_conds:def_norm_ZZZ}
		\norm{\Delta \tZ} = \norm{\talpha - a^*}_1 + \max_I \norm{\tolx_I - x^*_I}
		+ \norm{\tbeta - b^*}_1 + \max_J \norm{\toly_J - y^*_J}.
	\end{equation}
	Clearly $\norm{\cdot}$ defines a norm on $\ZZZ - Z^*$.
\end{itemize}
%	Importantly, $\HHH$ is the sum of an anti-symmetric and a positive semi-definite matrix.
%	So we can rewrite and lower-bound $\tgaplocs$ as
%	\begin{align*}
%		\tgaplocs(z; \hz) 
%		= (\hZ - Z)^\top \HHH \Delta \hZ
%		&= (\Delta \hZ - \Delta Z)^\top \HHH \Delta \hZ \\
%		&= -\Delta Z^\top \HHH \Delta \hZ 
%		+ (\Delta \hZ)^\top \HHH \Delta \hZ \\
%		&\geq -\Delta Z^\top \HHH \Delta \hZ.
%	\end{align*}

So far we showed that, for all $\hz$ with $V(\hz)$ less than some $r$ so that \autoref{lm:growth_conds:ulb_on_weights_locally} applies,
\begin{equation*}
	\max_{W, P} \tgap(z; \hz)
	% &= \max_{W, P} \tgaplocs(z; \hz) + O \left( V_1(\hz) \right) \\
	= \max_{\substack{
			A \in \Delta_{n^*}, X \in \XXX^{n^*} \\
			B \in \Delta_{m^*}, Y \in \YYY^{m^*}
	}}
	-\begin{pmatrix}
		\Delta A \\
		\Delta X \\
		\Delta B \\
		\Delta Y
	\end{pmatrix}^\top
	\HHH
	\begin{pmatrix}
		\tDelta \hola \\
		\Delta \holx \\
		\tDelta \holb \\
		\Delta \holy
	\end{pmatrix}
	+ O \left( V_1(\hz) \right).
\end{equation*}
To complete the proof of the lemma, it suffices to prove that there exist $r, C > 0$ only dependent on $(f, \XXX, \YYY)$ such that if $V_1(\hz) \leq r$, then
\begin{equation*}
	\max_{\substack{
			A \in \Delta_{n^*}, X \in \XXX^{n^*} \\
			B \in \Delta_{m^*}, Y \in \YYY^{m^*}
	}}
	-\Delta Z^\top \HHH \tDelta \hZ
	\geq 
	C \sqrt{
		\sum_I d_h(w^*_I, \holw_I)
		+ \sum_I \holw_I \norm{\Delta \holp_I}^2
	}.
\end{equation*}
We will do so by working with the aggregated variables $\Delta Z, \tDelta \hZ$.
So let us first clarify the relation between the norm on $\ZZZ - Z^*$ and the desired divergence.
Namely, we show that $\norm{\tDelta \hZ}$ is equivalent to 
% $\sqrt{V_1(\hz) - \holw_0}$ 
the square root of the desired divergence
up to additive terms in $O(\holw_0) \leq O(V_1(\hz))$.

\begin{claim} \label{claim:growth_conds:DeltaZ_equiv_etaveesigmaV}
    Suppose $\min_I \hola_I, \min_J \holb_J \geq \wulb$ for some constant $\wulb>0$. Then
	\begin{equation*}
		2\wulb \sum_I d_h(w^*_I, \holw_I)
		+ \sum_I \holw_I \norm{\Delta \holp_I}^2
		\,\leq\,
		\norm{\hZ - Z^*}^2
        \,\leq \,
		8 (n^* \wedge m^*) \sum_I d_h(w^*_I, \holw_I)
		+ \frac{4}{\wulb} \sum_I \holw_I \norm{\Delta \holp_I}^2
		% \lesssim V_1(\hz).
	\end{equation*}
    and $\abs{ \norm{\tDelta Z} - \norm{\hZ - Z^*} } \leq \holw_0$.
\end{claim}
% \TODO{do we need the second inequality anywhere?} --> I think so yes, in \autoref{sec:proof_rel_Lya_NI}

\begin{proof}
	For the first inequality, for the weight part:
	$h: x \mapsto x \log x - x + 1$ is $\frac{1}{\wulb}$-smooth over $[\wulb, 1]$ 
    (since $1 \leq h''(x) = \frac{1}{x} \leq \frac{1}{c}$), 
    so
    $\forall s, s' \in [\wulb, 1],~
    d_h(s, s') \leq \frac{1}{2\wulb} \abs{s-s'}^2$.
	So
	\begin{equation*}
		\sum_I d_h(a^*_I, \hola_I)
		\leq \frac{1}{2\wulb} \sum_I (a^*_I-\hola_I)^2 
		\leq \frac{1}{2\wulb} \left( \sum_I \abs{a^*_I-\hola_I} \right)^2.
	\end{equation*}
	For the position part: just write
    $\sum_I \hola_I \norm{\Delta \holx_I}^2
    \leq \max_I \norm{\Delta \holx_I}^2$.
	
	For the second inequality, for the weight part:
	$h$ is $1$-strongly concave over $[\wulb, 1]$, so
    $\forall s, s' \in [\wulb, 1],~
    d_h(s, s') \geq \frac{1}{2} \abs{s-s'}^2$.
	So
	\begin{equation*}
		\sum_I d_h(a^*_I, \hola_I)
		\geq \frac{1}{2} \sum_I (a^*_I-\hola_I)^2
		\geq \frac{1}{2 n^*} \left( \sum_I \abs{a^*_I-\hola_I} \right)^2.
	\end{equation*}
	For the position part: since $\wulb \leq \hola_I$, just write
	$
		\max_I \norm{\Delta \holx_I}^2
		\leq \sum_I \frac{\hola_I}{\wulb} \norm{\Delta \holx_I}^2.
	$

    The last inequality of the claim,
    $\abs{ \norm{\tDelta Z} - \norm{\hZ - Z^*} } \leq \holw_0$,
    follows directly from the definition of the norm $\norm{\cdot}$ in \eqref{eq:growth_conds:def_norm_ZZZ}
    and the definition of $\tDelta \hola = \Delta\hola + \ola_0 a^*$.
\end{proof}

\subsubsection{Step 2: ``Steepness'' of the reduced game}

The following claim is the crucial point of our analysis. It extends \cite[Lemma~14]{wei_linear_2021} to the case of continuous instead of finite games, using the Wasserstein-Fisher-Rao instead of the Fisher-Rao geometry. In particular, the proof crucially relies on the assumption that the MNE is unique.

\begin{claim} \label{claim:growth_conds:steepness_of_reduced_game}
	For all $\tDelta\hZ \in (\ZZZ-Z^*) \setminus \{ 0 \}$, 
	$\max_{\substack{
			A \in \Delta_{n^*}, X \in \XXX^{n^*} \\
			B \in \Delta_{m^*}, Y \in \YYY^{m^*}
	}}
	-\Delta Z^\top \HHH \tDelta \hZ
	> 0$.
\end{claim}

\begin{proof}
    Let $\tDelta \hZ = \begin{pmatrix}
    	\tDelta \halpha \\
    	\Delta \holx \\
    	\tDelta \hbeta \\
    	\Delta \holy
    \end{pmatrix}
    \in (\ZZZ-Z^*) \setminus \{ 0 \}$.
    Suppose by contradiction
    ${
    	% \max_{\substack{
    	% 	A \in \Delta_{n^*}, X \in \XXX^{n^*} \\
    	% 	B \in \Delta_{m^*}, Y \in \YYY^{m^*}
    	% }}
    	\max_{A,X,B,Y}\,
    	-\Delta Z^\top \HHH \tDelta \hZ \leq 0
    }$.
    Since the set
    $\left\lbrace
    	\Delta Z = \begin{pmatrix}
    		\Delta A \\
    		\Delta X \\
    		\Delta B \\
    		\Delta Y
    	\end{pmatrix} \!;
    	A \in \Delta_{n^*}, X \in \XXX^{n^*} \!,
    	B \in \Delta_{m^*}, Y \in \YYY^{m^*} \!
    \right\rbrace$
    contains a (relative) neighborhood of zero 
    since $a^*$ is in the interior of $\Delta_{n^*}$,
%	(since $a^*$ resp.\ $b^*$ is in the interior of $\Delta_{n^*}$ resp.\ $\Delta_{m^*}$),
    clearly the inequality to contradict is equivalent to
    \begin{equation}
    \label{eq:stronglocassumWS:tocontradict}
    	\forall 
    	A \in \Delta_{n^*}, X \in \XXX^{n^*},
    	B \in \Delta_{m^*}, Y \in \YYY^{m^*},~
    	\Delta Z^\top \HHH \tDelta \hZ = 0.
    \end{equation}
    
    Pose for some $\lambda \neq 0$ to be specified
    \begin{equation*}
    	A = a^* + \lambda \tDelta \halpha
    	~~ ~~~~\text{and}~~ ~~~~
    	B = b^* + \lambda \tDelta \hbeta.
    \end{equation*}
    It is straightforward to check that $\sum_I A_I = 1$ and $\sum_J B_J = 1$.
    Moreover since $a^*$ lies in the interior of $\Delta_{n^*}$, $\abs{\lambda}$ can be chosen small enough so that $A_I \geq 0$, and so $A \in \Delta_{n^*}$ (and likewise $B \in \Delta_{m^*}$).
    Further pose
    \begin{equation*}
    	X = x^* + \lambda \Delta \holx
    	~~ ~~~~\text{and}~~ ~~~~
    	Y = y^* + \lambda \Delta \holy.
    \end{equation*}
    Evaluating \eqref{eq:stronglocassumWS:tocontradict} at this $(A, X, B, Y)$ yields
    \begin{align*}
    	0 &= \lambda \begin{pmatrix}
    		\tDelta \halpha \\
    		\Delta \holx \\
    		\tDelta \hbeta \\
    		\Delta \holy
    	\end{pmatrix}^\top
    	\begin{bmatrix}
    		0 & 0 & \gmatss & \partial_y \gmatss b^* \\
    		0 & a^* H_x & a^* \partial_x \gmatss & a^* \partial_{xy}^2 \gmatss b^* \\
    		-(\gmatss)^\top & -(a^* \partial_x \gmatss)^\top & 0 & 0 \\
    		-(\partial_y \gmatss b^*)^\top & -(a^* \partial_{xy}^2 \gmatss b^*)^\top & 0 & b^* H_y
    	\end{bmatrix}
    	\begin{pmatrix}
    		\tDelta \halpha \\
    		\Delta \holx \\
    		\tDelta \hbeta \\
    		\Delta \holy
    	\end{pmatrix} \\
    	&= \tDelta \hZ^\top \HHH \tDelta \hZ \\
    	&= a^* \Delta \holx^\top H_x \Delta \holx
    	+ b^* \Delta \holy^\top H_y \Delta \holy \\
    	&= \sum_I a^*_I \norm{\Delta \holx_I}_{H_I}^2
    	+ \sum_J b^*_J \norm{\Delta \holy_J}_{H_J}^2.
    \end{align*}
    Since $H_I, H_J \succ 0$, this implies that $\Delta \holx = 0$ and $\Delta \holy = 0$.
    So the inequality to contradict reduces to
    \begin{equation}
    \label{eq:stronglocassumWS:tocontradict_reduced}
    	\forall 
    	A \in \Delta_{n^*}, X \in \XXX^{n^*}\!,
    	B \in \Delta_{m^*}, Y \in \YYY^{m^*}\!,~
    	\begin{pmatrix}
    		\Delta A \\
    		\Delta X \\
    		\Delta B \\
    		\Delta Y
    	\end{pmatrix}^\top
    	\begin{bmatrix}
    		0 & \gmatss \\
    		0 & a^* \partial_x \gmatss \\
    		-(\gmatss)^\top & 0 \\
    		-(\partial_y \gmatss b^*)^\top & 0
    	\end{bmatrix}
    	\begin{pmatrix}
    		\tDelta \halpha \\
    %			\Delta \holx \\
    		\tDelta \hbeta \\
    %			\Delta \holy
    	\end{pmatrix} 
    	= 0.
        % ~~~~ ~~
    \end{equation}
    
    % We can then proceed similarly as in the exact-parametrized case.
    Since $\tDelta \hZ \neq 0$ and $\Delta \holx, \Delta \holy = 0$, then w.l.o.g.\ $\tDelta \halpha \neq 0$. We want to show that there exists $\theta>0$ such that, denoting
	\begin{equation*}
		\halpha^\theta = a^* + \theta \tDelta \halpha
		~~~~\text{and}~~~~
		\hmu^\theta 
		= \sum_I \halpha^\theta_I \, \delta_{x^*_I}
		= \mu^* + \theta \sum_I \tDelta\halpha_I \, \delta_{x^*_I},
	\end{equation*}
	$(\hmu^\theta, \nu^*)$ is a MNE, which will contradict uniqueness of the MNE $(\mu^*, \nu^*)$.
	Equivalently, we want to show that the first variation $(\hmu^\theta)^\top F$ is everywhere upper-bounded by $\gval$:
	\begin{equation*}
		\forall y \in \YYY,~ ((\hmu^\theta)^\top F)(y) \leq \gval.
	\end{equation*}
	First remark that:
    \begin{enumerate}
%		\item By evaluating \eqref{eq:stronglocassumWS:tocontradict} at $\overline{\delta}$ such that $\overline{\delta a}, \overline{\delta b}, \overline{\delta x} = 0$ and $\overline{\delta y}$ arbitrary, we have that	
    	\item By the non-degeneracy Assumption~\ref{assum:6}, there exists $\tau>0$ such that
    	\begin{equation*}
    		\forall J, \forall y,~ \norm{y-y^*_J} \eqqcolon \norm{\delta y} \leq \tau \implies ((\mu^*)^\top F)(y) \leq \gval - \frac{1}{4} \sigma_{\min} \norm{\delta y}^2.
    	\end{equation*}
    	\item By 
%    	the third row of \eqref{eq:stronglocassumWS:tocontradict_reduced} evaluated at $B = e_J$,
		\eqref{eq:stronglocassumWS:tocontradict_reduced} evaluated at $A=a^*, X=x^*, Y=y^*$ and $B=e_J$,
    	\begin{align*}
    		\forall J,~
    		& [\tDelta \halpha \gmatss]_J - \tDelta \halpha \gmatss b^* = 0 \\
    		& [\tDelta \halpha \gmatss]_J = \tDelta \halpha \gmatss b^* = 0.
    	\end{align*}
    	\item By 
%    	the fourth row of \eqref{eq:stronglocassumWS:tocontradict_reduced},
		\eqref{eq:stronglocassumWS:tocontradict_reduced} evaluated at $A=a^*, X=x^*, B=b^*, Y_{J'} = y^*_{J'}$ for $J' \neq J$ and $Y_J$ arbitrary,
    	\begin{equation*}
    		\forall J,~
    		[\tDelta \halpha \partial_y \gmatss]_J = 0.
    	\end{equation*}
    \end{enumerate}
    Now,
    \begin{itemize}
        \item 
        Fix $J \in [m^*]$. Let us show that there exists $\theta_{0J}>0$ such that for all $\theta \leq \theta_{0J}$, we have
    	$\forall y \in B_{y^*_J, \tau},~ ((\hmu^\theta)^\top F)(y) \leq \gval$.
    	Indeed for all $y \in B_{y^*_J,\tau}$ and letting $\delta y = y-y^*_J$,
    	\begin{align*}
    		((\hmu^\theta)^\top F)(y)
    		&= ((\mu^*)^\top F)(y) + \theta \sum_I \tDelta \halpha_I f(x^*_I, y) \\
    		&\leq \gval - \frac{1}{4} \sigma_{\min} \norm{\delta y}^2 
    		+ \theta \sum_I \tDelta \halpha_I \left[
    			\gmatss_{IJ} + \partial_y \gmatss_{IJ} \cdot \delta y + O(\norm{\delta y}^2) 
    		\right]
    	\end{align*}
        by point 1. 
    	Now
    	\begin{align*}
    		\sum_I \tDelta \halpha_I \gmatss_{IJ} 
    		&= [\tDelta \halpha \gmatss]_J = 0
    	\end{align*}
    	by point 2 and
    	\begin{equation*}
    		\sum_I \tDelta \halpha_I \partial_y \gmatss_{IJ} 
    		= 0
    	\end{equation*}
    	by point 3.
    	So
    	\begin{align*}
    		((\hmu^\theta)^\top F)(y)
    		&\leq \gval - \frac{1}{4} \sigma_{\min} \norm{\delta y}^2 
    		+ \theta \left[ 
    			0 + 0 \cdot \delta y + O(\norm{\delta y}^2) 
    		\right] \\
    		&\leq \gval - \left( \frac{1}{4} \sigma_{\min} + O(\theta) \right) \norm{\delta y}^2.
    	\end{align*}
    	So clearly we can choose such a $\theta_{0J}$.
    	\item
    	By the non-degeneracy Assumption~\ref{assum:5} there exists $\gvsep_\tau>0$ such that
    	for any $y \in \YYY \setminus \left( \cup_j B_{y^*_J,\tau} \right)$, 
    	$((\mu^*)^\top F)(y) \leq \gval - \gvsep_\tau$.
    	So for all such $y$,
    	\begin{equation*}
    		((\hmu^\theta)^\top F)(y) \leq \gval - \gvsep_\tau + O \left( \theta \norm{\tDelta \halpha} \right).
    	\end{equation*}
    	So we can indeed choose $0< \theta \leq \min_J \theta_{0J}$ that satisfies the requirement.
    \qedhere
    \end{itemize}
\end{proof}

\subsubsection{Step 3: Leverage homogeneity}

\begin{claim} \label{claim:growth_conds:lb_max_bilin_over_WP}
	There exists a constant $C>0$ (dependent only on $(f, \XXX, \YYY)$) such that
	\begin{equation*}
		\forall \tDelta\hZ \in \ZZZ-Z^*,~
		\max_{\substack{
				A \in \Delta_{n^*}, X \in \XXX^{n^*} \\
				B \in \Delta_{m^*}, Y \in \YYY^{m^*}
		}}
		-\Delta Z^\top \HHH \tDelta \hZ
		\geq C \norm{\tDelta \hZ}.
	\end{equation*}
\end{claim}

\begin{proof}
	Let for concision 
	$g(\tDelta \hZ) = \max_{\substack{
			A \in \Delta_{n^*}, X \in \XXX^{n^*} \\
			B \in \Delta_{m^*}, Y \in \YYY^{m^*}
	}}
	-\Delta Z^\top \HHH \tDelta \hZ$.
	Note that $g$ is continuous and positive-homogeneous.
	
	Since $a^*$ resp.\ $b^*$ lies in the interior of its domain, it is not hard to check that there exists $r>0$ (only dependent on $a^*, b^*$) such that, for any $\tDelta\hZ \in \ZZZ-Z^*$, then $Z^* + r \frac{\tDelta \hZ}{\norm{\tDelta \hZ}} \in \ZZZ$.
	In other words, any $\tDelta\hZ \in \ZZZ-Z^*$ can be written as $\frac{\norm{\tDelta \hZ}}{r} \Delta \tZ$ for some $\tZ \in \SSS_{Z^*,r} \coloneqq \left\lbrace \tZ \in \ZZZ; \norm{\Delta \tZ} = r \right\rbrace$.
	
	Now by \autoref{claim:growth_conds:steepness_of_reduced_game},
    $g(\Delta \tZ) > 0$ for all $\tZ \in \SSS_{Z^*,r}$.
	Since $\SSS_{Z^*,r}$ is a compact set and $g$ is continuous, we have
	\begin{align*}
		\forall \tZ \in \ZZZ ~~\text{s.t.}~~ \norm{\Delta \tZ} = r,~ &
		g(\Delta \tZ) \geq \inf_{\SSS_{Z^*,r}} g > 0 \\
		\text{so by positive-homogeneity,} \quad\quad\quad
		\forall \tDelta\hZ \in \ZZZ-Z^*,~ &
		g(\tDelta \hZ) \geq \left( \inf_{\SSS_{Z^*,r}} g \right) \frac{\norm{\tDelta \hZ}}{r} \eqqcolon C \norm{\tDelta \hZ}.
    \rqedhere
	\end{align*}
\end{proof}

%\begin{remark}
%	For small enough $r$, $\frac{1}{r} \left( \inf_{\SSS_{Z^*,r}} g \right)$ does not depend on $r$.
%\end{remark}

\autoref{lm:growth_conds:error_bound} follows by
% combining the four above claims.
using \autoref{claim:growth_conds:DeltaZ_equiv_etaveesigmaV} to further lower-bound the result of \autoref{claim:growth_conds:lb_max_bilin_over_WP}, and by substituting into \autoref{claim:growth_conds:tgaplocs_proxyz}.

% \subsection{Proof of ``local star-convexity-concavity''} \label{subsec:growth_conds:starconvex}
\subsection{Proof of ``local star-convexity-concavity''} \label{subsec:growth_conds:starconvex}

\begin{lemma} \label{lm:growth_conds:starconvex}
    Consider the Lyapunov function $V_1$ as in \eqref{eq:cv_proof:gencase:def_V} with the partitions of unity $(\varphi_I)_I$ and $(\psi_J)_J$ as in \eqref{eq:cv_proof:gencase:def_varphi}.%
    \footnote{
        For this lemma, the precise choice of the partitions of unity $(\varphi_I)_I$ and $(\psi_J)_J$ does not matter, only the fact that $\support(\varphi_I) \subset B_{x^*_I, \lambda\tau}$ and $\support(\psi_J) \subset B_{y^*_J, \lambda\tau}$.
    }
    
    Consider any $\hz = (\ha, \hx, \hb, \hy) \in \Delta_n \times \XXX^n \times \Delta_m \times \YYY^m$,
    and let $\zps = (\aps, \xps, \bps, \yps)$ ``proxy solution particles'' similarly as in \eqref{eq:cv_proof:gencase:choice_of_proxy*}:
    \begin{equation*}
    	\xps_i \coloneqq \hx_i + \sum_{I \in [n^*]} \hvarphi_{Ii} (x^*_I - \hx_i)
    	~~~\text{and}~~~
    	\aps_i \coloneqq \sum_{I \in [n^*]} a^*_I \frac{\hvarphi_{Ii} \ha_i}{\hola_I},
    \end{equation*}
    and similarly for $\aps, \bps$.
    Suppose $\left( \min_I \hola_I \right) \wedge \left( \min_J \holb_J \right) \geq c$ for some $c>0$.
    Then, denoting $\hmu = \sum_{i=1}^n \ha_i \delta_{\hx_i}$ and $\hnu = \sum_{j=1}^m \hb_j \delta_{\hy_j}$,
    \begin{align*}
    	\tgap(\zps; \hz)
    	% &= F(\hmu, \nu^*) - F(\mu^*, \hnu) 
        &= \ha^\top \gmaths b^* - (a^*)^\top \gmatsh \hb \\
        &~~~~ 
		+ \frac{1}{2} \sum_I \sum_i \hvarphi_{Ii} \ha_i \norm{\hx_i-x^*_I}_{H_I}^2 
    	+ \frac{1}{2} \sum_J \sum_j \hpsi_{Jj} \hb_j \norm{\hy_j-y^*_J}_{H_J}^2 \\
    	&~~~~ + O \left( c^{-1} V_1(\hz)^{3/2} \right)
    	+ \bR
        \qquad\qquad\qquad\qquad
         \text{where}~~~~
         \abs{\bR} \leq 2 \smoothnessf_3 c^{-1} \cdot \lambda \tau \cdot V_\pos(\hz) \\
         &\geq F(\hmu, \nu^*) - F(\mu^*, \hnu) 
         + V_\pos(\hz) \left( \sigma_{\min} - 2 \smoothnessf_3 c^{-1} \cdot \lambda \tau \right)
         + O \left( c^{-1} V_1(\hz)^{3/2} \right).
    \end{align*}
\end{lemma}

Recall from \autoref{lm:growth_conds:ulb_on_weights_locally} that we can always ensure $\left( \min_I \hola_I \right) \wedge \left( \min_J \holb_J \right) \geq \wulb$ for some constant $\wulb>0$ by assuming $V_1(\hz) \leq r$ for some constant $r>0$.
In order to recover the informal statement of \autoref{subsubsec:cv_proof:proofingr_growthconds:starconvex},
note that if in addition $\lambda\tau \leq \frac{\wulb~ \sigma_{\min}}{4 \smoothnessf_3}$,
% i.e., $\sigma_{\min} - 2 \smoothnessf_3 \wulb^{-1} \cdot \lambda \tau \geq \frac{\sigma_{\min}}{2}$,
then
\begin{equation*}
    \tgap(\zps; \hz)
    \geq F(\hmu, \nu^*) - F(\mu^*, \hnu) 
    + \frac{\sigma_{\min}}{2} V_\pos(\hz) + O \left( V_1(\hz)^{3/2} \right).
\end{equation*}

The proof proceeds by Taylor expansions to estimate 
% $\tgap(\zps; \hz) - \left[ \ha^\top \gmaths b^* - (a^*)^\top \gmatsh \hb \right]$.
$\tgap(\zps; \hz)$.
This involves rather tedious calculations, 
\ifextended%
    so we defer it to 
    \autoref{subsec:tedious_calcs:proof__lm--growth_conds--starconvex}.
\else%
    so we place it as \AdditionalMaterial\ in Section~G.2 of the arXiv version of this paper.
\fi
In a nutshell, we do Taylor expansions of 
% $\gmathh_{ij} = f(\hx_i, \hy_j)$
$f(\hx_i, \hy_j)$
around $(x^*_I, \hy_j)$ or $(\hx_i, y^*_J)$.
In order to make 
% $a^* \gmatsh \hb$ 
% $\ha^\top \gmaths b^* - (a^*)^\top \gmatsh \hb$
$\ha^\top \gmaths b^*$ and $(a^*)^\top \gmatsh \hb$
appear, at first we only expand the side with $\zps - \hz$ --- e.g., when estimating the terms 
$\innerprod{\nabla_a F_{n,m}(\hz)}{\ha-\aps}$ 
and $\innerprod{\nabla_x F_{n,m}(\hz)}{\hx-\xps}$,
start by expanding only with respect to $x$ and keeping $\hy, \hb$ as is.
Once the expansion of $\tgap(\zps, \hz)$ (the equality) is proved, the lower bound follow straightforwardly.

% \end{document}

% !TEX root = ../main.tex
% \documentclass[../main]{subfiles}
% \begin{document}

\section{Proof of the relation between Lyapunov function and NI error}
% \section{Relation between Lyapunov function and NI error} 
\label{sec:rel_Lya_NI}

%The exact-parametrization case, \autoref{prop:cv_proof:exactparam:NI_equiv_V}, follows from exactly the same arguments as in the general case, except we apply \autoref{lm:growth_conds:quadr_exactparam} instead of \autoref{lm:growth_conds:quadr_gencase} to lower-bound $F(\hmu, \nu^*) - F(\mu^*, \hnu)$.
%The difference in exponents between the lower bounds in the two lemmas comes from \TODO{ref to a remark in the course of the proof for the general case}.

In this section, we show that the Lyapunov function can be used as a proxy for the NI error.
We present the proof for the general case (\autoref{prop:cv_proof:gencase:NI_equiv_V}), and describe in \autoref{subsec:rel_Lya_NI:exactparam} the necessary adaptations to prove the proposition for the exact-parametrization case (\autoref{prop:cv_proof:exactparam:NI_equiv_V}).
% follows as a limit case when $\lambda\tau$ is held constant and $\tau \to \infty$, along exactly the same lines as in the proof of \autoref{lm:growth_conds:quadr_exactparam}.

As announced in the main text, here is a more quantitative version of \autoref{prop:cv_proof:gencase:NI_equiv_V}.

\begin{proposition} \label{prop:rel_Lya_NI:NI_equiv_V_moreprecise}
	Define $V_1$ as in \eqref{eq:cv_proof:gencase:def_V} with the partitions of unity $(\varphi_I)_I$ and $(\psi_J)_J$ as in \eqref{eq:cv_proof:gencase:def_varphi}.
	Suppose that
    $\lambda\tau \leq \frac{\sigma_{\min}}{2 \smoothnessf_3}$.
	
	There exist constants $C_1, C_2$ dependent on $(f, \XXX, \YYY)$ and a constant $K$ dependent on $\lambda, \tau$ such that,
	for any $\hz = (\ha, \hx, \hb, \hy) \in \Delta_n \times \XXX^n \times \Delta_m \times \YYY^m$,
    denoting $\hmu = \sum_i \ha_i \delta_{\hx_i}$ and $\hnu = \sum_j \hb_j \delta_{\hy_j}$,
	\begin{equation*}
    	C_1 K \left[ \left( \min_I \hola_I \wedge \min_J \holb_J \right) V_1(\hz) \right]^{5/4} 
		\leq \NI(\hmu, \hnu)
		\leq C_2 \sqrt{V_1(\hz)}.
	\end{equation*}
	
	Moreover, there exists $r>0$ dependent only on $(f, \XXX, \YYY)$ such that, if $\NI(\hmu, \hnu) \leq r K$, then $\min_I \holw_I \geq \wulb = \frac{w^*_{\min}}{4}$, and so
	\begin{equation*}
		C_1 \wulb^{5/4} ~ K ~ V_1(\hz)^{5/4} \leq \NI(\hmu, \hnu).
	\end{equation*}
\end{proposition}

The expression of $K$ can be found in \eqref{eq:rel_Lya_NI:def_K}. In particular, if $\lambda, \tau$ are chosen as in the proof of convergence for the general case \eqref{eq:proof_gencase:choice_lambda_tau}, then $K \asymp \sqrt{\sigma}$.

The rest of this section is dedicated to proving the above proposition, with the exception of the last subsection where we deal with the exact-parametrization case.

\subsection{Proof of the first inequality: \texorpdfstring{$\NI \lesssim \sqrt{V_1}$}{NI lesssim sqrt(V1)}}

By bilinearity of $F(\mu, \nu)$, for any
$\hmu = \sum_{i=1}^n \ha_i \delta_{\hx_i}$, $\hnu = \sum_{j=1}^m \hb_j \delta_{\hy_j}$,
\begin{align*}
	\NI(\hmu, \hnu) 
	= \max_{\mu, \nu} F(\hmu, \nu) - F(\mu, \hnu)
	&= \max_{\mu, \nu} \int_\YYY \sum_i \ha_i f(\hx_i, \cdot) d\nu - \int_\XXX \sum_j \hb_j f(\cdot, \hy_j) d\mu \\
	&= \max_{x, y} \sum_i \ha_i f(\hx_i, y) - \sum_j \hb_j f(x, \hy_j).
% 	&= \max_{x, y} F(\hmu, \delta_y) - F(\delta_x, \hnu) \\
% 	&= \max_{x, y} F_{n, 1}(\ha, \hx, (1), y) - F_{1, m}((1), x, \hb, \hy).
\end{align*}
Now, denoting $\hvarphi_{Ii} = \varphi_I(\hx_i)$, for any $y \in \YYY$
\begin{align*}
% 	F_{n, 1}(\ha, \hx, (1), y) &= \sum_i \ha_i f(\hx_i, y) \\
    \sum_i \ha_i f(\hx_i, y)
	&= \sum_I \sum_i \hvarphi_{Ii} \ha_i f(\hx_i, y) + \sum_i \hvarphi_{0i} \ha_i f(\hx_i, y) \\
	&= \sum_I \sum_i \hvarphi_{Ii} \ha_i \left( f(x^*_I, y) + O(\norm{\hx_i-x^*_I}) \right) 
	+ O(\hola_0) \\
	&= \sum_I \hola_I f(x^*_I, y) 
    + O\left( 
        % \sqrt{ \sum_I \sum_i \hvarphi_{Ii} \ha_i \norm{\hx_i-x^*_I}^2 } 
        \sum_I \sum_i \hvarphi_{Ii} \ha_i \norm{\hx_i-x^*_I}
    \right) 
    + O(\hola_0) \\
	&\leq \gval + O\left( \norm{\Delta \hola}_1 \right) + O\left( \sqrt{V_\pos(\ha, \hx)} \right)+ O(\hola_0)
	= \gval + O\left( \sqrt{V_1(\ha, \hx)} \right)
\end{align*}
where we used Jensen's inequality on $s\mapsto s^2$ and \eqref{eq:aux_lemmas:explicit_Vposax} to bound $\sum_{I,i} \hvarphi_{Ii} \ha_i \norm{\hx_i-x^*_I}$.
Similarly, for any $x \in \XXX$,
$
% 	F_{1, m}((1), x, \hb, \hy) 
    \sum_j \hb_j f(x, \hy_j)
	\geq \gval + O\left( \sqrt{V_1(\hb, \hy)} \right)
$.
Hence
\begin{equation*}
	\NI(\hmu, \hnu)
% 	= \max_{x, y} F_{n, 1}(\ha, \hx, (1), y) - F_{1, m}((1), x, \hb, \hy)
	= \max_{x, y} \sum_i \ha_i f(\hx_i, y) - \sum_j \hb_j f(x, \hy_j)
	\lesssim \sqrt{V_1(\hz)}.
\end{equation*}
This shows the first inequality in \autoref{prop:rel_Lya_NI:NI_equiv_V_moreprecise}.

\subsection{Proof of the second inequality: \texorpdfstring{$\NI \gtrsim \left[ \left( \min_I \ola_I \wedge \min_J \olb_J \right) V_1 \right]^{5/4}$}{NI gtrsim [(min aI wedge min bJ) V1]5/4}}

% \paragraph{Lower bound on $F(\hmu, \nu^*) - F(\mu^*, \hnu)$.}
\paragraph{Lower bound on ``gap'' to solution $(\mu^*, \nu^*)$.}

\autoref{lm:growth_conds:quadr_gencase} directly implies a lower bound on 
$\NI(\hmu, \hnu) = \max_{\mu, \nu} F(\hmu, \nu) - F(\mu, \hnu) \geq F(\hmu, \nu^*) - F(\mu^*, \hnu)$.
% , namely
% \begin{align*}
%     F(\hmu, \nu^*) - F(\mu^*, \hnu) 
%     % = \ha \gmaths b^* - a^* \gmatsh \hb
% 	&\geq 
% 	\max \left\{
% 		\left[
% 		    \frac{\sigma_{\min}}{2} \wedge \frac{2 \gvsep}{(\lambda\tau)^2}
% 		\right]
% 		V_\pos(\ha, \hx, \hb, \hy), ~~
% 		\left[ 
% 			\frac{\sigma_{\min}}{4} \frac{3 (\lambda\tau)^2}{\lambda^3} \wedge \gvsep
% 		\right]
% 		(\hola_0 + \holb_0)
% 	\right\}
% \end{align*}
% for some constant $\gvsep>0$ only dependent on $(f, \XXX, \YYY)$.
For concision, within this section, denote $K$ the constant appearing in that lower bound:
\begin{equation} \label{eq:rel_Lya_NI:def_K}
    \NI(\hmu, \hnu) \geq K \left( \holw_0 + V_\pos(\hz) \right)
    ~~\text{where}~~
    K = \frac{1}{2} \left(
        \frac{\sigma_{\min}}{2} \wedge \frac{2 \gvsep}{(\lambda\tau)^2} 
        \wedge \frac{\sigma_{\min}}{4} \frac{3 (\lambda\tau)^2}{\lambda^3} \wedge \gvsep
    \right)
\end{equation}
and where $\gvsep>0$ is a constant only dependent on $(f, \XXX, \YYY)$.
It remains to lower-bound $\NI(\hmu, \hnu)$ by $\norm{\Delta \holw}_1$ with some exponent.

% \paragraph{Lower bound on $\max_{A,X,B,Y} F(\ha, \hx, B, Y) - F(A, X, \hb, \hy)$.}
\paragraph{Lower bound on maximum ``gap'' to perturbations of the solution.}

In the remainder of this section, we adopt again the notations of \autoref{subsec:growth_conds:error_bound} (Eq.~\eqref{eq:growth_conds:def_norm_ZZZ}) for the set $\ZZZ$, the norm on $\ZZZ - Z^*$, and the shorthand $\tDelta \hZ \in \ZZZ-Z^*$.

\begin{claim} \label{claim:rel_Lya_NI:Fnm-Fnm_bilinear}
    Suppose $\lambda\tau \leq \frac{\sigma_{\min}}{4 \smoothnessf_3}$.
	For any $A \in \Delta_{n^*}, X \in \XXX^{n^*}, B \in \Delta_{m^*}, Y \in \YYY^{m^*}$,
	\begin{multline*}
		F_{n,m^*}(\ha, \hx, B, Y) - F_{n^*,m}(A, X, \hb, \hy) \\
		\geq -\begin{pmatrix}
			\Delta A \\
			\Delta X \\
			\Delta B \\
			\Delta Y
		\end{pmatrix}^\top
    	\begin{bmatrix}
    		0 & 0 & \gmatss & \partial_y \gmatss b^* \\
			0 & a^* \frac{\sigma_{\min}}{2} \id & a^* \partial_x \gmatss & a^* \partial_{xy}^2 \gmatss b^* \\
    		-(\gmatss)^\top & -(a^* \partial_x \gmatss)^\top & 0 & 0 \\
    		-(\partial_y \gmatss b^*)^\top & -(a^* \partial_{xy}^2 \gmatss b^*)^\top & 0 & b^* \frac{\sigma_{\min}}{2} \id
    	\end{bmatrix}
		\begin{pmatrix}
			\tDelta \hola \\
			\Delta \holx \\
			\tDelta \holb \\
			\Delta \holy
		\end{pmatrix} \\
		+ O \left(
	        \norm{\Delta \hZ}^3
	       % + \norm{\Delta Z} \left( \holw_0 + V_\pos(\hz) \right)
	        + \left( \holw_0 + V_\pos(\hz) \right)^2
	        + \norm{\Delta Z}^2
	    \right)
	\end{multline*}
% 	where $\norm{\Delta Z} = \norm{\Delta A}_1 + \max_I \norm{\Delta X_I} + \norm{\Delta B}_1 + \max_J \norm{\Delta Y_J}$.
    where $\tDelta \hola = \Delta \hola + \ola_0 a^*$ and $\tDelta \holb = \Delta \holb + \holb_0 b^*$,
    and where we denote $\left[ a^* \frac{\sigma_{\min}}{2} \id \right]_{I I'} = \ind_{I=I'} a^*_I \frac{\sigma_{\min}}{2} \id_{\XXX}$ for each $I, I'$.
\end{claim}

The proof of this claim follows from
\ifextended%
    simple but tedious calculations, which we defer to
    \autoref{subsec:tedious_calcs:proof__claim--rel_Lya_NI--Fnm-Fnm_bilinear}.
\else%
    simple but tedious calculations. For full details, one can refer to Section~G.3 of the arXiv version of this paper as \AdditionalMaterial.
\fi
Essentially, we do Taylor expansions of $f$ around $(x^*_I, y^*_J)$, rearrange the terms so as to get an expression which is order-$2$-exact in $\Delta \hZ$ and order-$1$-exact in $\Delta Z$, and check that the remaining terms are non-negative or negligible.

%\TODO{idea: follow the same reasoning as Claim~1.5 of ``analysis of particle PPA, general case v2.3'', but restrict $\norm{\Delta Z}$ (the important thing is to restrict $\norm{\Delta X}$ and $\norm{\Delta Y}$) to be of order at most $\norm{\Delta \hZ}^{3/2}$, so that the error in the above Claim is of order 3; note that this will scale up the lower bound into $\norm{\Delta \hZ} * \norm{\Delta \hZ}^{3/2} = \norm{\Delta \hZ}^{5/2}$, which is still dominant compared to the order-3 error (but not great) (but that's the best I can do AFAICT). Note that some care will be required (in fact I'm still not 100\% sure it works, I really hope it does) when doing this "restrict the test $\Delta Z$ to have small norm, namely controlled by $\norm{\Delta \hZ}$ itself" argument.} --> done, yep it works

Denote $\HHH$ the 
% symbolic
block
matrix in the above claim.
We reuse the result of \autoref{claim:growth_conds:lb_max_bilin_over_WP} --- actually it was proved for a slightly different $\HHH$, with $\frac{\sigma_{\min}}{2} \id$ replaced by $H_x$ resp.\ $H_y$,
% resp.\ $\frac{\sigma_{\min}}{2} \id$, 
but the proof can be very easily adapted ---:
There exist a constant $C>0$ (dependent only on $(f, \XXX, \YYY)$) such that
\begin{equation*}
	\forall \tDelta\hZ \in \ZZZ-Z^*,~
	\max_{\substack{
			A \in \Delta_{n^*}, X \in \XXX^{n^*} \\
			B \in \Delta_{m^*}, Y \in \YYY^{m^*}
	}}
	-\Delta Z^\top \HHH \tDelta \hZ
	\geq C \norm{\tDelta \hZ}.
\end{equation*}
We further refine that result slightly by also exploiting the positive-homogeneity with respect to $\Delta Z$ instead of just $\tDelta \hZ$.

\begin{claim} \label{claim:rel_Lya_NI:lb_DeltaZleqq}
	There exist constants $q_1, C>0$ dependent only on $(f, \XXX, \YYY)$ such that for any $q \leq q_1$,
	\begin{equation*}
		\forall \tDelta\hZ \in \ZZZ-Z^*,~
		\max_{\norm{\Delta Z} \leq q} -\Delta Z^\top \HHH \tDelta \hZ \geq C q \norm{\tDelta \hZ}.
	\end{equation*}
\end{claim}

\begin{proof}
    Recall the notation $\ZZZ = \Delta_{n^*} \times \XXX^{n^*} \times \Delta_{m^*} \times \YYY^{m^*}$. Denote $\Delta \ZZZ = \ZZZ - Z^*$ and
	\begin{equation*}
		D = \left\lbrace
			\begin{pmatrix}
				\Delta A \\
				\Delta X \\
				\Delta B \\
				\Delta Y
			\end{pmatrix}
			\in \RR^{n^*} \times \XXX^{n^*} \times \RR^{m^*} \times \YYY^{m^*}
			~;~
			\sum_I \Delta A_I = \sum_J \Delta B_J = 0
			~\text{and}~
			\norm{\Delta Z} \leq 1
		\right\rbrace.
	\end{equation*}
%     ($\hZ \in \ZZZ$ but $Z \in \tZZZ$.)
% 
% 	We just rescale the inequality from \autoref{claim:growth_conds:lb_max_bilin_over_WP}.
% 	More precisely:
	Since $a^*, x^*, b^*, y^*$ lie in the (relative) interior of their domains, then there exist $q_1, q_2>0$ such that
	$
		q_1 D \subset \Delta \ZZZ \subset q_2 D
	$.
	So, using \autoref{claim:growth_conds:lb_max_bilin_over_WP} for the second inequality,
	\begin{align*}
		\max_{v \in q_2 D} -v^\top \HHH \tDelta \hZ
		&\geq
		\max_{v \in \Delta \ZZZ} -v^\top \HHH \tDelta \hZ 
		= \max_{Z \in \ZZZ} -\Delta Z^\top \HHH \tDelta \hZ
        \geq C \norm{\tDelta \hZ} \\
		\max_{v \in D} -v^\top \HHH \tDelta \hZ
		= \frac{1}{q_2} \max_{v \in q_2 D} -v^\top \HHH \tDelta \hZ
		&\geq \frac{1}{q_2} C \norm{\tDelta \hZ}
	\end{align*}
	and so, for any $q \leq q_1$, since $q D \subset q_1 D \subset \Delta \ZZZ$,
	\begin{align*}
		\max_{Z \in \ZZZ ~\text{s.t.}~ \norm{\Delta Z} \leq q} -\Delta Z^\top \HHH \tDelta \hZ
		= \max_{v \in q D \cap \Delta \ZZZ} -v^\top \HHH \tDelta \hZ 
		&= \max_{v \in q D} -v^\top \HHH \tDelta \hZ \\
		&= q \max_{v \in D} -v^\top \HHH \tDelta \hZ
        \geq \frac{q}{q_2} C \norm{\tDelta \hZ}.
    \rqedhere
	\end{align*}
\end{proof}

With this we can prove the following lower bound on $\NI(\hmu, \hnu)$:
\begin{claim} \label{lm:rel_Lya_NI:lb_maxgap_beta}
	There exists a constant $C$ dependent only on $(f, \XXX, \YYY)$ such that, for any $\beta \geq 1$,
	\begin{multline*}
		\max_{A, X, B, Y} F_{n,m^*}(\ha, \hx, B, Y) - F_{n^*,m}(A, X, \hb, \hy) \\
		\geq C \norm{\tDelta \hZ}^{1+\beta}
		+ O \left( 
			\norm{\tDelta \hZ}^{2\beta}
			+ \norm{\tDelta \hZ}^3 + \left( \holw_0 + V_\pos(\hz) \right)^2
		\right).
	\end{multline*}
	
    In particular for $\beta=3/2$, we have that
    \begin{align*}
    	\NI(\hmu, \hnu) &\geq C \norm{\tDelta \hZ}^{5/2} + O \left( \norm{\tDelta \hZ}^3 + \left( \holw_0 + V_\pos(\hz) \right)^2 \right) \\
        &\geq C' \norm{\Delta \hZ}^{5/2} + O \left( \norm{\Delta \hZ}^3 + \left( \holw_0 + V_\pos(\hz) \right)^2 \right)
    \end{align*}
    where $C' = 2^{1-\frac52}\, C$,
    and there exists $r>0$ such that
    \begin{equation*}
    	\forall (\ha, \hx, \hb, \hy) ~\text{s.t.}~
    	\norm{\Delta \hZ} \leq r,~
    	\NI(\hmu, \hnu)
    	\geq \frac{C'}{2} \norm{\Delta \hZ}^{5/2} + O \left( \left( \holw_0 + V_\pos(\hz) \right)^2 \right).
    \end{equation*}
\end{claim}

\begin{proof}
	Note that for any $(\ha, \hx, \hb, \hy)$, $\norm{\tDelta \hZ} \leq 4 + 2\diamXY \eqqcolon R'$.
	To prove the first part of the claim, for any fixed $(\ha, \hx, \hb, \hy)$, simply apply \autoref{claim:rel_Lya_NI:lb_DeltaZleqq} with $q = q_1 \frac{\norm{\tDelta \hZ}^\beta}{(R')^\beta}$ and substitute into \autoref{claim:rel_Lya_NI:Fnm-Fnm_bilinear}.
    Note that whether the error term in the $O(\cdot)$ is a power of $\norm{\Delta \hZ}$ or of $\norm{\tDelta \hZ}$ is irrelevant, since they differ by at most $\holw_0$ by \autoref{claim:growth_conds:DeltaZ_equiv_etaveesigmaV} and the $O(\cdot)$ already contains a term in $\holw_0^2$.
	
	The second part of the claim follows by substituting $\beta=3/2$; for the second inequality we used that
    $\forall A, B>0$, $\left( \frac{A+B}{2} \right)^{5/2} \leq \frac12 (A^{5/2} + B^{5/2})$
    by Jensen's inequality,
    and so
    \begin{equation*}
        \norm{\Delta \hZ}^{5/2} \leq \left( \norm{\tDelta \hZ} + \holw_0 \right)^{5/2} \leq 2^{-1+\frac52} \left( \norm{\tDelta \hZ}^{5/2} + \holw_0^{5/2} \right) = 2^{-1+\frac52} \norm{\tDelta \hZ}^{5/2} + O(\holw_0^2).
    \end{equation*}

    For the third part of the claim, let $C_1$ denote the constant hidden in the $O(\cdot)$ in the second inequality of the second part, and pick $r$ such that
	\begin{equation*}
	    C' \norm{\Delta \hZ}^{5/2} - C_1 \norm{\Delta \hZ}^3 
	    = \norm{\Delta \hZ}^{5/2} \left( C' - C_1 \norm{\Delta \hZ}^{1/2} \right) 
	    \geq \norm{\Delta \hZ}^{5/2} \frac{C'}{2}
	\end{equation*}
	for any $\norm{\Delta \hZ} \leq r$.
\end{proof}

This gives the desired bound on a neighborhood of $(\mu^*, \nu^*)$. Outside of that neighborhood, we can simply use that $\NI$ is non-zero and continuous and has a compact domain in the relevant topology; though some care is needed to ensure constants independent of $n$ and $m$.
% In the following claim, we let with abuse of notation $\NI(\ha, \hx, \hb, \hy) = \NI(\hmu, \hnu)$.
\begin{claim}
% 	For any $r>0$, denoting 
% 	$\SSS_{n,m,Z^*,r} = \left\lbrace
% 		(\ha, \hx, \hb, \hy); \norm{\Delta \hZ} \geq r
% 	\right\rbrace$,
% 	we have
% 	\begin{equation*}
% % 		\inf_{\SSS_{n,m,Z^*,r}} \NI(\ha, \hx, \hb, \hy) > 0.
% 		\inf_{\SSS_{n,m,Z^*,r}} \NI(\hmu, \hnu) > 0.
% 	\end{equation*}
    For any $r>0$, there exists a constant $C'>0$ dependent only on $(f, \XXX, \YYY)$ such that 
    $\forall (\ha, \hx, \hb, \hy) ~\text{s.t.}~ \norm{\Delta \hZ} > r,~ \NI(\hmu, \hnu) \geq C'$.
\end{claim}

\begin{proof}
    Fix $r>0$ and let us show that the set
    $\bigcup_{n, m \in \NN^*} \Big\{ \NI(\hmu, \hnu);~ (\ha, \hx, \hb, \hy) \in \Delta_n \times \XXX^n \times \Delta_m \times \YYY^m ~\text{and}~ \norm{\Delta \hZ} \geq r \Big\}$
    is bounded below by a positive constant $C'$. 
    Suppose the contrary, that is, that there exist sequences $(n_p)_p$, $(m_p)_p$, and $(\ha_p, \hx_p, \hb_p, \hy_p) \in \Delta_{n_p} \times \XXX^{n_p} \times \Delta_{m_p} \times \YYY^{m_p}$ for each $p$,
    such that $\norm{\Delta \hZ_p} \geq r$ and $\lim_{p \to \infty} \NI(\hmu_p, \hnu_p) = 0$.
    Since $\XXX$ and $\YYY$ are compact, then $\{\hmu_p\}_p$ and $\{\hnu_p\}_p$ are tight, so up to taking a subsequence, we may assume $\hmu_p$ and $\hnu_p$ converge weakly to some $\hmu_\infty$ resp.\ $\hnu_\infty$.
    Recall from their definition in \autoref{subsec:cv_proof:gencase} that the quantities $\ola, \olx$ are well-defined in terms of moments of $\mu$ (independently of any notion of particle parametrization), so the vectors $\hZ_p$ also converge to some $\hZ_\infty$, with in particular $\norm{\Delta \hZ_\infty} \geq r$.
    Now one can check that $\NI: (\hmu, \hnu) \mapsto \max_y (\hmu^\top F)(y) - \min_x (F \hnu)(x)$ is weakly continuous over $\PPP(\XXX) \times \PPP(\YYY)$.
    % , since each term is weakly continuous as the maximum (or minimum) or uniformly weakly continuous functions. \TODO{check this.} --> yes, this is valid, but long and annoying to detail. Can use essentially the same idea as in https://math.stackexchange.com/a/66462
    So 
    $\NI(\hmu_\infty, \hnu_\infty) = 0$, so $(\hmu_\infty, \hnu_\infty) = (\mu^*, \nu^*)$ since the MNE is assumed unique, so $\hZ_\infty = Z^*$ by considering the corresponding moments; but this contradicts the fact that $\norm{\Delta \hZ_\infty} \geq r >0$.
    % This concludes the proof by contradiction.
\end{proof}

% \begin{proof}
%     It suffices to show that $\inf_{\SSS_{n,m,Z^*,r}} \NI(\hmu, \hnu) > 0$ where $\SSS_{n,m,Z^*,r} = \left\lbrace
% 		(\ha, \hx, \hb, \hy); \norm{\Delta \hZ} \geq r
% 	\right\rbrace$.

% 	The set $\SSS_{n,m,Z^*,r}$ is closed as a preimage by the continuous function $(\ha, \hx, \hb, \hy) \mapsto \norm{\Delta \hZ}$, and compact as a closed subset of the compact $\Delta_n \times \XXX^n \times \Delta_m \times \YYY^m$.
% 	Furthermore, the mapping
% 	\begin{equation*}
% 		(\ha, \hx, \hb, \hy) \mapsto \NI(\hmu, \hnu) = \max_y \sum_i \ha_i f(\hx_i, y) - \min_x \sum_j \hb_j f(x, \hy_j)
% 	\end{equation*}
% 	is continuous, since each term is continuous as the maximum (or minimum) of uniformly continuous functions.
% % 	\footnote{\fnsurl{https://math.stackexchange.com/a/66462}}
% 	It just remains to check that $\NI$ is non-zero on $\SSS_{n,m,Z^*,r}$, which follows from uniqueness of the MNE since $\NI(\hmu, \hnu) = 0 \implies (\hmu, \hnu) = (\mu^*, \nu^*) \implies \hZ = Z^*$.
% \end{proof}

Putting together the two above claims, and using that $\norm{\Delta \hZ}$ is anyway bounded by $4+2\diamXY= O(1)$, we have shown that there exists a constant $C>0$ dependent only on $(f, \XXX, \YYY)$ such that
\begin{equation} \label{eq:rel_Lya_NI:lb_NI_DeltaZ_global}
	\NI(\hmu, \hnu) \geq C \norm{\Delta \hZ}^{5/2} + O \left( \left( \holw_0 + V_\pos(\hz) \right)^2 \right)
\end{equation}
for any $(\ha, \hx, \hb, \hy)$.

\paragraph{Conclusion.}

In the first paragraph we have shown ($K$ is given by Eq.~\eqref{eq:rel_Lya_NI:def_K})
\begin{equation*}
	\NI(\hmu, \hnu)
	\geq F(\hmu, \nu^*) - F(\mu^*, \hnu) 
	\geq K \left( \holw_0 + V_\pos(\hz) \right).
\end{equation*}
In the second paragraph we have shown (Eq.~\eqref{eq:rel_Lya_NI:lb_NI_DeltaZ_global})
\begin{equation*}
	\NI(\hmu, \hnu)
	\geq \max_{A,X,B,Y} F(\ha, \hx, B, Y) - F(A, X, \hb, \hy) 
	\geq C \norm{\Delta \hZ}^{5/2} - C_1 \left( \holw_0 + V_\pos(\hz) \right)^2
\end{equation*}
for some $C, C_1>0$ only dependent on $(f, \XXX, \YYY)$;
further, using that $\norm{\Delta \hZ}^2 \geq 2 (\min_I \holw_I) \cdot \sum_I d_h(w^*_I, \holw_I) + \sum_I \holw_I \norm{\Delta \holp_I}^2$ (\autoref{claim:growth_conds:DeltaZ_equiv_etaveesigmaV}),
we have
\begin{align*}
	\NI(\hmu, \hnu) 
	&\geq 
	C \left( 2 \,\big(\!\min_I \holw_I\big) \sum_I d_h(w^*_I, \holw_I) + \sum_I \holw_I \norm{\Delta \holp}^2 \right)^{5/4} 
	- C_1 \left( \left( \holw_0 + V_\pos(\hz) \right)^2 \right) \\
	&\geq 
	C \left( 2 \,\big(\!\min_I \holw_I\big) \sum_I d_h(w^*_I, \holw_I) + \sum_I \holw_I \norm{\Delta \holp}^2 \right)^{5/4} 
	- \underbrace{C_1 (2 + \diamXY)}_{\eqqcolon C_2} \left( \holw_0 + V_\pos(\hz) \right).
\end{align*}

Taking a convex combination of these two inequalities with ratio $\frac{C_2}{C_2 + \frac{K}{2}}$, and since $K = O(1)$ by definition, we get
\begin{equation*}
	\NI(\hmu, \hnu)
	\geq C_3 K \left( (\min_I \holw_I) \cdot V_1(\ha, \hx, \hb, \hy) \right)^{5/4}
\end{equation*}
for some $C_3>0$ only dependent on $(f, \XXX, \YYY)$.
% \TODO{this step hides a lot of easy considerations which I did on paper but it's very annoying to write down}
This concludes the proof of the second inequality of the proposition.

% \begin{remark}[The exponent in the lower bound]
%     \TODO{to check.}
% 	One could check that the following statement also holds: for any $(\ha, \hx, \hb, \hy)$, for any $0 < \eps \leq \frac{1}{4}$,
% 	\begin{equation*}
% 		K \left( (\min_I \holw_I) \cdot V_1(\hz) \right)^{1+\eps}
% 		\lesssim \NI(\hmu, \hnu)
% 		\lesssim \sqrt{V_1(\hz)}
% 	\end{equation*}
% 	where this time the $\lesssim$ on the left also depends on $\eps$. ($\eps$ is just, up to affine rescaling, the parameter we call $\beta$ in the proof.)
% \end{remark}

\subsection{Proof of the second part of the proposition}

% In the claim above, the dependency of $\min_I \hola_I, \min_J \holb_J$ could be made explicit, via Lemma~1.4 of ``analysis of particle PPA with exact parametrization'': $\wulb \geq \exp \left( -\frac{H(w^*) + \eta V(\ha, \hx, \hb, \hy)}{w^*_{\min}} \right)$.
% This would give the bound
% \begin{equation*}
% 	\NI \geq \sqrt{\sigma} C' e^{-c_1 \eta V} (\eta V)^{5/4}
% \end{equation*}
% for some constants $c_1$ and $C'$, which is not enough to conclude that $V$ is arbitrarily small when $\NI$ is small.
% However, we can still prove that fact using \eqref{eq:rel_Lya_NI:lb_NI_DeltaZ_global} directly.

We showed in Eq.~\eqref{eq:rel_Lya_NI:lb_NI_DeltaZ_global} that there exist constants $C, C_1>0$ such that
\begin{equation*}
	\NI(\hmu, \hnu)
	\geq C \norm{\Delta \hZ}^{5/2} - C_1 \left( \holw_0 + V_\pos(\hz) \right)^2
    ~~~\text{i.e.,}~~~
    \norm{\Delta \hZ}^{5/2} \leq 
	\frac{1}{C} \NI(\hmu, \hnu) + \frac{C_1}{C} \left( \holw_0 + V_\pos(\hz) \right)^2.
\end{equation*}
Now by Eq.~\eqref{eq:rel_Lya_NI:def_K}, we have
\begin{equation*}
    \NI(\hmu, \hnu) \geq K \left( \holw_0 + V_\pos(\hz) \right) 
	~~~\text{i.e.,}~~~
    \holw_0 + V_\pos(\hz) \leq \frac{1}{K} \NI(\hmu, \hnu).
\end{equation*}
Let $r \leq r_1 \wedge r_2$ where $r_1, r_2$ are defined by
$
    \frac{K}{C} r_1 = \frac{C_1}{C} r_2^2 = \frac{1}{2} \left( \frac{w^*_{\min}}{10} \right)^{5/2}
$.
Then, for any $(\ha, \hx, \hb, \hy)$ such that $\NI(\hmu, \hnu) \leq r K$,
\begin{align*}
    \frac{1}{C} \NI(\hmu, \hnu) 
    \leq \frac{1}{C} r_1 K
    &= \frac{1}{2} \left( \frac{w^*_{\min}}{10} \right)^{5/2}
    & &\text{and} &
    \frac{C_1}{C} \left( \holw_0 + V_\pos(\hz) \right)^2
    \leq \frac{C_1}{C} r_2^2
    &= \frac{1}{2} \left( \frac{w^*_{\min}}{10} \right)^{5/2},
\end{align*}
and so
$
    \norm{\Delta \hZ}^{5/2} \leq \left( \frac{w^*_{\min}}{10} \right)^{5/2}
$,
and in particular $\min_I \hola_I, \min_J \holb_J \geq \frac{w^*_{\min}}{4}$.
Note that by definition, $K = O(1)$, and so $r$ can be chosen independent of $\lambda, \tau$ and only dependent on $(f, \XXX, \YYY)$.

This concludes the proof of \autoref{prop:rel_Lya_NI:NI_equiv_V_moreprecise}, and so of \autoref{prop:cv_proof:gencase:NI_equiv_V}.

\subsection{Proof for the exact-parametrization case (\autoref{prop:cv_proof:exactparam:NI_equiv_V})}
\label{subsec:rel_Lya_NI:exactparam}

The first part of \autoref{prop:cv_proof:exactparam:NI_equiv_V} (the upper bound on $\NI$) follows from exactly the same computations as in the general case.

The second part of the proposition follows from the same considerations as for the general case; only Eq.~\eqref{eq:rel_Lya_NI:def_K} and \autoref{claim:rel_Lya_NI:Fnm-Fnm_bilinear} need to be adapted.
For the former, simply use \autoref{lm:growth_conds:quadr_exactparam} instead of \autoref{lm:growth_conds:quadr_gencase}.
For the latter, the same bound as in the general case holds; indeed this can be deduced from the general case using \autoref{claim:growth_conds:gencase_imply_exactparam}, by holding $\lambda\tau$ constant and letting $\tau \to \infty$.

% \begin{remark}
%     The second part of the proposition only holds when the solution particles are ordered such that $x^*_i$ is close to $\hx_i$.
% \end{remark}

% \end{document}

% !TEX root = ../main.tex
% \documentclass[../main]{subfiles}
% \begin{document}

\section{Proof of convergence in the general case} \label{sec:proof_gencase}

In this section we prove \autoref{thm:cv_proof:gencase:loc_exp_cv}.

\paragraph{Choice of the partitions of unity's parameters.}

Our specific choice for the parameters $\lambda, \tau$ appearing in the definition of $(\varphi_I)_I$, Eq.~\eqref{eq:cv_proof:gencase:def_varphi}, will not come into play until later in the proof, but to fix ideas we give their expressions right away.
We choose
\begin{equation} \label{eq:proof_gencase:choice_lambda_tau}
	\lambda^3 = \frac{1}{\sqrt{\sigma}}
	~~\text{and}~~
	\lambda\tau = \min\left\lbrace
		\sqrt{\frac{1}{2} \frac{\sigma}{\eta}},~
		\frac{\wulb~ \sigma_{\min}}{4 \smoothnessf_3},~
		\frac{\min_{I,I'} \norm{x^*_I-x^*_{I'}} \wedge \min_{J,J'} \norm{y^*_J-y^*_{J'}}}{4}
	\right\rbrace
	\asymp 1.
\end{equation}

Intuitively, in terms of the illustration \autoref{fig:misc:constru_lya}, the cut-off abscissa $\lambda\tau$ should be thought of as $\Theta(1)$, and the blue curve as being ``spiky'' with a scale of $\tau = \Theta(\sigma^{1/6})$, when $\eta, \sigma$ are small.

Note that
$\eps = e^{-\lambda^3/3} = e^{-1/(3 \sqrt{\sigma})}$
(the value of $\varphi_I$ at the cut-off)
is exponentially small for $\sigma \asymp \eta$ small, so that
%we have for example
%\begin{align*}
%% 	\begin{cases}
%% 		\frac{\lambda^3}{\alpha} \eps &\leq \frac{1}{10} \frac{\sigma_{\min}}{4} (\lambda\tau)^2 \eta \\
%% 		\eps &\leq \frac{1}{10} \gvsep \eta
%% 	\end{cases}
%%        --> pb: we haven't defined $\gvsep$ yet at this point...
%    \frac{\lambda^3}{3} \eps &\leq \frac{1}{10} \frac{\sigma_{\min}}{4} (\lambda\tau)^2 \eta
%\end{align*}
%whenever $\sigma$ and $\eta$ are small enough.
$\eps \cdot \mathit{poly}(\eta, \sigma, \lambda, \tau)$ is arbitrary small for $\eta, \sigma$ small enough, where $\mathit{poly}(...)$ is any polynomial expression of the arguments.
Essentially, any term where $\eps$ appears can be neglected (will be compensated by other terms), for $\eta$ and $\sigma$ small enough.

\subsection{Making the Lyapunov function appear in the characterizing inequality}

% \TODO{We could even include this subsection in the main text... A small writing obstacle would be that if we don't use the $w_i, p_i$ shorthands, then the final inequality will take twice as much space, and if we do, we must manage to introduce those shorthands concisely --- which is doable I think.}
% -> deal with this later

We start from the characterizing inequality \eqref{eq:cv_proof:ppa_ineq_foc}; for reference it reads
\begin{align*}
    \forall z,~
%	 \in \Delta_n \times \XXX^n \times \Delta_m \times \YYY^m,~ & \\
	\eta \tgap(z ; z^{k+1}) 
	&\leq 
	\sum_i (a_i-a^{k+1}_i) \log \frac{a^{k+1}_i}{a^k_i}
	+ \sum_j (b_j-b^{k+1}_j) \log \frac{b^{k+1}_j}{b^k_j} \\
	& ~~ + \frac{\eta}{\sigma} \sum_i a^k_i \innerprod{x^{k+1}_i-x^k_i}{x_i-x^{k+1}_i}
	+ \frac{\eta}{\sigma} \sum_j b^k_j \innerprod{y^{k+1}_j-y^k_j}{y_j-y^{k+1}_j}
\end{align*}
where 
$
	\tgap(z ; \hz) 
	= \innerprod{\begin{pmatrix}
		\nabla_a \\
		\nabla_x \\
		-\nabla_b \\
		-\nabla_y
	\end{pmatrix}
	F_{n,m}(\hz)}{\begin{pmatrix}
		\ha - a \\
		\hx - x \\
		\hb - b \\
		\hy - y
	\end{pmatrix}}
$.
As announced in \autoref{subsec:cv_proof:gencase}, \eqref{eq:cv_proof:gencase:choice_of_proxy*}, 
our first step is to evaluate it at ``proxy solution particles'' $(\aps, \xps, \bps, \yps) \in \Delta_n \times \XXX^n \times \Delta_m \times \YYY^m$ given by
\begin{equation*}
	\xps_i = x^{k+1}_i + \sum_{I \in [n^*]} \varphi^{k+1}_{Ii} (x^*_I - x^{k+1}_i) 
	~~~\text{and}~~~
	\aps_i = \sum_{I \in [n^*]} a^*_I \frac{\varphi^{k+1}_{Ii} a^{k+1}_i}{\ola^{k+1}_I}
\end{equation*}
and similarly for $\bps, \yps$.

\paragraph{Position terms.}

The term $\frac{\eta}{\sigma} \sum_i a^k_i \innerprod{x^{k+1}_i-x^k_i}{\xps_i-x^{k+1}_i}$ on the right-hand side of \eqref{eq:cv_proof:ppa_ineq_foc} becomes, by Pythagorean identity and Eq.~\eqref{eq:aux_lemmas:explicit_Vposax},
\begin{align*}
	& \sum_{I,i} a^k_i \, \varphi^{k+1}_{Ii} \innerprod{x^{k+1}_i-x^k_i}{x^*_I-x^{k+1}_i}
    = \frac{1}{2} \sum_{I, i} a^k_i \, \varphi^{k+1}_{Ii} \left( \norm{x^*_I-x^k_i}^2 - \norm{x^*_I-x^{k+1}_i}^2 - \norm{x^{k+1}_i-x^k_i}^2 \right) \\
	&= \left[ 
		\frac{1}{2} \sum_{I, i} a^k_i \, \varphi^{k+1}_{Ii} \norm{x^*_I-x^k_i}^2 
	\right] - \left[
		\frac{1}{2} \sum_{I, i} a^k_i \, \varphi^{k+1}_{Ii} \norm{x^*_I-x^{k+1}_i}^2
	\right] - \left[
		\frac{1}{2} \sum_{I, i} a^k_i \, \varphi^{k+1}_{Ii} \norm{x^{k+1}_i-x^k_i}^2
	\right] \\
	&= \left[ V_\pos(a^k, x^k) + \frac{1}{2} \sum_{I, i} a^k_i \left( \varphi^{k+1}_{Ii}-\varphi^k_{Ii} \right) \norm{x^*_I-x^k_i}^2 \right] \\
	& - \left[ V_\pos(a^{k+1}\!, x^{k+1}) + \frac{1}{2} \sum_{I, i} \left( a^k_i-a^{k+1}_i \right) \varphi^{k+1}_{Ii} \norm{x^*_I-x^{k+1}_i}^2 \right] 
	- \left[ \frac{1}{2} \sum_i a^k_i \left( 1-\varphi^{k+1}_{0i} \right) \norm{x^{k+1}_i-x^k_i}^2 \right].
\end{align*}

\paragraph{Weight terms.}

For all $I \in [0, n^*]$, let
$u^{k+1,I}_i = \frac{\varphi^{k+1}_{Ii} a^{k+1}_i}{\ola^{k+1}_I}$,
so that $u^{k+1,I} \in \Delta_n$ for each $I$.
Then since $\aps_i = \sum_I a^*_I u^{k+1,I}_i$,
\begin{align*}
	\sum_i \aps_i \log \frac{a^{k+1}_i}{a^k_i} 
	&= \sum_I \sum_i a^*_I u^{k+1,I}_i \log \frac{a^{k+1}_i}{a^k_i} \\
	&= \sum_I a^*_I \sum_i u^{k+1,I}_i \log \frac{\ola^{k+1}_I}{\ola^k_I} \cdot \frac{a^{k+1}_i / \ola^{k+1}_I}{a^k_i / \ola^k_I} \\
	&= \sum_I a^*_I \left( \sum_i u^{k+1,I}_i \right) \log \frac{\ola^{k+1}_I}{\ola^k_I}
	+ \sum_I a^*_I \sum_i u^{k+1,I}_i \log \frac{a^{k+1}_i / \ola^{k+1}_I}{a^k_i / \ola^k_I} \\
	&= \KLdiv(a^*, \ola^k) - \KLdiv(a^*, \ola^{k+1}) 
	+ \sum_I a^*_I \sum_i u^{k+1,I}_i \log \frac{a^{k+1}_i / \ola^{k+1}_I}{a^k_i / \ola^k_I}
\end{align*}
and, since $a^{k+1}_i = \sum_{I \in [0,n^*]} \ola^{k+1}_I u^{k+1,I}_i$,
by the same calculation
with $a^*_I$ replaced by $\ola^{k+1}_I$
\begin{align*}
	\sum_i a^{k+1}_i \log \frac{a^{k+1}_i}{a^k_i}
	&= \sum_{I \in [0,n^*]} \sum_i \ola^{k+1}_I u^{k+1,I}_i \log \frac{a^{k+1}_i}{a^k_i} \\
	&= \KLdiv(\ola^{k+1}, \ola^k) 
	+ \sum_{I \in [0,n^*]} \ola^{k+1}_I \sum_i u^{k+1,I}_i \log \frac{a^{k+1}_i / \ola^{k+1}_I}{a^k_i / \ola^k_I},
\end{align*}
and the term for $I=0$ in this last sum is equal to
\begin{align*}
	\sum_i \varphi^{k+1}_{0i} a^{k+1}_i \!\left[ \log \frac{a^{k+1}_i}{a^k_i} - \log \frac{\ola^{k+1}_0}{\ola^k_0} \right] 
	&= \sum_i \varphi^{k+1}_{0i} \left[ d_h(a^{k+1}_i, a^k_i) + a^{k+1}_i - a^k_i \right] \\
	&~~~~ - \left[ d_h(\ola^{k+1}_0, \ola^k_0) + \ola^{k+1}_0 - \ola^k_0 \right] \\[0.5em]
	&= \sum_i \varphi^{k+1}_{0i} d_h(a^{k+1}_i, a^k_i)
	- d_h(\ola^{k+1}_0, \ola^k_0)
	- \sum_i \left( \varphi^{k+1}_{0i} - \varphi^k_{0i} \right) a^k_i.
\end{align*}
So the weight term
on the right-hand side of \eqref{eq:cv_proof:ppa_ineq_foc} becomes
\begin{align}
	\sum_i (\aps_i-a^{k+1}_i) \log \frac{a^{k+1}_i}{a^k_i}
	&= V_\wei(a^k, x^k) - V_\wei(a^{k+1}, x^{k+1})
	- \sum_I d_h(\ola^{k+1}_I, \ola^k_I) \\
	&~~~~ + \sum_I (a^*_I - \ola^{k+1}_I) \sum_i u^{k+1,I}_i \log \frac{a^{k+1}_i / \ola^{k+1}_I}{a^k_i / \ola^k_I} \\
	&~~~~ - \sum_i \varphi^{k+1}_{0i} d_h(a^{k+1}_i, a^k_i)
	~+~ \sum_i (\varphi^{k+1}_{0i} - \varphi^k_{0i}) a^k_i.
\label{eq:proof_gencase:cq_choice_of_wei_pre}
\end{align}

All in all, evaluating the characterizing inequality \eqref{eq:cv_proof:ppa_ineq_foc} at $(\aps, \xps, \bps, \yps) = \zps$ yields
% \TODO{the equation numbering is buggy...} -> solved by following https://latex.org/forum/viewtopic.php?t=22540 (put the "eq:proof_gencase:ppa_ineq_star" label inside the align, even within the subequations, in order to get the equation number to be displayed)
\begin{subequations}
    \begin{align} 
    \label{eq:proof_gencase:ppa_ineq_star}
    	\eta \tgap(\zps; z^{k+1}) 
    	&\leq
    	V(z^k) - V(z^{k+1}) \\
    	&~~~~ - \sum_I d_h(\olw^{k+1}_I, \olw^k_I) - \frac{\eta}{2 \sigma} \sum_i w^k_i (1-\varphi_{0i}^{k+1}) \norm{p^{k+1}_i-p^k_i}^2 \\
    	&~~~~ + \sum_I (w^*_I-\olw^{k+1}_I)
    	\sum_i \frac{\varphi^{k+1}_{Ii} w^{k+1}_i}{\olw^{k+1}_I} \log \frac{w^{k+1}_i / \olw^{k+1}_I}{w^k_i / \olw^k_I} 
    \label{eq:proof_gencase:ppa_ineq_star:errterm1} \tag{err1}
    	\\
    	&~~~~ - \sum_i \varphi^{k+1}_{0i} d_h(w^{k+1}_i, w^k_i) \\
        &~~~~ + \sum_i (\varphi^{k+1}_{0i} - \varphi^k_{0i}) w^k_i
        ~+~ \frac{\eta}{2 \sigma} \sum_I \sum_i w^k_i \left( \varphi^{k+1}_{Ii}-\varphi^k_{Ii} \right) \norm{p^*_I-p^k_i}^2 ~~~
    \label{eq:proof_gencase:ppa_ineq_star:errterm2} \tag{err2}
    	\\
    	&~~~~ + \frac{\eta}{2 \sigma} \sum_I \sum_i (w^{k+1}_i-w^k_i) \varphi^{k+1}_{Ii} \norm{p^*_I-p^{k+1}_i}^2
    \label{eq:proof_gencase:ppa_ineq_star:errterm3} \tag{err3}
    \end{align}
\end{subequations}
where we let for concision 
$w^k = \begin{pmatrix}
    a^k \\
    b^k
\end{pmatrix} \in \Delta_n \times \Delta_m$
and
$p^k = \begin{pmatrix}
    x^k \\
    y^k
\end{pmatrix} \in \XXX^n \times \YYY^m$,
and similarly for $\olw^k, w^* \in \Delta_{[0, n^*]} \times \Delta_{[0, m^*]}$ and $\olp^k, p^* \in \XXX^{n^*} \times \YYY^{m^*}$.

The left-hand side looks like ``gap from MNE to iterates'' so morally non-negative, which we will show and quantify.
The second line consists of minus ``divergence from $(k+1)$ to $k$'' terms which we will lower-bound.
The last four lines consist of error terms which we will control.

% \eqref{eq:proof_gencase:ppa_ineq_star}
% \eqref{eq:proof_gencase:ppa_ineq_star:errterm1}
% \eqref{eq:proof_gencase:ppa_ineq_star:errterm2}
% \eqref{eq:proof_gencase:ppa_ineq_star:errterm3}

\subsection{Preliminary lemmas} \label{subsec:proof_gencase:prelim_lems}

As this phrase is used many times in the proof,
let us emphasize again that by ``for $\eta, \sigma$ small enough'' we always mean that a property holds for all $\eta \leq \eta_0, \sigma \leq \sigma_0$ for some $\eta_0, \sigma_0$ only dependent on $(f, \XXX, \YYY)$ and $\Gamma_0$ (the same constants that may be hidden in $O(\cdot)$'s).

Next, we state some useful elementary facts about the algorithm.
The following equations are clear from the update rule \eqref{eq:main_res:ppa_upd}.
Alternatively, they can be seen as a consequence of \eqref{eq:cv_proof:ppa_ineq_foc} --- holding with equality --- 
applied to $(a^{k+1} + \delta a, x^{k+1}, b^{k+1}, y^{k+1})$ for all $\delta a \in \{ \bmone_n \}^\top$,
resp.\ $(a^{k+1}, x^{k+1} + \delta x^{(i_0)}, b^{k+1}, y^{k+1})$ where $\delta x^{(i_0)}_i = \ind_{i=i_0}$.
\begin{align} \label{eq:proof_gencase:ppa_upd_eqs}
	a^{k+1}_i &=  a^k_i ~ e^{-\eta [(\gmat^{k+1} b^{k+1})_i - \gval]} / Z
	~~~\text{where}~~~
	Z = \sum_{i'} a^k_{i'} ~ e^{-\eta [(\gmat^{k+1} b^{k+1})_{i'} - \gval]} \\
	\text{and}~~~~
	x^{k+1}_i &= x^k_i - \sigma \frac{a^{k+1}_i}{a^k_i} \partial_x \gmat^{k+1}_{i \bullet} b^{k+1}.
\end{align}

We formalize the trivial fact that iterates move by no more than the step-size at each time-step.
\begin{lemma} \label{lm:proof_gencase:bound_movement}
	For all $k$, 
	\begin{align*} 
    % \label{eq:proof_gencase:coarsely_bound_movement_wei}
		\forall i,~ \abs{a^{k+1}_i-a^k_i} 
		&\leq \min\{ a^k_i, a^{k+1}_i \} (e^{2\eta L_0}-1) \\
		\norm{a^{k+1}-a^k}_1 &\leq e^{2\eta L_0}-1
		= O(\eta)
	\end{align*}
	and in particular $a^{k+1}_i = a^k_i (1+O(\eta))$,
	and similarly for $b$,
%	(Here $O(\cdot)$ hides $\smoothnessf_0$ and $\eta_0$.)
%	
%	Furthermore, for all $k$,
	and
	\begin{equation} \label{eq:proof_gencase:coarsely_bound_movement_pos}
		\forall i,~ \norm{x^{k+1}_i-x^k_i} 
		\leq \sigma e^{2 \eta \smoothnessf_0} \smoothnessf_1
		= O(\sigma)
	\end{equation}
	and similarly for $y$.
%	(Here $O(\cdot)$ hides $\smoothnessf_0$, $\smoothnessf_1$, $\eta_0$ and $\sigma_0$.)
\end{lemma}

\begin{proof}
	From \eqref{eq:proof_gencase:ppa_upd_eqs}, we have
	\begin{align*}
		&\forall i,~&
		a_i^{k+1} &= a^k_i \frac{e^{-\eta (\gmat^{k+1} b^{k+1})_i}}{\sum_{i'} a^k_{i'} e^{-\eta (\gmat^{k+1} b^{k+1})_{i'}}} && \\
		\text{and so}~~ & &
		e^{-2\eta L_0} =
		\frac{e^{-\eta L_0}}{\sum_{i'} a^k_{i'} e^{\eta L_0}}
		\leq \frac{a^{k+1}_i}{a^k_i}
		&\leq \frac{e^{\eta L_0}}{\sum_{i'} a^k_{i'} e^{-\eta L_0}} 
		= e^{2 \eta L_0} &&
	\end{align*}
	from which the first result follows.
	
	Furthermore, also from \eqref{eq:proof_gencase:ppa_upd_eqs},
	$
		\norm{x^{k+1}_i-x^k_i} 
		= \sigma \frac{a^{k+1}_i}{a^k_i} \norm{ \partial_x \gmat^{k+1}_{i \bullet} b^{k+1} } 
		\leq \sigma~ e^{2 \eta \smoothnessf_0} \smoothnessf_1.
	$
\end{proof}

Recall from \autoref{lm:growth_conds:ulb_on_weights_locally} that 
locally (i.e.\ if $V(z^k)$ is small enough), we have a constant lower bound on the iterates' aggregated weights $\ola^k_I, \olb^k_I$. 
That is, there exists $r>0$ (dependent only on $a^*, b^*$) such that
$
    V(z) \leq r \implies
    \left( \min_{I \neq 0} \ola_I \right) \wedge \left( \min_{J \neq 0} \olb_J \right) 
    \geq \frac{a^*_{\min} \wedge b^*_{\min}}{2}
$.
Thanks to the considerations above, we can show that locally, the aggregated weights are lower-bounded by a constant for iterates both at $k$ and at $k+1$ .
\begin{lemma} \label{lm:proof_gencase:ulb_on_weights_locally_k+1}
	There exists $r>0$ (only dependent on $a^*, b^*$) such that if $V(z^k) \leq r$, then for small enough $\eta, \sigma$,
	\begin{equation*}
	    \left( \min_{I \neq 0}~ \ola^k_I \right)
	    \wedge
	    \left( \min_{I \neq 0}~ \ola^{k+1}_I \right)
	    \wedge
	    \left( \min_{J \neq 0}~ \olb^k_J \right)
	    \wedge
	    \left( \min_{J \neq 0}~ \olb^{k+1}_J \right)
	    \geq 
        \frac{a^*_{\min} \wedge b^*_{\min}}{4}
        \eqqcolon \wulb.
	\end{equation*}
\end{lemma}

The proof is conceptually simple but annoyingly technical due to the fact that, to compare 
$\ola^k_I = \sum_i \varphi_I(x^k_i) a^k_i$
and
$\ola^{k+1}_I = \sum_i \varphi_I(x^{k+1}_i) a^{k+1}_i$,
we also need to control the variation between $\varphi_I(x^k_i)$ and $\varphi_I(x^{k+1}_i)$.
Namely we have the following bound, which will also be useful elsewhere in this appendix.
\begin{lemma} \label{lm:proof_gencase:control_varphi_t+1_t}
	For any $I \in [n^*]$,
	\begin{equation*}
		\abs{ \sum_i (\varphi^{k+1}_{Ii} - \varphi^k_{Ii}) a^k_i }
		\lesssim 
        \eps \left( \ola^k_0 + \ola^{k+1}_0 \right) 
        + \sqrt{\sigma} \min\left\{ 1, V(z^{k+1}) \right\}
		\lesssim 
        \eps \ola^k_0
        + \sqrt{\sigma} V(z^{k+1}).
	\end{equation*}
\end{lemma}
The proof of \autoref{lm:proof_gencase:control_varphi_t+1_t} is quite technical 
and is deferred 
\ifextended%
    to \autoref{subsec:proof_gencase:delayed_proofs}.
\else%
    to Section~E.7.2 of the arXiv version of this paper as \AdditionalMaterial.
\fi
In particular it relies on our specific choice of partitions of unity \eqref{eq:cv_proof:gencase:def_varphi} and of $\lambda, \tau$ \eqref{eq:proof_gencase:choice_lambda_tau}.

\begin{proof}[Proof of \autoref{lm:proof_gencase:ulb_on_weights_locally_k+1}]
	Let $r$ the constant from \autoref{lm:growth_conds:ulb_on_weights_locally}, so that
	$\eta V(z) \leq r \implies \min_I \ola_I, \min_J \olb_J \geq \frac{a^*_{\min} \wedge b^*_{\min}}{2} = 2\wulb$.
	This immediately ensures that $\min_I \ola^k_I, \min_I \olb^k_I \geq 2 \wulb$.
	
	For any $I \in [n^*]$,
	\begin{align*}
		\ola^{k+1}_I - \ola^k_I
		= \sum_i \varphi^{k+1}_{Ii} a^{k+1}_i - \varphi^k_{Ii} a^k_i
		&= \sum_i \varphi^{k+1}_{Ii} (a^{k+1}_i-a^k_i)
		+ \sum_i (\varphi^{k+1}_{Ii} - \varphi^k_{Ii}) a^k_i \\
		&\geq -O(\eta) \ola^{k+1}_I + \sum_i (\varphi^{k+1}_{Ii} - \varphi^k_{Ii}) a^k_i \\
		\ola^{k+1}_I &\geq (1-O(\eta)) \left[ \ola^k_I + \sum_i (\varphi^{k+1}_{Ii} - \varphi^k_{Ii}) a^k_i \right].
	\end{align*}
	Now $\ola^k_I \geq 2 \wulb$, and by \autoref{lm:proof_gencase:control_varphi_t+1_t},
	$
		\abs{ \sum_i (\varphi^{k+1}_{Ii} - \varphi^k_{Ii}) a^k_i }
        \lesssim 
        \eps + \sqrt{\sigma}
	$.
    % Hence,
    % $
    %     \ola^{k+1}_I 
    %     \geq
    %     (1-O(\eta)) \left[
    %         2 \wulb
    %         - O(\eps + \sqrt{\sigma)})
    %     \right]
    % $.
	So for $\eta, \sigma$ small enough, we indeed have $\ola^{k+1}_I \geq \wulb$, and similarly $\olb^{k+1}_J \geq \wulb$ for all $J \in [m^*]$.
\end{proof}

\emph{In the remainder of this section except \autoref{subsec:proof_gencase:proof_conclusion}, we assume that the conditions of this lemma are satisfied.}
As a first useful consequence, we have that both $V(a^k, x^k)$ and $V(a^{k+1}, x^{k+1})$ are uniformly bounded. 
Indeed, $V(z^k) \leq r$ and $V_\pos(a, x) = \sum_{I,i} \varphi_{Ii} a_i \norm{x_i-x^*_I}^2 \leq \diamXY^2 = O(1)$ for any $(a,x)$ anyway,
and
\begin{equation} \label{eq:proof_gencase:etaveesigma_V_bounded}
	V_\wei(a^{k+1}, x^{k+1}) = \ola^{k+1}_0 + \sum_I d_h(a^*_I, \ola^{k+1}_I)
	\leq 1 + \sum_I \frac{1}{\wulb} (a^*_I - \ola^{k+1}_I)^2
	\leq 1 + \frac{n^*}{\wulb} = O(1)
\end{equation}
% since $h'' \leq \frac{1}{c}$ over $[c,1]$.
by $\frac{1}{c}$-smoothness of $h$ over $[c,1]$.
A second useful consequence is that, by \autoref{claim:growth_conds:DeltaZ_equiv_etaveesigmaV},
for $(a, x) = (a^k, x^k)$ or $(a^{k+1}, x^{k+1})$,
\begin{equation*}
    \ola_0 + \norm{\ola - a^*}_1^2 \asymp V_\wei(a, x) 
    ~~~~\text{and}~~~~
    \max_I \norm{\olx_I-x^*_I}^2 \lesssim V_\pos(a, x).
\end{equation*}
% (In fact $V_\pos(a, x) \asymp \max_I \norm{\olx_I-x^*_I}^2 + \trace \Sigma_I$.)

\subsection{Lower-bounding \texorpdfstring{$\tgap(\zps; z^{k+1})$}{tildegap(zps; zk+1)}}

Let us give a quantitative lower bound on the term on the left-hand side of \eqref{eq:proof_gencase:ppa_ineq_star}:
$\tgap(\zps; z^{k+1})$.
It is here that we make use of the ``quadratic growth'' and ``star-convexity-concavity'' properties discussed in 
Sections~\ref{subsubsec:cv_proof:proofingr_growthconds:quadr} and \ref{subsubsec:cv_proof:proofingr_growthconds:starconvex}.

\begin{lemma} \label{lm:proof_gencase:lb_tgap_olw0}
    There exists a constant $\gvsep>0$ only dependent on $(f, \XXX, \YYY)$ such that
	\begin{equation*}
		\tgap(\zps; z^{k+1})
		\geq 
		\left[ 
			\frac{\sigma_{\min}}{4} \frac{3 (\lambda\tau)^3}{\lambda^3} \wedge \gvsep
		\right]
		(\ola^{k+1}_0 + \olb^{k+1}_0)
		+ \frac{\sigma_{\min}}{2} V_\pos(z^{k+1})
        + O \left( V(z^{k+1})^{3/2} \right).
	\end{equation*}
\end{lemma}

\begin{proof}
    As a direct consequence of \autoref{lm:growth_conds:starconvex},
    since $\left( \min_I \ola^{k+1}_I \right) \wedge \left( \min_J \olb^{k+1}_J \right) \geq \wulb$,
    then
    denoting $\mu^{k+1} = \sum_{i=1}^n a^{k+1}_i \delta_{x^{k+1}_i}$
    and $\nu^{k+1} = \sum_{j=1}^m b^{k+1}_j \delta_{y^{k+1}_j}$,
	\begin{equation*}
		\tgap(\zps; z^{k+1})
		\geq F(\mu^{k+1}, \nu^*) - F(\mu^*, \nu^{k+1}) 
        + V_\pos(z^{k+1}) 
        \underbrace{\left( 
            \sigma_{\min} - 2 \smoothnessf_3 \wulb^{-1} \cdot \lambda \tau 
        \right)}
        + O\left( V(z^{k+1})^{3/2} \right).
    \end{equation*}
    Note that by our choice of $\lambda\tau \leq \frac{\wulb \sigma_{\min}}{4 \smoothnessf_3}$, we have $\sigma_{\min} - 2 \smoothnessf_3 \wulb^{-1} \cdot \lambda \tau \geq \frac{\sigma_{\min}}{2}$.
    
	Furthermore, as a direct consequence of \autoref{lm:growth_conds:quadr_gencase}, we have
	\begin{equation*}
		F(\mu^{k+1}, \nu^*) - F(\mu^*, \nu^{k+1})
		\geq 
		\left[ 
			\frac{\sigma_{\min}}{4} \frac{3 (\lambda\tau)^2}{\lambda^3} \wedge \gvsep
		\right]
		(\ola^{k+1}_0 + \olb^{k+1}_0)
	\end{equation*}
	for some constant $\gvsep>0$ only dependent on $(f, \XXX, \YYY)$.
	The lemma follows by combining the two inequalities.
\end{proof}

\subsection{Controlling the error terms}

The proofs for the lemmas in this subsection are all technical, relying on our specific choice of partitions of unity $(\varphi_I)_I$ as well as of parameters $\lambda, \tau$.
We defer the proofs to 
\ifextended%
    \autoref{subsec:proof_gencase:delayed_proofs}.
\else%
    Section~E.7 of the arXiv version of this paper as \AdditionalMaterial.
\fi

\begin{lemma}[Bound for \eqref{eq:proof_gencase:ppa_ineq_star:errterm1}]
\label{lm:proof_gencase:errterm1__control_KLdiv_uI_vI_errterm}
    For $\eta, \sigma$ small enough,
	\begin{align*}
        \forall I,~
        \sum_i \frac{\varphi^{k+1}_{Ii} w^{k+1}_i}{\olw^{k+1}_I} \log \frac{w^{k+1}_i / \olw^{k+1}_I}{w^k_i / \olw^k_I}  
        &= O\left( \eps \olw^k_0 + \sqrt{\eta} V(z^{k+1}) \right). \\
        \text{In particular,}~~
		\sum_I (w^*_I-\olw^{k+1}_I)
    	\sum_i \frac{\varphi^{k+1}_{Ii} w^{k+1}_i}{\olw^{k+1}_I} \log \frac{w^{k+1}_i / \olw^{k+1}_I}{w^k_i / \olw^k_I}  
		% &\lesssim 
		% \norm{\Delta \olw^{k+1}}_1 \cdot 
		% \left[ \eps \olw^k_0 + \sqrt{\eta} V(z^{k+1}) \right] \\
		&\lesssim 
		\sqrt{ V(z^{k+1}) } \cdot
		\left[ \eps \olw^k_0 + \sqrt{\eta} V(z^{k+1}) \right].
	\end{align*}
\end{lemma}

%\begin{remark}
%	The proof of the following lemma crucially relies on the sign of $\varphi^k_{Ii} - \varphi^{k+1}_{Ii}$.
%	That is, it wouldn't be sufficient to just bound
%	$
%		\sum_I \sum_i \abs{\varphi^k_{Ii} - \varphi^{k+1}_{Ii}} w^k_i
%	$
%	using \autoref{lm:gencase:control_varphi_t+1_t}.
%	If we did, we would obtain an upper bound containing a term of the form
%	\begin{equation*}
%		\lesssim \lambda^3 \sigma V_\pos(a^{k+1}, x^{k+1}).
%	\end{equation*}
%	(See previous versions of this document for details.)
%	This would not be compensated by the $\frac{1}{2} \sigma_{\min} \sigma V_\pos(z^{k+1})$ in the lower bound of \autoref{lm:gencase:lb_tgap_olw0}, unless $\lambda^3$ is small, which would make $\eps$ large.
%	
%	Frustratingly, in numerical experiments, $\sum_I \sum_i \abs{\varphi^k_{Ii} - \varphi^{k+1}_{Ii}} w^k_i$ does in fact seem negligible compared to $V(z^{k+1})$... (But maybe it's because worst-case is hard to see in experiments.)
%\end{remark}
% 
\begin{lemma}[Bound for \eqref{eq:proof_gencase:ppa_ineq_star:errterm2}] 
\label{lm:proof_gencase:errterm2__control_tricky_errterms}
    For $\eta, \sigma$ small enough,
	\begin{equation*}
	    \sum_i (\varphi^{k+1}_{0i} - \varphi^k_{0i}) w^k_i
        +\frac{\eta}{2 \sigma} \sum_{I,i} w^k_i \left( \varphi^{k+1}_{Ii}-\varphi^k_{Ii} \right) \norm{p^*_I-p^k_i}^2
        \lesssim	
        \eps \olw^{k+1}_0
		+ \eta^{3/2} V_\pos(z^{k+1})
		+ \sqrt{\eta} V(z^{k+1})^{3/2}.
	\end{equation*}
%	and the constant hidden in the $\lesssim$ is $1$ (actually $1+O(\eta)$) for the first term.
\end{lemma}

\begin{lemma}[Bound for \eqref{eq:proof_gencase:ppa_ineq_star:errterm3}] \label{lm:proof_gencase:errterm3__control_nonBregmanness_errterm}
    For $\eta, \sigma$ small enough,
	\begin{equation*}
		\frac{\eta}{2 \sigma} \sum_I \sum_i (w^{k+1}_i-w^k_i) \varphi^{k+1}_{Ii} \norm{p^*_I-p^{k+1}_i}^2
		\lesssim \eta V(z^{k+1})^{3/2}.
	\end{equation*}
\end{lemma}

\paragraph{\texorpdfstring{$V(z^k)$}{V(zk)} does not grow too fast.}
At this point, we have all we need to show the following.
\begin{lemma} \label{lm:proof_gencase:Vk+1_small}
    For $\eta, \sigma$ small enough, 
    $
		V(z^{k+1}) \leq 2 V(z^k).
	$
\end{lemma}

\begin{proof}
	Starting from \eqref{eq:proof_gencase:ppa_ineq_star},
	upper-bound the second line by $0$, 
	lower-bound the left-hand side using \autoref{lm:proof_gencase:lb_tgap_olw0}
	and bound the error terms using
    Lemmas~\ref{lm:proof_gencase:errterm1__control_KLdiv_uI_vI_errterm}, \ref{lm:proof_gencase:errterm2__control_tricky_errterms} and~\ref{lm:proof_gencase:errterm3__control_nonBregmanness_errterm}.
	Simplify the obtained inequality using that 
	$\lambda\tau \asymp 1$,
    $\lambda^3 = \frac{1}{\sqrt{\sigma}}$,
    $\eta \asymp \sigma$,
	and $V(z^{k+1}) = O(1)$ as noted above \eqref{eq:proof_gencase:etaveesigma_V_bounded}.
	Rearranging, we get
	\begin{align*}
	    V(z^{k+1}) - V(z^k) 
		&\leq O(\eps) V(z^k) + O \left( \sqrt{\eta} \cdot V(z^{k+1}) \right) \\
		V(z^{k+1}) &\leq (1+O(\sqrt{\eta})) (1+O(\eps)) V(z^k)
	\end{align*}
	and so $V(z^{k+1}) \leq 2 V(z^k)$ for $\eta, \sigma$ small enough,
	as announced.
\end{proof}

\subsection{Lower-bounding the ``divergence from \texorpdfstring{$(k+1)$ to $k$}{k+1 to k}'' terms}

In this subsection, we lower-bound the quantity
\begin{equation*}
    D(k+1, k) \coloneqq \sum_I d_h(\olw^{k+1}_I, \olw^k_I) + \frac{\eta}{2 \sigma} \sum_i w^k_i (1-\varphi_{0i}^{k+1}) \norm{p^{k+1}_i-p^k_i}^2
\end{equation*}
appearing with a negative sign on the right-hand side of \eqref{eq:proof_gencase:ppa_ineq_star}.
The bound relies on the ``error bound''-type property discussed in \autoref{subsubsec:cv_proof:proofingr_growthconds:error_bound}.

\begin{lemma} \label{lm:proof_gencase:lb_Dk+1_k}
    There exist constants $r, C >0$ only dependent on $(f, \XXX, \YYY)$ and $\Gamma_0$ such that, for $V(z^k) \leq r$ and small enough $\eta, \sigma$,
    then $D(k+1, k)$ is lower-bounded by
    \begin{equation*}
        D(k+1, k) 
    	\geq C \eta^2 \left(
    		\sum_I d_h(w^*_I, \olw^{k+1}_I)
    		+ \sum_I \olw^{k+1}_I \norm{\Delta \olp^{k+1}_I}^2
    	\right)
    	+ O \left( \left( \eps \olw^k_0 \right)^2 + \eta V(z^{k+1})^2 \right).
    \end{equation*}
\end{lemma}

% Note that there is neither $\olw^{k+1}_0$ nor any term involving $\trace(H_I \Sigma^{k+1}_I)$ on the right-hand side, so our lower bound on $D(k+1, k)$ will not be proportional to $V(z^{k+1})$. But that's okay, because we already have ${ \olw^{k+1}_0 + \sigma V_\pos(z^{k+1}) }$ on the left-hand side of \eqref{eq:proof_gencase:ppa_ineq_star} thanks to \autoref{lm:proof_gencase:lb_tgap_olw0}.

The remainder of this subsection is dedicated to proving this lemma.

For any $(A_I)_{I \in [n^*]}, (B_J)_{J \in [m^*]}$ (and $A_0=B_0=0$) and $(X_I)_{I \in [n^*]}, (Y_J)_{J \in [m^*]}$,
% let analogously to \eqref{eq:cv_proof:gencase:choice_of_proxy*} the ``proxy particles''
define ``proxy particles'' as in \eqref{eq:growth_conds:encode_WP_hat}:
\begin{equation} \label{eq:proof_gencase:encode_WP}
	x_i = x^{k+1}_i + \sum_I \varphi^{k+1}_{Ii} (X_I - x^{k+1}_i)
	~~~\text{and}~~~
	a_i = \sum_I A_I \frac{\varphi^{k+1}_{Ii} a^{k+1}_i}{\ola^{k+1}_I},
\end{equation}
% \begin{equation} \label{eq:proof_gencase:encode_WP}
% 	x_i &= x^{k+1}_i + \sum_I \varphi^{k+1}_{Ii} (X_I - x^{k+1}_i) \\
% 	\text{and}\quad\quad\quad
% 	a_i &= \sum_I A_I u^{k+1,I}_i
% 	\quad\quad\quad \text{where} \quad\quad\quad
% 	\forall I,~ u^{k+1,I}_i = \frac{\varphi^{k+1}_{Ii} a^{k+1}_i}{\ola^{k+1}_I}.
% \end{equation}
similarly for $b$ and $y$ and let $z = (a, x, b, y)$, leaving the dependence on $A, X, B$ and $Y$ implicit to lighten notation.
For concision, let as usual
$w = \begin{pmatrix}
    a \\
    b
\end{pmatrix}$,
$p = \begin{pmatrix}
    x \\
    y
\end{pmatrix}$,
$W = \begin{pmatrix}
    A \\
    B
\end{pmatrix}$, and
$P = \begin{pmatrix}
    X \\
    Y
\end{pmatrix}$.
We proceed by upper- and lower-bounding the gradient-norm-like quantity
$\max_{A,X,B,Y} \tgap(z; z^{k+1})$.

% \paragraph{Upper part of the sandwich.}
\paragraph{Upper bound on the gradient-norm-like quantity.}

\begin{claim} \label{claim:gencase:sandwich_maxtgap_upper}
    % There exist constants $r, C>0$ only dependent on $(f, \XXX, \YYY)$ and $\Gamma_0$ such that for $V(z^k) \leq r$ and small enough $\eta, \sigma$,
    There exists $C>0$ only dependent on $(f, \XXX, \YYY)$ and $\Gamma_0$ such that
	\begin{equation*}
		\max_{A,X,B,Y} \tgap(z; z^{k+1}) 
		\leq 
        \frac{C}{\eta}
        \sqrt{D(k+1, k)}
        ~+~ \frac{1}{\eta} O \left( \eps \olw^k_0 + \sqrt{\eta} V(z^{k+1}) \right).
	\end{equation*}
\end{claim}

\begin{proof}
	Evaluate \eqref{eq:cv_proof:ppa_ineq_foc} at the proxy particles $z=(a, x, b, y)$:
	\begin{equation*}
		\forall A, X, B, Y,~
		\eta \tgap(z; z^{k+1}) 
		\leq \sum_i (w_i-w^{k+1}_i) \log \frac{w^{k+1}_i}{w^k_i}
		+ \frac{\eta}{\sigma} \sum_i w^k_i \innerprod{p^{k+1}_i-p^k_i}{p_i-p^{k+1}_i}.
	\end{equation*}
	
	On the right-hand side, we get for the position terms
	\begin{align*}
		\sum_i w^k_i \innerprod{p^{k+1}_i-p^k_i}{p_i-p^{k+1}_i} 
		&= \sum_I \sum_i w^k_i \varphi^{k+1}_{Ii} 
        \Big\langle 
            p^{k+1}-p^k_i,
            \underbrace{
                P_I-p^{k+1}_i
            }_{\norm{\cdot} \leq R}
        \Big\rangle \\
		&\leq \diamXY \sum_I \sum_i w^k_i \varphi^{k+1}_{Ii} \norm{p^{k+1}_i-p^k_i}
        = \diamXY \sum_i w^k_i (1-\varphi^{k+1}_{0i}) \norm{p^{k+1}_i-p^k_i} \\
		&\leq 
		\diamXY 
		\underbrace{
			\sqrt{\sum_i w^k_i (1-\varphi^{k+1}_{0i})}
		}_{\leq \sqrt{2}}
		\sqrt{
			\sum_i w^k_i (1-\varphi^{k+1}_{0i}) \norm{p^{k+1}_i-p^k_i}^2
		}.
	\end{align*}
    In the last inequality, we used Cauchy-Schwarz inequality.
	
	For the weight terms,
	by the same calculation as for \eqref{eq:proof_gencase:cq_choice_of_wei_pre} with $w^*$ replaced by $W$,
	\begin{align*}
		\sum_i (w_i-w^{k+1}_i) \log \frac{w^{k+1}_i}{w^k_i}
		&= \sum_I W_I \log \frac{\olw^{k+1}_I}{\olw^k_I}
		- \sum_I d_h(\olw^{k+1}_I, \olw^k_I) \\
		&~~~~ + \sum_I (W_I - \olw^{k+1}_I) \sum_i \frac{\varphi^{k+1}_{Ii} a^{k+1}_i}{\ola^{k+1}_I} \log \frac{w^{k+1}_i / \olw^{k+1}_I}{w^k_i / \olw^k_I} \\
		&~~~~ - \sum_i \varphi^{k+1}_{0i} d_h(w^{k+1}_i, w^k_i)
		+ \sum_i (\varphi^{k+1}_{0i} - \varphi^k_{0i}) w^k_i.
	\end{align*}
	The second and fourth terms are non-positive, the third term is bounded by
	$
		O \left( \eps \olw^k_0 + \sqrt{\eta} V(z^{k+1}) \right)
    $
	by \autoref{lm:proof_gencase:errterm1__control_KLdiv_uI_vI_errterm},
	and so is the last term
	by \autoref{lm:proof_gencase:control_varphi_t+1_t}.
	Further upper-bound the first term, using that $\olw^k_I \geq \wulb$ (\autoref{lm:proof_gencase:ulb_on_weights_locally_k+1}), by
	\begin{align*}
		\sum_I W_I \log \frac{\olw^{k+1}_I}{\olw^k_I} 
		&\leq \sum_I \frac{W_I}{\olw^k_I} \left( \olw^{k+1}_I-\olw^k_I \right) 
        \leq \frac{1}{\wulb} \sum_I \abs{\olw^{k+1}_I-\olw^k_I} \\
		&\leq \frac{\sqrt{n^* + m^*}}{\wulb} \sqrt{ \sum_I \abs{\olw^{k+1}_I-\olw^k_I}^2 }
		\leq \frac{\sqrt{2 (n^* + m^*)}}{\wulb} \sqrt{\sum_I d_h(\olw^{k+1}_I, \olw^k_I)}
	\end{align*}
	since $h$ is $1$-strongly convex over $[0,1]$.
    Thus,
	\begin{equation*}
		\sum_i (w_i-w^{k+1}_i) \log \frac{w^{k+1}_i}{w^k_i}
		\leq 
		\frac{\sqrt{2 (n^* + m^*)}}{c} \sqrt{\sum_I d_h(\olw^{k+1}_I, \olw^k_I)}
		+ O \left( \eps \olw^k_0 + \sqrt{\eta} V(z^{k+1}) \right).
	\end{equation*}

	By putting the two parts together and using that $\sqrt{A} + \sqrt{B} \leq \sqrt{2} \sqrt{A+B}$, we obtain
	\begin{equation*}
		\eta \max_{A,X,B,Y} \tgap(z; z^{k+1}) 
		\leq 
        2 \left( \frac{\sqrt{n^* + m^*}}{\wulb} \vee \diamXY \sqrt{\frac{\eta}{\sigma}} \right)
        \sqrt{D(k+1, k)}
        ~+~ O \left( \eps \olw^k_0 + \sqrt{\eta} V(z^{k+1}) \right)
	\end{equation*}
    and the claim follows directly.
\end{proof}

% \paragraph{Lower part of the sandwich.}
\paragraph{Lower bound on the gradient-norm-like quantity.}

As a direct application of \autoref{lm:growth_conds:error_bound},
there exist $r', C' > 0$ only dependent on $(f, \XXX, \YYY)$ and $\Gamma_0$ such that if $V(z^{k+1}) \leq r'$, then
\begin{equation*}
	\max_{A,X,B,Y} \tgap(z; z^{k+1}) 
	\geq 
	C' \sqrt{
		\sum_I d_h(w^*_I, \olw^{k+1}_I)
		+ \sum_I \olw^{k+1}_I \norm{\Delta \olp^{k+1}_I}^2
	}
	+ O \left( V(z^{k+1}) \right).
\end{equation*}
Note that thanks to \autoref{lm:proof_gencase:Vk+1_small}, we can indeed assume $V(z^{k+1}) \leq r'$ by choosing $r$ small enough in the statement of \autoref{lm:proof_gencase:lb_Dk+1_k}.

\paragraph{Putting the two bounds together.}
All in all, we showed that (assuming $V(z^k) \leq r$ for some $r$ small enough)
\begin{equation*}
    \frac{C}{\eta}
    \sqrt{D(k+1, k)}
    + \frac{1}{\eta} O \left( \eps \olw^k_0 + \sqrt{\eta} V(z^{k+1}) \right) 
	\geq C' \sqrt{
		\sum_I d_h(w^*_I, \olw^{k+1}_I)
		+ \sum_I \olw^{k+1}_I \norm{\Delta \olp^{k+1}_I}^2
	}
	+ O \left( V(z^{k+1}) \right)
\end{equation*}
for some $C, C'$ only dependent on $(f, \XXX, \YYY)$ and $\Gamma_0$.
Rearranging and taking squares,
\begin{equation*}
    D(k+1, k) \geq \left(\frac{C'}{C}\right)^2 \eta^2 \left(
        \sum_I d_h(w^*_I, \olw^{k+1}_I)
		+ \sum_I \olw^{k+1}_I \norm{\Delta \olp^{k+1}_I}^2
    \right)
    + O\left( (\eps \olw^k_0)^2 + \eta V(z^{k+1})^2 \right).
\end{equation*}
This concludes the proof of \autoref{lm:proof_gencase:lb_Dk+1_k}.

\subsection{Proof conclusion} \label{subsec:proof_gencase:proof_conclusion}

It just remains to put everything together by substituting the terms by their lower bound in \eqref{eq:proof_gencase:ppa_ineq_star}.
Our choice of $\lambda, \tau$ and our assumption that $\sigma \asymp \eta$ simplify things considerably.
The only subtlety is that some of the terms in the upper bound of \autoref{lm:proof_gencase:errterm1__control_KLdiv_uI_vI_errterm} and \autoref{lm:proof_gencase:errterm2__control_tricky_errterms} need to be compensated by the lower bound of \autoref{lm:proof_gencase:lb_tgap_olw0}, which can be done 
by assuming $\eta, \sigma$ small enough.

In the remainder of this subsection, $C_1, C_2, ... >0$ will denote constants dependent only on $(f, \XXX, \YYY)$ and $\Gamma_0$
(the same things as what we hide in $O(\cdot)$).

\paragraph{Putting all the bounds together.}
Assume that $V(z^k) \leq r$ and that $r, \eta, \sigma$ are small enough so that all of the lemmas apply.
Just substitute the terms in \eqref{eq:proof_gencase:ppa_ineq_star} by their bounds:
\begin{equation*}
    V(z^{k+1}) - V(z^k) 
    \leq 
    -\eta \tgap(\zps; z^{k+1})
    - D(k+1, k)
    + \mathrm{(err1)}
    + \mathrm{(err2)}
    + \mathrm{(err3)}.
\end{equation*}
By \autoref{lm:proof_gencase:lb_tgap_olw0}, there exist $C_1, C_2$ such that
\begin{equation*}
	\eta \tgap(\zps; z^{k+1}) \geq \eta \cdot C_1 \sqrt{\sigma} \olw^{k+1}_0 + \eta \cdot C_2 V_\pos(z^{k+1})
	+ O \left( \eta V(z^{k+1})^{3/2} \right).
\end{equation*}
By \autoref{lm:proof_gencase:lb_Dk+1_k}, there exists $C_3$ such that
\begin{equation*}
	D(k+1, k) 
	\geq C_3 \eta^2 \left(
		\sum_I d_h(w^*_I, \olw^{k+1}_I)
		+ \sum_I \olw^{k+1}_I \norm{\Delta \olp^{k+1}_I}^2
	\right)
	+ O \left( \left( \eps \olw^k_0 \right)^2 + \eta V(z^{k+1})^2 \right).
\end{equation*}
By \autoref{lm:proof_gencase:errterm1__control_KLdiv_uI_vI_errterm},
\begin{align*}
    \mathrm{(err1)} 
	&\lesssim
    \sqrt{\eta} V(z^{k+1})^{3/2}
	+ \eps \olw^k_0 \cdot \sqrt{V(z^{k+1})} \\
	&\lesssim
	\sqrt{\eta} V(z^{k+1})^{3/2}
	+ \eps \left( \olw^k_0 \right)^2
	+ \eps V(z^{k+1}).
\end{align*}
By \autoref{lm:proof_gencase:errterm2__control_tricky_errterms},
\begin{equation*}
    \mathrm{(err2)} 
	\lesssim 
	\eps \olw^{k+1}_0
	+ \eta^{3/2} V_\pos(z^{k+1})
	+ \sqrt{\eta} V(z^{k+1})^{3/2}.
\end{equation*}
By \autoref{lm:proof_gencase:errterm3__control_nonBregmanness_errterm},
\begin{equation*}
    \mathrm{(err3)} 
	\lesssim \eta V(z^{k+1})^{3/2}.
\end{equation*}
All in all,
since $\olw^k_0 = O(V(z^k))$ by Eq.~\eqref{eq:aux_lemmas:decomp_Vwei}, 
we get 
\begin{align*}
	V(z^{k+1}) - V(z^k)
	&\leq 
	- \left( C_1 \eta^{3/2} - O(\eps) \right) \olw^{k+1}_0 
	- \left( C_2 - O(\sqrt{\eta}) \right) \eta V_\pos(z^{k+1}) \\
	&~~~~ - C_3 \eta^2 \left(
		\sum_I d_h(w^*_I, \olw^{k+1}_I)
		+ \sum_I \olw^{k+1}_I \norm{\Delta \olp^{k+1}_I}^2
	\right) \\
	&~~~~ + O \left( \sqrt{\eta} V(z^{k+1})^{3/2} \right)
	+ O \left( \eps V(z^{k+1}) \right)
	+ O \left( \eps V(z^k)^2 \right).
\end{align*}

To be explicit, this means that there exists a constant $M>0$ dependent only on $(f, \XXX, \YYY)$ and $\Gamma_0$ such that
\begin{align*}
	V(z^{k+1}) - V(z^k)
	&\leq 
	- \left( C_1 \eta^{3/2} -  M \eps \right) \olw^{k+1}_0 
	- \left( C_2 - M \sqrt{\eta} \right) \eta V_\pos(z^{k+1}) \\
	&~~~~ - C_3 \eta^2 \left(
    	\sum_I d_h(w^*_I, \olw^{k+1}_I)
    	+ \sum_I \olw^{k+1}_I \norm{\Delta \olp^{k+1}_I}^2
	\right) \\
	&~~~~ + M \left( 
		\sqrt{\eta} V(z^{k+1})^{3/2}
		+ \eps V(z^{k+1})
		+ \eps V(z^k)^2
	\right).
\end{align*}
Since $\eps = e^{-\lambda^3/3} = e^{-\frac{1}{3 \sqrt{\sigma}}}$ and $\sigma \asymp \eta$, then for small enough $\eta, \sigma$, we have $C_1 \eta^{3/2} -  M \eps \geq \frac{C_1}{2} \eta^{3/2}$.
Furthermore, for small enough $\eta$, we have $C_2 - M \sqrt{\eta} \geq \frac{C_2}{2}$.
Thus, there exists $C_4>0$ such that
\begin{equation*}
	V(z^{k+1}) - V(z^k)
	\leq 
	-C_4 \eta^2 V(z^{k+1})
	+ M \left( 
		\sqrt{\eta} V(z^{k+1})^{3/2}
		+ \eps V(z^{k+1})
		+ \eps V(z^k)^2
	\right).
\end{equation*}
Moreover for small enough $\eta, \sigma$, we have $C_4 \eta^2 - M \eps \geq \frac{C_4}{2} \eta^2$.
Then,
\begin{gather*}
	V(z^{k+1}) - V(z^k)
	\leq 
	-(C_4/2) \eta^2 V(z^{k+1})
	+ M \left( 
		\sqrt{\eta} V(z^{k+1})^{3/2}
		+ \eps V(z^k)^2
	\right) \\
	V(z^{k+1}) \left[ 1 + (C_4/2) \eta^2 - M \sqrt{\eta} \sqrt{V(z^{k+1})} \right]
	\leq V(z^k) \left[ 1 + M \eps V(z^k) \right].
\end{gather*}

\paragraph{Sufficient decrease of the Lyapunov function.}
For a fixed $r_0>0$ to be chosen (small enough so that all of the lemmas apply), assume that $V(z^k) \leq r_0$.
By \autoref{lm:proof_gencase:Vk+1_small}, we can assume $\eta, \sigma$ small enough such that 
$V(z^{k+1}) \leq 2 r_0$.
We then have that
\begin{align}
	V(z^{k+1}) \left[ 1 + (C_4/2) \eta^2 - M \sqrt{2 \eta r_0} \right]
	&\leq V(z^k) \left[ 1 + M \eps r_0 \right] \\
	\frac{V(z^{k+1})}{V(z^k)}
	&\leq \frac{
		1 + M \eps r_0
	}{
		1 + (C_4/2) \eta^2 - M \sqrt{2 \eta r_0}
	}
    ~\eqqcolon 1-\kappa.
\label{eq:proof_gencase:pre_kappa}
\end{align}
Clearly $r_0$ can be chosen (dependent on $\eta$)
such that the right-hand side is strictly less than $1$, i.e., $\kappa>0$.
By induction, if $V(z^0) \leq r_0$, then for all $k$, $V(z^{k+1}) \leq r_0$ and 
\begin{align*}
	\forall k,~
	\frac{V(z^{k+1})}{V(z^k)} &\leq 1 - \kappa \\
	V(z^k) &\leq V(z^0) (1-\kappa)^k.
\end{align*}
This concludes the proof of \autoref{thm:cv_proof:gencase:loc_exp_cv}.

\begin{remark}
    More precisely, for the right-hand side of \eqref{eq:proof_gencase:pre_kappa} to be less than $1$, $r_0$ needs to be chosen less than $\eta^3$ times a constant (dependent on $(f, \XXX, \YYY)$ and $\Gamma_0$).
    The rate $\kappa$ can be seen to be of order $\eta^2$, for any admissible choice of $r_0$.
\end{remark}

% \end{document}

\ifextended%
    % !TEX root = ../main.tex
% \documentclass[../main]{subfiles}
% \begin{document}

\subsection{Delayed technical proofs} \label{subsec:proof_gencase:delayed_proofs}

In some of the proofs of this subsection we use the expressions and a priori bounds for $V_\pos$ and $V_\wei$ from \autoref{subsec:aux_lemmas:useful_exprs_V} without explicit mention.

\subsubsection{Auxiliary claims}

As a consequence of the fact that $\norm{x^{k+1}_i-x^k_i} = O(\sigma)$ (\autoref{lm:proof_gencase:bound_movement}),
we can meaningfully classify the particles according to which $\support(\varphi_I)$ they belong to, both at $k$ and at $k+1$.
For a fixed $k$, denote
\begin{equation} \label{eq:proof_gencase:def_NNNI}
    \forall I \in [n^*],~ 
    \NNN(I) = \left\lbrace
        i 
        % ~~\text{s.t}~~
        ;~~
        x^{k+1}_i \in B_{x^*_I, \lambda\tau} 
        ~\text{or}~
        x^k_i \in B_{x^*_I, \lambda\tau} 
	\right\rbrace
	~~~~\text{and}~~~~
	\NNN(0) = [n] \setminus \left( \cup_I \NNN(I) \right).
\end{equation}
Since $x^{k+1}_i-x^k_i = O(\sigma)$, for $\sigma$ chosen small enough compared to $\min_{I \neq I'} \norm{x^*_I-x^*_{I'}}$ we have that
$x^k_i \in B_{x^*_I, \lambda\tau} \implies \forall I' \neq I,~ x^k_i, x^{k+1}_i \not\in B_{x^*_{I'}, \lambda\tau}$, and so the $\NNN(I)$ are pair-wise disjoint.
In other words, 
$\bigsqcup_{I \in [0,n^*]} \NNN(I)$ then forms a partition of $[n]$; and similarly for the $(y_j)_{j \in [m]}$.
\emph{In the remainder of this section, we assume $\sigma$ small enough so that this is the case.}

% \TODO{the notation $i \in \NNN(I)$ can be confusing because it looks symmetric, it should be $i \in \mathrm{cl}(I)$, but that's already much heavier notation and we need to use it often... Maybe $i \rightsquigarrow I$?}
% \TODO{Can write $\NNN(I) = \left\lbrace i; i \in \NNN(I) \right\rbrace$ and we show that $[n] = \sqcup_I \NNN(I)$}

Let 
\begin{equation*}
	\forall x \in \XXX,~ \tvarphi_I(x) = \exp\left( -\frac{\norm{x-x^*_I}^3}{3 \tau^3} \right)
\end{equation*}
so that $\varphi_I(x)$ coincides with $\tvarphi_I(x)$ if and only if
$\norm{x-x^*_I} \leq \lambda\tau$.

\begin{claim} \label{claim:proof_gencase:delayed_proofs:bound_varphi_tvarphi}
	For any $I \in [n^*]$,
	\begin{equation*}
		\sum_i (\varphi^{k+1}_{Ii} -\varphi^k_{Ii}) a^k_i
		= \sum_{i \in \NNN(I)} (\tvarphi^{k+1}_{Ii} -\tvarphi^k_{Ii}) a^k_i
		+ O \left( \eps (\ola^{k+1}_0 + \ola^k_0) \right).
	\end{equation*}
% 	\TODO{what we prove is stronger: that $\sum_i \varphi^k_{Ii} a^k_i = \sum_{i \in \NNN(I)} \tvarphi^k_{Ii} a^k_i + O(\eps \hola^k_0)$, and similarly (with an extra $\eta$ in the big-O) for $\varphi^{k+1}$} --> hmm not sure actually...
\end{claim}

\begin{proof}
	Since $\varphi_I(x)$ coincides with $\tvarphi_I(x)$ if and only if 
	$\norm{x-x^*_I} \leq \lambda\tau$,
	\begin{itemize}
	    \item if $i \not\in \NNN(I)$, i.e., if both $x^{k+1}_i, x^k_i \not\in B_{x^*_I, \lambda\tau}$, then
	    $\varphi^{k+1}_{Ii} - \varphi^k_{Ii} = 0$;
	    \item if $i \in \NNN(I)$,
	    \begin{multline*}
    		\abs{
    			( \varphi^{k+1}_{Ii} - \varphi^k_{Ii} )
    			- ( \tvarphi^{k+1}_{Ii} - \tvarphi^k_{Ii} )
    		}
    		\leq \abs{ \varphi^{k+1}_{Ii} - \tvarphi^{k+1}_{Ii} }
    		+ \abs{ \varphi^k_{Ii} - \tvarphi^k_{Ii} } \\
    		\leq \eps \cdot \ind\left[ x^{k+1}_i \not\in B_{x^*_I, \lambda\tau} ~\wedge~ x^k_i \in B_{x^*_I, \lambda\tau} \right]
    		~+~ \eps \cdot \ind\left[ x^{k+1}_i \in B_{x^*_I, \lambda\tau} ~\wedge~ x^k_i \not\in B_{x^*_I, \lambda\tau} \right].
	    \end{multline*}
	\end{itemize}
	Further note that, by definition,
	\begin{align*}
		\sum_I \sum_{i \in \NNN(I)} \eps \ind\left[ x^{k+1}_i \not\in B_{x^*_I, \lambda\tau} ~\wedge~ x^k_i \in B_{x^*_I, \lambda\tau} \right] a^k_i
		&\leq \eps \sum_i \varphi_0(x^{k+1}_i) a^k_i
		\leq \eps (1+O(\eta)) \ola^{k+1}_0 \\
		\text{and}~~~
		\sum_I \sum_{i \in \NNN(I)} \eps \ind\left[ x^{k+1}_i \in B_{x^*_I, \lambda\tau} ~\wedge~ x^k_i \not\in B_{x^*_I, \lambda\tau} \right] a^k_i
		&\leq \eps \sum_i \varphi_0(x^k_i) a^k_i
		= \eps \ola^k_0.
	\end{align*}
	Thus
	\begin{equation*}
	    \abs{
    	    \sum_i (\varphi^{k+1}_{Ii} -\varphi^k_{Ii}) a^k_i
    		- \sum_{i \in \NNN(I)} (\tvarphi^{k+1}_{Ii} -\tvarphi^k_{Ii}) a^k_i
    	}
		\leq \eps \left( (1+O(\eta)) \ola^{k+1}_0 + \ola^k_0 \right)
	\end{equation*}
	and hence the announced estimate.
\end{proof}

The following claim follows from a Taylor expansion of $x \mapsto \norm{x-x^*_I}^3$.
\begin{claim} \label{claim:proof_gencase:delayed_proofs:taylor_expan_norm_pow3}
	For any $i, I$, we have
% 	\begin{align*}
% % 		\abs{ \norm{x^{k+1}_i-x^*_I}^3 - \norm{x^k_i-x^*_I}^3 }
% % 		&\leq 3 \norm{x^{k+1}_i-x^k_i} \norm{x^{k+1}_i-x^*_I}^2 
% % 		+ 3 \norm{x^{k+1}_i-x^k_i}^2 \norm{x^{k+1}_i-x^*_I} \\
% % 		&~~~~ + 
% % 		O \left( \norm{x^{k+1}_i-x^k_i}^3 \right)
% 		\abs{ \norm{x^{k+1}_i-x^*_I}^3 - \norm{x^k_i-x^*_I}^3 }
% 		&\leq 3 \norm{x^{k+1}_i-x^k_i} \norm{x^{k+1}_i-x^*_I}^2 
% 		+ O \left( \norm{x^{k+1}_i-x^k_i}^2 \norm{x^{k+1}_i-x^*_I} \right)
% 	\end{align*}
% 	and 
	\begin{align*}
% 		\norm{x^{k+1}_i-x^*_I}^3 - \norm{x^k_i-x^*_I}^3
% 		&\leq 3 \innerprod{x^{k+1}_i-x^k_i}{x^{k+1}_i-x^*_I} \norm{x^{k+1}_i-x^*_I}
% 		+ 3 \norm{x^{k+1}_i-x^k_i}^2 \norm{x^{k+1}_i-x^*_I} \\
% 		&~~~~ + 
% 		%		3 (3-1) (3-2) 
% 		O \left( \norm{x^{k+1}_i-x^k_i}^3 \right).
		\norm{x^{k+1}_i-x^*_I}^3 - \norm{x^k_i-x^*_I}^3
		&= 3 \innerprod{x^{k+1}_i-x^k_i}{x^{k+1}_i-x^*_I} \norm{x^{k+1}_i-x^*_I} \\
% 		+ O \left( \norm{x^{k+1}_i-x^k_i}^3 \right).
		&~~~~ + O \left( \norm{x^{k+1}_i-x^k_i}^2 \norm{x^{k+1}_i-x^*_I} + \norm{x^{k+1}_i-x^k_i}^3 \right).
	\end{align*}
% 
% 	\TODO{do we use the order-3 expansions anywhere?...} --> indeed I don't think so. I commented it out
\end{claim}

\begin{proof}
	The first and second derivatives of $\norm{\cdot-x^*_I}^3$ are given, up to translation, by
	\begin{align*}
		(\nabla \norm{\cdot}^3)(x) 
		&= 3 \norm{x} x
		& &\text{and} &
		0 \preceq
		(\nabla^2 \norm{\cdot}^3)(x) 
		&=
		3 \norm{x} \id + 3 \norm{x} \frac{x x^\top}{\norm{x}^2}
		\preceq 6 \norm{x} \id.
	\end{align*}
	By Taylor expansion of $\norm{\cdot-x^*_I}^3$ centered at $x^{k+1}_i$ with remainder in Lagrange form,
	there exists $\theta \in [0,1]$ such that
	\begin{equation*}
		\norm{x^{k+1}_i-x^*_I}^3 - \norm{x^k_i-x^*_I}^3
		= 3 \innerprod{x^{k+1}_i-x^k_i}{x^{k+1}_i-x^*_I} \norm{x^{k+1}_i-x^*_I}
		- \bR
	\end{equation*}
	where
	\begin{align*}
	    \bR &= \frac{1}{2} (x^{k+1}_i-x^k_i)^\top 
	    \left[
	        \left( \nabla^2 \norm{\cdot-x^*_I}^3 \right) 
	        \left( \theta x^{k+1}_i + (1-\theta) x^k_i \right)
	    \right] 
	    (x^{k+1}_i-x^k_i) \\
	    2 \abs{\bR} &\leq 6 \norm{x^{k+1}_i-x^k_i}^2 \norm{\theta (x^{k+1}_i-x^*_I) + (1-\theta) (x^k_i-x^*_I) } \\
	    &\leq 6 \norm{x^{k+1}_i-x^k_i}^2 \Big( 
	        \norm{x^{k+1}_i-x^*_I} 
	        + \underbrace{ \norm{x^k_i-x^*_I} }_{\leq \norm{x^{k+1}_i-x^*_I} + \norm{x^{k+1}_i-x^k_i}}
        \Big) \\
	    &\leq 12 \norm{x^{k+1}_i-x^k_i}^2 \norm{x^{k+1}_i-x^*_I}
	    + 6 \norm{x^{k+1}_i-x^k_i}^3.
     \rqedhere
    \end{align*}
\end{proof}

We will repeatedly use the following Taylor expansions of the local payoff matrices.
\begin{claim}
	For any $i, I$,
	\begin{equation}
	\label{eq:proof_gencase:delayed_proofs:expan_Mb_i}
		(\gmat^{k+1} b^{k+1})_i
		= \gval 
		+ \frac{1}{2} \norm{x^{k+1}_i-x^*_I}_{H_I}^2
		+ O(\norm{x^{k+1}_i-x^*_I}^3)
		+ O\left( \sqrt{V(b^{k+1}, y^{k+1})} \right)
	\end{equation}
	and more precisely if $\norm{x^{k+1}_i-x^*_I} \leq \frac{\sigma_{\min}}{2 \smoothnessf_3}$, then
	\begin{equation}
	\label{eq:proof_gencase:delayed_proofs:expan_Mb_i_pos}
		(\gmat^{k+1} b^{k+1})_i
		\geq \gval + \frac{1}{4} \sigma_{\min} \norm{x^{k+1}_i-x^*_I}^2 
		+ O\left( \sqrt{V(b^{k+1}, y^{k+1})} \right).
	\end{equation}
	Furthermore,
	\begin{equation}
	\label{eq:proof_gencase:delayed_proofs:expan_partial_x_Mb_i}
		\partial_x \gmat^{k+1}_{i \bullet} b^{k+1}
		= 
		(x^{k+1}_i-x^*_I)^\top H_I
		+ O(\norm{x^{k+1}_i-x^*_I}^2)
		+ O\left( \min\left\{ 1, \sqrt{V(b^{k+1}, y^{k+1})} \right\} \right) 
	\end{equation}
	and more precisely if $\norm{x^{k+1}_i-x^*_I} \leq \frac{\sigma_{\min}}{2 \smoothnessf_3}$, then
	\begin{equation}
	\label{eq:proof_gencase:delayed_proofs:expan_partial_x_Mb_i_pos}
		(x^{k+1}_i-x^*_I)^\top \partial_x \gmat^{k+1}_{i \bullet} b^{k+1}
		\geq \frac{\sigma_{\min}}{2} \norm{x^{k+1}_i-x^*_I}^2
		+ O\left( \norm{x^{k+1}_i-x^*_I} \sqrt{V(b^{k+1}, y^{k+1})} \right).
	\end{equation}
\end{claim}

Note that our choice of $\lambda, \tau$ implies $\lambda\tau \leq \frac{\sigma_{\min}}{4 \smoothnessf_3}$,
and that for all $i \in \NNN(I)$, $\norm{x^{k+1}_i-x^*_I} \leq \lambda\tau + O(\sigma)$.
\emph{In the remainder of this section, we assume $\sigma$ small enough so that $\norm{x^{k+1}_i-x^*_I} \leq \frac{\sigma_{\min}}{2 \smoothnessf_3}$ holds for all $i \in \NNN(I)$} and similarly for the $y_j, y^*_J$.

\begin{proof}
	To lighten notation in the calculations,
    denote $\hx = x^{k+1}$, $\hy = y^{k+1}$ and $\hb = b^{k+1}$.
    
    By Taylor expansion, for all $i, I$,
	\begin{align*}
		&
        (\gmat^{k+1} b^{k+1})_i = (\gmathh \hb)_i
		= \sum_J \sum_j \hpsi_{Jj} \left[ \gmaths_{iJ} + \gmathh_{ij} - \gmaths_{iJ} \right] \hb_j
		+ \sum_j \hpsi_{0j} \gmathh_{ij} \hb_j \\
		&= \sum_J \gmaths_{iJ} \holb_J 
		+ \sum_J \sum_j \hpsi_{Jj} O(\norm{\hy_j-y^*_J}) \hb_j
		+ O(\holb_0) \\
		&= \sum_J \left[ \gmatss_{IJ} + (\hx_i-x^*_I)^\top \partial_x \gmatss_{IJ} + \frac{1}{2} ((\hx_i-x^*_I)^2)^\top \partial_{xx}^2 \gmatss_{IJ} + O(\norm{\hx_i-x^*_I}^3) \right] \holb_J 
        ~ + O \left( \sqrt{V_\pos(\hb, \hy)}  \,+ \holb_0 \right) \\
		&= \gval + \sum_J \gmatss_{IJ} \Delta \holb_J 
		+ \sum_J (\hx_i-x^*_I)^\top \partial_x \gmatss_{IJ} \Delta \holb_J
		+ \frac{1}{2} \norm{\hx_i-x^*_I}_{H_I}^2
		+ O(\norm{\hx_i-x^*_I}^3) 
        ~ + O\left( \sqrt{V(\hb, \hy)} \right) \\
		&= \gval 
		+ \frac{1}{2} \norm{\hx_i-x^*_I}_{H_I}^2
		+ O(\norm{\hx_i-x^*_I}^3)
		+ O\left( \sqrt{V(\hb, \hy)} \right).
	\end{align*}
	More precisely, 
	\begin{align*}
		(\gmathh \hb)_i
		&\geq \gval + \frac{1}{2} \norm{\hx_i-x^*_I}_{H_I}^2
		- \frac{\smoothnessf_3}{2} \norm{\hx_i-x^*_I}^3
		+ O\left( \sqrt{V(\hb, \hy)} \right) \\
        &\geq \gval + \frac{1}{4} \sigma_{\min} \norm{\hx_i-x^*_I}^2 
		+ O\left( \sqrt{V(\hb, \hy)} \right)
        ~~~~\text{if}~~ 
        \norm{\hx_i-x^*_I} \leq \frac{\sigma_{\min}}{2 \smoothnessf_3}.
	\end{align*}

	Also by Taylor expansion, for all $i, I$,
	\begin{align*}
		\MoveEqLeft
		\partial_x \gmat^{k+1}_{i \bullet} b^{k+1} = \partial_x \gmathh_{i \bullet} \hb
		= \sum_J \sum_j \hpsi_{Jj} \left[ \partial_x \gmaths_{iJ} + \partial_x \gmathh_{ij} - \partial_x \gmaths_{iJ} \right] \hb_j
		+ \sum_j \hpsi_{0j} \partial_x \gmathh_{ij} \hb_j \\
		&= \sum_J \partial_x \gmaths_{iJ} \holb_J
		+ \sum_J \sum_j \hpsi_{Jj} O(\norm{\hy_j-y^*_J}) \hb_j
		+ O(\holb_0) \\
		&= \sum_J \left[ \partial_x \gmatss_{IJ} + (\hx_i-x^*_I)^\top \partial_{xx}^2 \gmatss_{IJ} + O(\norm{\hx_i-x^*_I}^2) \right] \holb_J
		+ O\left( \sqrt{V_\pos(\hb, \hy)} \right)
		+ O(\holb_0) \\
		&= \partial_x \gmatss_{I\bullet} \Delta \holb
		+ \sum_J (\hx_i-x^*_I)^\top \partial_{xx}^2 \gmatss_{IJ} \holb_J
		+ O(\norm{\hx_i-x^*_I}^2)
		+ O\left( \min\left\{ 1, \sqrt{V(\hb, \hy)} \right\} \right) \\
		&= (\hx_i-x^*_I)^\top H_I
		+ O(\norm{\hx_i-x^*_I}^2)
		+ O\left( \min\left\{ 1, \sqrt{V(\hb, \hy)} \right\} \right).
	\end{align*}
    On lines 4 and 5, the fact that the last error term is $O(1)$ can be checked by noting that $\holb_0, \norm{\Delta \holb} \leq 2$ and $V_\pos(\hb, \hy) = \sum_J \sum_j \hpsi_{Jj} \hb_j \norm{\hy_j-y^*_J}^2 \leq \diamXY^2$.
	More precisely, for any $\delta x$,
	\begin{equation*}
		\innerprod{\delta x}{\partial_x \gmathh_{i \bullet} \hb}
		\geq \innerprod{\delta x}{H_I (\hx_i-x^*_I)}
		- \smoothnessf_3 \norm{\delta x} \norm{\hx_i-x^*_I}^2
		+ O\left(  \norm{\delta x} \sqrt{V(\hb, \hy)} \right).
	\end{equation*}
	So if $\norm{\hx_i-x^*_I} \leq \frac{\sigma_{\min}}{2 \smoothnessf_3}$, then
	\begin{align*}
		(\hx_i-x^*_I)^\top \partial_x \gmathh_{i \bullet} \hb
		&\geq \sigma_{\min} \norm{\hx_i-x^*_I}^2 - \smoothnessf_3 \norm{\hx_i-x^*_I}^3
		+ O\left( \norm{\hx_i-x^*_I} \sqrt{V(\hb, \hy)} \right) \\
		&\geq \frac{\sigma_{\min}}{2} \norm{\hx_i-x^*_I}^2
		+ O\left( \norm{\hx_i-x^*_I} \sqrt{V(\hb, \hy)} \right).
    \rqedhere
	\end{align*}
\end{proof}

\subsubsection{Proof of \autoref{lm:proof_gencase:control_varphi_t+1_t}}

% In the following proof we use, for the first time, our choice of $\lambda^3 = \frac{1}{\sqrt{\sigma}}$ and $\lambda\tau \asymp 1$.  --> remarked under the lemma statement in Sec. "Preliminary lemmas"
Note that the proof does not make use of the fact that $V(z^{k+1}) = O(1)$
(inequality~\eqref{eq:proof_gencase:etaveesigma_V_bounded}), so as to avoid circular reasoning since we showed that fact as a consequence of \autoref{lm:proof_gencase:control_varphi_t+1_t}.

\begin{proof}
	Fix $I \in [n^*]$.
	We showed in \autoref{claim:proof_gencase:delayed_proofs:bound_varphi_tvarphi} that
	\begin{equation*}
		\sum_i (\varphi^{k+1}_{Ii} -\varphi^k_{Ii}) a^k_i
		= \sum_{i \in \NNN(I)} (\tvarphi^{k+1}_{Ii} -\tvarphi^k_{Ii}) a^k_i
		+ O \left( \eps (\ola^{k+1}_0 + \ola^k_0) \right)
	\end{equation*}
	and it remains to upper- and lower-bound the first term.

	First note that
	\begin{equation} \label{eq:proof_gencase:delayed_proofs:bound_psialpha_t+1_t}
		\tvarphi^k_{Ii} - \tvarphi^{k+1}_{Ii}
		= \tvarphi^k_{Ii} \cdot \left(\! 1 - \frac{\tvarphi^{k+1}_{Ii}}{\tvarphi^k_{Ii}} \right)
		\leq \tvarphi^k_{Ii} \cdot \left( \log \tvarphi^k_{Ii} - \log \tvarphi^{k+1}_{Ii} \right)
		= \tvarphi^k_{Ii} \cdot
		\frac{1}{3 \tau^3}
		\left(
			\norm{x^{k+1}_i-x^*_I}^3 - \norm{x^k_i-x^*_I}^3
		\right)
	\end{equation}
	and that
	\begin{equation*}
		\tvarphi^k_{Ii} - \tvarphi^{k+1}_{Ii}
		= \tvarphi^{k+1}_{Ii} \cdot \left( \frac{\tvarphi^k_{Ii}}{\tvarphi^{k+1}_{Ii}} - 1 \!\right)
		\geq \tvarphi^{k+1}_{Ii} \cdot \left( \log \tvarphi^k_{Ii} - \log \tvarphi^{k+1}_{Ii} \right)
		= \tvarphi^{k+1}_{Ii} \cdot
		\frac{1}{3 \tau^3}
		\left(
			\norm{x^{k+1}_i-x^*_I}^3 - \norm{x^k_i-x^*_I}^3
		\right).
	\end{equation*}
	Furthermore, by \autoref{claim:proof_gencase:delayed_proofs:taylor_expan_norm_pow3},
	\begin{align*}
		\abs{ \norm{x^{k+1}_i-x^*_I}^3 - \norm{x^k_i-x^*_I}^3 }
		&\leq 3 \norm{x^{k+1}_i-x^k_i} \norm{x^{k+1}_i-x^*_I}^2 \\
        &~~~~ + O \left( \norm{x^{k+1}_i-x^k_i}^2 \norm{x^{k+1}_i-x^*_I} + \norm{x^{k+1}_i-x^k_i}^3 \right) \\
        &\lesssim \norm{x^{k+1}_i-x^k_i} \norm{x^{k+1}_i-x^*_I} + \norm{x^{k+1}_i-x^k_i}^2
	\end{align*}
    and by
	the update equation \eqref{eq:proof_gencase:ppa_upd_eqs}
	and the expansion \eqref{eq:proof_gencase:delayed_proofs:expan_partial_x_Mb_i},
	since $a^k_i = (1+O(\eta)) a^{k+1}_i \asymp a^{k+1}_i$ by \autoref{lm:proof_gencase:bound_movement},
	\begin{equation*}
	    \norm{x^{k+1}_i-x^k_i}
	    = \sigma \frac{a^{k+1}_i}{a^k_i} \norm{\partial_x \gmat^{k+1}_{i \bullet} b^{k+1}} 
        \lesssim \sigma
	    \left( 
    		\norm{x^{k+1}_i-x^*_I}
    		~+~ 1 \wedge \sqrt{V(b^{k+1}, y^{k+1})}
		\right),
	\end{equation*}
	and so
	\begin{align*}
		\MoveEqLeft \abs{ \norm{x^{k+1}_i-x^*_I}^3 - \norm{x^k_i-x^*_I}^3 } \\
		&\lesssim \sigma \norm{x^{k+1}_i-x^*_I}^2 + \sigma 
        \left[ 1 \wedge \sqrt{V(b^{k+1}, y^{k+1})} \right]
        \norm{x^{k+1}_i-x^*_I} 
        + \sigma^2 \left[ 1 \wedge V(b^{k+1}, y^{k+1}) \right] \\
		&\lesssim \sigma \norm{x^{k+1}_i-x^*_I}^2 + \sigma \left[ 1 \wedge V(b^{k+1}, y^{k+1}) \right]
	\end{align*}
    where we used that $ab \leq \frac{a^2}{2} + \frac{b^2}{2}$ to bound the second term of the first line.
	So
% 	since $a^k_i = (1+O(\eta)) a^{k+1}_i \lesssim a^{k+1}_i$,
% 	and $\norm{x^{k+1}_i-x^*_I} = O(\lambda\tau+\sigma) = O(\lambda\tau)$ for all $i \in \NNN(I)$
% 	by \autoref{lm:proof_gencase:bound_movement},
	\begin{align*}
		\sum_{i \in \NNN(I)} (\tvarphi^k_{Ii} -\tvarphi^{k+1}_{Ii}) a^k_i
		&\leq \frac{1}{3 \tau^3} \sum_{i \in \NNN(I)} \tvarphi^k_{Ii} a^k_i \cdot 
		\left( \norm{x^{k+1}_i-x^*_I}^3 - \norm{x^k_i-x^*_I}^3 \right) \\
		&\lesssim \sigma \frac{\lambda^3}{(\lambda\tau)^3} 
		\sum_{i \in \NNN(I)} \tvarphi^k_{Ii} a^{k+1}_i 
		\left( 
		    \norm{x^{k+1}_i-x^*_I}^2
		    + \left[ 1 \wedge V(b^{k+1}, y^{k+1}) \right]
		\right) \\
		&\lesssim \sqrt{\sigma}
		\left(
		    \sum_{i \in \NNN(I)} \tvarphi^k_{Ii} a^{k+1}_i \norm{x^{k+1}_i-x^*_I}^2
		    + \left[ 1 \wedge V(b^{k+1}, y^{k+1}) \right]
		\right)
	\end{align*}
	where the last line follows from our choice of $\lambda^3 = \frac{1}{\sqrt{\sigma}}$ and $\lambda\tau \asymp 1$.
	Similarly, on the other side,
	\begin{equation*}
		\sum_{i \in \NNN(I)} (\tvarphi^{k+1}_{Ii} -\tvarphi^k_{Ii}) a^k_i
		\lesssim \sqrt{\sigma}
		\left(
		    \sum_{i \in \NNN(I)} \tvarphi^{k+1}_{Ii} a^{k+1}_i \norm{x^{k+1}_i-x^*_I}^2
		    + \left[ 1 \wedge V(b^{k+1}, y^{k+1}) \right]
		\right)
	\end{equation*}

	Finally, it remains to bound $\sum_{i \in \NNN(I)} \!\tvarphi^k_{Ii} a^{k+1}_i \!\norm{x^{k+1}_i-x^*_I}^2$
    % as well as 
    and
    $\sum_{i \in \NNN(I)} \!\tvarphi^{k+1}_{Ii} a^{k+1}_i \!\norm{x^{k+1}_i-x^*_I}^2$
	in terms of
	$\sum_i \varphi^{k+1}_{Ii} a^{k+1}_i \norm{x^{k+1}_i-x^*_I}^2 = O\left( \sigma V_\pos(a^{k+1}, x^{k+1}) \right)$.
	This is done in the following \autoref{claim:gencase:bound_tvarphit_tvarphit+1_Vpos}.

    By putting everything together, we obtain that
	\begin{equation*}
		\abs{ \sum_i (\varphi^{k+1}_{Ii} - \varphi^k_{Ii}) a^k_i }
		\lesssim 
        \eps \left( \ola^k_0 + \ola^{k+1}_0 \right) 
        + \sqrt{\sigma} \left[ 1 \wedge V(z^{k+1}) \right],
	\end{equation*}
    which is the first inequality of the lemma.
    The second inequality of the lemma follows by noting that
    $\eps \ola^{k+1}_0 = O\left( \sqrt{\sigma} V(z^{k+1}) \right)$,
    since $\eps = e^{-1/(3 \sqrt{\sigma})} = O(\sqrt{\sigma})$.
\end{proof}

\begin{claim} \label{claim:gencase:bound_tvarphit_tvarphit+1_Vpos}
	For small enough $\eta$ and $\sigma$,
	for any $I \in [n^*]$,
	\begin{align*}
		\sum_{i \in \NNN(I)} \tvarphi^k_{Ii} a^{k+1}_i \norm{x^{k+1}_i-x^*_I}^2
		&\lesssim
		\sum_{i \in \NNN(I)} \tvarphi^{k+1}_{Ii} a^{k+1}_i \norm{x^{k+1}_i-x^*_I}^2 \\
		&\lesssim
		\eps \ola^{k+1}_0 + \sum_{i \in \NNN(I)} \varphi^{k+1}_{Ii} a^{k+1}_i \norm{x^{k+1}_i-x^*_I}^2.
	\end{align*}
\end{claim}

\begin{proof}
    By the same reasoning as in the proof of \autoref{claim:proof_gencase:delayed_proofs:bound_varphi_tvarphi}, one can show that
    \begin{equation*}
		\sum_{i \in \NNN(I)} \left( \tvarphi^{k+1}_{Ii} - \varphi^{k+1}_{Ii} \right) a^{k+1}_i \norm{x^{k+1}_i-x^*_I}^2
		\leq
		\eps \ola^{k+1}_0 \diamXY^2.
	\end{equation*}
    Hence the second inequality.
    
	For the first inequality: As we saw in 
    \eqref{eq:proof_gencase:delayed_proofs:bound_psialpha_t+1_t},
	\begin{align*}
		&\sum_{i \in \NNN(I)} \left( \tvarphi^k_{Ii}-\tvarphi^{k+1}_{Ii} \right) a^{k+1}_i \norm{x^{k+1}_i-x^*_I}^2 \\
		&\qquad \leq \frac{1}{3 \tau^3} \sum_{i \in \NNN(I)} \tvarphi^k_{Ii} a^{k+1}_i \norm{x^{k+1}_i-x^*_I}^2 
		\left( \norm{x^{k+1}_i-x^*_I}^3 - \norm{x^k_i-x^*_I}^3 \right),
	\end{align*}
	and by \autoref{claim:proof_gencase:delayed_proofs:taylor_expan_norm_pow3} and \autoref{lm:proof_gencase:bound_movement},
    $\norm{x^{k+1}_i-x^*_I}^3 - \norm{x^k_i-x^*_I}^3 \lesssim \norm{x^{k+1}_i-x^k_i} = O(\sigma)$.
	So, by our choice of $\lambda^3 = \frac{1}{\sqrt{\sigma}}$ and $\lambda\tau \asymp 1$,
	\begin{equation*}
		\sum_{i \in \NNN(I)} \left( \tvarphi^k_{Ii}-\tvarphi^{k+1}_{Ii} \right) a^{k+1}_i \norm{x^{k+1}_i-x^*_I}^2 
		\lesssim
        \underbrace{
		      \sigma \frac{\lambda^3}{(\lambda\tau)^3} 
        }_{\asymp \sqrt{\sigma}}
        \cdot
		\sum_{i \in \NNN(I)} \tvarphi^k_{Ii} a^{k+1}_i \norm{x^{k+1}_i-x^*_I}^2.
	\end{equation*}
	Thus,
	\begin{align*}
		(1 - O(\sqrt{\sigma})) \sum_{i \in \NNN(I)} \tvarphi^k_{Ii} a^{k+1}_i \norm{x^{k+1}_i-x^*_I}^2
		&\leq \sum_{i \in \NNN(I)} \tvarphi^{k+1}_{Ii} a^{k+1}_i \norm{x^{k+1}_i-x^*_I}^2 \\
		\sum_{i \in \NNN(I)} \tvarphi^k_{Ii} a^{k+1}_i \norm{x^{k+1}_i-x^*_I}^2
		&\leq (1 + O(\sqrt{\sigma}))
		\sum_{i \in \NNN(I)} \tvarphi^{k+1}_{Ii} a^{k+1}_i \norm{x^{k+1}_i-x^*_I}^2.
    \rqedhere
	\end{align*}
\end{proof}

\subsubsection{Proof of \autoref{lm:proof_gencase:errterm1__control_KLdiv_uI_vI_errterm}
(bound on \texorpdfstring{$\mathrm{(err1)}$}{(err1)}) }

\begin{claim}
	For any $I \in [n^*]$,
	\begin{equation*}
	    \sum_i \frac{\varphi^{k+1}_{Ii} a^{k+1}_i}{\ola^{k+1}_I} \log \frac{a^{k+1}_i / \ola^{k+1}_I}{a^k_i / \ola^k_I} 
		= 
		O\left(
			\sum_i (\varphi^{k+1}_{Ii} - \varphi^k_{Ii}) a^k_i
		\right)
		+ \eta^2 O \left( V(z^{k+1}) \right).
	\end{equation*}
\end{claim}

\autoref{lm:proof_gencase:errterm1__control_KLdiv_uI_vI_errterm} follows straightforwardly from the claim and from \autoref{lm:proof_gencase:control_varphi_t+1_t}.

\begin{proof}[Proof of the claim]
    Fix $I \in [n^*]$.
	Let $u^{k+1,I}_i = \frac{\varphi^{k+1}_{Ii} a^{k+1}_i}{\ola^{k+1}_I}$,
	and we want to bound
	\begin{equation*}
%		\KLdiv(u^{k+1,I}, u^{k,I}) + \sum_i u^{k+1,I} \log \frac{\varphi^k_{Ii}}{\varphi^{k+1}_{Ii}}
%		&= \sum_i u^{k+1,I}_i \log \frac{u^{k+1,I}_i}{u^{k,I}_i} \frac{\varphi^k_{Ii}}{\varphi^{k+1}_{Ii}} \\
%		&= 
		\sum_i u^{k+1,I}_i \log \frac{a^{k+1}_i/\ola^{k+1}_I}{a^k_i/\ola^k_I}
		= \sum_i u^{k+1,I}_i \log \frac{Z a^{k+1}_i}{a^k_i} 
		~ + ~ \log \frac{\ola^k_I}{Z \ola^{k+1}_I}.
	\end{equation*}
	By \eqref{eq:proof_gencase:ppa_upd_eqs}, we have
	$
		\log \frac{a^{k+1}_i}{a^k_i} = -\eta [(\gmat^{k+1} b^{k+1})_i - \gval] - \log Z
	$ where $
		Z = \sum_{i'} a^k_{i'} ~ e^{-\eta [(\gmat^{k+1} b^{k+1})_{i'} - \gval]},
	$
	so
	\begin{equation*}
		\sum_i u^{k+1,I}_i \log \frac{Z a^{k+1}_i}{a^k_i} 
		= \sum_i u^{k+1,I}_i (-\eta) [(\gmat^{k+1} b^{k+1})_i - \gval]
	\end{equation*}
	and
	\begin{align*}
		\log \frac{\ola^k_I}{Z \ola^{k+1}_I}
		&= \log \frac{
			\sum_i \varphi^k_{Ii} a^k_i
		}{
			Z \sum_i \varphi^{k+1}_{Ii} a^{k+1}_i
		} \\
		&= \log \frac{
			\sum_i \varphi^k_{Ii} a^k_i
		}{
			\sum_i \varphi^{k+1}_{Ii} a^k_i
		} 
		+ \log \frac{
			\sum_i \varphi^{k+1}_{Ii} a^k_i
		}{
			Z \sum_i \varphi^{k+1}_{Ii} a^{k+1}_i
		} \\
		&= \log \frac{
			\sum_i \varphi^k_{Ii} a^k_i
		}{
			\sum_i \varphi^{k+1}_{Ii} a^k_i
		} 
		+ \log \frac{
			\sum_i \varphi^{k+1}_{Ii} a^{k+1}_i e^{\eta [(\gmat^{k+1} b^{k+1})_i - \gval]}
		}{
			\sum_i \varphi^{k+1}_{Ii} a^{k+1}_i
		} \\
		&= \log \frac{
			\sum_i \varphi^k_{Ii} a^k_i
		}{
			\sum_i \varphi^{k+1}_{Ii} a^k_i
		} 
		+ \log \sum_i u^{k+1,I}_i e^{\eta [(\gmat^{k+1} b^{k+1})_i - \gval]}.
	\end{align*}

	Now by Jensen inequality on concavity of $\log$, since $\sum_i u^{k+1,I}_i = 1$,
	\begin{equation*}
		\log \left[ \sum_i u^{k+1,I}_i e^{\eta [(\gmat^{k+1} b^{k+1})_i - \gval]} \right]
		+ \sum_i u^{k+1,I}_i (-\eta) [(\gmat^{k+1} b^{k+1})_i - \gval] \geq 0.
	\end{equation*}
	Furthermore, since $\log x \leq x-1$
	and $e^x = 1 + x + O(x^2)$
	and using \eqref{eq:proof_gencase:delayed_proofs:expan_Mb_i},
%	\footnote{
%		The fact that the terms linear in $\eta \left[ (\gmat^{k+1} b^{k+1})_i - \gval \right]$ cancel each other out like this seems miraculous, to the point that I wonder if I didn't miss something trivial...
%		
%		If they didn't cancel each other out, it would still be ok, as we can also give a satisfactory (I think) bound on 
%		$\eta \abs{ \sum_i \frac{\varphi_{Ii} a^{k+1}_i}{\ola^{k+1}_I} \left[ (\gmat^{k+1} b^{k+1})_i - \gval \right] }$
%		by going one order higher in the Taylor expansion of $\gmat^{k+1}_{ij} - \gmatss_{IJ}$.
%	}
	\begin{align*}
		& \log \left[ \sum_i u^{k+1,I}_i e^{\eta [(\gmat^{k+1} b^{k+1})_i - \gval]} \right]
		+ \sum_i u^{k+1,I}_i (-\eta) [(\gmat^{k+1} b^{k+1})_i - \gval] \\
		&\leq \sum_i u^{k+1,I}_i \left( e^{\eta [(\gmat^{k+1} b^{k+1})_i - \gval]} - 1 - \eta [(\gmat^{k+1} b^{k+1})_i - \gval] \right) \\
		&= \sum_i u^{k+1,I}_i O \left( \eta^2 [(\gmat^{k+1} b^{k+1})_i - \gval]^2 \right) \\
		&= \eta^2 \sum_i \frac{\varphi^{k+1}_{Ii} a^{k+1}_i}{\ola^{k+1}_I} 
		O \left( \norm{x^{k+1}_i-x^*_I}^2 + V(b^{k+1}, y^{k+1}) \right) \\
		&\lesssim \eta^2 O \left( V(z^{k+1}) \right).
	\end{align*}
	Thus,
	\begin{equation*}
		\sum_i u^{k+1,I}_i \log \frac{a^{k+1}_i / \ola^{k+1}_I}{a^k_i / \ola^k_I}
		= \log \frac{
			\sum_i \varphi^k_{Ii} a^k_i
		}{
			\sum_i \varphi^{k+1}_{Ii} a^k_i
		} + \eta^2 O \left( V(z^{k+1}) \right).
	\end{equation*}
	
	To upper- and lower-bound the first term, note that by \autoref{lm:proof_gencase:ulb_on_weights_locally_k+1} 
% 	and \autoref{lm:proof_gencase:bound_movement}
	\begin{equation*}
		\sum_i \varphi^k_{Ii} a^k_i = \ola^k_I \geq \wulb
		~~~~\text{and}~~~~ 
		\sum_i \varphi^{k+1}_{Ii} a^k_i \geq \sum_i \varphi^{k+1}_{Ii} a^{k+1}_i (1 - O(\eta)) = \ola^{k+1}_I (1-O(\eta)) \geq \wulb (1 - O(\eta)).
	\end{equation*}
	So just by bounding the derivative of $\log$ we have
	\begin{equation*}
		\abs{\log \sum_i \varphi^k_{Ii} a^k_i - \log \sum_i \varphi^{k+1}_{Ii} a^k_i}
		\leq \underbrace{
			\frac{1}{\wulb (1-O(\eta))}
		}_{= O(1)}
		\abs{ \sum_i (\varphi^{k+1}_{Ii} - \varphi^k_{Ii}) a^k_i }.
    \qedhere
	\end{equation*}
\end{proof}

\subsubsection{Proof of \autoref{lm:proof_gencase:errterm2__control_tricky_errterms}
(bound on \texorpdfstring{$\mathrm{(err2)}$}{(err2)}) }

\begin{proof}
	Focus on the $a$ terms.
	The quantity we want to upper-bound is
	\begin{align*}
		&~~ \sum_i (\varphi^{k+1}_{0i} - \varphi^k_{0i}) a^k_i
		+ \frac{\eta}{2 \sigma} \sum_I \sum_i a^k_i \left( \varphi^{k+1}_{Ii}-\varphi^k_{Ii} \right) \norm{x^*_I-x^k_i}^2 \\
		&= \sum_I \sum_i \left( \varphi^k_{Ii} - \varphi^{k+1}_{Ii} \right) a^k_i \left[ 1 - \frac{\eta}{2 \sigma} \norm{x^*_I-x^k_i}^2 \right].
	\end{align*}
	Fix $I \in [n^*]$.
	Note that the sum is only over indices $i \in \NNN(I)$ (otherwise $\varphi^{k+1}_{Ii} = \varphi^k_{Ii} = 0$) and that we have for all such $i$
	\begin{align*}
		\norm{x^*_I-x^k_i}^2 
		&\leq (\lambda \tau + O(\sigma))^2
		\leq 2 (\lambda \tau)^2 + O(\sigma^2)
		\leq \frac{\sigma}{\eta} + O(\sigma^2) \\
		\frac{\eta}{2 \sigma} \norm{x^*_I-x^k_i}^2
		&\leq \frac{1}{2} + O(\eta \sigma)
	\end{align*}
	due to our choice of $\lambda, \tau$ that ensures that
	$(\lambda \tau)^2 \leq \frac{1}{2} \frac{\sigma}{\eta}$.
	So for small enough $\eta, \sigma$, we have for all $i \in \NNN(I)$ that
	${0 \leq 1 - \frac{\eta}{2 \sigma} \norm{x^*_I-x^k_i}^2 \leq 1}$.
	
	Following a similar reasoning as for \autoref{claim:proof_gencase:delayed_proofs:bound_varphi_tvarphi},
    but taking into account that we know the sign of the objects involved and that we are only interested in an upper bound,
	one can check that
	\begin{equation*}
		\sum_i \! \left( \varphi^k_{Ii} - \varphi^{k+1}_{Ii} \right) \!a^k_i \!\left[ 1 \!-\! \frac{\eta}{\sigma} \norm{x^*_I-x^k_i}^2 \right] 
		\leq \!\! \sum_{i \in \NNN(I)} \!\! \left( \tvarphi^k_{Ii} - \tvarphi^{k+1}_{Ii} \right) \!a^k_i \!\left[ 1 \!-\! \frac{\eta}{\sigma} \norm{x^*_I-x^k_i}^2 \right] 
        \,+\, \eps (1+O(\eta)) \ola^{k+1}_0
	\end{equation*}
    --- note that $\ola^k_0$ does not appear on the right-hand side.
	It remains to bound the first term.
	
	As we already saw in the proof of \autoref{lm:proof_gencase:control_varphi_t+1_t} (Eq.~\eqref{eq:proof_gencase:delayed_proofs:bound_psialpha_t+1_t}), we have
	\begin{equation*}
		\tvarphi^k_{Ii} - \tvarphi^{k+1}_{Ii}
		= \tvarphi^k_{Ii} \cdot \left(\! 1 - \frac{\tvarphi^{k+1}_{Ii}}{\tvarphi^k_{Ii}} \right) 
		\leq \tvarphi^k_{Ii} \cdot \left( \log \tvarphi^k_{Ii} - \log \tvarphi^{k+1}_{Ii} \right)
		= \tvarphi^k_{Ii} \cdot
		\frac{1}{3 \tau^3}
		\left(
			\norm{x^{k+1}_i-x^*_I}^3 - \norm{x^k_i-x^*_I}^3
		\right)
	\end{equation*}
	and by \autoref{claim:proof_gencase:delayed_proofs:taylor_expan_norm_pow3}, we have the Taylor expansion 
	\begin{align*}
		\norm{x^{k+1}_i-x^*_I}^3 - \norm{x^k_i-x^*_I}^3
		&= 3 \innerprod{x^{k+1}_i-x^k_i}{x^{k+1}_i-x^*_I} \norm{x^{k+1}_i-x^*_I} \\
		&~~~~
		+ O \left( \norm{x^{k+1}_i-x^k_i}^2 \norm{x^{k+1}_i-x^*_I} + \norm{x^{k+1}_i-x^k_i}^3 \right).
	\end{align*}
    Now by the update equation \eqref{eq:proof_gencase:ppa_upd_eqs} and the expansion \eqref{eq:proof_gencase:delayed_proofs:expan_partial_x_Mb_i_pos}, for any $i \in \NNN(I)$,
	\begin{align*}
		\innerprod{x^{k+1}_i-x^k_i}{x^{k+1}_i-x^*_I}
		&= \innerprod{-\sigma \frac{a^{k+1}_i}{a^k_i} \partial_x \gmat^{k+1}_{i \bullet} b^{k+1}}{x^{k+1}_i-x^*_I} \\
        % &= -\sigma \frac{a^{k+1}_i}{a^k_i} (x^{k+1}_i-x^*_I)^\top \partial_x \gmat^{k+1}_{i \bullet} b^{k+1} \\
		&\leq -\sigma \frac{a^{k+1}_i}{a^k_i} \frac{\sigma_{\min}}{2} \norm{x^{k+1}_i-x^*_I}^2
		+ \sigma \frac{a^{k+1}_i}{a^k_i} O\left( \norm{x^{k+1}_i-x^*_I} \sqrt{V(b^{k+1}, y^{k+1})} \right) \\
		&\leq \sigma (1+O(\eta)) O\left( \norm{x^{k+1}_i-x^*_I} \sqrt{V(b^{k+1}, y^{k+1})} \right)
	\end{align*}
	and so the terms arising from the order-1 terms in the Taylor expansion of $\norm{x^{k+1}_i-x^*_I}^3 - \norm{x^k_i-x^*_I}^3$ are upper-bounded by
	\begin{align}
		&
		\sum_{i \in \NNN(I)} \tvarphi^k_{Ii} a^k_i \left[ 1 - \frac{\eta}{\sigma} \norm{x^*_I-x^k_i}^2 \right] \cdot
		\frac{1}{\tau^3} \innerprod{x^{k+1}_i-x^k_i}{x^{k+1}_i-x^*_I} \norm{x^{k+1}_i-x^*_I} \\
		&\lesssim 
		\frac{\sigma}{\tau^3} \sum_{i \in \NNN(I)} \tvarphi^k_{Ii} a^{k+1}_i \left[ 1 - \frac{\eta}{\sigma} \norm{x^*_I-x^k_i}^2 \right] \cdot
		\norm{x^{k+1}_i-x^*_I}^2 \cdot
		\sqrt{V(b^{k+1}, y^{k+1})} \\
		&\lesssim 
		\sigma \frac{\lambda^3}{(\lambda\tau)^3}
		\left( \sum_{i \in \NNN(I)} \tvarphi^k_{Ii} a^{k+1}_i \norm{x^{k+1}_i-x^*_I}^2 \right)
		\sqrt{V(b^{k+1}, y^{k+1})}.
	\label{eq:proof_gencase:delayed_proofs:prooferr2_order1terms}
	\end{align}
	Here in the last line we just bounded $\abs{ 1 - \frac{\eta}{\sigma} \norm{x^*_I-x^k_i}^2 }$ by $1$.
	Further, by the update equation \eqref{eq:proof_gencase:ppa_upd_eqs} and the expansion \eqref{eq:proof_gencase:delayed_proofs:expan_partial_x_Mb_i},
	\begin{equation*}
	    \norm{x^{k+1}_i-x^k_i}
	    = \sigma \frac{a^{k+1}_i}{a^k_i} \norm{\partial_x \gmat^{k+1}_{i \bullet} b^{k+1}}\lesssim \sigma (1+O(\eta)) 
	    \left( 
    		\norm{x^{k+1}_i-x^*_I}
    		+ \sqrt{V(b^{k+1}, y^{k+1})}
		\right)
	\end{equation*}
	so the terms arising from higher-order terms ($O(\cdot)$ terms) in the Taylor expansion of $\norm{x^{k+1}_i-x^*_I}^3 - \norm{x^k_i-x^*_I}^3$
	are upper-bounded by
	\begin{align}
	    &
	    \sum_{i \in \NNN(I)} \tvarphi^k_{Ii} a^k_i \left[ 1 - \frac{\eta}{\sigma} \norm{x^*_I-x^k_i}^2 \right] \cdot \frac{1}{\tau^3} O \left( \norm{x^{k+1}_i-x^k_i}^2 \norm{x^{k+1}_i-x^*_I} + \norm{x^{k+1}_i-x^k_i}^3 \right) 
        \nonumber \\
		&\lesssim \frac{1}{\tau^3} \sum_{i \in \NNN(I)} \tvarphi^k_{Ii} a^{k+1}_i 
        \!\cdot\!
		\left( 
			\sigma^2 \norm{x^{k+1}_i-x^*_I}^3
			+ \sigma^2 \!\cdot\! V(b^{k+1}, y^{k+1}) \!\cdot\! \norm{x^{k+1}_i-x^*_I}
			+ \sigma^3 \left[ V(b^{k+1}, y^{k+1}) \right]^{3/2}
		\right)
        \nonumber \\
		&\lesssim \frac{1}{\tau^3} \sum_{i \in \NNN(I)} \tvarphi^k_{Ii} a^{k+1}_i 
        \cdot 
		\left( 
			\sigma^2 \norm{x^{k+1}_i-x^*_I}^3
			+ \sigma^2 \left[ V(b^{k+1}, y^{k+1}) \right]^{3/2}
		\right)
		\nonumber \\
		&\lesssim 
		\sigma^2 \frac{\lambda^3}{(\lambda\tau)^2}
		\left( \sum_{i \in \NNN(I)} \tvarphi^k_{Ii} a^{k+1}_i \norm{x^{k+1}_i-x^*_I}^2 \right)
		+ \sigma^2 \frac{\lambda^3}{(\lambda \tau)^3} 
		\left[ V(b^{k+1}, y^{k+1}) \right]^{3/2}
	\label{eq:proof_gencase:delayed_proofs:prooferr2_order2terms}
	\end{align}
	where in the second line we just bounded $\abs{ 1 - \frac{\eta}{\sigma} \norm{x^*_I-x^k_i}^2 }$ by $1$ again,
	the third line follows from Young's inequality,
	and the last line uses that $\norm{x^{k+1}_i-x^*_I} = O(\lambda\tau + \sigma) = O(\lambda\tau)$ for $i \in \NNN(I)$.

	Putting everything together, by summing \eqref{eq:proof_gencase:delayed_proofs:prooferr2_order1terms} and \eqref{eq:proof_gencase:delayed_proofs:prooferr2_order2terms}
	and by using \autoref{claim:gencase:bound_tvarphit_tvarphit+1_Vpos} to bound
	$\sum_{i \in \NNN(I)} \tvarphi^k_{Ii} a^{k+1}_i \norm{x^{k+1}_i-x^*_I}^2$,
	we get
	\begin{align*}
		\MoveEqLeft
		\sum_I \sum_i \left( \varphi^k_{Ii} - \varphi^{k+1}_{Ii} \right) a^k_i \left[ 1 - \frac{\eta}{2 \sigma} \norm{x^*_I-x^k_i}^2 \right] \\
		&\lesssim
		\eps (\ola^{k+1}_0 + \ola^k_0)
		+ \sigma^2 \frac{\lambda^3}{(\lambda \tau)^3}
		\left[ V(b^{k+1}, y^{k+1}) \right]^{3/2} \\
		&\qquad +
		\sigma \frac{\lambda^3}{(\lambda\tau)^3}
		\left( \sum_I \sum_{i \in \NNN(I)} \tvarphi^k_{Ii} a^{k+1}_i \norm{x^{k+1}_i-x^*_I}^2 \right)
		\left(
		    \sqrt{V(b^{k+1}, y^{k+1})} 
		    + \sigma \lambda\tau
	    \right)
		\\
		&\lesssim 
		\eps (\ola^{k+1}_0 + \ola^k_0)
		+ \sigma^2 \frac{\lambda^3}{(\lambda \tau)^3}
		\left[ V(b^{k+1}, y^{k+1}) \right]^{3/2} \\
		&\qquad +
		\sigma \frac{\lambda^3}{(\lambda\tau)^3}
		\left( \eps \ola^{k+1}_0 + V_\pos(a^{k+1}, x^{k+1}) \right)
		\left(
		    \sqrt{V(b^{k+1}, y^{k+1})} 
		    + \sigma \lambda\tau
	    \right).
	\end{align*}
	Finally, we use that $\lambda\tau \asymp 1$ 
	and that $\lambda^3 = \frac{1}{\sqrt{\sigma}}$ 
	and that $\eta \asymp \sigma$ to simplify the bound,
	and we obtain as announced that the above is upper-bounded up to a constant factor by
	\begin{equation*}
        % \sum_I \sum_i \left( \varphi^k_{Ii} - \varphi^{k+1}_{Ii} \right) a^k_i \left[ 1 - \frac{\eta}{2 \sigma} \norm{x^*_I-x^k_i}^2 \right]
        % \lesssim 
		\eps (\ola^{k+1}_0 + \ola^k_0)
		+ \sqrt{\sigma} \cdot \sigma V_\pos(a^{k+1}, x^{k+1})
		+ \sqrt{\sigma} V(z^{k+1})^{3/2}.
    \qedhere
	\end{equation*}
\end{proof}

\subsubsection{Proof of \autoref{lm:proof_gencase:errterm3__control_nonBregmanness_errterm} (bound on \texorpdfstring{$\mathrm{(err3)}$}{(err3)}) }

\begin{proof}
	Focus on the $a$ terms. We want to bound
	$\frac{\eta}{2 \sigma} \sum_I \sum_i (a^{k+1}_i - a^k_i) \varphi^{k+1}_{Ii} \norm{x^*_I-x^{k+1}_i}^2$.
	By the update equation \eqref{eq:proof_gencase:ppa_upd_eqs},
	\begin{align*}
		a^{k+1}_i &= a^k_i ~ e^{-\eta [(\gmat^{k+1} b^{k+1})_i - \gval]} / Z
		~~~\text{where}~~~
		Z = \sum_{i'} a^k_{i'} ~ e^{-\eta [(\gmat^{k+1} b^{k+1})_{i'} - \gval]} \\
		\text{i.e.}~~~~
		a^{k+1}_i - a^k_i &= a^{k+1}_i \left[ 1 - e^{\eta [(\gmat^{k+1} b^{k+1})_i - \gval]} Z \right] \\
		&= a^{k+1}_i \left[ 1 - e^{\eta [(\gmat^{k+1} b^{k+1})_i - \gval]} + e^{\eta [(\gmat^{k+1} b^{k+1})_i - \gval]} (1-Z) \right].
	\end{align*}
	So
	\begin{align*}
		\sum_I \sum_i (a^{k+1}_i - a^k_i) \varphi^{k+1}_{Ii} \norm{x^*_I-x^{k+1}_i}^2
		&= \sum_I \sum_i \varphi^{k+1}_{Ii} a^{k+1}_i \norm{x^*_I-x^{k+1}_i}^2 \left( 1 - e^{\eta [(\gmat^{k+1} b^{k+1})_i - \gval]} \right) \\
		&~ + \sum_I \sum_i \varphi^{k+1}_{Ii} a^{k+1}_i \norm{x^*_I-x^{k+1}_i}^2 e^{\eta [(\gmat^{k+1} b^{k+1})_i - \gval]} ~ (1-Z).
	\end{align*}
	
	For the first term, 
	since $1-e^x \leq -x$, then
	by the expansion \eqref{eq:proof_gencase:delayed_proofs:expan_Mb_i_pos},
	\begin{align*}
		1 - e^{\eta [(\gmat^{k+1} b^{k+1})_i - \gval]}
		\leq -\eta [(\gmat^{k+1} b^{k+1})_i - \gval]
		&\leq -\eta \left[ 
			\frac{1}{4} \sigma_{\min} \norm{x^{k+1}_i-x^*_I}^2 
			+ O\left( \!\sqrt{V(b^{k+1}, y^{k+1})} \right)
		\right] \\
		&\lesssim \eta \sqrt{V(b^{k+1}, y^{k+1})} .
	\end{align*}
	So we get
	\begin{align*}
		\MoveEqLeft \sum_I \sum_i \varphi^{k+1}_{Ii} a^{k+1}_i \norm{x^*_I-x^{k+1}_i}^2 
		\left( 1 - e^{\eta [(\gmat^{k+1} b^{k+1})_i - \gval]} \right) \\
		&\lesssim 
		\sum_I \sum_i \varphi^{k+1}_{Ii} a^{k+1}_i \norm{x^*_I-x^{k+1}_i}^2 \cdot
		\eta \sqrt{V(b^{k+1}, y^{k+1}) } \\
		&\lesssim V_\pos(a^{k+1}, x^{k+1}) \cdot \eta \sqrt{ V(b^{k+1}, y^{k+1}) }.
	\end{align*}
	
	For the second term, write $Z$ as
	\begin{equation*}
		\forall i \in [n],
        a^{k+1}_i = a^k_i e^{-\eta [(\gmat^{k+1} b^{k+1})_i - \gval]} / Z 
        ~~\implies~~
		1/Z = \sum_i a^{k+1}_i e^{\eta [(\gmat^{k+1} b^{k+1})_i - \gval]}.
	\end{equation*}
	So using the expansion \eqref{eq:proof_gencase:delayed_proofs:expan_Mb_i},
	\begin{align*}
		\frac{1}{Z}-1
		&= \sum_i a^{k+1}_i \left( e^{\eta [(\gmat^{k+1} b^{k+1})_i - \gval]} - 1 \right)
        = \sum_i a^{k+1}_i O \left( \eta [(\gmat^{k+1} b^{k+1})_i - \gval] \right) \\
		&= \sum_I \sum_i \varphi^{k+1}_{Ii} a^{k+1}_i O \left( \eta [(\gmat^{k+1} b^{k+1})_i - \gval] \right)
		~+~ \sum_i \varphi^{k+1}_{0i} a^{k+1}_i O \left( \eta [(\gmat^{k+1} b^{k+1})_i - \gval] \right) \\
		&= \eta \sum_I \sum_i \varphi^{k+1}_{Ii} a^{k+1}_i O \left( \norm{x^{k+1}_i-x^*_I}^2 + \sqrt{V(b^{k+1}, y^{k+1}) } \right)
		~+~ O(\eta \ola^{k+1}_0) \\
		&= \eta O \left(
			\sigma V_\pos(a^{k+1}, x^{k+1})
			+ \sqrt{V(b^{k+1}, y^{k+1}) }
			+ \ola^{k+1}_0
		\right)
		= \eta O \left( \sqrt{V(z^{k+1}) } \right).
	\end{align*}
	So, since $e^{\eta [(\gmat^{k+1} b^{k+1})_i - \gval]} \leq e^{2 \eta_0 \smoothnessf_0}$ 
	and $Z = \sum_{i'} a^k_{i'} ~ e^{-\eta [(\gmat^{k+1} b^{k+1})_{i'} - \gval]} \leq e^{2 \eta_0 \smoothnessf_0} = O(1)$,
	we get
	\begin{align*}
		\MoveEqLeft  \sum_I \sum_i \varphi^{k+1}_{Ii} a^{k+1}_i \norm{x^*_I-x^{k+1}_i}^2 
		e^{\eta [(\gmat^{k+1} b^{k+1})_i - \gval]} ~ (1-Z) \\
		&\leq
		\sum_I \sum_i \varphi^{k+1}_{Ii} a^{k+1}_i \norm{x^*_I-x^{k+1}_i}^2 \cdot
		e^{4 \eta_0 \smoothnessf_0} \abs{\frac{1}{Z} - 1} \\
		&\lesssim V_\pos(a^{k+1}, x^{k+1}) \abs{\frac{1}{Z}-1} \\
		&\lesssim V_\pos(a^{k+1}, x^{k+1}) \cdot \eta \sqrt{V(z^{k+1})}.
	\end{align*}
 
	Putting the two bounds together gives the announced inequality.
\end{proof}

% \end{document}
%
\fi
% !TEX root = ../main.tex
% \documentclass[../main]{subfiles}
% \begin{document}

\section{Auxiliary lemmas} \label{sec:aux_lemmas}

\begin{lemma}[KL- vs.\ $\chi^2$-divergence comparison] \label{lm:aux_lemmas:KL_chi2_comparison}
	For any $a, \ha \in \Delta_n$, denoting $\KLdiv(a, \ha)$ the KL-divergence and $\chi^2(a, \ha) = \sum_i \frac{(\ha_i-a_i)^2}{\ha_i}$ the $\chi^2$-divergence,
	\begin{equation*}
		\KLdiv(a, \ha)
		\leq \log \left( 1+\chi^2(a, \ha) \right)
		\leq \chi^2(a, \ha)
    \qquad \text{and} \qquad
		\KLdiv(a, \ha) \geq \left( \max_i \frac{a_i}{\ha_i} \right)^{-1} \chi^2(a, \ha).
	\end{equation*}
% 	\footnote{\fnsurl{https://cstheory.stackexchange.com/a/45523}}
\end{lemma}

\begin{proof}
	For the first inequality, use Jensen's inequality on $\log$, and that $\log(1+x) \leq x$ for all $x$.
	
	For the second inequality, recall that the $f$-divergence of $a$ \wrt\ $\ha$ is defined as
	$
		\sum_{i=1}^n \ha_i f\left( \frac{a_i}{\ha_i} \right)
	$.
	The KL-divergence $\KLdiv(a, \ha)$ is the $f$-divergence for $f_{\KLdiv}(t) = t \log(t)$
	and the $\chi^2$-divergence $\chi^2(a, \ha)$ is the $f$-divergence for $f_{\chi^2}(t) = t^2-t$.
	Note that for any $c > 0$,
	\begin{equation*}
	    % \min_{t \leq c} \frac{f_{\KLdiv}(t)}{f_{\chi^2}(t)} \geq \frac{1}{c}.
	    \forall t \leq c,~~
        f_{\KLdiv}(t) \geq \frac{1}{c} f_{\chi^2}(t).
	\end{equation*}
    The claimed inequality follows by evaluating at $t = \frac{a_i}{\ha_i}$ and taking the sum weighted by $\ha_i$.
\end{proof}

\begin{lemma} \label{lm:aux_lemmas:bregman_divgce}
	Let $D_\Phi$ denote the Bregman divergence associated to an arbitrary differentiable function $\Phi: X \to \RR$, that is, 
	$D_\Phi(x, y) = \Phi(x) - \Phi(y) - \innerprod{\nabla\Phi(y)}{x-y}$.
	\begin{itemize}
		\item $D_\Phi$ is non-negative if and only if $\Phi$ is convex, and $D_\Phi$ is zero if and only if $\Phi$ is linear.
		\item If $\Phi$ is convex, then $D_\Phi(x,y)$ is convex in $x$ (but not in $y$ in general).
% 		\footnote{\fnsurl{https://mathoverflow.net/a/18484}}
		\item $D_\Phi$ is linear in $\Phi$, and $D_{D_\Phi(\cdot, z)}(x, y) = D_\Phi(x, y)$.
		\item We have the three-point identity (or Pythagorean identity)
		\begin{equation} \label{eq:bregman_3pt}
			\forall x, y, z \in X,~
			D_\Phi(x, z) = D_\Phi(x, y) + D_\Phi(y, z) - \innerprod{\nabla\Phi(z) - \nabla\Phi(y)}{x-y}.
		\end{equation}
		In particular, 
		\begin{equation*}
			\forall x, y, z \in X,~
			\nabla_y D_\Phi(y, z)^\top (x-y) = D_\Phi(x, z) - D_\Phi(x, y) - D_\Phi(y, z).
		\end{equation*}
		\item For any fixed $y_0$ and $y_1$, $x \mapsto D_\Phi(x, y_1) - D_\Phi(x, y_0)$ is affine in $x$.
	\end{itemize}
\end{lemma}

The following lemma is just a rewriting of the facts remarked in \autoref{subsec:main_res:pb_stmt} about the structure of the problem. 

\begin{lemma} \label{lm:aux_lemmas:variat_char_MNE}
    Under the Assumptions~\ref{assum:1}-\ref{assum:6}, letting $\gval = F(\mu^*, \nu^*)$,
    % the min-max value,
	\begin{align*}
		\forall x \in \XXX,~ & (F \nu^*)(x) \geq \gval
		&
		\forall y \in \YYY,~ & ((\mu^*)^\top F)(y) \leq \gval
		\\
		\text{and}~~~
		\forall x \in \support(\mu^*),~ & (F \nu^*)(x) = \gval
		&
		\forall y \in \support(\nu^*),~ & ((\mu^*)^\top F)(y) = \gval
	\end{align*}
	and we have the first- and second-order conditions
	\begin{align*}
		\forall x \in \support(\mu^*),~ \partial_x (F \nu^*)(x) &= 0
		&
		\forall y \in \support(\nu^*),~ \partial_y ((\mu^*)^\top F)(y) &= 0
		\\
		\text{and}~~~ \partial^2_{xx} (F \nu^*)(x) &\succ 0
		&
		\partial^2_{yy} ((\mu^*)^\top F)(y) &\prec 0.
	\end{align*}
	
	Using the shorthand notations detailed in \autoref{sec:notations_for_proofs}, this means that
% 	\begin{align*}
% 		\forall I \in [n^*],~ \left( \gmatss b^* \right)_I &= \gval
% 		&
% 		\forall J \in [m^*],~ \left( a^* \gmatss \right)_J &= \gval
% 		\\
% 		\text{and}~~~  \left( \partial_x \gmatss b^* \right)_I &= 0 
% 		&
% 		\left( a^* \partial_y \gmatss \right)_J &= 0 
% 		\\
% 		\text{and}~~~ H_I = \left( \partial^2_{xx} \gmatss b^* \right)_I &\succeq 0
% 		&
% 		-H_J = \left( a^* \partial^2_{yy} \gmatss \right)_J &\preceq 0.
% 	\end{align*}
	\begin{align*}
		\forall I \in [n^*],~ \gmatss_{I \bullet} b^* &= \gval
		&
		\forall J \in [m^*],~ (a^*)^\top \gmatss_{\bullet J} &= \gval
		\\
		\text{and}~~~  \partial_x \gmatss_{I \bullet} b^* &= 0 
		&
		(a^*)^\top  \partial_y \gmatss_{\bullet J} &= 0 
		\\
		\text{and}~~~ \partial^2_{xx} \gmatss_{I \bullet} b^* &\eqqcolon H_I \succ 0
		&
		(a^*)^\top  \partial^2_{yy} \gmatss_{\bullet J} &\eqqcolon -H_J \prec 0.
	\end{align*}
\end{lemma}

\subsection{Useful expressions and a priori bounds for \texorpdfstring{$V_\wei$}{Vwei} and \texorpdfstring{$V_\pos$}{Vpos}} \label{subsec:aux_lemmas:useful_exprs_V}

In many technical proofs, it will be helpful to keep in mind the following decomposition of the Lyapunov function:
% $V(a, x) = V_\wei(a, x) + V_\pos(a, x)$:
\begin{itemize}
    \item By Pythagorean identity, we have the bias-variance decomposition
    \begin{align} 
    	2 V_\pos(a, x)
    	&= \sum_I \ola_I \left( 
    		\norm{x^*_I-\olx_I}^2 
    		+ \trace(\Sigma_I)
    	\right) 
    %	= \sum_I \sum_i a_i \varphi_{Ii} \norm{x^*_I-\olx_I}^2
    %	+ \sum_I \sum_i a_i \varphi_{Ii} \norm{x_i-\olx_I}^2 \\
    	= \sum_I \sum_i a_i \varphi_{Ii} \left( \norm{x^*_I-\olx_I}^2 + \norm{x_i-\olx_I}^2 \right) \\
    	&= \sum_I \sum_i a_i \varphi_{Ii} \left( \norm{x^*_I-x_i}^2 + 2 \innerprod{x^*_I-\olx_I}{x_i-\olx_I} \right) \\
    	&= \sum_I \sum_i a_i \varphi_{Ii} \norm{x^*_I-x_i}^2.
    \label{eq:aux_lemmas:explicit_Vposax}
    \end{align}
    In particular, note that by Jensen's inequality,
    $
        \sum_I \norm{\left( \ola \odot \Delta \olx \right)_I}
        = \sum_I \norm{\sum_i \varphi_{Ii} a_i (x_i-x^*_I)}
        \leq \sum_{I,i} \varphi_{Ii} a_i \norm{x_i-x^*_I} 
        \leq \sqrt{ \sum_{I,i} \varphi_{Ii} a_i \norm{x_i-x^*_I}^2 }
        = \sqrt{ 2 V_\pos(a,x) }
    $.
    \item The stray weights $\ola_0, \olb_0$ play a special role.
    % in the Lyapunov function.
    Indeed since $a^*_0 = 0$, then $d_h(a^*_0, \ola_0) = \ola_0$, so
    \begin{equation} \label{eq:aux_lemmas:decomp_Vwei}
    	V_\wei(a, x)
    	= \KLdiv(a^*, \ola) 
    	= \sum_I d_h(a^*_I, \ola_I) + \ola_0.
    \end{equation}
    Now $d_h(s, s') \geq 0$.
    So $V_\wei$ can be viewed as a sum of two terms, both positive: The first one measures the (unnormalized entropic Bregman) distance between $a^*$ and $(\ola_I)_{I \in [n^*]}$ in $\RR_+^{n^*}$, and the second one accounts for the stray weights $\ola_0$.
    
    In particular,
    $\ola_0 = O( V_\wei(a, x) )$, while for $I \neq 0$ we only have 
    $\abs{\ola_I-a^*_I} \leq \norm{\ola - a^*}_1 = O( \sqrt{V_\wei(a, x)} )$ (by Pinsker's inequality or by $1$-strong convexity of $h$).
\end{itemize}

% \end{document}

\ifextended%
    % !TEX root = ../main.tex
% \documentclass[../main]{subfiles}
% \begin{document}

\section{Calculations} \label{sec:tedious_calcs}

In this section we present the simple but tedious calculations 
that constitute the proofs of
\autoref{claim:growth_conds:tgaplocs_proxyz},
\autoref{lm:growth_conds:starconvex},
and \autoref{claim:rel_Lya_NI:Fnm-Fnm_bilinear}.
They all consist in writing Taylor expansions of $f$ around $(x^*_I, y^*_J)$ or $(\hx_i, y^*_J)$ or $(x^*_I, \hy_j)$,
% and using the equalities of \autoref{lm:aux_lemmas:variat_char_MNE}.
and applying 
the facts collected in \autoref{lm:aux_lemmas:variat_char_MNE} and \autoref{subsec:aux_lemmas:useful_exprs_V},
which we will use without explicit mention throughout this section.

We will also repeatedly use that as a consequence of \eqref{eq:aux_lemmas:explicit_Vposax} and \eqref{eq:aux_lemmas:decomp_Vwei},
\begin{align*}
    \sum_I \norm{\left( \Delta \ola \odot \Delta \olx \right)_I}
    &= \sum_I \abs{\Delta \ola_I} \norm{\Delta \olx_I}
    = \sum_I \frac{1}{\ola_I} \abs{\Delta \ola_I} \norm{\ola_I \Delta \olx_I} \\
    &\leq \left( \min_I \ola_I \right)^{-1} \norm{\Delta \ola}_1 ~ \max_I \norm{\ola_I \Delta \olx_I} \\
    &= O \left(
        \left( \min_I \ola_I \right)^{-1} \sqrt{V_\wei(a, x)} \sqrt{V_\pos(a, x)}
    \right)
    = O \left( \left( \min_I \ola_I \right)^{-1} V_1(a,x) \right).
\end{align*}

% as well as the estimates of \autoref{claim:growth_conds:DeltaZ_equiv_etaveesigmaV}.
% Unfortunately, it does not seem possible to present a unified calculation that implies all three desired results, due to the fact that $z$ and $\hz$ have symmetric roles in $\gap(z, \hz)$ and not in $\tgap(z; \hz)$. (Basically we cannot avoid having an $O(\norm{\Delta Z}^3)$ term in the expansion of $\gap$, so we cannot just use "expansion of $\gap$" + "$\tgap \approx \gap$" to get \autoref{claim:growth_conds:tgaplocs_proxyz}. cf "tedious-calcs-bak.tex".)

\subsection{Proof of \autoref{claim:growth_conds:tgaplocs_proxyz}}
\label{subsec:tedious_calcs:proof__claim--growth_conds--tgaplocs_proxyz}

\begin{proof}[Proof of \autoref{claim:growth_conds:tgaplocs_proxyz}]
    % Let us compute separately the four terms of
    We wish to compute
    $
    	\tgap(z; \hz) 
    	= \innerprod{
    		\begin{pmatrix}
    			\nabla_a \\
    			\nabla_x \\
    			-\nabla_b \\
    			-\nabla_y
    		\end{pmatrix}
    		F_{n,m}(\hz) 
    	}{
    		\begin{pmatrix}
    			\ha - a \\
    			\hx - x \\
    			\hb - b \\
    			\hy - y \\
    		\end{pmatrix}
    	}
    $,
    where for ease of reference we recall that
    $x_i \coloneqq \hx_i + \sum_I \hvarphi_{Ii} (X_I - \hx_i)$
    and
    $a_i \coloneqq \sum_I A_I \frac{\hvarphi_{Ii} \ha_i}{\hola_I}$.
    Let us compute separately the four terms in this expression.
    
    Focus on the $\nabla_a F_{n,m}$ term (and the $-\nabla_b F_{n,m}$ term is dealt with analogously).
    By definition $\innerprod{\delta a}{\nabla_a F_{n,m}(\hz)} = (\delta a)^\top \gmathh \hb$, so
    \begin{equation*}
    	\innerprod{\nabla_a F_{n,m}(\hz)}{\ha-a}
    	= \ha^\top \gmathh \hb - a^\top \gmathh \hb
    \end{equation*}
    and by Taylor expansions of $f$ around $(x^*_I, y^*_J)$,
    \begin{align*}
        & - a^\top \gmathh \hb
    	= - \sum_I A_I \sum_{ij} \frac{\hvarphi_{Ii} \ha_i}{\hola_I} \gmathh_{ij} \hb_j \\
    	&= - \sum_{IJ} \frac{A_I}{\hola_I} \sum_{ij} \hvarphi_{Ii} \psi_{Jj} \cdot \ha_i
    	\left[
        	\gmatss_{IJ}
        	+ (\hx_i-x^*_I)^\top \partial_x \gmatss_{IJ}
        	+ \partial_y \gmatss_{IJ} (\hy_j-y^*_J)
        	+ O(\norm{\hx_i-x^*_I}^2 + \norm{\hy_j-y^*_J}^2)
    	\right] \hb_j
        + O \!\left( \holb_0 \!\right) \\
    	&= - A^\top \gmatss \holb
    	- (A \odot (\holx-x^*))^\top \partial_x \gmatss \holb
    	- A^\top \partial_y \gmatss ((\holy-y^*) \odot \holb) 
        + O \left( 
            \holb_0
            + (\min_I \hola_I)^{-1} V_\pos(\ha, \hx)
            + V_\pos(\hb, \hy) 
        \right) \\
    	&= - A^\top \gmatss \holb
    	- (A \odot \Delta \holx)^\top \partial_x \gmatss \Delta \holb
    	- \Delta A^\top \partial_y \gmatss (\Delta \holy \odot \holb)
    	+ O \left(\left( \min_I \holw_I \right)^{-1} V_1(\hz) \right) \\
    	&= - A^\top \gmatss \holb
    	- \Delta A^\top \partial_y \gmatss (\Delta \holy \odot b^*)
    	+ O \left(\left( \min_I \holw_I \right)^{-1} V_1(\hz) \right).
    \end{align*}
    
    For the $\nabla_x F_{n,m}$ term (and the $-\nabla_y F_{n,m}$ term is dealt with analogously),
    for any $\delta x$ we have by definition $\innerprod{\delta x}{\nabla_x F_{n,m}(\hz)} = (\ha \odot \delta x)^\top \partial_x \gmathh \hb$, so
    \begin{align*}
    	& \innerprod{\nabla_x F_{n,m}(\hz)}{\hx - x}
    	= \sum_i \ha_i (\hx_i - x_i)^\top \partial_x \gmathh_{i \bullet} \hb
        = \sum_I \sum_i \hvarphi_{Ii} \ha_i (\hx_i - X_I)^\top \partial_x \gmathh_{i \bullet} \hb \\
    	&\!= \sum_{IJ} \sum_{ij} \hvarphi_{Ii} \hpsi_{Jj} \cdot \ha_i (\hx_i - X_I)^\top \!
    	\left[
        	\partial_x \gmatss_{IJ}
        	+ (\hx_i - x^*_I)^\top \partial_{xx}^2 \gmatss_{IJ}
        	+ \partial_{xy}^2 \gmatss_{IJ} (\hy_j - y^*_J)
        	+ O \!\left( \norm{\hx_i\!-\!x^*_I}^2 \!+\! \norm{\hy_j\!-\!y^*_J}^2 \right)
    	\right] \hb_j \\
    	&~~~~ + O \left( \holb_0 \right) \\
    	&\!= (\hola \odot (\holx-X))^\top \partial_x \gmatss \holb
    	+ \sum_{I,i} \hvarphi_{Ii} \ha_i (\hx_i-X_I, \hx_i-x^*_I)^\top \partial_{xx}^2 \gmatss_{I \bullet} \holb 
    	+ (\hola \odot (\holx-X))^\top \partial_{xy}^2 \gmatss ((\holy-y^*) \odot \holb) \\
    	&~~~~ 
    	+ O \left( \holb_0 + V_\pos(\hz) \right) \\
    	&\!= (\hola \odot (\holx-X))^\top \partial_x \gmatss \Delta \holb
    	+ \sum_{I,i} \hvarphi_{Ii} \ha_i \innerprod{\hx_i-X_I}{\hx_i-x^*_I}_{H_I}\!
    	+ (\hola \odot (\holx-X))^\top \partial_{xy}^2 \gmatss ((\holy-y^*) \odot \holb)
        ~ + O\left( V_1(\hz) \right) \\
    	&\!= -(\hola \odot \Delta X)^\top \partial_x \gmatss \Delta \holb
    	+ \sum_{I,i} \hvarphi_{Ii} \ha_i \innerprod{\hx_i-X_I}{\hx_i-x^*_I}_{H_I}\!
    	- (\hola \odot \Delta X) \partial_{xy}^2 \gmatss (\Delta \holy \odot \holb) 
        ~ + O\left( V_1(\hz) \right) \\
    	&\!= -(a^* \odot \Delta X)^\top \partial_x \gmatss \Delta \holb
    	+ \sum_{I,i} \hvarphi_{Ii} \ha_i \innerprod{\hx_i-X_I}{\hx_i-x^*_I}_{H_I}\!
    	- (a^* \odot \Delta X)^\top \partial_{xy}^2 \gmatss (\Delta \holy \odot b^*)
    	+ O\left( \!\big(\! \min_I \holw_I \big)^{-1} V_1(\hz) \right).
    \end{align*}
    Further transform the second term as
    \begin{align*}
        \sum_{I,i} \hvarphi_{Ii} \ha_i \innerprod{\hx_i-X_I}{\hx_i-x^*_I}_{H_I}
        &= \sum_{I,i} \hvarphi_{Ii} \ha_i \norm{\hx_i-x^*_I}_{H_I}^2 
        + \sum_{I,i} \hvarphi_{Ii} \ha_i \innerprod{x^*_I-X_I}{\hx_i-x^*_I}_{H_I} \\
        &= O \left( V_\pos(\ha, \hx) \right)
        - \sum_I \hola_I \innerprod{\Delta X_I}{\Delta \holx_I}_{H_I}.
    \end{align*}
    
    Putting everything together, and using that
    \begin{equation*}
        -A^\top \gmatss \holb + \hola^\top \gmatss B
        = -A^\top \gmatss \Delta \holb + \Delta \hola^\top \gmatss B
        = -\Delta A^\top \gmatss \Delta \holb + \Delta \hola^\top \gmatss \Delta B + O \left( \hola_0 + \holb_0 \right),
    \end{equation*}
    we get the estimate
	\begin{align*} 
% \label{eq:gencase:tgaplocs}
		&\tgap(z; \hz)
		= - \Delta A^\top \gmatss \Delta \holb ~+~ \Delta \hola^\top \gmatss \Delta B \\
		&~~~~ + \Delta \hola^\top \partial_y \gmatss (\Delta Y \odot b^*)
		~-~ (a^* \odot \Delta X)^\top \partial_x \gmatss \Delta \holb \\
		&~~~~ - \Delta A^\top \partial_y \gmatss (\Delta \holy \odot b^*)
		~+~ (a^* \odot \Delta \holx)^\top \partial_x \gmatss \Delta B \\
		&~~~~
% 		- \sum_I w^*_I \innerprod{\Delta P_I}{\Delta \holp_I}_{H_I} \\
		- \sum\nolimits_I a^*_I \innerprod{\Delta X_I}{\Delta \holx_I}_{H_I}
		- \sum\nolimits_J b^*_J \innerprod{\Delta Y_I}{\Delta \holy_J}_{H_J} \\
		&~~~~ - (a^* \odot \Delta X)^\top \partial_{xy}^2 \gmatss (\Delta \holy \odot b^*)
		~+~ (a^* \odot \Delta \holx)^\top \partial_{xy}^2 \gmatss (\Delta Y \odot b^*)
        + O \left( (\min_I \holw_I)^{-1} V_1(\hz) \right).
	\end{align*}
    This expression differs from the one in the Claim statement only in the first two lines: it involves $\Delta \hola$ instead of $\tDelta \hola = \Delta \hola + \hola_0 a^*$, and likewise for the $\holb$.
    The difference between the above expression and the one in the Claim statement is thus bounded by $O(\ola_0 + \olb_0) = O(V_1(\hz))$, and can thus be absorbed in the $O \left( (\min_I \holw_I)^{-1} V_1(\hz) \right)$ term.
\end{proof}

\subsection{Proof of \autoref{lm:growth_conds:starconvex}}
\label{subsec:tedious_calcs:proof__lm--growth_conds--starconvex}

\begin{proof}[Proof of \autoref{lm:growth_conds:starconvex}]
    % Let us compute separately the four terms of
    We wish to compute
    $
    	\tgap(\zps; \hz) 
    	= \innerprod{
    		\begin{pmatrix}
    			\nabla_a \\
    			\nabla_x \\
    			-\nabla_b \\
    			-\nabla_y
    		\end{pmatrix}
    		F_{n,m}(\hz) 
    	}{
    		\begin{pmatrix}
    			\ha - \aps \\
    			\hx - \xps \\
    			\hb - \bps \\
    			\hy - \yps \\
    		\end{pmatrix}
    	}
    $,
    where for ease of reference we recall that
    $\xps_i \coloneqq \hx_i + \sum_I \hvarphi_{Ii} (x^*_I - \hx_i)$
    and
    $\aps_i \coloneqq \sum_I a^*_I \frac{\hvarphi_{Ii} \ha_i}{\hola_I}$.
    Let us compute separately the four terms in this expression.

    In the calculations below, we write $\beps, \beps_{Iij} \in [-1,1]$ or $\in B_{0,1}$ to denote quantities possibly dependent on summation indices, and that may change from line to line.
    This is done in order to track error terms with more precision than using $O(\cdot)$'s.
    
    Focus on the $\nabla_a F_{n,m}$ term (and the $-\nabla_b F_{n,m}$ term is dealt with analogously).
    By definition $\innerprod{\delta a}{\nabla_a F_{n,m}(\hz)} = (\delta a)^\top \gmathh \hb$, so
    \begin{equation*}
    	\innerprod{\nabla_a F_{n,m}(\hz)}{\ha-a}
    	= \ha^\top \gmathh \hb - (\aps)^\top \gmathh \hb,
    \end{equation*}
    and by Taylor expansions of $f$ around $(x^*_I, \hy_j)$,
    \begin{align*}
        & - (\aps)^\top \gmathh \hb 
    	= - \sum_I a^*_I \sum_{ij} \frac{\hvarphi_{Ii} \ha_i}{\hola_I} \gmathh_{ij} \hb_j \\
    	&= - \sum_I a^*_I \sum_{ij} \frac{\hvarphi_{Ii} \ha_i}{\hola_I}
    	\left[
        	\gmatsh_{Ij}
        	+ (\hx_i - x^*_I)^\top \partial_x \gmatsh_{Ij}
        	+ \frac{1}{2} ((\hx_i - x^*_I)^2)^\top \partial_{xx}^2 \gmatsh_{Ij}
        	+ L_3 \beps_{Iij} \norm{\hx_i-x^*_I}^3
    	\right] \hb_j \\
    	&= - (a^*)^\top \gmatsh \hb
    	- (a^* \odot \Delta \holx)^\top \partial_x \gmatsh \hb 
    	- \frac{1}{2} \sum_I a^*_I \sum_i \frac{\hvarphi_{Ii} \ha_i}{\hola_I} ((\hx_i-x^*_I)^2)^\top \partial_{xx}^2 \gmatsh_{I \bullet} \hb \\
    	&~~~~ + L_3 \beps \sum_I a^*_I \sum_i \frac{\hvarphi_{Ii} \ha_i}{\hola_I} \norm{\hx_i-x^*_I}^3.
    \end{align*}
    Now for any $I, J, j$,
    $
        \partial_{xx}^2 \gmatsh_{Ij} = \partial_{xx}^2 \gmatss_{IJ} + O \left( \norm{\hy_j-y^*_J} \right)
    $,
    so
    \begin{align*}
        & \sum_I a^*_I \sum_i \frac{\hvarphi_{Ii} \ha_i}{\hola_I} ((\hx_i-x^*_I)^2)^\top \partial_{xx}^2 \gmatsh_{I \bullet} \hb \\
        &= \sum_{IJ} \frac{a^*_I}{\hola_I} \sum_{ij} \hvarphi_{Ii} \psi_{Jj} \cdot \ha_i ((\hx_i-x^*_I)^2)^\top \partial_{xx}^2 \gmatsh_{I \bullet} \hb_j
        + \sum_I \frac{a^*_I}{\hola_I} \sum_{ij} \hvarphi_{Ii} \psi_{0j} \cdot \ha_i ((\hx_i-x^*_I)^2)^\top \partial_{xx}^2 \gmatsh_{I \bullet} \hb_j \\
        &= \sum_{IJ} \frac{a^*_I}{\hola_I} \sum_{ij} \hvarphi_{Ii} \psi_{Jj} \cdot \ha_i ((\hx_i-x^*_I)^2)^\top \left[ \partial_{xx}^2 \gmatss_{IJ} + O \left( \norm{\hy_j-y^*_J} \right) \right] \hb_j
        + O \left( 
            \sum_I \frac{a^*_I}{\hola_I} \sum_i \hvarphi_{Ii} \ha_i \norm{\hx_i-x^*_I}^2 
            \cdot \holb_0 
        \right) \\
        &= \sum_{IJ} \frac{a^*_I}{\hola_I} \sum_i \hvarphi_{Ii} \ha_i ((\hx_i-x^*_I)^2)^\top \partial_{xx}^2 \gmatss_{IJ} \holb_J
        + O \left( 
            \sum_I \frac{a^*_I}{\hola_I} \sum_i \hvarphi_{Ii} \ha_i \norm{\hx_i-x^*_I}^2 
            \left( \sqrt{V_\pos(\hb, \hy)} + \holb_0 \right)
        \right) \\
        &= \sum_I \frac{a^*_I}{\hola_I} \sum_i \hvarphi_{Ii} \ha_i \norm{\hx_i-x^*_I}_{H_I}^2
        + O \left( 
            % \sum_I \frac{a^*_I}{\hola_I} \sum_i \hvarphi_{Ii} \ha_i \norm{\hx_i-x^*_I}^2 
            \left( \min_I \hola_I \right)^{-1} V_\pos(\ha, \hx)
            \left( \sqrt{V_\pos(\hb, \hy)} + \norm{\Delta \holb}_1 \right) 
        \right).
    \end{align*}
    
    For the $\nabla_x F_{n,m}$ term (and the $-\nabla_y F_{n,m}$ term is dealt with analogously),
    for any $\delta x$ we have
    by definition $\innerprod{\delta x}{\nabla_x F_{n,m}(\hz)} = (\ha \odot \delta x)^\top \partial_x \gmathh \hb$, so
    \begin{align*}
        & \innerprod{\nabla_x F_{n,m}(\hz)}{\hx - \xps}
    	= \sum_i \ha_i (\hx_i - \xps_i)^\top \partial_x \gmathh_{i \bullet} \hb
        = \sum_I \sum_i \hvarphi_{Ii} \ha_i (\hx_i - x^*_I)^\top \partial_x \gmathh_{i \bullet} \hb \\
    	&~ = \sum_I \sum_{ij} \hvarphi_{Ii} \ha_i (\hx_i - x^*_I)^\top 
    	\left[
        	\partial_x \gmatsh_{Ij}
        	+ (\hx_i - x^*_I)^\top \partial_{xx}^2 \gmatsh_{Ij} 
        	+ L_3 \beps_{Iij} \norm{\hx_i-x^*_I}^2
    	\right] \hb_j \\
    	&~ = (\hola \odot \Delta \holx)^\top \partial_x \gmatsh \hb
    	+ \sum_I \sum_i \hvarphi_{Ii} \ha_i ((\hx_i-x^*_I)^2)^\top \partial_{xx}^2 \gmatsh_{I \bullet} \hb
    	+ L_3 \beps \sum_I \sum_i \hvarphi_{Ii} \ha_i \norm{\hx_i-x^*_I}^3.
    \end{align*}
    Now by the same calculation as previously,
    \begin{equation*}
    	\sum_I \sum_i \hvarphi_{Ii} \ha_i ((\hx_i-x^*_I)^2)^\top \partial_{xx}^2 \gmatsh_{I \bullet} \hb 
        = \sum_I \sum_i \hvarphi_{Ii} \ha_i \norm{\hx_i-x^*_I}_{H_I}^2
        + O \left( V_\pos(\ha, \hx) \sqrt{V_1(\hb, \hy)} \right).
    \end{equation*}
    
    Putting everything together we get
    \begin{align*}
    	\tgap(z; \hz)
    	&= \ha^\top \gmaths b^* - (a^*)^\top \gmatsh \hb \\
    	&~~ + ((\hola-a^*) \odot \Delta \holx)^\top \partial_x \gmatsh \hb
    	- \ha^\top \partial_y \gmaths (\Delta \holy \odot (\holb-b^*)) \\
    	&~~ + \sum_I \left( 1 - \frac{1}{2} \frac{a^*_I}{\hola_I} \right) \sum_i \hvarphi_{Ii} \ha_i \norm{\hx_i-x^*_I}_{H_I}^2 
    	+ \sum_J \left( 1 - \frac{1}{2} \frac{b^*_J}{\holb_J} \right) \sum_j \hpsi_{Jj} \hb_j \norm{\hy_j-y^*_J}_{H_J}^2 \\
        &~~ + L_3 \beps
        \sum_I \left( 1 + \frac{w^*_I}{\holw_I} \right) \sum_i \hvarphi_{Ii} \hw_i \norm{\hp_i-p^*_I}^3
        + O \left( \left( \min_I \hola_I \right)^{-1} V_1(\hz)^{3/2} \right).
    \end{align*}
    Now,
    \begin{itemize}
        \item The terms on the second line are negligible, as
        \begin{align*}
        	\partial_x \gmatsh_{I \bullet} \hb
        	&= \sum_J \sum_j \hpsi_{Jj} \left[ \partial_x \gmatss_{IJ} + O(\norm{\hy_j-y^*_J}) \right] \hb_j + O \left( \holb_0 \right) \\
        	&= \partial_x \gmatss_{I \bullet} \Delta \holb
        	+ O \left( \sqrt{V_\pos(\hb, \hy) } \right) + O \left( \holb_0 \right) 
            = O \left( \sqrt{V_1(\hb, \hy)} \right) \\
        	\text{so}~~
        	((\hola-a^*) \odot \Delta \holx)^\top \partial_x \gmatsh \hb 
        	&= O \left( 
                \left( \min_I \hola_I \right)^{-1} V_1(\ha, \hx) \cdot \sqrt{V_1(\hb, \hy)}
            \right).
        \end{align*}
        \item Part of the terms on the third line turn out to be negligible, as
        \begin{align*}
            \sum_I \left( \frac{1}{2} - \frac{1}{2} \frac{a^*_I}{\hola_I} \right) \sum_i \hvarphi_{Ii} \ha_i \norm{\hx_i-x^*_I}_{H_I}^2
            &= \frac{1}{2} \sum_I \frac{\hola_I-a^*_I}{\hola_I} \sum_i \hvarphi_{Ii} \ha_i \norm{\hx_i-x^*_I}^2 \\
            &= O \left(
                \left( \min_I \hola_I \right)^{-1} \norm{\hola-a^*}_1 V_\pos(\ha, \hx)
            \right).
        \end{align*}
    	\item On the last line, the term in $L_3 \beps$ is absolutely bounded by
        \begin{equation*}
            L_3 \cdot 1 \cdot \sum_I 2 \left( \min_I \holw_I \right)^{-1} \hvarphi_{Ii} \hw_i \norm{\hp_i-p^*_I}^2 \cdot \lambda\tau
            \leq 2 L_3 \left( \min_I \holw_I \right)^{-1} \cdot V_\pos(\hz) \cdot \lambda\tau.
        \end{equation*}
    \end{itemize}
    So as announced,
    \begin{align*}
    	\tgap(z; \hz) 
    	% &= F_{n,m^*}(\ha, \hx, b^*, y^*) - F_{n^*, m}(a^*, x^*, \hb, \hy) \\
        &= \ha^\top \gmaths b^* - (a^*)^\top \gmatsh \hb
        + \frac{1}{2} \sum_{I,i} \hvarphi_{Ii} \ha_i \norm{\hx_i-x^*_I}_{H_I}^2 
    	+ \frac{1}{2} \sum_{J,j} \sum_j \hpsi_{Jj} \hb_j \norm{\hy_j-y^*_J}_{H_J}^2 \\
    	&~~ + O \left( \left( \min_I \holw_I \right)^{-1} V_1(\hz)^{3/2} \right)
        + \beps \cdot \left[ 
            2 \smoothnessf_3 \left( \min_I \holw_I \right)^{-1} \cdot V_\pos(\hz) \cdot \lambda \tau 
        \right].
    \rqedhere
    \end{align*}
\end{proof}

\subsection{Proof of \autoref{claim:rel_Lya_NI:Fnm-Fnm_bilinear}}
\label{subsec:tedious_calcs:proof__claim--rel_Lya_NI--Fnm-Fnm_bilinear}

\begin{proof}[Proof of \autoref{claim:rel_Lya_NI:Fnm-Fnm_bilinear}]
    Let any $\hz = (\ha, \hx, \hb, \hy) \in \Delta_n \times \XXX^n \times \Delta_m \times \YYY^m$
    and $Z = (A, X, B, Y) \in \Delta_{n^*} \times \XXX^{n^*} \times \Delta_{m^*} \times \YYY^{m^*}$,
    and denote $\hZ = (\hola, \holx, \holb, \holy)$.
    Recall that, as defined in \eqref{eq:growth_conds:def_norm_ZZZ}, for any $\tX \in \XXX^{n^*}$, $\norm{\Delta \tX} = \max_I \norm{\Delta \tX_I}$,
    and for any $\tZ = (\tA, \tX, \tB, \tY)$, $\norm{\Delta \tZ} = \norm{\Delta \tA}_1 + \norm{\Delta \tX} + \norm{\Delta \tB}_1 + \norm{\Delta \tY}$.

    In the calculations below, we write $\beps, \beps_{Ii}, \beps_{IJi} \in [-1,1]$ to denote quantities possibly dependent on summation indices, and that may change from line to line.
    This is done in order to track error terms with more precision than using $O(\cdot)$'s.
    
	By Taylor expansions of $f$ around $(x^*_I, y^*_J)$, we have
	\begin{align*}
		& F_{n,m^*}(\ha, \hx, B, Y) 
        = \sum_{i=1}^n \sum_{J \in [m^*]} \ha_i f(\hx_i, Y_J) B_J \\
		&= \sum_{IJ} \sum_i \hvarphi_{Ii} \ha_i B_J \Bigg[
    		\gmatss_{IJ} + (\hx_i-x^*_I)^\top \partial_x \gmatss_{IJ} + \partial_y \gmatss_{IJ} (Y_J-y^*_J) + (\hx_i-x^*_I)^\top \partial_{xy}^2 \gmatss_{IJ} (Y_J-y^*_J) \\
    		&~~~~ ~~~~ ~~
    		+ \frac{1}{2} ((\hx_i-x^*_I)^2)^\top \partial_{xx}^2 \gmatss_{IJ}
            + L_3 \beps_{IJi} \norm{\hx_i-x^*_I}^3
            + O \left( \norm{Y_J-y^*_J}^2 \right)
		\Bigg] \\
		&~~ + \sum_J \sum_i \hvarphi_{0i} \ha_i B_J \left[
		    \gmaths_{iJ} + O \left( \norm{Y_J-y^*_J} \right)
		\right] \\
% \displaybreak \\  % manual pagebreak https://tex.stackexchange.com/questions/151122/custom-pagebreak-in-align-equation
		&= \hola^\top \gmatss B + (\hola \odot \Delta \holx)^\top \partial_x \gmatss B + \hola^\top \partial_y \gmatss (\Delta Y \odot B) + (\hola \odot \Delta \holx)^\top \partial_{xy}^2 \gmatss (\Delta Y \odot B) \\
        &~~ + \frac{1}{2} \sum_I \sum_i \hvarphi_{Ii} \ha_i \norm{\hx_i \!-\! x^*_I}_{H_I}^2 
        + \sum_{IJ} \sum_i \hvarphi_{Ii} \ha_i \Delta B_J O \left( \norm{\hx_i \!-\! x^*_I}^2 \right) 
        + \sum_I \sum_i \hvarphi_{Ii} \ha_i L_3 \beps_{Ii} \norm{\hx_i-x^*_I}^3 \\
        &~~ + (\hvarphi_0 \odot \ha)^\top \gmaths B 
        + O \left( \hola_0 \norm{\Delta Y} \right) 
        + O \left( \norm{\Delta Y}^2 \right) \\
		&= \hola^\top \gmatss B + (\hvarphi_0 \odot \ha)^\top \gmaths B \\
		&~~ + (\hola \odot \Delta \holx)^\top \partial_x \gmatss B + \hola^\top \partial_y \gmatss (\Delta Y \odot B) + (\hola \odot \Delta \holx)^\top \partial_{xy}^2 \gmatss (\Delta Y \odot B) \\
		&~~ + \frac{1}{2} \sum_I \sum_i \hvarphi_{Ii} \ha_i \left( \norm{\hx_i-x^*_I}_{H_I}^2 + 2 L_3 \beps_{Ii} \norm{\hx_i-x^*_I}^3 \right) \\
		&~~ + O \left( \hola_0 \norm{\Delta Y} \right) 
        + O \left( \norm{\Delta Y}^2 \right)
		+ O \left( \norm{\Delta B}_1 \cdot V_\pos(\ha, \hx) \right).
	\end{align*}
    Here,
    \begin{itemize}
        \item The first line can be rewritten as
    	\begin{align*}
    		& ~~ \hola^\top \gmatss B + (\hvarphi_0 \odot \ha)^\top \gmaths B \\
    		&= (1-\hola_0) \gval + \hola^\top \gmatss \Delta B + (\hvarphi_0 \odot \ha)^\top \gmaths B \\
    		&= (1-\hola_0) \gval + \Delta \hola^\top \gmatss \Delta B + (\hvarphi_0 \odot \ha)^\top \gmaths B \\
    		&= \gval + \Delta \hola^\top \gmatss \Delta B + (\hvarphi_0 \odot \ha)^\top \left[ \gmaths B - \gval \bmone \right] \\
    		&= \gval 
            + \Delta \hola^\top \gmatss \Delta B 
            + (\hvarphi_0 \odot \ha)^\top
            \underbrace{
                \left[ \gmaths b^* - \gval \bmone \right] 
            }_{\geq 0}
            + O \left( \hola_0 \norm{\Delta B}_1 \right).
    	\end{align*}
    	\item On the second line, it is not hard to check that
    	\begin{equation*}
    	    (\hola \odot \Delta \holx)^\top \partial_x \gmatss B
    	    = (\hola \odot \Delta \holx)^\top \partial_x \gmatss \Delta B
    	    = (a^* \odot \Delta \holx)^\top \partial_x \gmatss \Delta B + O \left( \norm{\Delta \hola}_1 \norm{\Delta \holx} \norm{\Delta B}_1 \right)
    	\end{equation*}
    	and likewise for the other terms, and so
    	\begin{align*}
    	    &~~ (\hola \odot \Delta \holx)^\top \partial_x \gmatss B + \hola^\top \partial_y \gmatss (\Delta Y \odot B) + (\hola \odot \Delta \holx)^\top \partial_{xy}^2 \gmatss (\Delta Y \odot B) \\
    	    &= (a^* \odot \Delta \holx)^\top \partial_x \gmatss \Delta B + \Delta \hola^\top \partial_y \gmatss (\Delta Y \odot b^*) + (a^* \odot \Delta \holx)^\top \partial_{xy}^2 \gmatss (\Delta Y \odot b^*) \\
    	    & \hspace{24em} + O \left(
    	        \norm{\Delta \hZ}^3 + \norm{\Delta Z}^3
    	    \right).
    	\end{align*}
    	\item We can lower-bound the third line as
        \begin{equation*}
            \frac{1}{2} \sum_{I,i} \hvarphi_{Ii} \ha_i \left( \norm{\hx_i-x^*_I}_{H_I}^2 + 2 L_3 \beps_{Ii} \norm{\hx_i-x^*_I}^3 \right)
            \geq \frac{1}{2} \sum_{I,i} \hvarphi_{Ii} \ha_i \norm{\hx_i-x^*_I}^2 
            \underbrace{\left( 
                \sigma_{\min} - 2 L_3 \beps_{Ii} \cdot \lambda\tau 
            \right)}_{\geq \frac{\sigma_{\min}}{2}}
        \end{equation*}
    	by our assumption that $\lambda\tau \leq \frac{\sigma_{\min}}{4 \smoothnessf_3}$.
    	Further, we can decompose this lower bound as
    	\begin{equation*}
    	    \sum_{I,i} \hvarphi_{Ii} \ha_i \norm{\hx_i \!-\! x^*_I}^2 
    	    = \sum_I \hola_I \Bigg( 
                \norm{\holx_I \!-\! x^*_I}^2 
                + \underbrace{
                    \trace(\hSigma_I) 
                }_{\geq 0}
            \Bigg)
    	    \geq \sum_I a^*_I \norm{\holx_I \!-\! x^*_I}^2 
            + O \left( \norm{\Delta \hola}_1 \norm{\Delta \holx}^2 \right)
        \end{equation*}
        and transform the remaining quadratic term, using that $2\innerprod{a}{b} \leq 2 \norm{a} \norm{b} \leq \norm{a}^2 + \norm{b}^2$, as
        \begin{equation*}
            \sum_I a^*_I \norm{\holx_I-x^*_I}^2 
            \geq \sum_I a^*_I \left( -2\innerprod{\Delta \holx_I}{\Delta X_I} - \norm{\Delta X_I}^2 \right)
            = -2 \sum_I a^*_I \innerprod{\Delta \holx_I}{\Delta X_I} 
            + O \left( \norm{\Delta X}^2 \right).
    	\end{equation*}
    \end{itemize}
	Thus, combining the above bounds, we have
	\begin{align*}
		F_{n,m^*}(\ha, \hx, B, Y)
		&\geq \gval + \Delta \hola^\top \gmatss \Delta B \\
		&~ + (a^* \odot \Delta \holx)^\top \partial_x \gmatss \Delta B 
        + \Delta \hola^\top \partial_y \gmatss (\Delta Y \odot b^*) 
        + (a^* \odot \Delta \holx)^\top \partial_{xy}^2 \gmatss (\Delta Y \odot b^*) \\
		&~ - \frac{\sigma_{\min}}{2} \sum_I a^*_I \innerprod{\Delta \holx_I}{\Delta X_I} 
        \, + O \left(
	        \norm{\Delta \hZ}^3 \!+ \norm{\Delta Z}^2
	        \!+ \norm{\Delta Z} \left( \hola_0 + V_\pos(\ha, \hx) \right) \!
	    \right).
	\end{align*}
 
    One can derive the analogous upper bound for $F_{n^*, m}(A, X, \hb, \hy)$. Combining the two, we obtain
    \begin{align*}
        & F_{n, m^*}(\ha, \hx, B, Y) - F_{n^*, m}(A, X, \hb, \hy) \\
        &\geq -\begin{pmatrix}
		    \Delta A \\
		    \Delta X \\
		    \Delta B \\
		    \Delta Y
		\end{pmatrix}^\top
		\begin{bmatrix}
            0 & 0 & \gmatss & \partial_y \gmatss b^* \\
            0 & a^* \frac{\sigma_{\min}}{2} \id & a^* \partial_x \gmatss & a^* \partial_{xy}^2 \gmatss b^* \\
            -(\gmatss)^\top & -(a^* \partial_x \gmatss)^\top & 0 & 0 \\
            -(\partial_y \gmatss b^*)^\top & -(a^* \partial_{xy}^2 \gmatss b^*)^\top & 0 & b^* \frac{\sigma_{\min}}{2} \id
		\end{bmatrix}
        \begin{pmatrix}
            \Delta \hola \\
            \Delta \holx \\
            \Delta \holb \\
            \Delta \holy
        \end{pmatrix} \\
		&~~~~ + O \left(
	        \norm{\Delta \hZ}^3
	        + \norm{\Delta Z} \left( \holw_0 + V_\pos(\hz) \right)
	        + \norm{\Delta Z}^2
	    \right) \\
        &= -\begin{pmatrix}
		    \Delta A \\
		    \Delta X \\
		    \Delta B \\
		    \Delta Y
		\end{pmatrix}^\top
		\begin{bmatrix}
            0 & 0 & \gmatss & \partial_y \gmatss b^* \\
            0 & a^* \frac{\sigma_{\min}}{2} \id & a^* \partial_x \gmatss & a^* \partial_{xy}^2 \gmatss b^* \\
            -(\gmatss)^\top & -(a^* \partial_x \gmatss)^\top & 0 & 0 \\
            -(\partial_y \gmatss b^*)^\top & -(a^* \partial_{xy}^2 \gmatss b^*)^\top & 0 & b^* \frac{\sigma_{\min}}{2} \id
		\end{bmatrix}
        \begin{pmatrix}
            \tDelta \hola \\
            \Delta \holx \\
            \tDelta \holb \\
            \Delta \holy
        \end{pmatrix} \\
		&~~~~ + O \left(
	        \norm{\Delta \hZ}^3
	        + \left( \holw_0 + V_\pos(\hz) \right)^2
	        + \norm{\Delta Z}^2
	    \right)
    \end{align*}
    where on the last line we used that $\norm{\Delta \hola - \tDelta \hola}_1 \leq \hola_0$ by definition of $\tDelta \hola = \Delta \hola + \hola_0 a^*$, and the fact that $\forall A, B >0, 2AB \leq A^2 + B^2$.
\end{proof}

\fi
% !TEX root = ../main.tex
% \documentclass[../main]{subfiles}
% \begin{document}

\section{Proof of convergence of CP-MP in the exact-\hspace{0em}parametrization case} \label{sec:mp_exactparam}

In this section we prove \autoref{prop:cv_proof:exactparam:mp_cv}, which states that CP-MP converges under the same conditions and with the same rate as CP-PP. The proof essentially combines the convergence result for CP-PP
with the general fact that the Mirror Prox and Proximal Point updates coincide up to order-3 terms,
% in the step-size, 
a consequence of the two following lemmas.
% 
% We start by deriving order-2 expressions for the mirror prox and proximal point updates.
% Due to the presence of constraints and to the fact that the divergence used is (non-Euclidean and) not even a Bregman divergence in our case, the results we state are quite general, and they may be of independent interest.
% (In the lemma statements and in their proofs we do not use the notations from \autoref{sec:notations_for_proofs}.)

\begin{lemma} \label{lm:mp_exactparam:mp_expan}
    Let $A \in \RR^{m \times d}$ and $b \in \RR^m$ for some $m < d$, with $A$ having full rank, and denote $\ZZZ = \left\lbrace z \in \RR^d; Az = b \right\rbrace$.
    Define the semi-norm
    \begin{equation*}
        \forall v \in \RR^d,~ 
        \norm{v}_{*\ZZZ} 
        = \max_{\substack{\norm{\delta} \leq 1 \\ A \delta = 0}}
        \innerprod{\delta}{v}
    \end{equation*}
    where $\norm{\delta}$ is the usual Euclidean norm.
    Consider some function $F: \RR^d \to \RR$ with Lipschitz-continuous second-order differentials.
    
    Let $D: \ZZZ \times \ZZZ \to \RR_+$ such that for any $z^0$ in some subset $\ZZZ_0 \subset \ZZZ$,
    \begin{itemize}
    	\item $D^0(z) \coloneqq D(z, z^0)$ is strongly convex and smooth over $z \in \ZZZ$, and the constants do not depend on $z^0$.
    	\item There exist $H = H_{ij} \succ 0$ and $K = K_{ijk}$ an order-3 symmetric tensor (that depend on $z^0$) such that
    	$\nabla D^0(z) = H (z-z^0) + \frac{K}{2} (z-z^0)^2 + O \left( \norm{z-z^0}^3 \right)$.
    	That is, using Einstein's summation notation,
    	\begin{equation*}
    		[\nabla D^0(z)]_i = H_{ij} (z-z^0)^j + \frac{1}{2} K_{ijk} (z-z^0)^j (z-z^0)^k + O \left( \norm{z-z^0}^3 \right).
    	\end{equation*}
    \end{itemize}
    To be clear, we assume that $\sigma_{\min}(H)^{-1}$ and the norms of $H$ and $K$, as well as the constant hidden in the $O(\cdot)$ in the above equation, are all bounded by a constant that does not depend on $z^0$.
    
    Consider the Mirror Prox (MP) update $z_{\mathrm{MP}}^{k+1} = \mathrm{MP}(z^k; \eta)$ defined by
    \begin{align*}
    	\hz &= \argmin_{Az=b} \innerprod{\nabla F(z^k)}{z} + \frac{1}{\eta} D(z, z^k) \\
    	z_{\mathrm{MP}}^{k+1} &= \argmin_{Az=b} \innerprod{\nabla F(\hz)}{z} + \frac{1}{\eta} D(z, z^k).
    \end{align*}
    Then, if $z^k \in \ZZZ_0$ (and for $\eta$ small enough),
    \begin{align*}
        z_{\mathrm{MP}}^{k+1} - z^k
        &= -\eta H^{-1} P \nabla F(z^k)
        + \eta^2 H^{-1} P \cdot \left( \nabla^2 F(z^k) H^{-1} P \nabla F(z^k) - \frac{1}{2} K \left[ H^{-1} P \nabla F(z^k) \right]^2 \right) \\
        &~~~ + O \left( \eta^3 \norm{\nabla F(z^k)}_{*\ZZZ}^2 \right)
    \end{align*}
    where $O(\cdot)$ hides only the aforementioned constant and the smoothness constants of $F$, and where
    \begin{equation*}
        P = I - A^\top \left[ A H^{-1} A^\top \right]^{-1} A H^{-1}.
    \end{equation*}
\end{lemma}

\begin{lemma} \label{lm:mp_exactparam:pp_expan}
    Under the same conditions 
    and using the same notations
    as in the previous lemma, the Proximal Point (PP) update
    $z_{\mathrm{PP}}^{k+1} = \mathrm{PP}(z^k; \eta)$ defined by
    \begin{equation*}
    	z_{\mathrm{PP}}^{k+1} = \argmin_{Az=b} F(z) + \frac{1}{\eta} D(z, z^k)
    \end{equation*}
    satisfies
    \begin{align*}
        z_{\mathrm{PP}}^{k+1} - z^k
        &= -\eta H^{-1} P \nabla F(z^k)
        + \eta^2 H^{-1} P \cdot \left( \nabla^2 F(z^k) H^{-1} P \nabla F(z^k) - \frac{1}{2} K \left[ H^{-1} P \nabla F(z^k) \right]^2 \right) \\
        &~~~ + O \left( \eta^3 \norm{\nabla F(z^k)}_{*\ZZZ} \right).
    \end{align*}
    That is,
    $z_{\mathrm{PP}}^{k+1} - z_{\mathrm{MP}}^{k+1} = O \left( \eta^3 \norm{\nabla F(z^k)}_{*\ZZZ} \right)$.
\end{lemma}

\begin{remark} \label{rk:mp_exacparam:norm_P_starZZZ}
    One can check that $P^2 = P$ and that $H^{-1} P$ is symmetric, i.e., $P$ is a projection which is orthogonal for $\innerprod{\cdot}{\cdot}_{H^{-1}}$.
    Furthermore, $P A^\top = 0$, i.e., $P^\top$ projects onto the kernel of $A$, and so the semi-norm 
    $\norm{P v} = \max_{\norm{\delta} \leq 1} \innerprod{\delta}{P v}$
    is dominated by $\norm{v}_{*\ZZZ}$ 
    (the operator norm of $P^\top$ being bounded by a constant).
\end{remark}

\begin{remark}
    In the Euclidean case where $D(\cdot, \cdot) = \frac{1}{2} \norm{\cdot-\cdot}^2$, we recover the formulas from
    % \cite[proof of Corollary~1]{lu_osr-resolution_2021}.
    \cite[Prop.~2]{mokhtari_unified_2019}.
    In the case where the divergence function $D(\cdot, \cdot)$ is a Bregman divergence, we provide a finer (order-2) expansion than \cite[Prop.~1]{amid_reparameterizing_2020} (which was order-1).
    
% 	We also note that in the Bregman divergence case, the $K_{ijk}$ are nothing else than the Christoffel symbols (up to lowering of an index) of the divergence function $D(\cdot, \cdot)$ at the point $z^0$ \cite{amari_information_2010}.
% 	In particular it is likely that these results are contained in some form in the Riemannian geometry literature, but we have not been able to locate a statement usable for our purpose.
\end{remark}

In the next subsection we show how to prove \autoref{prop:cv_proof:exactparam:mp_cv} using (a min-max version of) the two above lemmas, and
\ifextended%
    in the two following subsections we prove \autoref{lm:mp_exactparam:mp_expan} and \autoref{lm:mp_exactparam:pp_expan} respectively.
\else%
    the proofs of \Autoref{lm:mp_exactparam:mp_expan,lm:mp_exactparam:pp_expan} are placed in 
    Sections~H.2, H.3 of the arXiv version of this paper
    as \AdditionalMaterial.
\fi

\subsection{Proof of \autoref{prop:cv_proof:exactparam:mp_cv}}

In this subsection we assume the exact-parametrization setting, i.e.\ $n=n^*, m=m^*$, and we use the notations introduced in \autoref{sec:notations_for_proofs}.

\paragraph{Preliminaries.}
For all $z, \hz \in \Delta_n \times \XXX^n \times \Delta_m \times \YYY^m$, denote
\begin{equation*}
    D\left( (a, x), (\ha, \hx) \right) = \KLdiv(a, \ha) + \frac{\eta}{2 \sigma} \sum_i \ha_i \norm{x_i-\hx_i}^2,
\end{equation*}
similarly for $D\left( (b, y), (\hb, \hy) \right)$, 
and $D(z, \hz) = D\left( (a, x), (\ha, \hx) \right) + D\left( (b, y), (\hb, \hy) \right)$; 
in particular we have $V(\hz) = D(z^*, \hz)$.
Also let
\begin{equation*}
    \norm{z-\hz}^2 = \norm{a-\ha}_1^2 + \max_i \norm{x_i-\hx_i}^2 + \norm{b-\hb}_1^2 + \max_j \norm{y_j-\hy_j}^2
\end{equation*}
and recall from \autoref{claim:growth_conds:DeltaZ_equiv_etaveesigmaV} that divergence and squared norm are equivalent in the sense that, if
$a_i, \ha_i, b_j, \hb_j \geq \wulb = \Theta(1)$ for all $i, j$, then
$
    \norm{z-\hz}^2 \lesssim D(z, \hz) \lesssim \norm{z-\hz}^2
$.

Furthermore, denote
$
    g(z) = \begin{pmatrix}
        \nabla_a \\
        \nabla_x \\
        -\nabla_b \\
        -\nabla_y
    \end{pmatrix} F_{n,m}(z)
$
and define the semi-norm
\begin{equation*}
    \norm{v}_{*\ZZZ}
    = \max_{\substack{ \delta z \in \RR^n \times \XXX^n \times \RR^m \times \YYY^m \\ 
    \norm{\delta z} \leq 1 \\
    \bmone^\top \delta a = \bmone^\top \delta b = 0 }}
    \innerprod{\delta z}{v}.
\end{equation*}
By definition, $z^*$ is a stationary point of the vector flow $g(z)$ under the constraint $z \in \Delta_n \times \XXX^n \times \Delta_m \times \YYY^m$, and $z^*$ belongs to the relative interior of that domain. So by smoothness of $F_{n,m}$,
\begin{equation*}
    \norm{g(z)}_{*\ZZZ} = \norm{g(z) - g(z^*)}_{*\ZZZ} \lesssim \norm{z-z^*}.
\end{equation*}

\paragraph{Comparing the CP-MP and CP-PP updates.}
Starting from $z^k$,
% the CP-MP update $\tz^{k+1} = \text{CP-MP}(z^k; \eta, \sigma)$
the CP-MP update $\tz^{k+1}$
is given by
\begin{align*}
	\hz &= \argmin_{\substack{a \in \Delta_n \\ x \in \XXX^n}} 
	\argmax_{\substack{b \in \Delta_m \\ y \in \YYY^m}}~
	\innerprod{g(z^k)}{z} 
	+ \frac{1}{\eta} \left[ 
	    D\left( (a, x), (a^k, x^k) \right) - D\left( (b, y), (b^k, y^k) \right) 
    \right] \\
	\tz^{k+1} &= \argmin_{\substack{a \in \Delta_n \\ x \in \XXX^n}} 
	\argmax_{\substack{b \in \Delta_m \\ y \in \YYY^m}}~
	\innerprod{g(\hz)}{z} 
	+ \frac{1}{\eta} \left[ 
	    D\left( (a, x), (a^k, x^k) \right) - D\left( (b, y), (b^k, y^k) \right) 
    \right]
\end{align*}
% and the CP-PP update $z^{k+1} = \text{PP-MP}(z^k; \eta, \sigma)$ is equivalent to
and the CP-PP update $z^{k+1}$ by
\begin{equation*}
	z^{k+1} = \argmin_{\substack{a \in \Delta_n \\ x \in \XXX^n}} 
	\argmax_{\substack{b \in \Delta_m \\ y \in \YYY^m}}~
	F_{n,m}(z) 
	+ \frac{1}{\eta} \left[ 
	    D\left( (a, x), (a^k, x^k) \right) - D\left( (b, y), (b^k, y^k) \right) 
    \right].
\end{equation*}
It is not hard to adapt the proofs of \autoref{lm:mp_exactparam:mp_expan} and \autoref{lm:mp_exactparam:pp_expan} to cover min-max updates of these forms, as $\KLdiv(a, \ha) = +\infty$ for $a$ on the relative boundary of $\Delta_n$ so that the constraints reduce to $\bmone^\top a = \bmone^\top b = 1$.
Furthermore, it is not hard to show that $\norm{\Delta \tz^{k+1}}, \norm{\Delta z^{k+1}} \lesssim \norm{\Delta z^k} + \eta$, so in particular by choosing $r_0$ and $\eta$ small enough, we may assume 
$z^k, z^{k+1}, \tz^{k+1} \in \ZZZ_0 
= \left\lbrace z ;~ \min_i a_i, \min_j b_j \geq \wulb \right\rbrace$,
and the assumptions of Lemmas~\ref{lm:mp_exactparam:mp_expan} and \ref{lm:mp_exactparam:pp_expan} on $D(\cdot, z^k)$ are satisfied.
Thus we have
\begin{equation*}
    \norm{\tz^{k+1} - z^{k+1}} 
    \lesssim \eta^3 \norm{g(z^k)}_{*\ZZZ}
    \lesssim \eta^3 \norm{\Delta z^k}.
\end{equation*}

Let us convert the above bound on $\norm{\tz^{k+1} - z^{k+1}}$ into a bound on $\abs{ V(z^{k+1}) - V(\tz^{k+1}) }$.
% , i.e., their divergences to $z^*$.
% rather than squared distances to $z^*$.
Since we can assume $z^{k+1}, \tz^{k+1} \in \ZZZ_0$, by using that $h: s \mapsto s \log s - s + 1$ is $\wulb$-smooth over $[\wulb, 1]$ one easily checks that,
denoting $w = (a, b)$ and $p=(x, y)$,
\begin{align*}
    \KLdiv(w^*, w^{k+1}) - \KLdiv(w^*, \tw^{k+1})
    &= \KLdiv(\tw^{k+1}, w^{k+1}) - \sum_i \left( h'(w^{k+1}_i)-h'(\tw^{k+1}_i) \right) \left( w^*_i-\tw^{k+1}_i\right) \\
    &= O \left( \norm{w^{k+1}-\tw^{k+1}}^2 + \norm{w^{k+1}-\tw^{k+1}} \norm{\Delta \tw^{k+1}} \right)
\end{align*}
by Bregman three-point identity;
and similarly, for each $i$
\begin{equation*}
    \norm{p^*_i-p^{k+1}_i}^2 - \norm{p^*_i-\tp^{k+1}_i}^2
    = O \left( \norm{p^{k+1}_i-\tp^{k+1}_i}^2 + \norm{p^{k+1}_i-\tp^{k+1}_i} \norm{\Delta \tp^{k+1}_i} \right).
\end{equation*}
Now it is not hard to show 
\ifextended%
    --- in fact this is just \eqref{eq:mp_exactparam:MP:firstestim} below --- 
\fi
that $\norm{\tz^{k+1}-z^k} \lesssim \eta \norm{\nabla g(z^k)}$
and so $\norm{\Delta \tz^{k+1}} \lesssim \norm{\Delta z^k}$.
So
\begin{equation*}
    \abs{ V(z^{k+1}) - V(\tz^{k+1}) }
    = \abs{ D(z^*, z^{k+1}) - D(z^*, \tz^{k+1}) }
    \lesssim 
    \eta^3 \norm{\Delta z^k}^2
    \lesssim \eta^3 D(z^*, z^k)
    = \eta^3 V(z^k).
\end{equation*}

\paragraph{Proof conclusion.}
In the proof of \autoref{thm:cv_proof:exactparam:loc_exp_cv} we showed that
\begin{equation*}
	V(z^{k+1}) 
	\leq V(z^k) - (C/2) \eta^2 V(z^{k+1})
\end{equation*}
for some $C$ dependent only on $(f, \XXX, \YYY)$ and $\Gamma_0$, for $\eta, \sigma$ small enough and $r_0$ small enough (depending on $\eta, \sigma$).
So
\begin{equation*}
	V(\tz^{k+1}) 
	\leq V(z^k) - (C/2) \eta^2 V(\tz^{k+1})
	+ O \left( \eta^3 V(z^k) \right),
\end{equation*}
and we can conclude to the local exponential convergence of the sequence of CP-MP iterates in exactly the same way as for \autoref{thm:cv_proof:exactparam:loc_exp_cv}.

% \end{document}

\ifextended%
    % !TEX root = ../main.tex
% \documentclass[../main]{subfiles}
% \begin{document}

\subsection{Proof of \autoref{lm:mp_exactparam:mp_expan}} \label{subsec:mp_exactparam:MP}

To lighten notation and since we focus on a single iteration, instead of ``$z^{k+1}_{\mathrm{MP}} = \mathrm{MP}(z^k; \eta)$'' we will consider $\mathrm{MP}(z^0; \eta) = z^2$ with
\begin{align}
	z^1 = \argmin_{Az=b} \innerprod{\nabla F(z^0)}{z} + \frac{1}{\eta} D(z, z^0) 
	\label{eq:mp_exactparam:MP:def1} \tag{U1}
	\\
	z^2 = \argmin_{Az=b} \innerprod{\nabla F(z^1)}{z} + \frac{1}{\eta} D(z, z^0).
	\label{eq:mp_exactparam:MP:def2} \tag{U2}
\end{align}
The goal is to get an order-2 expansion for $\delta z \coloneqq z^2 - z^0$.

Also to lighten notation, we will write $\norm{\nabla F(z^0)}$ for $\norm{\nabla F(z^0)}_{*\ZZZ}$ (and similarly for $\nabla F(z^1)$).

\paragraph{First estimates.}
By Lagrangian duality, there exist $\lambda^1, \lambda^2 \in \RR^m$ such that
\begin{align}
	\nabla F(z^0) + \frac{1}{\eta} \nabla D^0(z^1) - A^\top \lambda^1 &= 0 
	\label{eq:mp_exactparam:MP:statio1} \tag{S1}
	\\
	\text{and}~~~~ \nabla F(z^1) + \frac{1}{\eta} \nabla D^0(z^2) - A^\top \lambda^2 &= 0.
	\label{eq:mp_exactparam:MP:statio2} \tag{S2}
\end{align}
As a first consequence, since $(z^1-z^0)^\top A^\top = 0$, we get that
\begin{align*}
	(z^1-z^0)^\top \left[ \nabla F(z^0) + \frac{1}{\eta} \nabla D^0(z^1) \right] &= 0 \\
	\frac{\mu}{2} \norm{z^1-z^0}^2 \leq (z^1-z^0)^\top \left( \nabla D^0(z^1) - \nabla D^0(z^0) \right) &\leq \eta \norm{z^1-z^0} \norm{\nabla F(z^0)} \\
	\norm{z^1-z^0} &\lesssim \eta \norm{\nabla F(z^0)}
\end{align*}
and also consequently $\norm{\nabla F(z^1)} \leq \norm{\nabla F(z^0)} + O (\norm{z^1-z^0}) \lesssim \norm{\nabla F(z^0)}$.
Similarly since $(z^2-z^0)^\top A^\top = 0$,
\begin{equation} \label{eq:mp_exactparam:MP:firstestim}
	\norm{z^2-z^0} \lesssim \eta \norm{\nabla F(z^1)} \lesssim \eta \norm{\nabla F(z^0)}.
\end{equation}

\paragraph{An order-1 expansion for the first update \eqref{eq:mp_exactparam:MP:def1}.}
Next we want to get an explicit approximate expression for $\nabla F(z^1)$ only in terms of $z^0$, based on the expansion
\begin{equation*}
	\nabla F(z^1) = \nabla F(z^0) + \nabla^2 F(z^0) (z^1-z^0) + O \left( \norm{z^1-z^0})^2 \right).
\end{equation*}
For this we want to get an explicit approximate expression for $z^1-z^0$.

From \eqref{eq:mp_exactparam:MP:statio1} and an order-1 expansion of $\nabla D^0$, we have
\begin{align*}
	\eta \nabla F(z^0) + H (z^1-z^0) - \eta A^\top \lambda^1 &= O(\norm{z^1-z^0}^2) \\
	z^1-z^0 &= \eta H^{-1} \left( -\nabla F(z^0) + A^\top \lambda^1 \right) + O(\norm{z^1-z^0}^2)
\end{align*}
and this will get us an expression of $z^1-z^0$ of the correct order for this paragraph's purpose.
It remains to identify $A^\top \lambda^1$. An approximate expression of it
% (of sufficient order for this paragraph's purpose) 
can be obtained simply by
\begin{align*}
	\frac{1}{\eta} A (z^1-z^0) = 0 &= A H^{-1} \left( -\nabla F(z^0) + A^\top \lambda^1 \right) + \frac{1}{\eta} O(\norm{z^1-z^0}^2) \\
	A H^{-1} A^\top \lambda^1 &= A H^{-1} \nabla F(z^0) + \frac{1}{\eta} O(\norm{z^1-z^0}^2) \\
	\lambda^1 &= \left[ A H^{-1} A^\top \right]^{-1} A H^{-1} \nabla F(z^0) + \frac{1}{\eta} O(\norm{z^1-z^0}^2)
\end{align*}
since $A H^{-1} A^\top$ is invertible as an $m \times m$ product of full-rank matrices.
Thus we get
\begin{equation}
\label{eq:mp_exactparam:MP:def_P} \tag{P}
	z^1 - z^0 = - \eta H^{-1} 
	\underbrace{ 
		\left( I - A^\top \left[ A H^{-1} A^\top \right]^{-1} A H^{-1} \right) 
	}_{= P}
	\nabla F(z^0) + O(\norm{z^1-z^0}^2).
\end{equation}
%\begin{gather}
%	z^1 - z^0 = - \eta H^{-1} P \nabla F(z^0) + O(\norm{z^1-z^0}^2) \\
%	\text{where}~~~~
%	P \coloneqq I - A^\top \left[ A H^{-1} A^\top \right]^{-1} A H^{-1}.
%	\label{eq:mp_exactparam:MP:def_P} \tag{P}
%\end{gather}
% (Note that $P^2 = P$ and that $H^{-1} P$ is symmetric, i.e., $P$ is a projection which is orthogonal for $\innerprod{\cdot}{\cdot}_{H^{-1}}$). -> written in lemma statement
% This is exactly the expression from \cite[Proposition~1]{amid_reparameterizing_2020}.
To recap, we showed that
\begin{equation} \label{eq:mp_exactparam:MP:nablaFz1}
	\nabla F(z^1) = \nabla F(z^0) - \eta \nabla^2 F(z^0) H^{-1} P \nabla F(z^0) + O(\norm{z^1-z^0}^2).
\end{equation}

\paragraph{An order-2 expansion of $\delta z$ (the second update \eqref{eq:mp_exactparam:MP:def2}).}
Recall that we denote $\delta z = z^2 - z^0$.
From \eqref{eq:mp_exactparam:MP:statio2} and an order-1 expansion of $\nabla D^0$, by exactly the same calculations as in the previous paragraph,
\begin{align} \label{eq:mp_exactparam:MP:ord2_lambda2_deltaz}
	\lambda^2 &= \left[ A H^{-1} A^\top \right]^{-1} A H^{-1} \nabla F(z^1) + \frac{1}{\eta} O(\norm{\delta z}^2) \\
	\text{and}~~~~
	\delta z &= -\eta H^{-1} P \nabla F(z^1) + O(\norm{\delta z}^2).
\end{align}
However this is not precise enough for our goal, as the error is order-2 in $\eta$.

From \eqref{eq:mp_exactparam:MP:statio2} and an order-2 expansion of $\nabla D^0$, we have
\begin{align*}
	\eta \nabla F(z^1) + H_{ij} \delta z^j + \frac{1}{2} K_{ijk} \delta z^j \delta z^k - \eta A^\top \lambda^2 &= O(\norm{\delta z}^3) \\
	\left( H_{ij} + \frac{1}{2} K_{ijk} \delta z^k \right) \delta z^j &= -\eta \nabla F(z^1) + \eta A^\top \lambda^2 + O(\norm{\delta z}^3)
\end{align*}
where unmarked vectors are implicitly indexed by subscript $i$.
Denoting for concision
\begin{align*}
	G_{ij} &\coloneqq \frac{1}{2} K_{ijk} \delta z^k
	& & \text{and} &
	v &\coloneqq H^{-1} \left( \nabla F(z^1) - A^\top \lambda^2 \right)
\end{align*}
(note that $G$ is symmetric since $K$ is),
the above equation writes
\begin{align*}
	(H+G) \delta z &= -\eta Hv + O(\norm{\delta z}^3) \\
	(I + H^{-1} G) \delta z &= -\eta v + O(\norm{\delta z}^3).
\end{align*}
Now
$\norm{G} \lesssim \norm{\delta z} \lesssim \eta \norm{\nabla F(z^0)}$
so $(I + H^{-1} G)^{-1} = I - H^{-1} G + O\left( \eta^2 \norm{\nabla F(z^0)}^2 \right)$,
and by our order-1 estimates from \eqref{eq:mp_exactparam:MP:ord2_lambda2_deltaz} we have
$\norm{v} \lesssim \norm{P \nabla F(z^1)} + \frac{1}{\eta} \norm{\delta z}^2
\lesssim \norm{\nabla F(z^0)}$
using that $\norm{P \bullet} \lesssim \norm{\bullet}_{*\ZZZ}$ by \autoref{rk:mp_exacparam:norm_P_starZZZ}.
So
\begin{equation} \label{eq:mp_exactparam:MP:ord3_deltaz_prelim}
	\delta z = -\eta (I - H^{-1} G) v + O(\eta^3 \norm{\nabla F(z^0)}^3).
\end{equation}

It remains to estimate $v$, and namely the $A^\top \lambda^2$ term, up to $O \left( \eta^2 \norm{\nabla F(z^0)} \right)$ error terms.
To do this, write
\begin{align*}
	A \delta z = 0 &= -\eta A (I-H^{-1}G) v + O \left( \eta^3 \norm{\nabla F(z^0)}^3 \right) \\
	A (I-H^{-1}G) v &= O \left( \eta^2 \norm{\nabla F(z^0)}^3 \right).
\end{align*}
We already have an estimate of
% $\delta z$ and so of 
$G$ of the correct order thanks to \eqref{eq:mp_exactparam:MP:ord2_lambda2_deltaz}:
\begin{align}
%	\delta z &= -\eta H^{-1} P \nabla F(z^1) + O(\norm{\delta z}^2) \\
%	&= -\eta H^{-1} P \nabla F(z^0) + O \left( \eta^2 \norm{\nabla F(z^0)} \right) \\
%	\text{so}~~~
	G_{ij} &= \frac{1}{2} K_{ijk} \delta z^k \\
	&= \underbrace{
		\frac{1}{2} K_{ijk} \left[ -\eta H^{-1} P \nabla F(z^1) \right]^k 
	}_{\eqqcolon \tG_{ij}}
	+ O(\norm{\delta z}^2).
	\label{eq:mp_exactparam:MP:def_tG}
\end{align}
Since $\norm{v} \lesssim \norm{\nabla F(z^0)}$ we just need to solve for $A^\top \lambda^2$ in
$A (I - H^{-1} \tG) v = O \left( \eta^2 \norm{\nabla F(z^0)}^3 \right)$:
\begin{align*}
	A (I-H^{-1}\tG) H^{-1} \left( \nabla F(z^1) - A^\top \lambda^2 \right) &= O \left( \eta^2 \norm{\nabla F(z^0)}^3 \right) \\
	A (I-H^{-1}\tG) H^{-1} A^\top \lambda^2 &= A (I-H^{-1}\tG) H^{-1} \nabla F(z^1) + O \left( \eta^2 \norm{\nabla F(z^0)}^3 \right) \\
	\lambda^2 &= \left[ A (I-H^{-1}\tG) H^{-1} A^\top \right]^{-1} A (I-H^{-1}\tG) H^{-1} \nabla F(z^1) \\
    &\hspace{15em} + O \left( \eta^2 \norm{\nabla F(z^0)}^3 \right).
\end{align*}
Since $\norm{\tG} \lesssim \eta \norm{\nabla F(z^0)}$, we have the expansion
\begin{align*}
	A (I-H^{-1}\tG) H^{-1} A^\top
	&= A H^{-1} A^\top - A H^{-1} \tG H^{-1} A^\top \\
	&= A H^{-1} A^\top \left( I - \left[ A H^{-1} A^\top \right]^{-1} A H^{-1} \tG H^{-1} A^\top \right) \\
	\left[ A (I-H^{-1}\tG) H^{-1} A^\top \right]^{-1}
	&= \left( I - \left[ A H^{-1} A^\top \right]^{-1} A H^{-1} \tG H^{-1} A^\top \right)^{-1} \left[ A H^{-1} A^\top \right]^{-1} \\
	&= \left( I + \left[ A H^{-1} A^\top \right]^{-1} A H^{-1} \tG H^{-1} A^\top \right) \left[ A H^{-1} A^\top \right]^{-1} \\
	&\hspace{20em} + O \left( \eta^2 \norm{\nabla F(z^0)}^2 \right) \\
	&= \left( \left[ A H^{-1} A^\top \right]^{-1} + \left[ A H^{-1} A^\top \right]^{-1} A H^{-1} \tG H^{-1} A^\top \left[ A H^{-1} A^\top \right]^{-1} \right) \\
	&\hspace{20em} + O \left( \eta^2 \norm{\nabla F(z^0)}^2 \right).
%	\left[ A (I-\tG) H^{-1} A^\top \right]^{-1}
%	&= \left[ A H^{-1} A^\top - A \tG H^{-1} A^\top \right]^{-1} \\
%	&= \left( \left[ A H^{-1} A^\top \right]^{-1} + \left[ A H^{-1} A^\top \right]^{-1} A \tG H^{-1} A^\top \left[ A H^{-1} A^\top \right]^{-1} \right)
%	+ O \left( \eta^2 \norm{\nabla F(z^0)}^2 \right).
\end{align*}
Substituting and expanding the product,
and neglecting the terms in $\norm{\tG}^2 \norm{\nabla F(z^0)}$,
we get the following expression for $\lambda^2$; as a sanity-check, when we neglect the terms in $\norm{\tG} \norm{\nabla F(z^0)}$, we recover the estimate from \eqref{eq:mp_exactparam:MP:ord2_lambda2_deltaz}.
\begin{align*}
	& \lambda^2 
	+ O \left( \eta^2 \norm{\nabla F(z^0)}^3 \right) \\
	&= \left( \left[ A H^{-1} A^\top \right]^{-1} + \left[ A H^{-1} A^\top \right]^{-1} A H^{-1} \tG H^{-1} A^\top \left[ A H^{-1} A^\top \right]^{-1} \right)
	A (I-H^{-1}\tG) H^{-1} \nabla F(z^1) \\
	&= \left( \left[ A H^{-1} A^\top \right]^{-1} A + \left[ A H^{-1} A^\top \right]^{-1} A H^{-1} \tG H^{-1} A^\top \left[ A H^{-1} A^\top \right]^{-1} A - \left[ A H^{-1} A^\top \right]^{-1} A H^{-1} \tG \right)
	H^{-1} \nabla F(z^1) \\
	&= \left( \left[ A H^{-1} A^\top \right]^{-1} A - \left[ A H^{-1} A^\top \right]^{-1} A H^{-1} \tG \left( I - H^{-1} A^\top \left[ A H^{-1} A^\top \right]^{-1} A \right) \right)
	H^{-1} \nabla F(z^1) \\
	&= \left( \left[ A H^{-1} A^\top \right]^{-1} A H^{-1} - \left[ A H^{-1} A^\top \right]^{-1} A H^{-1} \tG H^{-1} \left( I - A^\top \left[ A H^{-1} A^\top \right]^{-1} A H^{-1} \right) \right)
	\nabla F(z^1) \\
	&= \left( \left[ A H^{-1} A^\top \right]^{-1} A H^{-1} - \left[ A H^{-1} A^\top \right]^{-1} A H^{-1} \tG H^{-1} P \right)
	\nabla F(z^1).
\end{align*}
Substituting, we get the following expression for $v = H^{-1} \left( \nabla F(z^1) - A^\top \lambda^2 \right)$:
\begin{align*}
%	v &= H^{-1} \left( \nabla F(z^1) - A^\top \lambda^2 \right) \\
	-Hv &= A^\top \lambda^2 - \nabla F(z^1) \\
	&= \left( A^\top \left[ A H^{-1} A^\top \right]^{-1} A H^{-1} - A^\top \left[ A H^{-1} A^\top \right]^{-1} A H^{-1} \tG H^{-1} P - I \right)
	\nabla F(z^1)
	+ O \left( \eta^2 \norm{\nabla F(z^0)}^3 \right) \\
	&= \left( -P - A^\top \left[ A H^{-1} A^\top \right]^{-1} A H^{-1} \tG H^{-1} P \right) \nabla F(z^1) 
	+ O \left( \eta^2 \norm{\nabla F(z^0)}^3 \right) \\
	v &= H^{-1} \left( I + A^\top \left[ A H^{-1} A^\top \right]^{-1} A H^{-1} \tG H^{-1} \right) P \nabla F(z^1)
	+ O \left( \eta^2 \norm{\nabla F(z^0)}^3 \right).
\end{align*}
Substituting into \eqref{eq:mp_exactparam:MP:ord3_deltaz_prelim} and using $\norm{G-\tG} \lesssim \eta^2 \norm{\nabla F(z^0)}^2$,
we get the following expression for $\delta z$:
\begin{align*}
	\delta z &= -\eta (I-H^{-1}\tG) v + O \left( \eta^3 \norm{\nabla F(z^0)}^3 \right) \\
	&= -\eta (I-H^{-1}\tG) H^{-1} \left( I + A^\top \left[ A H^{-1} A^\top \right]^{-1} A H^{-1} \tG H^{-1} \right) 
	P \nabla F(z^1)
	+ O \left( \eta^3 \norm{\nabla F(z^0)}^3 \right) \\
	&= -\eta H^{-1} (I-\tG H^{-1}) \left( I + A^\top \left[ A H^{-1} A^\top \right]^{-1} A H^{-1} \tG H^{-1} \right) 
	P \nabla F(z^1)
	+ O \left( \eta^3 \norm{\nabla F(z^0)}^3 \right) \\
	&= -\eta H^{-1} \left( I - \left( I- A^\top \left[ A H^{-1} A^\top \right]^{-1} A H^{-1} \right) \tG H^{-1} \right) 
	P \nabla F(z^1)
	+ O \left( \eta^3 \norm{\nabla F(z^0)}^3 \right) \\
	&= -\eta H^{-1} \left( I - P \tG H^{-1} \right) 
	P \nabla F(z^1)
	+ O \left( \eta^3 \norm{\nabla F(z^0)}^3 \right).
\end{align*}
Expanding and recalling the definition of $\tG$ \eqref{eq:mp_exactparam:MP:def_tG}, we get that
\begin{align*}
	\left[ \tG H^{-1} P \nabla F(z^1) \right]_i
	&= \tG_{ij} \left[ H^{-1} P \nabla F(z^1) \right]^j \\
	&= \frac{1}{2} K_{ijk} \left[ -\eta H^{-1} P \nabla F(z^1) \right]^k \left[ H^{-1} P \nabla F(z^1) \right]^j \\
	&= -\eta \frac{1}{2} K_{ijk} \left[ H^{-1} P \nabla F(z^1) \right]^k \left[ H^{-1} P \nabla F(z^1) \right]^j \\
	\text{or in shorthand,}~~~
	\tG H^{-1} P \nabla F(z^1)
	&= -\eta \frac{1}{2} K \left[ H^{-1} P \nabla F(z^1) \right]^2.
\end{align*}
So finally, we can write $\delta z$ as
\begin{equation} \label{eq:mp_exactparam:MP:ord3_deltaz_nablaFz1}
	\delta z = -\eta H^{-1} P \nabla F(z^1)
	- \eta^2 H^{-1} P \cdot \frac{1}{2} K \left[ H^{-1} P \nabla F(z^1) \right]^2
	+ O \left( \eta^3 \norm{\nabla F(z^0)}^3 \right).
\end{equation}
% (Incidentally, this calculation would directly give a precise expression for the MD update $\mathrm{MD}(z^0; \eta) = z^1$, simply replacing $z^2$ by $z^1$ and $\nabla F(z^1)$ by $\nabla F(z^0)$ everywhere.)

To make the expression of $\delta z$ fully explicit and conclude the analysis, let us substitute the expression of $\nabla F(z^1)$ from \eqref{eq:mp_exactparam:MP:nablaFz1}. Note that doing so makes us lose an order of precision in $\norm{\nabla F(z^0)}$ for the first term.
% , but that is perfectly okay for our goal.
\begin{align}
	\delta z &= -\eta H^{-1} P 
	\left[ \nabla F(z^0) - \eta \nabla^2 F(z^0) H^{-1} P \nabla F(z^0) \right]
	+ \eta^2 H^{-1} P \cdot \frac{1}{2} K \left[ H^{-1} P \nabla F(z^0) \right]^2
	+ O \left( \eta^3 \norm{\nabla F(z^0)}^2 \right)  \nonumber \\
	&= -\eta H^{-1} P \nabla F(z^0)
	+ \eta^2 H^{-1} P \cdot \left( \nabla^2 F(z^0) H^{-1} P \nabla F(z^0) - \frac{1}{2} K \left[ H^{-1} P \nabla F(z^0) \right]^2 \right)
	+ O \left( \eta^3 \norm{\nabla F(z^0)}^2 \right).
	\label{eq:mp_exactparam:MP:ord3_deltaz_nablaFz0}
\end{align}

\subsection{Proof of \autoref{lm:mp_exactparam:pp_expan}} \label{subsec:mp_exactparam:PP}

We keep the notations of the previous section, and this time we are interested in getting a similar Taylor expansion for the Proximal Point (PP) update $\mathrm{PP}(z^0; \eta) = z^\infty$ defined by
\begin{equation}
\label{eq:mp_exactparam:PP:definfty} \tag{U$\infty$}
	z^\infty = \argmin_{Az=b} F(z) + \frac{1}{\eta} D(z, z^0).
\end{equation}
The goal is to get an order-2 expansion for $\delta z \coloneqq z^\infty - z^0$.

Also again to lighten notation, we will write $\norm{\nabla F(z^0)}$ for $\norm{\nabla F(z^0)}_{*\ZZZ}$ and similarly for $\nabla F(z^\infty)$.

By Lagrangian duality, there exists $\lambda \in \RR^m$ such that
\begin{equation}
\label{eq:mp_exactparam:PP:statioinfty} \tag{S$\infty$}
	\nabla F(z^\infty) + \frac{1}{\eta} \nabla D^0(z^\infty) - A^\top \lambda = 0.
\end{equation}
By similar calculations as for MP, we get that
\begin{equation*}
	\norm{\delta z} 
	\lesssim \eta \norm{\nabla F(z^\infty)} 
	\lesssim \eta \norm{\nabla F(z^0)}.
\end{equation*}

\paragraph{An order-1 expansion of $\delta z$.}
From \eqref{eq:mp_exactparam:PP:statioinfty} and an order-1 expansion of $\nabla D^0$, we have
\begin{align*}
	\eta \nabla F(z^\infty) + H \delta z - \eta A^\top \lambda &= O(\norm{\delta z}^2) \\
	\delta z &= \eta H^{-1} \left( -\nabla F(z^\infty) + A^\top \lambda \right) + O(\norm{\delta z}^2).
\end{align*}
We can get an approximate expression of $A^\top \lambda$ by
\begin{align*}
	\frac{1}{\eta} A \delta z = 0 &= A H^{-1} \left( -\nabla F(z^\infty) + A^\top \lambda \right) + \frac{1}{\eta} O(\norm{\delta z}^2) \\
	A H^{-1} A^\top \lambda &= A H^{-1} \nabla F(z^\infty) + \frac{1}{\eta} O(\norm{\delta z}) \\
	\lambda &= \left[ A H^{-1} A^\top \right]^{-1} A H^{-1} \nabla F(z^\infty) + \frac{1}{\eta} O(\norm{\delta z}^2).
% \label{eq:mp_exactparam:PP:ord2_lambdainfty}
\end{align*}
Thus we get
\begin{equation*}
	\delta z = - \eta H^{-1} 
	\underbrace{ 
		\left( I - A^\top \left[ A H^{-1} A^\top \right]^{-1} A H^{-1} \right) 
	}_{= P}
	\nabla F(z^\infty) + O(\norm{\delta z}^2).
\end{equation*}

This limited-order expansion is sufficient for us to get an approximate expression of $\nabla F(z^\infty)$ in terms of $z^0$. Indeed,
\begin{align} \label{eq:mp_exactparam:PP:nablaFzinfty}
	\nabla F(z^\infty) &= \nabla F(z^0) + \nabla^2 F(z^0) \delta z + O(\norm{\delta z}^2) \\
	&= \nabla F(z^0) - \eta \nabla^2 F(z^0) H^{-1} P \nabla F(z^\infty) + O(\norm{\delta z}^2) \\
	\left( I + \eta \nabla^2 F(z^0) H^{-1} P \right) \nabla F(z^\infty) &= \nabla F(z^0) + O(\norm{\delta z}^2) \\
	\nabla F(z^\infty) &= \left( I - \eta \nabla^2 F(z^0) H^{-1} P \right) \nabla F(z^0) + O\left( \eta^2 \norm{\nabla F(z^0)} \right)
\end{align}
using \autoref{rk:mp_exacparam:norm_P_starZZZ} to control the error in the last line.
% Note already the similarity with \eqref{eq:mp_exactparam:MP:nablaFz1}. The only minor difference is that the error term is less precise by one order in $\norm{\nabla F(z^0)}$.

\paragraph{An order-2 expansion of $\delta z$.}
From \eqref{eq:mp_exactparam:PP:statioinfty} and this time an order-2 expansion of $\nabla D^0$, we have more precisely
\begin{align*}
% \label{eq:mp_exactparam:PP:statioinfty_ord2}
	\eta \nabla F(z^\infty) + H_{ij} \delta z^j + \frac{1}{2} K_{ijk} \delta z^j \delta z^k - \eta A^\top \lambda &= O(\norm{\delta z}^3) \\
	\left( H_{ij} + \frac{1}{2} K_{ijk} \delta z^k \right) \delta z^j &= -\eta \nabla F(z^\infty) + \eta A^\top \lambda + O(\norm{\delta z}^3)
\end{align*}
where unmarked vectors are implicitly indexed by subscript $i$.
Denote for concision
\begin{align*}
	G_{ij} &\coloneqq \frac{1}{2} K_{ijk} \delta z^k
	& & \text{and} &
	v &\coloneqq H^{-1} \left( \nabla F(z^\infty) - A^\top \lambda \right).
\end{align*}

We can unroll the exact same calculations as in the last paragraph of \autoref{subsec:mp_exactparam:MP} with $\nabla F(z^1)$ replaced by $\nabla F(z^\infty)$, and we obtain an equivalent of \eqref{eq:mp_exactparam:MP:ord3_deltaz_nablaFz1} for PP:
\begin{equation*}
	\delta z = -\eta H^{-1} P \nabla F(z^\infty)
	- \eta^2 H^{-1} P \cdot \frac{1}{2} K \left[ H^{-1} P \nabla F(z^\infty) \right]^2
	+ O \left( \eta^3 \norm{\nabla F(z^0)}^3 \right).
\end{equation*}

Since the expression of $\nabla F(z^\infty)$ in terms of $z^0$ for PP \eqref{eq:mp_exactparam:PP:nablaFzinfty}
is exactly the same as the one of $\nabla F(z^1)$ for MD \eqref{eq:mp_exactparam:MP:nablaFz1}, 
the very last step of the calculations is also the same, and we get
\begin{equation}
\label{eq:mp_exactparam:PP:ord3_deltaz_nablaFz0}
	\delta z
	= -\eta H^{-1} P \nabla F(z^0)
	+ \eta^2 H^{-1} P \cdot \left( \nabla^2 F(z^0) H^{-1} P \nabla F(z^0) - \frac{1}{2} K \left[ H^{-1} P \nabla F(z^0) \right]^2 \right)
	+ O \left( \eta^3 \norm{\nabla F(z^0)} \right).
\end{equation}
The only difference is that we lose yet another order of precision in $\norm{\nabla F(z^0)}$ in the error term compared to \eqref{eq:mp_exactparam:MP:ord3_deltaz_nablaFz0}.
% , but it is still precise enough for our purpose.

% Thus, these two subsections prove that in the conditions described at the beginning, the MP and the PP updates differ only by $O \left( \eta^3 \norm{\nabla F(z^0)} \right)$.

% \begin{remark}
% 	If $D(\cdot, \cdot)$ was known to be a Bregman divergence, then we could use that $\nabla D(z, \hz) = -\nabla D(\hz, z)$ 
% 	(where $\nabla$ denotes derivative w.r.t.\ the first argument),
% 	so that the MP \eqref{eq:mp_exactparam:MP:def1} and PP updates \eqref{eq:mp_exactparam:PP:definfty} are dual in the sense that
% 	$\hz = \mathrm{MP}(z; \eta) \iff z = \mathrm{PP}(\hz; -\eta)$.
% 	We could then derive an approximate expression for $\nabla F(z^\infty)$ in terms of $z^0$ by inverting the relation \eqref{eq:mp_exactparam:MP:nablaFz1}, taking care to first establish that $(H^{-1} P)_{z^\infty} = (H^{-1} P)_{z^0} + O(\norm{\delta z})$.
% 	We would directly obtain the desired expansion by substituting into the MD version of \eqref{eq:mp_exactparam:MP:ord3_deltaz_nablaFz1}.
% 
% 	In fact one can check from \eqref{eq:mp_exactparam:PP:statioinfty_ord2} that under our assumptions, the MD and PP updates are still dual ``up to'' $O(\norm{\delta z}^3)$ terms, so the above proof strategy would also work. Since it would 
% % 	anyway require performing the same computations (inverting the relation between $\nabla F(z^\infty)$ and $\nabla F(z^0)$), 
%     require performing no less computations,
% 	we chose to present a direct proof for ease of reading.
% \end{remark}

% \end{document}
%
\fi
% !TEX root = ../main.tex
% \documentclass[../main]{subfiles}
% \begin{document}

\section{Proof of the main result} \label{sec:proof_mainres}

In this section we show in detail how our main result \autoref{thm:main_res:loc_exp_cv_NI} follows from combining \autoref{prop:cv_proof:gencase:NI_equiv_V}, \autoref{prop:cv_proof:gencase:WFR_equiv_V} and \autoref{thm:cv_proof:gencase:loc_exp_cv}.

\begin{proof}
    Fix $\Gamma_0 \geq 1$.
    Choose $\eta_0, \sigma_0$ as in \autoref{thm:cv_proof:gencase:loc_exp_cv}.
    Fix any $\eta \leq \eta_0, \sigma \leq \sigma_0$ with $\Gamma_0^{-1} \leq \frac{\sigma}{\eta} \leq \Gamma_0$.
    
    Let $\lambda, \tau$ as in \eqref{eq:proof_gencase:choice_lambda_tau},
    let $(\varphi_I)_I, (\psi_J)_J$ as in \eqref{eq:cv_proof:gencase:def_varphi},
    and let $V_1$ and $V$ as in \eqref{eq:cv_proof:gencase:def_V}.
    
    Let $\tC, \tC', \tr$ as in \autoref{prop:cv_proof:gencase:NI_equiv_V}.
    Let $\tilde{\tilde{C}} ,\tilde{\tr}$ as in \autoref{prop:cv_proof:gencase:WFR_equiv_V}.
    Let $K$ resp.\ $R_0$ the quantities denoted $\kappa$ resp.\ $r_0$ in \autoref{thm:cv_proof:gencase:loc_exp_cv}.

    Let $r_0 = \min\left\{ \tr, \tilde{\tr}, \tC' \left( \frac{R_0}{\Gamma_0} \right)^{5/4} \right\}$.
    Denote $(z^k)_k$ the CP-PP iterates and $\mu^k = \sum_{i=1}^n a^k_i \delta_{x^k_i}$, $\nu^k = \sum_{j=1}^m b^k_j \delta_{y^k_j}$.
    Suppose $\NI(\mu^0, \nu^0) \leq r_0$, then by the second part of \autoref{prop:cv_proof:gencase:NI_equiv_V},
    \begin{align*}
        \tC' V_1(z^0)^{5/4} &\leq \NI(\mu^0, \nu^0) \leq r_0 \\
        \implies~
        V_1(z^0) &\leq \left( \frac{r_0}{\tC'} \right)^{4/5} \\
        \implies~
        V(z^0) &\leq \Gamma_0 V_1(z^0) \leq \Gamma_0 \left( \frac{r_0}{\tC'} \right)^{4/5} \leq R_0.
    \end{align*}
    So by \autoref{thm:cv_proof:gencase:loc_exp_cv},
    $V(z^k) \leq V(z^0) (1-K)^k$, 
    and so by the first part of \autoref{prop:cv_proof:gencase:NI_equiv_V},
    \begin{equation*}
        \NI(\mu^k, \nu^k) \leq \tC \sqrt{V(z^k)} 
        \leq \tC \sqrt{V(z^0)} \left( \sqrt{1-K} \right)^k
        \leq \tC \sqrt{V(z^0)} \left( 1 - \frac{K}{2} \right)^k.
    \end{equation*}
    This proves the first inequality of \autoref{thm:main_res:loc_exp_cv_NI} by letting $C = \tC \sqrt{R_0}$ and $\kappa = \frac{K}{2}$.

    Moreover, by \autoref{prop:cv_proof:gencase:WFR_equiv_V},
    \begin{align*}
        \WFR_2^2(\mu^k, \mu^*) + \WFR_2^2(\nu^k, \nu^*) 
        &\leq 2 V(z^k) \left( 1 + \tilde{\tilde{C}} \frac{\eta}{\sigma} (\lambda \tau)^2 \right) \\
        &\leq 2 V(z^0) (1-K)^k \left( 1 + \tilde{\tilde{C}}~ \Gamma_0 (\lambda \tau)^2 \right) \\
        &\leq \underbrace{2 R_0 \left( 1 + \tilde{\tilde{C}}~ \Gamma_0 (\lambda \tau)^2 \right)} (1-\kappa)^k,
    \end{align*}
    which proves the second inequality of \autoref{thm:main_res:loc_exp_cv_NI} by letting $C'$ be the underbraced expression.
\end{proof}

% \end{document}

\end{document}